\documentclass[12pt]{article}
\usepackage{graphicx}
\usepackage{amsmath}
\usepackage{amsfonts}
\usepackage{amssymb}
\usepackage{ifthen}
\usepackage{bm,bbm}
\usepackage{longtable}

\newtheorem{theorem}{Theorem}[section]
\newtheorem{lemma}[theorem]{Lemma}
\newtheorem{corollary}[theorem]{Corollary}

\newtheorem{definition}[theorem]{Definition}
\newtheorem{example}[theorem]{Example}

\newcommand{\DEFINED}[1]{{\bf #1}}

\newcommand{\SET}[1]{\mathcal{#1}}

\newcommand{\Ii}{{\bf i}}

\newcommand{\ZEROvect}{\mathbf{0}}

\newcommand{\ONESvect}{\mathbf{e}}

\newcommand{\PHASE}[1]{e^{\Ii #1}}

\newcommand{\STbasis}[1]{e_{#1}}

\newcommand{\REALS}{\mathbbm{R}}

\newcommand{\COMPLEX}{\mathbbm{C}}

\newcommand{\BISTOCHSPACE}{\mathcal{B}}
\newcommand{\UNITARY}{\mathcal{U}}

\newcommand{\COMPLEXmatricesNxN}[1]{\mathbbm{M}_{#1}}

\newcommand{\columnVECTORS}[1]{\mathbbm{M}_{{#1},1}}

\newcommand{\REALcolumnVECTORS}[1]{\mathbbm{M}_{{#1},1}(\mathbbm{R})}

\newcommand{\FIXEDmoduliMATRICES}[1]{\mathcal{M}_{#1}}

\newcommand{\ENPHASEDmatrices}[1]{\mathcal{T}_{#1}}

\newcommand{\IDENTITY}[1]{\mathbf{id}_{#1}}

\newcommand{\CONJ}[1]{\overline{#1}}

\newcommand{\ABSOLUTEvalue}[1]{\left| #1 \right|}

\newcommand{\diag}{\mathbf{diag}}

\newcommand{\trace}{\mathbf{tr}}

\newcommand{\RELofEQUI}{\simeq}

\renewcommand{\Re}{\mathbf{Re}}
\renewcommand{\Im}{\mathbf{Im}}

\newcommand{\DIFFERENTIAL}[2]{{D \! #1}_{#2}}

\newcommand{\TANGENTspace}[2]{\mathbf{T}_{#2} #1}

\newcommand{\SPANNEDspace}[1]{{\mathbf{span}}_{#1}}

\newcommand{\SPANC}{\SPANNEDspace{\mathbbm{C}}}

\newcommand{\NULLspace}[1]{{\mathbf{null}}_{#1}}
\newcommand{\NULL}{\NULLspace{}}

\newcommand{\DIMR}{\dim_{\REALS}}
\newcommand{\DIMC}{\dim_{\COMPLEX}}

\newcommand{\SPACE}[1]{\mathbbm{#1}}

\newcommand{\DEFECT}{\mathbf{d}}

\newcommand{\undephasedDEFECT}{\mathbf{D}}

\newcommand{\FEASIBLEspace}[1]{\mathbbm{D}_{#1}}

\newcommand{\FEASIBLEspaceDEPHASED}[1]{\mathbbm{D}'_{#1}}

\newcommand{\TRIVIALspace}[1]{\mathbbm{T}_{#1}}

\newcommand{\RANK}{\mathbf{rank}}

\newcommand{\lcm}{\mathbf{lcm}}

\newcommand{\HADprod}{\circ}

\newcommand{\matrixEXP}[1]{e^{#1}}

\newcommand{\ELEMENTof}[3]{{\left[ #1 \right]}_{#2,#3}}

\newcommand{\ROWof}[2]{{\left[ #1 \right]}_{#2,:}}

\newcommand{\COLUMNof}[2]{{\left[ #1 \right]}_{:,#2}}

\newcommand{\TRANSPOSE}[1]{{\left( #1 \right)}^T}

\newcommand{\hermTRANSPOSE}[1]{{\left( #1 \right)}^{*}}

\newcommand{\DELTA}[2]
{
  \ifthenelse{\equal{#1}{} \or \equal{#2}{}}
    {\Delta^{#1#2}}
    {\Delta^{#1,#2}}
}

\newcommand{\ithMATRIX}[2]{{#1}^{(#2)}}

\newcommand{\vectorINDEX}[1]{\left( #1 \right)}

\newcommand{\PROOFstart}{{\bf Proof}\\}
\newcommand{\PROOFend}{$\blacksquare$}

\newcommand{\Cu}[1]{\mathcal{C}_{#1}}
\newcommand{\Du}[1]{\mathcal{D}_{#1}}
\newcommand{\Iu}[1]{\mathcal{I}_{#1}}
\newcommand{\Ju}[1]{\mathcal{J}_{#1}}

\newcommand{\Vu}[2]{\mathbbm{V}^{#1}_{#2}}

\newcommand{\VuOLD}[1]{\mathbbm{V}_{#1}}

\newcommand{\REspaceOF}[1]{\left( {#1} \right)^{\Re}}
\newcommand{\IMspaceOF}[1]{\left( {#1} \right)^{\Im}}

\newcommand{\STinnerPRODUCT}[2]{\left\langle {#1}, {#2} \right\rangle}

\newcommand{\UinnerPRODUCT}[3]{\left\langle {#2}, {#3} \right\rangle_{#1}}

\newcommand{\DIAGwith}[1]{D_{#1}}

\newcommand{\leftPERM}[2]{P'_{{#1},{#2}}}
\newcommand{\rightPERM}[2]{P''_{{#1},{#2}}}
\newcommand{\leftDIAG}[2]{D'_{{#1},{#2}}}
\newcommand{\rightDIAG}[2]{D''_{{#1},{#2}}}

\newcommand{\pairOPERATOR}[2]{\left[\!\!\left[ {#1}, {#2} \right]\!\!\right]}

\newcommand{\RESTRICTEDto}[2]{{#1}_{\left| {#2} \right.}}

\newcommand{\UorthCOMPLEMENT}[2]{{#2}^{\bot_{\UinnerPRODUCT{{#1}}{}{}}}}

\newcommand{\FOURIER}[1]{F_{#1}}

\newcommand{\Igroup}[1]{\mathrm{I}_{#1}}

\newcommand{\enphasedPERMS}{\mathbf{P}}       
\newcommand{\PP}{\enphasedPERMS^2}              
\newcommand{\PERMS}{\mathbf{P}_\mathbf{S}}    
\newcommand{\PsPs}{\PERMS^2}                           
\newcommand{\enphasedSHIFTS}[1]{\mathbf{P}_{\mathbf{Z}^{#1}}}     
\newcommand{\PzPz}[1]{{\enphasedSHIFTS{#1}}^{\!\!\!\!\!\!\!\!\! 2\mbox{   }}}                            

\newcommand{\ZN}[1]{\mathbbm{Z}_{#1}}      
\newcommand{\ZNclass}[1]{\mbox{\boldmath $[$}#1\mbox{\boldmath $]$}}    
\newcommand{\REPR}[1]{\tilde{#1}}   

\newcommand{\ZNprimePOWER}[5]{\ZN{{#1}_{#2}^{{#3}_{#4}^{(#5)}}}}

\newcommand{\uniMORPHISM}[1]{#1}       
\newcommand{\ISOM}[1]{#1}         

\newcommand{\ISOMORPHISMS}[2]{\mathit{Isom}\left( \Igroup{#1}, \Igroup{#2} \right)}

\newcommand{\GROUP}[1]{\mathbf{#1}}     
\newcommand{\GsymmetryGROUP}[1]{\GROUP{G}^{#1}}    
\newcommand{\STAB}[2]{Stab_{#1}\left(#2\right)}      
\newcommand{\MAP}[3]{Map_{#1}\left(#2,#3\right)}   

\newcommand{\hMAP}[3]{h^{{#1},{#2}}_{#3}}

\newcommand{\INVERSE}[1]{\left( #1 \right)^{-1}}

\newcommand{\INTEGER}{\mathbbm{Z}}

\newcommand{\Wk}[2]{W^{#1}_{#2}}

\newcommand{\Zk}[2]{Z^{#1}_{#2}}

\newcommand{\Pro}[1]{P_{#1}}

\newcommand{\Smap}[1]{\mathcal{S}^{#1}}    
\newcommand{\Tmap}[1]{\mathcal{T}^{#1}}    

\newcommand{\Sformula}[4]{{#4} \cdot \Wk{#1}{#2} \Zk{#1}{#3}}
\newcommand{\Tformula}[4]{{#4}  \CONJ{{#1}_{{#2},{#3}}} \cdot  \Wk{#1}{#3} \Zk{#1}{-{#2}}} 

\newcommand{\SUPP}{\mathbf{supp}}

\newcommand{\kxBELOW}[1]{{\tilde{m}}_{#1}}

\newcommand{\ROWtype}[2]{%
                          \ifthenelse{\equal{#2}{0}}
                             {1}
                             {\ifthenelse{\equal{#2}{1}}
                                 {1/(#1)}
                                 {1/({#1}^{#2})}%
                             }%
                        }

\newcommand{\ROWStype}[2]{p_{\ROWtype{#1}{#2}}}

\newcommand{\ORDERinGROUP}[2]{\mathbf{ord}_{#1}\left( {#2} \right)}

\newcommand{\GROUPzero}[1]{0_{#1}}

\newcommand{\defectPOLY}[2]{\mathcal{P}\left({#1},{#2}\right)}

\newcommand{\VALUEat}[2]{{#1}_{\left| {#2} \right.}}

\begin{document}

\title{Defect and equivalence of unitary matrices.\\ The Fourier case.}

\author{
  Wojciech Tadej   \\
  \smallskip
  {\small Center for Theoretical Physics, Polish Academy of Sciences, Warsaw, Poland}\\
 \smallskip
 {\small e-mail: wtadej@wp.pl}
%
}
\date{\today}

\maketitle

\begin{abstract}
Consider the real space $\FEASIBLEspace{U}$ of directions moving into which from a unitary $N \times N$ matrix $U$ we do not disturb 
its unitarity and the moduli  of its entries in the first order. $\DIMR\left( \FEASIBLEspace{U} \right)$ is called the defect of $U$
and denoted $\undephasedDEFECT(U)$. We give an account of Alexander Karabegov's theory where $\FEASIBLEspace{U}$
is parametrized by the imaginary subspace of the eigenspace, associated with $\lambda = 1$, of a certain unitary operator $\Iu{U}$
on $\COMPLEXmatricesNxN{N}$, and where $\undephasedDEFECT(U)$ is the multiplicity of $1$ in the spectrum of $\Iu{U}$.
This characterization allows us to establish dependence of 
	$\undephasedDEFECT\left( \ithMATRIX{U}{1} \otimes \ldots \otimes  \ithMATRIX{U}{r} \right)$ on 
$\undephasedDEFECT\left( \ithMATRIX{U}{k} \right)$'s,
to derive formulas expressing $\undephasedDEFECT(F)$ for a Fourier matrix $F$ of the size being a power of a prime number,
as well as to show the multiplicativity of $\undephasedDEFECT(F)$ with respect to Kronecker factors of $F$  if their
sizes are pairwise relatively prime. Also partly due to the role of symmetries of $U$ in the determination of the 
eigenspaces of $\Iu{U}$ we study the 'permute and enphase' symmetries and equivalence of Fourier matrices,
associated with arbitrary finite abelian groups.
\end{abstract}

%
%
\section{Introduction}
	\label{sec_introduction}

The notion of defect of an $N \times N$ unitary matrix is a fruit of our interest in the structure of the set of  $N \times N$ unistochatic matrices,
within the set of all $N \times N$ doubly stochastic matrices (being a convex polytope of the full dimension within $\BISTOCHSPACE$, the 
$(N-1)^2$ dimensional real manifold -- a hyperplane -- formed by real matrices with the row and column sums equal to $1$).
By unistochastic matrices we mean those doubly stochastic matrices which are obtained from unitary matrices (belonging to the real $N^2$ dimensional
unitary manifold $\UNITARY$) by squaring the moduli of their entries. For a given $U \in \UNITARY$ and the map $f$ performing this operation:\ \ 
$\ELEMENTof{f(U)}{i}{j}\ =\ \ABSOLUTEvalue{U_{i,j}}^2$\ \ we considered the tangent map (the differential) of $f$ at $U$ restricted to the space tangent
to $\UNITARY$ at $U$:\ \ $\RESTRICTEDto{\DIFFERENTIAL{f}{U}}{\TANGENTspace{\UNITARY}{U}}$.\ \ 
The hyperplane $f(U)\ +\ \DIFFERENTIAL{f}{U}\left( \TANGENTspace{\UNITARY}{U} \right)$ approximates, around $f(U)$, only the part of the unistochastic set 
$f\left( \UNITARY \right)$ coming from a small neighbourhood of $U$  in $\UNITARY$, for there may be other $V \in \UNITARY$ mapped by $f$ into $f(U)$.

We were wondering whether $f(U) + \DIFFERENTIAL{f}{U}\left( \SET{V} \subset \TANGENTspace{\UNITARY}{U} \right)$, for $\SET{V}$ being
a small open $N^2$ dimensional ball around $\ZEROvect$, is an open $(N-1)^2$ dimensional  ellipsoid around $f(U)$ in $\BISTOCHSPACE$. In particular we
asked this question about the 'cyclic' Fourier matrix 
$\frac{1}{\sqrt{N}} \FOURIER{N}\ =\ \left[ \frac{1}{\sqrt{N}} \PHASE{\frac{2\pi}{N} i j} \right]_{i,j=0\ldots (N-1)}$ because we wanted to know
what the unistochastic set looks like around the flat doubly stochastic matrix $f\left( \frac{1}{\sqrt{N}} \FOURIER{N} \right)$ filled all with $\frac{1}{N}$'s.
The difference between $\TANGENTspace{\BISTOCHSPACE}{f(U)}\ =\ (N-1)^2$ and 
the dimension of the above ellipsoid was given the name of defect. It was more formally defined as
\begin{equation}
	\label{eq_defect_original_definition}
	\DEFECT(U)\ \ \stackrel{def}{=}\ \ (N-1)^2\ -\ \DIMR\left( \DIFFERENTIAL{f}{U}\left( \TANGENTspace{\UNITARY}{U} \right) \right)
\end{equation}
and investigated in \cite{Defect}. The restricted differential (a real linear map), as any other linear map between finite dimensional spaces,
satisfies the law saying that the sum of  dimensions of its kernel and its image is equal to the dimension of the space of arguments, hence
\begin{eqnarray}
	\label{eq_dim_ker_image_sum_equality_for_differential_of_f}
	\DIMR\left( \TANGENTspace{\UNITARY}{U} \right)
	&   =    &
	\DIMR\left( \ker \RESTRICTEDto{\DIFFERENTIAL{f}{U}}{\TANGENTspace{\UNITARY}{U}} \right)
		\ \ \ \ \ \ \ \ +\ \ \ \ \ \ \ \ 
	\DIMR\left(  \DIFFERENTIAL{f}{U}\left( \TANGENTspace{\UNITARY}{U} \right)  \right)         \\
	\nonumber
	 \parallel\ \ \ \ \ \ \ \ \ \ \ \ \      &             &     
	 \         \parallel   \ \ \ \ \ \ \ \ \ \ \ \ \ \ \ \ \ \ \ \ \ \ \ \ \ \ \ \ \ \ \  \ \ \ \ \ \ \ \ \ \ \ \ \ \ \   \parallel      \\
	\nonumber
	N^2\ \ \ \ \ \ \ \ \ \ \ 	              &      &    
	\undephasedDEFECT(U) \stackrel{def}{=} \DEFECT(U) + (2N-1)     \ \ \ \ \ \ \ \ \ \ \ \  (N-1)^2 - \DEFECT(U)    \ \ .
\end{eqnarray}

In the present work we define the defect with the use of the left component of the right hand side of the above equality.
The kernel there is the intersection of $\TANGENTspace{\UNITARY}{U}$ and $\ker \DIFFERENTIAL{f}{U}$. The latter space
\underline{is} equal to $\TANGENTspace{\FIXEDmoduliMATRICES{U}}{U}$, the space tangent at $U$ to $\FIXEDmoduliMATRICES{U}$,
the real $N^2$ dimensional manifold (considered in the next section) of matrices with the moduli of their entries fixed and identical with those of $U$,
\underline{under} the crucial assumption that $U$ has no zero entries. The reason for making this assumption is that it is required by
the characterization of the defect described in section \ref{sec_Karabegovs_theory}, which is later used in the calculation of the defect of
(rescaled unitary) Fourier matrices whose entries are unimodular. Thus the defect $\DEFECT(U)$ is defined
here as
\begin{equation}
	\label{eq_defect_new_definition}
	\DEFECT(U) 
			\ \ \stackrel{def}{=}\ \ 
	\DIMR\left(  \TANGENTspace{\UNITARY}{U} \ \cap\ \TANGENTspace{\FIXEDmoduliMATRICES{U}}{U}   \right)
		\ \ -\ \ 
	(2N-1)
\end{equation}
and this change of perspective is associated with the fact later realized by us that $\DEFECT(U)$ can serve as an upper bound for the
dimension of a dephased (e.g. with the $1$st row and column fixed) manifold stemming from $U$ and contained in 
$\UNITARY \cap \FIXEDmoduliMATRICES{U}$. In particular manifolds of complex Hadamard matrices are quite of interest, 
due to their applications in mathematics and quantum physics, consisting
of (rescaled) unitary matrices with all the moduli of their entries  identical (see for example \cite{Catalogue_printed},\cite{Catalogue}).
A more detailed description of the two types of defect, the undephased $\undephasedDEFECT(U)$ and the dephased $\DEFECT(U)$, the reader 
will find in section \ref{sec_defect}. We think that this presentation is more readable than that in \cite{Defect}, and it is
the result of our growing experience and discussions with our colleagues and students.

Calculating the defect $\undephasedDEFECT(U)$ amounts to finding the solution to some linear system associated with $U$,
which is equivalent to finding all the directions moving into which from $U$ we do not disturb the unitarity and the  moduli in the first order,
and calculating the dimension of this solution. The considered linear system can be reshaped to form an eigenvalue problem for a certain operator
determined by $U$, and this approach, taken by Alexander Karabegov in \cite{Karabegov} and \cite{Karabegov_old}, was described, in our
modified notation and with additional explanations,  in section \ref{sec_Karabegovs_theory}. Most of this section is based on the hints and information
obtained from Alexander in one way or another, except for the application of his theory to the Kronecker products of unitary matrices
in subsection \ref{subsec_Kronecker_products}. This last item supplements our previous work on the subject \cite{KroneckerDefect}.
In the final part of subsection \ref{subsec_equivalence}, as well as only at some alternative proof in section \ref{sec_fourier_defect}
we expect that the reader has competence in the theory of complex representations of finite groups. This knowledge makes it possible to 
fully appreciate the way in which Alexander Karabegov exploits potential symmetries of $U$ to discover the eigenspaces of the 
operator he uses and calculate the defect.

By the symmetries of $U$ we mean the equivalences operations mapping $U$ back into $U$, where an equivalence operation 
comprises permuting the rows and columns and/or multiplying  the rows and columns by unimodular numbers. 
As it is explained in subsection \ref{subsec_equivalence}, the symmetries of $U$
preserve the space of directions not violating the moduli and the unitarity in the first order. Moreover, the symmetries map smooth manifolds
stemming from $U$  and consisting of unitary matrices with fixed moduli (that is manifolds contained in $\UNITARY \cap \FIXEDmoduliMATRICES{U}$,
e.g. those comprised of complex Hadamard matrices) into the very same or similar objects also stemming from $U$. The same applies to 
the $N \times N$ affine unitary  families stemming from $U$ (e.g. affine Hadamard families, see \cite{Catalogue}), a subtype of the previously mentioned ones,
whose members additionally have their matrices of phases forming affine subspaces of $\REALS^{N \times N}$. Therefore we think that
on the one hand the symmetries of $U$ are worth studying as they, to some degree of course, determine the aforementioned objects, since
some of these symmetries are just the symmetries of these very objects. On the other hand, the symmetry and in general the equivalence 
operations map these objects (stemming from $U$) onto similar ones (stemming from $U$ or an equivalent matrix) and as such
can be used in their classification, much like they separate the unitary matrices (or for example complex Hadamard matrices only) into the equivalence
classes. This was the motivation behind the creation of section \ref{sec_Fourier_matrices}, where the symmetries of Fourier matrices and equivalence between them 
are studied. 

Finally, in section \ref{sec_fourier_defect} we take advantage of the Alexander Karabegov's characterization of the defect $\undephasedDEFECT(U)$
as the multiplicity of $1$ in the spectrum of the appropriate operator to derive the formulas expressing the undephased defect of a Fourier matrix
and certain properties of this quantity in the context of Fourier matrices. 
The results in this section were announced in \cite{Mattriad2011}, and they were obtained independently by Teodor Banica in \cite{Banica_def_gen_four_matr},
who used different techniques. 
By an accident, we and Teodor conceived the same names, dephased and
undephased, for the two types of defect. Probably there is no better alternative.

In section \ref{sec_conclusions} we summarize the contents of this article.

%
%
\section{The defect of a (rescaled) unitary matrix.}
	\label{sec_defect}

Denoting by $\COMPLEXmatricesNxN{N}$ the $N^2$ dimensional complex space of all $N \times N$ complex matrices,
which is also a $2N^2$ dimensional real space, let
\begin{itemize}
	\item
		$g:\ \COMPLEXmatricesNxN{N} \longrightarrow \COMPLEXmatricesNxN{N}$
		denote the map returning the Gram matrix of its argument,
		\begin{equation}
			\label{eq_g_Gram_matrix_map}
			g(U)\ \ \stackrel{def}{=}\ \       \hermTRANSPOSE{U} U               \ \        ,
		\end{equation}

	\item
		$f:\ \COMPLEXmatricesNxN{N} \longrightarrow \COMPLEXmatricesNxN{N}$
		denote the map squaring the moduli of entries of its argument,
		\begin{equation}
			\label{eq_f_squaring_moduli_map}
			f(U)\ \ \stackrel{def}{=}\ \              U \HADprod \CONJ{U}                    \ \ ,
		\end{equation}
\end{itemize}
where $\HADprod$ denotes the entrywise (Hadamard) product.

The differentials of these maps at $U$, i.e. the linear parts of the maps\ \  $D\ \longrightarrow\ g(U + D)$,\ \ $D\ \longrightarrow\ f(U+D)$,
understood as linear maps between $\COMPLEXmatricesNxN{N}$ and $\COMPLEXmatricesNxN{N}$ treated as real spaces, are easily
computed as:
\begin{eqnarray}
	\label{eq_g_differential_at_U}
	\DIFFERENTIAL{g}{U}(D)
	&   =   &
    \hermTRANSPOSE{D} U      \ +\      \hermTRANSPOSE{U} D                                         \\
	\label{eq_f_differential_at_U}
	\DIFFERENTIAL{f}{U}(D)
	&    =    &
      D \HADprod \CONJ{U}      \ +\        U \HADprod \CONJ{D}     
     \ \ =\ \    2 \Re\left( D \HADprod \CONJ{U} \right)
\end{eqnarray}

From now on we assume that $U$ has no zero entries. Let $U$, our reference point, be a unitary matrix.

Next we consider the set $\UNITARY \subset \COMPLEXmatricesNxN{N}$ of unitary matrices and 
the set $\FIXEDmoduliMATRICES{U} \subset \COMPLEXmatricesNxN{N}$ of matrices with the moduli of their entries identical with those in $U$:
\begin{eqnarray}
	\label{eq_unitary_matrices}
	\UNITARY
	&   \stackrel{def}{=}   &
	\left\{\ V\ :\ \ g(V)\ \ =\ \ I\ \right\}             \ \ ,                     \\
	\label{eq_fixed_moduli_matrices}
	\FIXEDmoduliMATRICES{U}
	&   \stackrel{def}{=}    &
	\left\{\ V\ :\ \ f(V)\ \ =\ \ f(U)\ \right\}          \ \ .
\end{eqnarray}
Both these sets are known to be $N^2$ dimensional real manifolds in $\COMPLEXmatricesNxN{N}$, so they have well defined
tangent spaces at $U$, $\TANGENTspace{\UNITARY}{U}$ and $\TANGENTspace{\FIXEDmoduliMATRICES{U}}{U}$, of dimension $N^2$.
We will explicitly construct them for the purpose of this paper.

By a vector tangent to $\UNITARY$\ \ ( to $\FIXEDmoduliMATRICES{U}$ )\ \  at $U$ we understand the derivative, at $0$, of any differentiable curve 
$\gamma:\ \REALS \longrightarrow \UNITARY$\ \ ( $\gamma:\ \REALS \longrightarrow \FIXEDmoduliMATRICES{U}$ ),\ \ such that
$\gamma(0)\ =\ U$. All such vectors form the tangent space $\TANGENTspace{\UNITARY}{U}$\ \ ( $\TANGENTspace{\FIXEDmoduliMATRICES{U}}{U}$ ),
which is a real linear space. This will be justified in the next two paragraphs. 

Consider a differentiable curve $\gamma:\ \REALS \longrightarrow \UNITARY$. $g(\gamma(t)) - I\ =\ \ZEROvect$ is satisfied in a certain
range of parameter $t$ to which $0$ (corresponding to $\gamma(0)\ =\ U$) belongs. Differentiating this equality with respect to $t$ at $0$ we get
$\DIFFERENTIAL{g}{U}(\gamma'(0))\ =\ \ZEROvect$, hence the tangent vector $\gamma'(0)$ belongs to $\NULL \DIFFERENTIAL{g}{U}$.
The space $\TANGENTspace{\UNITARY}{U}$ is thus contained in this nullspace, but in fact these spaces are equal.  
$\NULL \DIFFERENTIAL{g}{U}$ is easily computed as $\{EU:\ E\ \mbox{antihermitian}\}$ and has dimension $N^2$. On the other hand,
$\TANGENTspace{\UNITARY}{U}$ contains the derivatives at $0$ of the smooth unitary curves $\matrixEXP{tE} U$, with
an arbitrary antihermitian $E$, which are equal to $EU$. Simultaneously we have convinced ourselves that $\TANGENTspace{\UNITARY}{U}$
is indeed a real linear space.

If curve $\gamma$ is contained in $\FIXEDmoduliMATRICES{U}$, it satisfies $f(\gamma(t)) - f(U)\ =\ \ZEROvect$, hence its derivative at $0$
satisfies $\DIFFERENTIAL{f}{U}(\gamma'(0))\ =\ \ZEROvect$ which means that $\TANGENTspace{\FIXEDmoduliMATRICES{U}}{U}$, formed  by all possible
$\gamma'(0)$ allowed here, is contained in $\NULL \DIFFERENTIAL{f}{U}$. This nullspace, due to our assumption  $\forall i,j\ U_{i,j} \neq 0$,
is the space $\{\Ii R \HADprod U:\ R\ \mbox{real}\}$ of dimension $N^2$, composed of matrices whose entries are perpendicular on the complex plane
to the corresponding entries in $U$. Every element $\Ii R \HADprod U$ is attained as the derivative at $0$ of the smooth curve
$\ELEMENTof{\gamma(t)}{i}{j}\ =\ U_{i,j} \PHASE{R_{i,j} t}$ contained in $\FIXEDmoduliMATRICES{U}$. Again we have shown that
$\TANGENTspace{\FIXEDmoduliMATRICES{U}}{U}$ is a real linear space.

Having both tangent spaces in an explicit form we consider their intersection. This intersection, as well as its dimension, are helpful tools
in the search for submanifolds of $\UNITARY$  containing $U$ and composed of matrices with the pattern of moduli identical with that of $U$.
Note that smooth manifolds of complex Hadamard matrices 
($N \times N$ unitary, with all the moduli of their entries equal, necessarily to $1/\sqrt{N}$) 
stemming from a given complex Hadamard matrix fall into this class.

Analyzing the tangent space $\TANGENTspace{\SET{F}}{U}$  to such a manifold $\SET{F}$ at $U$ 
(as before for $\UNITARY$ and $\FIXEDmoduliMATRICES{U}$)
we see that since $\gamma$, corresponding to any $\gamma'(0) \in \TANGENTspace{\SET{F}}{U}$, is contained both in 
$\UNITARY$ and $\FIXEDmoduliMATRICES{U}$, the tangent vector $\gamma'(0)$ lies in the intersection 
$\TANGENTspace{\UNITARY}{U}\ \          \cap\ \            \TANGENTspace{\FIXEDmoduliMATRICES{U}}{U}$.
This intersection will be distinguished by a separate symbol:
\begin{equation}
	\label{eq_feasible_space}
	\FEASIBLEspace{U}
					\ \ \stackrel{def}{=}\ \ 
	\TANGENTspace{\UNITARY}{U}                  \ \ \cap\ \                   \TANGENTspace{\FIXEDmoduliMATRICES{U}}{U}
\end{equation}
and called \DEFINED{the feasible space of $U$}. Its (real) dimension, denoted by:
\begin{equation}
	\label{eq_undepased_defect}
	\undephasedDEFECT(U)
						\ \ \stackrel{def}{=}\ \ 
	\DIMR \left( \FEASIBLEspace{U} \right)                   \ \ ,
\end{equation}
will be called \DEFINED{the undephased defect of $U$}.
The meaning of the word 'undephased' will be made clear later.

Having computed $\TANGENTspace{\UNITARY}{U}$ and $\TANGENTspace{\FIXEDmoduliMATRICES{U}}{U}$ above
we can give a more detailed description of $\FEASIBLEspace{U}$:
\begin{eqnarray}
	\label{eq_feasible_space_characterization}
	\lefteqn{\FEASIBLEspace{U}\ \ =}
	&   &      \\
	\nonumber
	&  =   &
	\left\{
		\Ii R \HADprod U\ :\ \ 
		R\ \ \mbox{real}\ \ \ \ \ \ \mbox{and}\ \ \ \ \ \ \mbox{$\exists E$ antihermitian s.t. $\Ii R \HADprod U\ =\ EU$}
	\right\}  
	\\
	\nonumber
	&    =   &
	\left\{
		\Ii R \HADprod U\ :\ \
		R\ \ \mbox{real}\ \ \ \ \ \ \mbox{and}\ \ \ \ \ \ \mbox{$(\Ii R \HADprod U) \hermTRANSPOSE{U}$\ antihermitian}
	\right\}
	\\
	\nonumber
	&     =   &
	\Ii  
		        \ \ \cdot\ \ 
	\left\{
		\mbox{$R$ real}\ :\ \ 
		\mbox{$(\Ii R \HADprod U) \hermTRANSPOSE{U}$\ antihermitian}
	\right\}
	              \ \ \HADprod\ \ 
	U                                                         \ \ ,
\end{eqnarray}
where in the last expression an image of the real isomorphism\ \ 
	$(X \rightarrow \Ii X \HADprod U):\ \COMPLEXmatricesNxN{N} \rightarrow \COMPLEXmatricesNxN{N}$
has been used.

The dimension of $\FEASIBLEspace{U}$ thus reads:
\begin{equation}
	\label{eq_old_definition_of_undephased_defect}
	\undephasedDEFECT(U)
		\ \ =\ \ 
	\DIMR\left(
		\left\{
			\mbox{$R$ real}\ :\ \ 
			\mbox{$(\Ii R \HADprod U) \hermTRANSPOSE{U}$\ antihermitian}
		\right\}
	\right)                              \ \ ,
\end{equation}      
where the real system defining the set on the right hand side takes the more explicit form:
\begin{equation}
	\label{eq_iRU_Uhansposed_antihermitian_system_in_matrix_form}
	\left( \Ii R \HADprod U \right)    \hermTRANSPOSE{U}
		\ \ -\ \ 
	U    \hermTRANSPOSE{  \CONJ{\Ii R} \HADprod U }
				\ \ \ =\ \ \
	\ZEROvect                                                                       \ \ ,
\end{equation}
or:
\begin{equation}
	\label{eq_iRU_Uhansposed_antihermitian_system_in_entrywise_form}
	\forall\ i < j  \ \ \ \  
	\sum_{k=1}^N  
		U_{i,k}  \CONJ{U}_{j,k}  \left( R_{i,k}\ -\ R_{j,k} \right)
				\ \ \ =\ \ \ 
	0                                                                       \ \ .
\end{equation}

The dimension $\undephasedDEFECT(U)$ of the solution space of the above system is an upper bound on the
dimension of a manifold $\SET{F}$ of the type described in the second paragraph preceding (\ref{eq_feasible_space}),
because:
\begin{equation}
	\label{eq_undephased_manifold_dim_estimation}
	\dim(\SET{F})
		\ \ =\ \ 
	\DIMR\left(   \TANGENTspace{\SET{F}}{U}  \right)
		\ \ \leq\ \ 
	\undephasedDEFECT(U)
		\ \ =\ \ 
	\DIMR\left(  \TANGENTspace{\UNITARY}{U}\ \cap\ \TANGENTspace{\FIXEDmoduliMATRICES{U}}{U}  \right)
	\ \ ,
\end{equation}
since, as it was said above, $\TANGENTspace{\SET{F}}{U}$ is contained in the intersection on the right.

To introduce another type of defect we need the notion of equivalence. Unitary matrices $U$ and $V$ are \DEFINED{equivalent},
written $U\ \RELofEQUI\ V$, if one can be obtained from the other by permuting the rows and/or columns and 
multiplying them by unimodular numbers. That is to say, 
\begin{equation}
	\label{eq_equivalence}
	U\ \RELofEQUI\ V
						\ \ \ \ \ \ \stackrel{def}{\Longleftrightarrow}\ \ \ \ \ \ 
	\exists\ 
		\begin{array}{c}
			P_r,\ P_c\ \ \ \mbox{permutation matrices}  \\
			D_r,\ D_c\ \ \ \mbox{unitary diagonal matrices}
		\end{array}
		\ \ \mbox{s.t.}\ \ \
			 V\ \ =\ \ P_r\ D_r\ \cdot U\ \cdot D_c\ P_c                    \ \ .
\end{equation}
This is an equivalence relation in the usual sense. The order of permutation and unitary diagonal matrices on both sides of $U$
in the above definition is of course unimportant.

Looking for unitaries with the pattern of moduli identical with that of $U$ we are mostly interested in inequivalent matrices.
In our earlier paper \cite{Defect} we introduced the notion of a dephased manifold (Definition 4.3) which cannot contain an 
'enphased' matrix $D_r V D_c$ if it already contains $V$ and $V \neq D_r V D_c$ (for $D_r$, $D_c$ unitary diagonal).
This definition we solely apply to submanifolds of $\UNITARY$ containing matrices with a fixed pattern of the moduli of entries.
Every member $V$ of a dephased manifold posesses a neighbourhood deprived of matrices equivalent to $V$ (Remark 4.4 in \cite{Defect}).

A class of dephased manifolds which is easy to analyze is the one consisting of manifolds composed of unitaries with 
some entries (not only their moduli) fixed, so that these manifolds have the 'dephased' property. For example, 
if a manifold, containing $U$ (which determines the moduli throughout) with no zero entries, is composed of matrices with
their entries in the $1$-st row and column fixed, then it is a dephased manifold 
(see also Lemma 4.6 in \cite{Defect} and the paragraph that follows).
In this case $2N-1$ entries are fixed, but of course one could choose any other $k$-th row and $l$-th column to fix, instead of
$k=1$ and $l=1$. Such a manifold can be left- and right-multiplied by such unitary diagonal matrices for which the resulting  dephased manifold
has its fixed entries real, that is deprived of their nonzero phase. That is where the term 'dephased' comes from.

The dimension of such a manifold, call it $\SET{F}$ again, cannot exceed the value $\undephasedDEFECT(U)\ -\ (2N-1)$, and
we have proved that in Theorem 4.7 in \cite{Defect}. Below we provide a less general, but simpler explanation of this bound.
Additionally, in Theorem 4.2 in \cite{Defect} we have proved that in the case when   $\undephasedDEFECT(U)\ -\ (2N-1)\ \ =\ \ 0$
matrix $U$ is isolated, i.e. the only unitaries with the same pattern of moduli as in $U$ belonging to a small neighbourhood  of $U$
are those of the form $D_r U D_c$, where $D_r,\ D_c$ are unitary diagonal. This means that $\SET{F}$ with $k$-th row and $l$-th column
fixed is a single point, $\SET{F}\ =\ \{U\}$, in such a situation. 
	(The proof goes along this line: $\undephasedDEFECT(U)\ =\ (2N-1)$ allows us to choose an $N^2 - (2N-1)$ equation subsystem $h(R)\ =\ \ZEROvect$
from the unitarity conditions on $\left[ U_{i,j} \PHASE{R_{i,j}} \right]_{i,j = 1,\ldots, N}$, where $R = \left[ R_{i,j} \right]_{i,j = 1,\ldots, N}$ is the 
real variable, with the full rank $\RANK\left(  \DIFFERENTIAL{h}{\ZEROvect} \right) \ =\ N^2 - (2N-1)$ of the differential at $\ZEROvect$.
Thus $h(R)\ =\ \ZEROvect$, and hence the unitarity conditions, define a $(2N-1)$ dimensional manifold around $\ZEROvect$, which must be equal
to $\left\{ \left[ \alpha_i  +  \beta_j \right]_{i,j = 1,\ldots, N}\ :\ \ \alpha_i, \beta_j \in \REALS \right\}$.  
So $\left[ U_{i,j} \PHASE{R_{i,j}} \right]_{i,j = 1,\ldots, N}$ close to $U$ is  unitary iff  it is equal to some
$\left[ U_{i,j} \PHASE{\left(\alpha_i +\beta_j\right)} \right]_{i,j = 1,\ldots, N}
\ \ =\ \ 
\diag\left(\PHASE{\alpha_1},\ldots,\PHASE{\alpha_N}\right)    
	\cdot    
U    
	\cdot    
\diag\left(\PHASE{\beta_1},\ldots,\PHASE{\beta_N}\right)$.)

The aforementioned bound\ \  $\dim(\SET{F})\  \leq\ \undephasedDEFECT(U)\ -\ (2N-1)$\ \  is implied by the fact  that $\TANGENTspace{\SET{F}}{U}$,
contained in $\FEASIBLEspace{U}$ (see the paragraph preceding (\ref{eq_feasible_space})),
cannot contain any nonzero element of the $(2N-1)$ dimensional real space spanned by the $2N$ dependent directions:
\begin{eqnarray} 
	\label{eq_enphasing_directions}
	\Ii \ \diag\left( \STbasis{i} \right) U
		\ \ =\ \ 
	\Ii \left( \STbasis{i} \ONESvect^T \right)  \HADprod U
			\ \ \ \ \in\ \ 
	 \TANGENTspace{\UNITARY}{U}\ \cap\ \TANGENTspace{\FIXEDmoduliMATRICES{U}}{U}
	\ \ \ \ \ \ i\ =\ 1,2,\ldots,N
	\ \ ,
	&    &     \\
	\nonumber
	U \cdot \Ii\  \diag\left( \STbasis{j} \right)
		\ \ =\ \ 
	\Ii \left(  \ONESvect  \STbasis{j}^T \right)   \HADprod U
	\ \ \ \ \in\ \ 
	 \TANGENTspace{\UNITARY}{U}\ \cap\ \TANGENTspace{\FIXEDmoduliMATRICES{U}}{U}
	\ \ \ \ \ \ j\ =\ 1,2,\ldots,N
	\ \ ,
	&     &    \\
	\nonumber
	\mbox{where}\ \ \ \STbasis{k}\ \ \mbox{is the $k$-th vertical vector of the standard basis},
		\ \ \ \ 
	\ONESvect\ \ =\ \ [1,1,1,\ldots,1]^T                         \ \ ,
	&    &  
\end{eqnarray}
because any nonzero element of this space has a nonzero entry in the $l$-th column or $k$-th row
(this can be seen by considering the basis (\ref{eq_enphasing_directions}) with $j=l$ excluded), 
which violates the constraint on $\SET{F}$ (consequently on $\TANGENTspace{\SET{F}}{U}$)
that its members  have the $k$-row and $l$-th column fixed (filled with zeros).

It can be shown that the $(2N-1)$ dimensional real space spanned by vectors (\ref{eq_enphasing_directions}) is the space
$\TANGENTspace{\ENPHASEDmatrices{U}}{U}$ tangent to the \DEFINED{enphasing manifold}:
\begin{equation}
	\label{eq_enphasing_manifold}
	\ENPHASEDmatrices{U}
					\ \ \stackrel{def}{=}\ \ 
	\left\{
		D_r U D_c\ :\ \ D_r,\ D_c\ \ \ \mbox{unitary diagonal}
	\right\}
\end{equation}
composed of matrices equivalent to $U$. For a proof see the end of this section.
Obviously $\ENPHASEDmatrices{U}$, as well as the above $\SET{F}$, are contained in 
$\UNITARY \cap \FIXEDmoduliMATRICES{U}$, 
so both $\TANGENTspace{\SET{F}}{U}$ and $\TANGENTspace{\ENPHASEDmatrices{U}}{U}$ are contained in 
$\FEASIBLEspace{U}\ =\ \TANGENTspace{\UNITARY}{U}\ \cap\ \TANGENTspace{\FIXEDmoduliMATRICES{U}}{U}$.
We do not however have to consider  
the space spanned by   (\ref{eq_enphasing_directions}) as a tangent space, 
the fact that it is contained in $\FEASIBLEspace{U}$ can be deduced directly from the form of the spanning vectors  (\ref{eq_enphasing_directions}).

As it has been said in the paragraph preceding the last one,  
	$\TANGENTspace{\SET{F}}{U}$ and $\TANGENTspace{\ENPHASEDmatrices{U}}{U}$
are independent spaces i.e.  
	$\TANGENTspace{\SET{F}}{U}\ \cap\ \TANGENTspace{\ENPHASEDmatrices{U}}{U}\ =\ \{ \ZEROvect \}$. 
Therefore their algebraic sum is a direct sum and 
\begin{eqnarray}
	\label{eq_dephased_manifold_tangent_space_dim_estimation}
	\lefteqn{ 
		\TANGENTspace{\SET{F}}{U}\ \oplus\ \TANGENTspace{\ENPHASEDmatrices{U}}{U}
		\ \ \subset\ \ \FEASIBLEspace{U}
	}
	&     &         \\
	\nonumber
	&      \Longrightarrow      &
	\DIMR\left( 
		\TANGENTspace{\SET{F}}{U}\ \oplus\ \TANGENTspace{\ENPHASEDmatrices{U}}{U}
	\right)
			\ \ =\ \ 
	\DIMR\left(      \TANGENTspace{\SET{F}}{U}      \right)
			\ \ +\ \ 
	\DIMR\left(           \TANGENTspace{\ENPHASEDmatrices{U}}{U}        \right)
			\ \ \leq\ \ 
	\undephasedDEFECT(U)\ =\ \DIMR\left(  \FEASIBLEspace{U} \right)
	\\
	\nonumber
	&    \Longrightarrow     &
	\DIMR\left(     \TANGENTspace{\SET{F}}{U}       \right)
			\ \ \leq\ \ 
	\undephasedDEFECT(U)\ -\ (2N-1)                                     \ \ ,
\end{eqnarray}
which completes our argument for the bound $\dim(\SET{F})\ \leq\ \undephasedDEFECT(U)\ -\ (2N-1)$.

The last quantity in (\ref{eq_dephased_manifold_tangent_space_dim_estimation}) will be distinguished
by a separate symbol (where $U$ is unitary $N \times N$):
\begin{equation}
	\label{eq_dephased_defect}
	\DEFECT(U)      \ \ =\ \       \undephasedDEFECT(U)\ \ -\ \ (2N - 1)
\end{equation}
and called \DEFINED{the dephased defect of $U$}.

  There exists a $\DEFECT(U)$ dimensional subspace of $\FEASIBLEspace{U}$, 
$\FEASIBLEspaceDEPHASED{U}$, 
which contains $\TANGENTspace{\SET{F}}{U}$  ( $\TANGENTspace{\SET{F}}{U}$ is a space composed of vectors-matrices
having zeros in the $k$-th row and $l$-th column). $\FEASIBLEspaceDEPHASED{U}$ is defined by additional independent complex
conditions on the elements of $\FEASIBLEspace{U}$:
\begin{equation}
	\label{eq_dephasing_conditions}
	X_{1,l} = 0,\ 
	\ldots,\ 
	X_{N,l} = 0,\ 
	X_{k,1} = 0,\ 
	\ldots,\ 
	X_{k,l-1} = 0,\ 
	X_{k,l+1} = 0,\ 
	\ldots,\ X_{k,N} = 0                          \ \ .
\end{equation}  
Since the elements of $\FEASIBLEspace{U}$ have the form $\Ii R \HADprod U$, where $R$ is a real matrix, condition
$X_{a,b}\ =\ 0$ can be replaced by $\Re(X_{a,b})\ =\ 0$ if $\Im(U_{a,b})\ \neq\ 0$, or by $\Im(X_{a,b})\ =\ 0$ in the other case.
Thus any of the complex conditions (\ref{eq_dephasing_conditions}) is equivalent to an appropriate real condition, imposing which
on $\FEASIBLEspace{U}$ means intersecting it with a $(2N^2 - 1)$ dimensional real subspace of $\COMPLEXmatricesNxN{N}$.
Intersection with such a subspace cannot reduce the dimension of the intersected space by more than $1$. Subsequent 
intersection of $\FEASIBLEspace{U}$ with the subspaces corresponding to conditions (\ref{eq_dephasing_conditions}), in the
order given there, deprives $\FEASIBLEspace{U}$ of vectors:
\begin{equation}
	\label{eq_vectors_remomoved_from_D_U}
	\Ii \STbasis{1} \ONESvect^T \HADprod U,\ 
	\ldots,\ 
	\Ii \STbasis{N} \ONESvect^T \HADprod U,\ \ 
	\Ii \ONESvect \STbasis{1}^T \HADprod U,\ 
	\ldots,\ 
	\Ii \ONESvect \STbasis{l-1}^T \HADprod U,\ \ 
	\Ii \ONESvect \STbasis{l+1}^T \HADprod U,\ 
	\ldots,\ 
	\Ii \ONESvect \STbasis{N}^T \HADprod U                         \ \ ,
\end{equation}
one at a time. Thus $\FEASIBLEspace{U}$ faces reduction of its dimension at each step, 
to become a $\undephasedDEFECT(U)\ -\ (2N-1)$ dimensional space $\FEASIBLEspaceDEPHASED{U}$ after $(2N-1)$ steps.
The obtained space $\FEASIBLEspaceDEPHASED{U}$ cannot contain any nonzero element of 
$\TANGENTspace{\ENPHASEDmatrices{U}}{U}$, that is 
	$\FEASIBLEspaceDEPHASED{U}\ \cap\ \TANGENTspace{\ENPHASEDmatrices{U}}{U}\ =\ \{\ZEROvect\}$, so
\begin{equation}
	\label{eq_decomposition_of_feasible_space}
	\FEASIBLEspace{U}\ \ =\ \ \FEASIBLEspaceDEPHASED{U}\ \ \oplus\ \ \TANGENTspace{\ENPHASEDmatrices{U}}{U}
	\ \ .
\end{equation}

The term 'dephased' in the context of complex Hadamard matrices, or unitary matrices with prescribed moduli, means that
a matrix is left- and right-multiplied by unitary diagonal matrices so that its entries in the $k$-th row and $l$-th column
are real nonnegative (deprived of a nonzero phase). 
A manifold with a fixed pattern of moduli $\SET{F} \ni U $ built from such matrices has members with this row and column fixed, so it is a manifold
of the type discussed above. We could consider another manifold
\begin{equation}
	\label{eq_Fprim_as_enphased_F}
	\SET{F'}      \ \ =\ \             \{D_r V D_c:\  V \in \SET{F},\ \ D_r,\ D_c\ \mbox{unitary diagonal} \}
\end{equation}
which can be smoothly parametrized around $U$ by 
\begin{equation}
	\label{eq_Fprim_parametrization}
	\diag\left(\PHASE{\alpha_1},\ldots,\PHASE{\alpha_N}\right)    
		\cdot    
	\SET{F}\left(\gamma_1,\ldots,\gamma_K\right)    
		\cdot    
	\diag\left(1,\PHASE{\beta_2},\ldots,\PHASE{\beta_N}\right)              \ \ ,
\end{equation}
where $\SET{F}\left(\gamma_1,\ldots,\gamma_K\right)$ parametrizes $\SET{F}$, with independent parameters, around $U$.
The parameters in  (\ref{eq_Fprim_parametrization}) are independent
because the derivatives of (\ref{eq_Fprim_parametrization}) wrt $\alpha_{...}, \beta_{...}$ belonging to $\TANGENTspace{\ENPHASEDmatrices{U}}{U}$
and those wrt $\gamma_{...}$ belonging to $\TANGENTspace{\SET{F}}{U}$ are all independent since
	$\TANGENTspace{\SET{F}}{U}\ \cap\ \TANGENTspace{\ENPHASEDmatrices{U}}{U}\ =\ \{ \ZEROvect \}$
(see the comment preceding (\ref{eq_dephased_manifold_tangent_space_dim_estimation}); 
derivatives are taken at point $\alpha_{...}, \beta_{...},\gamma_{...}$ corresponding, through (\ref{eq_Fprim_parametrization}),  to $U$).
Manifold $\SET{F}$ can be obtained from manifold $\SET{F}'$ by dephasing  its elements,
that is left- and right-multiplying them by appropriate unitary diagonal matrices, the same for every element of $\SET{F}'$.
Originally by a 'dephased manifold' we had understood a manifold obtained like $\SET{F}$ from $\SET{F}'$, but later
we decided to define it using its properties, like in the paragraph following the one containing (\ref{eq_equivalence}).
Now, while the dimension of $\SET{F}'$, yet 'undephased', 
is bound by the undephased defect $\undephasedDEFECT(U)$ because $\TANGENTspace{\SET{F}'}{U} \subset \FEASIBLEspace{U}$
(see the paragraph preceding (\ref{eq_feasible_space})), 
the dimension of the dephased $\SET{F}$  is bound by the dephased defect $\DEFECT(U)$, 
because  $\TANGENTspace{\SET{F}}{U} \oplus \TANGENTspace{\ENPHASEDmatrices{U}}{U} \subset \FEASIBLEspace{U}$  
(see (\ref{eq_dephased_manifold_tangent_space_dim_estimation}) and before).

Let us summarize. In this section we have introduced two quantities for a unitary matrix $U$, 
the undephased defect $\undephasedDEFECT(U)$ and the dephased defect $\DEFECT(U)$.
The first one is more convenient to define and it depends on the Kronecker product structure of $U$ in a simple way
(see Corollary \ref{cor_V1_spaces_Kron_multiplied_are_contained_in_V1_space} {\bf c)}\ \ and\ \ 
Lemma \ref{lem_multiplicativity_of_D_of_F}).
The second one is used in applications, as a bound for the dimension of manifolds composed
of (locally) inequivalent unitary matrices with prescribed moduli.

Both these quantities can be defined for a rescaled unitary matrix (still with no zero entries).
Let $\alpha\ \in\ \COMPLEX \setminus \{0\}$,\ \ 
$U \in \UNITARY\ \ s.t.\ \ \forall i,j\ U_{i,j} \neq 0$,\ \ 
and let $V\ =\ \alpha U\  \in\ \alpha \UNITARY\ =\ \{W:\ g(W)\ =\ |\alpha|^2 I \}$,\ 
where $g$ defined in (\ref{eq_g_Gram_matrix_map}). Then, similarly as before,
$\TANGENTspace{\alpha \UNITARY}{V}\ =\ \{EV:\ E\ \mbox{antihermitian} \}$ and we define:
\begin{eqnarray}
	\label{eq_feasible_space_for_rescaled_unitary}
	\lefteqn{
		\FEASIBLEspace{V}              
				\ \ \stackrel{def}{=}\ \  
		\TANGENTspace{\alpha \UNITARY}{V}\ \ \cap\ \ \TANGENTspace{\FIXEDmoduliMATRICES{V}}{V}
	}
	&    &      \\
	\nonumber
	&      =       &
	\left\{\ 
		\Ii R \HADprod V\ :\ \ 
		R\ \mbox{real}     
				\ \ \ \mbox{and}\ \ \       
		\left( \Ii R \HADprod V \right) \hermTRANSPOSE{V}\ \mbox{antihermitian}
	\right\}
	\\
	\nonumber
	&       =        &
	\alpha \cdot
	\left\{\ 
		\Ii R \HADprod U\ :\ \ 
		R\ \mbox{real}     
				\ \ \ \mbox{and}\ \ \       
		\left( \Ii R \HADprod U \right) \hermTRANSPOSE{U}\ \mbox{antihermitian}
	\right\}
	\\
	\nonumber
	&       =       &
	\alpha \cdot  \FEASIBLEspace{U}                    \ \ ,
\end{eqnarray}
so all the calculations leading to the undephased defect can be done with $U$ representing $V$, and vice versa:
\begin{equation}
	\label{eq_undephased_defect_for_rescaled_unitary}
	\undephasedDEFECT(V)    
			\ \ \stackrel{def}{=}\ \    
	\DIMR\left( \FEASIBLEspace{V} \right)
			\ \ =\ \ 
	\DIMR\left( \FEASIBLEspace{U} \right)
			\ \ =\ \ 
	\undephasedDEFECT(U)               \ \ ,
\end{equation}

The dephased defect is defined in a similar way as for the representing unitary: 
\begin{equation}
	\label{eq_dephased_defect_for_rescaled_unitary}
	\DEFECT(V)      \ \ \stackrel{def}{=}\ \      \undephasedDEFECT(V)\ \ -\ \ (2N - 1)               \ \ ,
\end{equation}
where $N$ is the size of $V$.
It is the dimension of space $\FEASIBLEspaceDEPHASED{V}$ obtained from $\FEASIBLEspace{V}$ by imposing
conditions (\ref{eq_dephasing_conditions}) on it, and of course (see (\ref{eq_decomposition_of_feasible_space})):
 \begin{equation}
	\label{eq_decomposition_of_rescaled_feasible_space}
	\FEASIBLEspace{V}\ =\ \FEASIBLEspaceDEPHASED{V}\ \oplus\ \TANGENTspace{\ENPHASEDmatrices{V}}{V}
\end{equation}
where	$\TANGENTspace{\ENPHASEDmatrices{V}}{V}\ =\ \alpha \TANGENTspace{\ENPHASEDmatrices{U}}{U}$
is spanned by (\ref{eq_enphasing_directions}) with $U$ replaced by $V$, and where 
$\FEASIBLEspaceDEPHASED{V}\ =\ \alpha \FEASIBLEspaceDEPHASED{U}$.

$\undephasedDEFECT(V)$ and $\DEFECT(V)$ play the role of bounds on the dimensions of manifolds containing $V$
and contained in $\alpha \UNITARY\ \cap \FIXEDmoduliMATRICES{V}$, analogous to the role of 
$\undephasedDEFECT(U)$ and $\DEFECT(U)$.

We will meet rescaled unitary matrices later in this paper, in sections \ref{sec_Fourier_matrices} and \ref{sec_fourier_defect}.
These will be Fourier matrices, often conveniently defined to have unimodular entries.

We end this section with a proof of the statement which begins the paragraph containing the definition 
(\ref{eq_enphasing_manifold}) of $\ENPHASEDmatrices{U}$. Namely, that the tangent space
$\TANGENTspace{\ENPHASEDmatrices{U}}{U}$ is spanned by vectors (\ref{eq_enphasing_directions}), 
i.e. it is equal to $\{\Ii (a \ONESvect^T + \ONESvect b^T) \HADprod U \ :\ \ a,b \in \REALcolumnVECTORS{N} \}$,
where $\REALcolumnVECTORS{N}$ represents the set of real $N \times 1$ matrices and 
$\ONESvect\ =\ [1\ 1\ \ldots\ 1]^T$.

Let $\gamma \subset \ENPHASEDmatrices{U}$ be a differentiable curve such that $\gamma(0)\ =\ U$. For each $t$ from its domain
$\gamma(t)\ =\ \diag(\PHASE{\alpha_1(t)},\ldots,\PHASE{\alpha_N(t)}) 
       \cdot      U      \cdot        \diag(\PHASE{\beta_1(t)},\ldots,\PHASE{\beta_N(t)})$, 
that is $\gamma_{i,j}(t)\ =\ U_{i,j} \PHASE{(\alpha_i(t) + \beta_j(t))}$, where we can assume that $\beta_1(t)\ =\ 0$.
Because $\gamma(t) \stackrel{t \rightarrow 0}{\longrightarrow} U$, in a sufficiently small neighbourhood of $0$ the phase
$\alpha_i(t)$ in $\gamma_{i,1}(t)\ =\ U_{i,1} \PHASE{\alpha_i(t)}$ can be chosen to be close to $0$. Having done that for 
all $i$ we look at the phase $(\alpha_1(t)\ +\ \beta_j(t))$ in $\gamma_{1,j}(t)\ =\ U_{1,j} \PHASE{(\alpha_1(t) + \beta_j(t))}$,
choosing $\beta_j(t)$ to be close to zero. Thanks to all these choices we have guaranteed that 
$(\alpha_i(t) + \beta_j(t))     \stackrel{t \rightarrow 0}{\longrightarrow}    0$, as 
$\gamma_{i,j}(t)      \stackrel{t \rightarrow 0}{\longrightarrow}    U_{i,j}$.

The tangent vector $\gamma'(0)$, the existence of which we assume, is a matrix with entries:
\begin{equation}
	\label{eq_gamma_prim_entry}
	\ELEMENTof{ \gamma'(0) }{i}{j}
			\ \ =\ \ 
	\lim_{t \rightarrow 0}
		U_{i,j}    \frac{1}{t-0}    \left( \PHASE{(\alpha_i(t) + \beta_j(t))}\ -\ 1 \right)
			\ \ =\ \ 
	\lim_{k \rightarrow \infty}
		U_{i,j}    \frac{1}{t_k-0}    \left( \PHASE{(\alpha_i(t_k) + \beta_j(t_k))}\ -\ 1 \right)              \ \ ,
\end{equation}
where $(t_k)_{k=1,2,\ldots}$ is an arbitrary sequence converging to $0$ in the domain of $\gamma$, containing no zeros.
By successive taking subsequences from it we can construct such a subsequence $(s_l)_{l=1,2,\ldots}$ of $(t_k)_{k=1,2,\ldots}$
that for every $i,j$ there holds one of the two:
\begin{itemize}
	\item  
		$\forall l\ \ (\alpha_i(s_l) + \beta_j(s_l))\ \neq\ 0$, in this case $\gamma'_{i,j}(0)$ takes the form
		\begin{equation}
			\label{eq_gamma_prim_entry_split_limes}
			\ELEMENTof{ \gamma'(0) }{i}{j}
					\ \ =\ \ 
			\lim_{l \rightarrow \infty}
				U_{i,j}   \frac{   \PHASE{(\alpha_i(s_l) + \beta_j(s_l))}\ -\ 1   }
							    {   (\alpha_i(s_l) + \beta_j(s_l))\ -\ 0                 }
				\cdot
				\frac{   (\alpha_i(s_l) + \beta_j(s_l))\ -\ 0  }
					  {   s_l\ -\ 0 }                                                               \ \ ,
		\end{equation}
		so the right factor converges because the left factor converges to $\Ii U_{i,j}$. We can then write
		\begin{equation}
			\label{eq_gamma_prim_entry_simple_limes}
			\ELEMENTof{ \gamma'(0) }{i}{j}
					\ \ =\ \ 
			\lim_{l \rightarrow \infty}
				\Ii U_{i,j} \cdot \frac{   (\alpha_i(s_l) + \beta_j(s_l))   }{  s_l  }             \ \ .
		\end{equation}
		
	\item 
		$\forall l\ \ (\alpha_i(s_l) + \beta_j(s_l))  \ =\ 0$, then there must be $\gamma'_{i,j}(0)\ =\ 0$, so
		in this case (\ref{eq_gamma_prim_entry_simple_limes}) also holds.

\end{itemize}
Matrix
\begin{eqnarray}
	\label{eq_gamma_prim_approximating_point}
	\lefteqn{
		\gamma^{(l)}
				\ \ =\ \ 
		\left[      \Ii    U_{i,j}    \frac{  (\alpha_i(s_l) + \beta_j(s_l))  }{  s_l  }     \right]_{i,j=1,\ldots,N}
				\ \ =
	}   &   &   \\
	\nonumber
	&   & 
	\Ii  
			\cdot
	\frac{1}{s_l}
	\left(
		\left[    
			\begin{array}{c}  
				\alpha_1(s_l)    \\   
				\vdots   \\   
				\alpha_N(s_l)
			\end{array}
		\right]
			\cdot    
		\ONESvect^T
				\ \ +\ \ 
		\ONESvect    
			\cdot    
		\left[    
			\begin{array}{ccc} 
				\beta_1(s_l)  &  \cdots  &  \beta_N(s_l)
			\end{array}
		\right]
	\right)
			\HADprod
	U
\end{eqnarray}
belongs to the $(2N-1)$ dimensional real space 
	$\{\Ii (a \ONESvect^T + \ONESvect b^T) \HADprod U \ :\ \ a,b \in \REALcolumnVECTORS{N} \}$.
Because this space is a closed set, 
	$\gamma'(0)\    \stackrel{(\ref{eq_gamma_prim_entry_simple_limes})}{=}\ \lim_{l \rightarrow \infty} \gamma^{(l)}$ 
also belongs to this space.
Every element $\Ii \left( a \ONESvect^T \ +\ \ONESvect b^T \right)  \HADprod U$ of this space, 
for\ \  $a\ = \left[ a_1, \ldots, a_N \right]^T$\ \ and\ \  $b\ = \left[ b_1, \ldots, b_N \right]^T$,\ \ 
is attained as a tangent vector in $\TANGENTspace{\ENPHASEDmatrices{U}}{U}$, 
for example as the derivative at $0$ of the curve 
	$\diag\left( \PHASE{a_1 t}, \ldots, \PHASE{a_N t} \right)      \cdot
	 U   \cdot
	\diag\left( \PHASE{b_1 t}, \ldots, \PHASE{b_N t} \right) \ \subset\ \ENPHASEDmatrices{U}$.

%
%
\section{Karabegov's theory}
	\label{sec_Karabegovs_theory}
 In this section we give an account of the theory presented by Alexander Karabegov in either of his articles \cite{Karabegov}, 
\cite{Karabegov_old}. This theory leads to an interesting characterization of the undephased defect of $U$ as the 
multiplicity of the eigenvalue $1$ of a certain operator on $\COMPLEXmatricesNxN{N}$ associated with $U$. The theory heavily depends 
on the assumption that a unitary $U$ (of size $N$) has no zero entries. Our presentation is different from that of Karabegov, 
but more convenient for us to the point that we have decided to prove all the facts in our own style. The Karabegov's results are extended
with consideration of equivalent unitary matrices and invariance of the undephased defect under equivalence operations.
The special case of self-equivalence (symmetry) is treated in a way slightly similar to the one used in the older preprint \cite{Karabegov_old},
but a connection to the newer presentation in \cite{Karabegov} is also made. We end with an application of the theory to 
Kronecker products of unitary matrices, which should be treated as a supplement to our previous work on the defect of these \cite{KroneckerDefect}.
Surely our exposition misses the quantum mechanical  context which motivated Karabegov, as can be seen in his papers.

\subsection{Berezin transform}
	\label{subsec_Berezin_transform}

Let us return to the linear system (\ref{eq_iRU_Uhansposed_antihermitian_system_in_matrix_form}), whose solution space parametrizes
space  $\FEASIBLEspace{U}$. A unitary matrix $U$ is assumed to have no zero entries. The matrix equality in (\ref{eq_iRU_Uhansposed_antihermitian_system_in_matrix_form}) can be rewritten as:
\begin{equation}
	\label{eq_Cu_Du_difference_system}
	\Cu{U}( \Ii R )\     -    \ \Du{U}( \Ii R )       \ \ =\ \  \ZEROvect                \ \ ,
\end{equation}
where $\Cu{U},\ \Du{U}:\ \COMPLEXmatricesNxN{N} \longrightarrow \COMPLEXmatricesNxN{N}$ are linear isomorphisms
introduced by Karabegov, given by the formulas:
\begin{eqnarray}
	\label{eq_Cu_definition}
	\Cu{U}(F)\ \ \stackrel{def}{=}\ \ (F \HADprod U) U^{*}\ \ ,
	\ \ \ \ \ \mbox{i.e.}\ \ 
	\ELEMENTof{\Cu{U}(F)}{i}{j}
			\ \ =\ \ 
	\sum_{k=1}^{N} 
		U_{i,k} \CONJ{U}_{j,k} F_{i,k}
	&   &          \\
	\label{eq_Du_definition}
	\Du{U}(F)\ \ \stackrel{def}{=}\ \ U (\CONJ{F} \HADprod U)^{*}\ \ ,
	\ \ \ \ \ \mbox{i.e.}\ \ 
	\ELEMENTof{\Du{U}(F)}{i}{j}
			\ \ =\ \ 
	\sum_{k=1}^{N} 
		U_{i,k} \CONJ{U}_{j,k} F_{j,k}       \ \ .
	&   &
\end{eqnarray}
(\ref{eq_Cu_Du_difference_system}) is further equivalent to
\begin{equation}
	\label{eq_invCu_Du_eigenvector_for_1_system}
	\Ii R          \ \ =\ \        \Cu{U}^{-1}  \Du{U}  (\Ii R)                   \ \ ,
\end{equation}
which is nothing but a statement that an imaginary matrix $\Ii R$ is an eigenvector of $\Cu{U}^{-1} \Du{U}$ associated 
with eigenvalue $1$, if it exists. Since the matrices listed in (\ref{eq_enphasing_directions}) all belong to $\FEASIBLEspace{U}$,
their parametrizing matrices $\STbasis{i} \ONESvect^T,\ \ONESvect \STbasis{j}^T,\ i,j=1,\ldots,N$ 
satisfy system (\ref{eq_iRU_Uhansposed_antihermitian_system_in_matrix_form}), hence multiplied by $\Ii$ they are imaginary
eigenvectors of $\Cu{U}^{-1} \Du{U}$ satisfying (\ref{eq_invCu_Du_eigenvector_for_1_system}). We can thus be sure that
$\Vu{U}{1}$, the complex eigenspace of $\Cu{U}^{-1} \Du{U}$ associated with eigenvalue 1, is never a zero space.

$\Cu{U}^{-1} \Du{U}$ is a complex operator and one could ask whether the dimension of its complex eigenspace
$\Vu{U}{1}$ determines the value of the undephased defect $\undephasedDEFECT(U)$, the dimension of the real space $\FEASIBLEspace{U}$.  
Through our journey to a positive answer to this question let us follow Karabegov, 
inspecting properties of $\Cu{U}$, $\Du{U}$ and $\Cu{U}^{-1} \Du{U}$. From now on, the latter operator will be denoted by:
\begin{equation}
	\label{eq_Iu_definition}
	\Iu{U}\ \ \stackrel{def}{=}\ \ \Cu{U}^{-1} \Du{U}
\end{equation}
and called, after Karabegov, the Berezin transform.

Next let us introduce a notation for the standard inner product and a weighted inner product on $\COMPLEXmatricesNxN{N}$:
\begin{eqnarray}
	\label{eq_standard_inner_product}
	\STinnerPRODUCT{F}{G}
	&		\ \ \stackrel{def}{=}\ \        &
	\sum_{i,j} 
		F_{i,j} \CONJ{G}_{i,j}
			\ \ =\ \ 
	\trace G^* F\ \ ,
	\\
	\label{eq_U_inner_product}
	\UinnerPRODUCT{U}{F}{G}
	&		\ \ \stackrel{def}{=}\ \      &
	\sum_{i,j} 
		\ABSOLUTEvalue{U_{i,j}}^2 F_{i,j} \CONJ{G}_{i,j}
			\ \ =\ \ 
	\trace (G \HADprod U)^* (F \HADprod U)
			\ \ =\ \
	\STinnerPRODUCT{F \HADprod U}{G \HADprod U}                 \ \ .
\end{eqnarray}

Most of the below properties can be found in \cite{Karabegov} or \cite{Karabegov_old}.


\begin{lemma}
	\label{lem_properties_of_Cu_Du_Iu}
	Let $U$ be a unitary $N \times N$ matrix with no zero entries. Then:
	\begin{description}
		\item[a)]
			$\UinnerPRODUCT{U}{F}{G}
						\ \ =\ \ 
				\STinnerPRODUCT{\Cu{U}(F)}{\Cu{U}(G)}
						\ \ =\ \ 
				\STinnerPRODUCT{\Du{U}(F)}{\Du{U}(G)}$
		\item[b)]
			$\hermTRANSPOSE{\Cu{U}(F)}\ \ =\ \ \Du{U}\left(\CONJ{F}\right)$,
			\ \ \ \ \ \ \ \ equivalently\ \ \ \ \ 
			$\hermTRANSPOSE{\Du{U}(F)}\ \ =\ \ \Cu{U}\left(\CONJ{F}\right)$.
		\item[c)]
			\begin{description}
				\item    $\Cu{U}(F)$ is antihermitian \ \ \ \ if and only if \ \ \ \ $F\ \ =\ \ - \Iu{U}\left( \CONJ{F} \right)$.
				\item    $\Cu{U}(F)$ is hermitian\ \ \ \  if and only if \ \ \ \ $F\ \ =\ \  \Iu{U}\left( \CONJ{F} \right)$.
			\end{description}
		\item[d)]
			\begin{description}
				\item     $\Du{U}(F)$ is antihermitian\ \ \  \ if and only if \ \  \ \ $\CONJ{F}\ \ =\ \ - \Iu{U}( F )$.
				\item     $\Du{U}(F)$ is hermitian\ \ \ \  if and only if\ \ \ \  $\CONJ{F}\ \ =\ \  \Iu{U}( F )$.
			\end{description}
		\item[e)]
			If\ \ \  $F\ \ =\ \ \pm\ \Iu{U}(G)$\ \ \ then
			\begin{quote}
				\begin{description}
					\item        $\Cu{U}(F)\ =\ \pm\ \Du{U}(G)$ is antihermitian\ \ \ \  if and only if\ \ \ \  $F\ =\ \mp\ \CONJ{G}$.
					\item        $\Cu{U}(F)\ =\ \pm\ \Du{U}(G)$ is hermitian \ \ \ \ if and only if\ \ \ \  $F\ =\ \pm\ \CONJ{G}$.
				\end{description}
			\end{quote}
		\item[f)]
			$F\ \ =\ \ \Iu{U}(G)$\ \ \ \ if and only if\ \ \ \  $\CONJ{G}\ \ =\ \ \Iu{U}\left( \CONJ{F} \right)$.
		\item[g)]
			Denoting by $A^{\ominus 1}$ the entrywise (Hadamard) inverse of matrix $A$,
			\begin{eqnarray}
				\label{eq_Iu_formula}
				\Iu{U}(F)
				&		\ \ =\ \      &
				U^{\ominus 1}  
					\HADprod  
				\left(      U      \left(    \CONJ{F} \HADprod U    \right)^{*}    U  \right)              \ \ ,
				  \\
				\label{eq_Iu_ij_entry_formula}
				\ELEMENTof{\Iu{U}(F)}{i}{j}
				&		\ \ =\ \     &
				\frac{1}{U_{i,j}} 
				\sum_{k,l}
					U_{i,l}  U_{k,j} \CONJ{U}_{k,l}  F_{k,l}
						\ \ =\ \ 
				\sum_{k,l}
					\frac{U_{i,l} U_{k,j}}{U_{i,j} U_{k,l}}   \ABSOLUTEvalue{U_{k,l}}^2    F_{k,l}            \ \ ,
			\end{eqnarray}
			so the trace of $\Iu{U}$ is equal to $N$.
		\item[h)]
			Denoting by $\DIAGwith{a}$ a diagonal matrix in $\COMPLEXmatricesNxN{N}$ whose consecutive diagonal 
			entries are the consecutive entries of a column matrix $a \in \columnVECTORS{N}$, 
			\begin{description}
				\item     $\Cu{U}\left(  a b^T  \right)     \ \ =\ \      \DIAGwith{a}  U  \DIAGwith{b}   U^*$            \ \ ,
				\item     $\Du{U}\left(  a b^T  \right)     \ \ =\ \       U  \DIAGwith{b}   U^*   \DIAGwith{a} $         \ \ .
			\end{description}
	\end{description}
\end{lemma}

\PROOFstart 
\begin{description}
	\item[a)]
		\begin{description}
			\item
				$\UinnerPRODUCT{U}{F}{G}
					\ \stackrel{(\ref{eq_U_inner_product})}{=}\ 
				\trace (G \HADprod U)^* (F \HADprod U)
					\ =\ 
				\trace\left( (G \HADprod U) U^* \right)^*     \cdot     (F \HADprod U) U^*
					\ =\ 
				\STinnerPRODUCT{\Cu{U}(F)}{\Cu{U}(G)}$\ ,  
			\item
				$\!\!\!\!\!\UinnerPRODUCT{U}{F}{G}
					\ \stackrel{(\ref{eq_U_inner_product})}{=}\
				\UinnerPRODUCT{U}{\CONJ{G}}{\CONJ{F}}
					\ \stackrel{(\ref{eq_U_inner_product})}{=}\
				\trace \left( \CONJ{F} \HADprod U \right)^* \left( \CONJ{G} \HADprod U \right)
					\ =\ 
				\trace \left(  U \left(  \CONJ{G} \HADprod U \right)^*  \right)^*    \cdot     U (\CONJ{F} \HADprod U)^*
					\ =\ 
				\STinnerPRODUCT{\Du{U}(F)}{\Du{U}(G)}$\ \ .
		\end{description}
	\item[b)]
		$\hermTRANSPOSE{ \Cu{U}(F) } 
			\ =\ 
		\hermTRANSPOSE{ (U \HADprod F) U^* }
			\ =\ 
		U \left( U \HADprod \CONJ{\CONJ{F}} \right)^*
			\ =\ 
		\Du{U}\left( \CONJ{F} \right)$\ \ .
	\item[c)]
		$\hermTRANSPOSE{ \Cu{U}(F) }
			\ =\ 
		\mp\ \Cu{U}(F)
				\ \ \stackrel{\mathbf{b)}}{\Longleftrightarrow}\ \ 
		\Du{U}\left( \CONJ{F} \right)
			\ =\ 
		\mp\ \Cu{U}(F)
				\ \ \Longleftrightarrow\ \ 
		\Cu{U}^{-1} \Du{U} \left( \CONJ{F} \right)
			\ =\ 
		\mp\ F$   \ \ .
	\item[d)]
		$\hermTRANSPOSE{ \Du{U}(F) }
			\ =\ 
		\mp\ \Du{U}(F)
				\ \ \stackrel{\mathbf{b)}}{\Longleftrightarrow}\ \ 
		\Cu{U}\left( \CONJ{F} \right)
			\ =\ 
		\mp\ \Du{U}(F)
				\ \ \Longleftrightarrow\ \ 
		\CONJ{F}   
			\ =\ 
		\mp\ \Cu{U}^{-1} \Du{U} ( F )$            \ \ .
	\item[e)]
		If\ \ $F\ \ =\ \ \pm\ \Cu{U}^{-1} \Du{U} ( G )$\ \ then\ \ $\Cu{U}(F)\ \ =\ \ \pm\ \Du{U}(G)$\ \ and \\
		$\hermTRANSPOSE{ \Cu{U}(F) }
			\ =\ 
		\mp\ \Cu{U}(F)
				\ \ \stackrel{\mathbf{b)}}{\Longleftrightarrow}\ \ 
		\Du{U}\left( \CONJ{F} \right)
			\ =\ 
		\mp\ \Cu{U}(F)
				\ \ \stackrel{\Cu{U}(F)\ \ =\ \ \pm\ \Du{U}(G)}{\Longleftrightarrow}\ \ 
		\Du{U}\left( \CONJ{F} \right)
			\ =\ 
		\mp\ \left( \pm\ \Du{U}(G) \right)
				\ \ \Longleftrightarrow\ \ 
		\CONJ{F}
			\ =\ 
		\mp\ \left( \pm\ G \right)
				\ \ \Longleftrightarrow\ \ 
		F
			\ =\ 
		\mp\ \left( \pm\ \CONJ{G} \right)$       \ \ .
	\item[f)]
		$F
			\ =\ 
		\Cu{U}^{-1} \Du{U} ( G )
				\Longleftrightarrow
		\Cu{U}\left( \CONJ{\CONJ{F}} \right)
			\ =\ 
		\Du{U}\left( \CONJ{\CONJ{G}} \right)
				\ \ \stackrel{\mathbf{b)}}{\Longleftrightarrow}\ \ 
		\hermTRANSPOSE{ \Du{U}\left( \CONJ{F} \right) }
			\ =\ 
		\hermTRANSPOSE{ \Cu{U}\left( \CONJ{G} \right) }
				\ \ \Longleftrightarrow\ \ 
		\CONJ{G}
			\ =\ 
		\Cu{U}^{-1} \Du{U} \left( \CONJ{F} \right)$         \ \ .
	\item[g)]
		If\ \ $H\ \ =\ \ \Cu{U}(G)\ \ =\ \ (G \HADprod U)U^*$\ \ then\ \ $G\ \ =\ \ U^{\ominus 1} \HADprod (HU)\ \ =\ \ \Cu{U}^{-1}(H)$, so
		to get $\Cu{U}^{-1} \Du{U} (F)$ we have to replace $H$ with $\Du{U} (F)$. 

		The $i,j$th entry of the right hand side of (\ref{eq_Iu_formula}) is equal to:\\
		$\frac{1}{U_{i,j}}\   U_{i,:} \cdot \left( F^T \HADprod U^* \right) \cdot U_{:,j}
			\ =\ 
		\frac{1}{U_{i,j}}\ \sum_{l,k}    U_{i,l} U_{k,j} \ELEMENTof{F^T \HADprod \CONJ{U}^T}{l}{k}
			\ =\ 
		\frac{1}{U_{i,j}}\ \sum_{k,l}    U_{i,l} U_{k,j} \CONJ{U}_{k,l}  F_{k,l}$,\\
		where  $U_{i,:}$  and  $U_{:,j}$ stand for the $i$th row and $j$th column of $U$, respectively.

		From (\ref{eq_Iu_ij_entry_formula}) we can extract the matrix of $\Iu{U}$ in a standard basis whatever ordered,   
		indexed by the elements of $\{1,\ldots,N\}\times\{1,\ldots,N\}$.  Let $\ELEMENTof{\Iu{U}}{(i,j)}{(k,l)}$ denote an element of this
		matrix, the $(i,j)$th coordinate of the image of the $(k,l)$th basis vector, with $1$ at position $k,l$ and zeros elsewhere, under $\Iu{U}$. 
		This element is equal to
		the number which multiplies $F_{k,l}$ in formula (\ref{eq_Iu_ij_entry_formula}). The trace is then equal to :
		\begin{equation}
			\label{eq_Iu_trace_calculation}
			\trace\left( \Iu{U} \right)
					\ \ =\ \ 
			\sum_{(i,j)}    \ELEMENTof{\Iu{U}}{(i,j)}{(i,j)}
					\ \ =\ \ 
			\sum_{(i,j)}    \ABSOLUTEvalue{U_{i,j}}^2
					\ \ =\ \ 
			N
		\end{equation}

	\item[h)]
		$\Cu{U}(ab^T)
			\ =\ 
		\left( ab^T \HADprod U \right)U^*
			\ =\ 
		\left( \DIAGwith{a} U \DIAGwith{b} \right)  U^*$    
		\ \ ,\\
		$\Du{U}(ab^T)
			\ =\ 
		U \left( \CONJ{ab^T} \HADprod U \right)^*
			\ =\ 
		U  \left( \CONJ{\DIAGwith{a}}  U  \CONJ{\DIAGwith{b}} \right)^*
			\ =\ 
		U  \left( \DIAGwith{b}  U^*  \DIAGwith{a} \right)$       \ \ .

\end{description}
\PROOFend  

These properties lead to the most important for us:


\begin{theorem}
	\label{theor_unitarity_and_eigenspaces_of_Iu}
Let U be a unitary $N \times N$ matrix with no zero entries. Then:
\begin{description}
	\item[a)]
		Operator $\Iu{U}:\ \COMPLEXmatricesNxN{N} \longrightarrow \COMPLEXmatricesNxN{N}$ is unitary with respect to the 
		weighted inner product $\UinnerPRODUCT{U}{}{}$, equivalently it has unimodular eigenvalues and an eigenbasis (of
		$\COMPLEXmatricesNxN{N}$) orthonormal wrt $\UinnerPRODUCT{U}{}{}$\ .
	\item[b)]
		If $V \in \Vu{U}{\lambda}$, the eigenspace of $\Iu{U}$ associated with eigenvalue $\lambda$, then also 
		$\CONJ{V},\ \Re(V),\ \Ii \cdot \Im(V) \in \Vu{U}{\lambda}$. As a consequence $\Vu{U}{\lambda}$, treated as a real space,
		decomposes into the real direct sum:
		\begin{equation}
			\label{eq_Vu_direct_sum_decomposition}
			\Vu{U}{\lambda}\ \ =\ \ \REspaceOF{\Vu{U}{\lambda}}\ \ \oplus\ \ \IMspaceOF{\Vu{U}{\lambda}}         \ \ ,
		\end{equation}
		where
		\begin{eqnarray}
			\label{eq_Vu_Re_definition}
			\REspaceOF{\Vu{U}{\lambda}}
			&   \stackrel{def}{=}   &
			\left\{\ 
				V \in \Vu{U}{\lambda}\ :\ \ V\ \ \mbox{real}\ 
			\right\}     
				\ \ =\ \ 
			\Re\left(  \Vu{U}{\lambda}  \right)  
			\ \ ,            \\
			\label{eq_Vu_Im_definition}
			\IMspaceOF{\Vu{U}{\lambda}}
			&   \stackrel{def}{=}   &
			\left\{\ 
				V \in \Vu{U}{\lambda}\ :\ \ V\ \ \mbox{imaginary}\
			\right\}
				 = 
			\Ii \cdot \Im\left( \Vu{U}{\lambda} \right)    
				 =
			\Ii \REspaceOF{\Vu{U}{\lambda}}    ,
		\end{eqnarray}
		so there holds
		\begin{equation}
			\label{eq_Vu_dimension}
			\DIMR\left(  \Vu{U}{\lambda}  \right)
				\ \ =\ \ 
			2 \DIMC\left(  \Vu{U}{\lambda}  \right)
					\ \ \ \ \ =\ \ \ \ \ 
			\DIMR\left(\     \REspaceOF{\Vu{U}{\lambda}}  \ \right)
				\ \ +\ \ 
			\DIMR\left(\     \IMspaceOF{\Vu{U}{\lambda}}  \ \right)                  \ \ ,
		\end{equation}
		where
		\begin{equation}
			\label{eq_Vu_and_Vu_Re_and_Vu_Im_dimension}
			\DIMC\left(  \Vu{U}{\lambda}  \right)
				\ \ =\ \ 
			\DIMR\left(\     \REspaceOF{\Vu{U}{\lambda}}  \ \right)
				\ \ =\ \ 
			\DIMR\left(\     \IMspaceOF{\Vu{U}{\lambda}}  \ \right)                \ \ .
		\end{equation}
		
	\item[c)]
		The $2N-1$ dimensional space 
		\begin{eqnarray}
			\label{eq_trivial_space}
			\TRIVIALspace{U}
			&
					\ \ \stackrel{def}{=}\ \ 
			&
			\left\{a \ONESvect^T + \ONESvect b^T\ :\ \ a,b \in \columnVECTORS{N},\ \   \ONESvect\ =\ [1\ 1\ \ldots\ 1]^T \ \right\}      \\
			\nonumber
			&
					\ \ =\ \ 
			&
			\left\{ \left[ a_i + b_j \right]_{i,j=1,\ldots,N}\ :\ \ a_i, b_j \in \COMPLEX \right\}
		\end{eqnarray}
		is contained in $\Vu{U}{1}$, hence the multiplicity of $1$ as the eigenvalue of $\Iu{U}$ is at least $2N-1$. 

		As in the case of $\Vu{U}{1}$,\ \ the real space $\TRIVIALspace{U}$ decomposes\ \ 
			$\TRIVIALspace{U}\ \ =\ \ \REspaceOF{ \TRIVIALspace{U} }      \ \ \oplus\ \    \IMspaceOF{ \TRIVIALspace{U} }$,
		where 
		\begin{eqnarray}
			\label{eq_Tu_Re_definition}
			\REspaceOF{ \TRIVIALspace{U} }			
			&   \stackrel{def}{=}   &
			\left\{\ T \in \TRIVIALspace{U}\ :\ \ T\ \ \mbox{real}\ \right\}
				\ \ =\ \ 
			\Re\left(  \TRIVIALspace{U}  \right)              \ \ ,             \\
			\label{eq_Tu_Im_definition}
			\IMspaceOF{ \TRIVIALspace{U} }
			&  \stackrel{def}{=}   &
			\left\{\ T \in \TRIVIALspace{U}\ :\ \ T\  \mbox{imaginary}\ \right\}
			\  =\  
			\Ii \cdot \Im\left(  \TRIVIALspace{U}  \right)      
			\  =\ 
			\Ii   \REspaceOF{ \TRIVIALspace{U} }        
		\end{eqnarray}
	
	\item[d)]
		\begin{eqnarray}
			\label{eq_V1_Im_characterization}
		\IMspaceOF{ \Vu{U}{1} }
		&   =   &
		\left\{\ V \in \Vu{U}{1}\ :\ \ V\ \mbox{imaginary}\ \right\}          \\
		\nonumber
		&    =    &
		\left\{\ V \in \Vu{U}{1}\ :\ \ \Cu{U}(V)\ \ \mbox{antihermitian}\ \right\}          \\
		\nonumber
		&   =    &
		\left\{\ V\ :\ \ V\ \ \mbox{imaginary}\ \ \ \ \ \mbox{and}\ \ \ \ \ \Cu{U}(V)\ \ \mbox{antihermitian}\ \right\}             \ \ ,
		\end{eqnarray}
		where $\Cu{U}$ can be replaced by $\Du{U}$.

	\item[e)]
		Space $\FEASIBLEspace{U}$ 
										(see (\ref{eq_feasible_space}),(\ref{eq_feasible_space_characterization}))
		and space  $\TANGENTspace{\ENPHASEDmatrices{U}}{U}\ \subset\ \FEASIBLEspace{U}$
																		(see the paragraph containing (\ref{eq_enphasing_manifold}))
		are isomorphically parametrized by $\IMspaceOF{\Vu{U}{1}}$ and $\IMspaceOF{\TRIVIALspace{U}} \subset \IMspaceOF{\Vu{U}{1}}$, respectively:
		\begin{eqnarray}
			\label{eq_feasible_space_parametrization}
			\FEASIBLEspace{U}
			&    =    &
			\IMspaceOF{\Vu{U}{1}}  \HADprod  U           \ \ \ =\ \ \      \Cu{U}\left(  \IMspaceOF{\Vu{U}{1}}  \right) \cdot U       \ \ ,    \\
			\label{eq_enphasing_manifold_tangent_space_parametrization}
			\TANGENTspace{\ENPHASEDmatrices{U}}{U}
			&    =    &
			\IMspaceOF{\TRIVIALspace{U}}  \HADprod  U           \ \ \ =\ \ \       \Cu{U}\left(   \IMspaceOF{\TRIVIALspace{U}}  \right)  \cdot U    \ \ ,
		\end{eqnarray}
		hence the dimension of $\FEASIBLEspace{U}$ (the undephased defect of $U$) equals:
		\begin{equation}
			\label{eq_undephased_defect_and_V1_dimension_equality}
			\undephasedDEFECT(U)\ \ =\ \ \DIMR\left( \IMspaceOF{\Vu{U}{1}} \right)  \ \ =\ \ \DIMC\left(  \Vu{U}{1} \right)       \ \ .
		\end{equation}

\end{description}
\end{theorem}

\PROOFstart  
\begin{description}
	\item[a)]
		$\Iu{U}$ preserves the weighted inner product $\UinnerPRODUCT{U}{}{}$, as can be seen from: \\
		$\UinnerPRODUCT{U}{F}{G}
				\  \stackrel{L. \ref{lem_properties_of_Cu_Du_Iu} \mathbf{a)}}{=}\ 
		\STinnerPRODUCT{\Du{U}(F)}{\Du{U}(G)}
				 \ =\  
		\STinnerPRODUCT{\Cu{U} \Cu{U}^{-1} \Du{U}(F)}{\Cu{U} \Cu{U}^{-1} \Du{U}(G)}
				\  \stackrel{L. \ref{lem_properties_of_Cu_Du_Iu}\mathbf{a)}}{=}\                       \\
		\UinnerPRODUCT{U}{ \Cu{U}^{-1} \Du{U}(F) }{ \Cu{U}^{-1} \Du{U}(G) }$  
	
	\item[b)]
		Let $F \in \Vu{U}{\lambda}$. Then, using the fact that $\lambda$, an eigenvalue of $\Iu{U}$, is unimodular:\\
		$\lambda F \ =\ \Iu{U}(F)  
				\ \ \stackrel{L. \ref{lem_properties_of_Cu_Du_Iu}\mathbf{f)}}{\Longrightarrow}\ \ 
		\CONJ{F}\ =\ \Iu{U}\left( \CONJ{\lambda} \CONJ{F} \right)
				\ \ \Longrightarrow\ \ 
		\Iu{U}\left(\CONJ{F}\right)\ =\ 1/\CONJ{\lambda} \cdot \CONJ{F} \ \stackrel{\mathbf{a)}}{=}\ \lambda \CONJ{F}$.\\
		Thus $\CONJ{F}$, $\left(F + \CONJ{F}\right)/2 \ =\ \Re(F)$ and $-\Ii\left(F - \CONJ{F}\right)/2\ =\ \Im(F)$ all belong to
		$\Vu{U}{\lambda}$. This justifies the second characterizations of $\REspaceOF{\Vu{U}{\lambda}}$ and $\IMspaceOF{\Vu{U}{\lambda}}$  
		in (\ref{eq_Vu_Re_definition}), (\ref{eq_Vu_Im_definition}).

		$\REspaceOF{\Vu{U}{\lambda}}$ and $\IMspaceOF{\Vu{U}{\lambda}}$ are real subspaces of $\Vu{U}{\lambda}$ treated as a real space.
		$\Vu{U}{\lambda}\ \ =\ \ \REspaceOF{\Vu{U}{\lambda}}\ +\ \IMspaceOF{\Vu{U}{\lambda}}$, an algebraic sum
		(because $\Vu{U}{\lambda}\ \ni\ V\ =\ \Re(V) + \Ii \cdot \Im(V)$), and $\ZEROvect$ can be decomposed into the elements of
		$\REspaceOF{\Vu{U}{\lambda}},\ \IMspaceOF{\Vu{U}{\lambda}}$ only in one way: $\ZEROvect\ =\ \ZEROvect + \ZEROvect$.
		Hence $\Vu{U}{\lambda}\ \ =\ \ \REspaceOF{\Vu{U}{\lambda}}\ \oplus\ \IMspaceOF{\Vu{U}{\lambda}}$.

		The calculation on dimensions is obvious.

	\item[c)]
		We use the notation introduced in Lemma \ref{lem_properties_of_Cu_Du_Iu} {\bf h)} and, 
		in consequence, the identity $\DIAGwith{\ONESvect}\ =\ I$
		to compute the values of $\Cu{U}$, $\Du{U}$ on an element of $\TRIVIALspace{U}$:\\
		$\Cu{U}\left(a \ONESvect^T\ +\ONESvect b^T\right) 
					\ \ \stackrel{L.\ref{lem_properties_of_Cu_Du_Iu} \mathbf{h)}}{=}\ \ 
		\DIAGwith{a}\ +\ U \DIAGwith{b} U^{*} 
					\ \ \stackrel{L.\ref{lem_properties_of_Cu_Du_Iu} \mathbf{h)}}{=}\ \ 
		\Du{U}\left(a \ONESvect^T\ +\ONESvect b^T\right)$,\ \ 
		from which 
		$\left(a \ONESvect^T\ +\ONESvect b^T\right)
				\ \ =\ \ 
		\Iu{U}\left(a \ONESvect^T\ +\ONESvect b^T\right)$. Thus $\TRIVIALspace{U}$ is contained in the eigenspace of $\Iu{U}$ 
		corresponding to eigenvalue $1$.

		$\Re\left(\TRIVIALspace{U}\right)$ is contained in $\REspaceOF{ \TRIVIALspace{U} }$ because
		$\Re\left(a \ONESvect^T + \ONESvect b^T\right)\ =\ \Re(a) \ONESvect^T +\ \ONESvect (\Re(b))^T$.
		Also $\REspaceOF{ \TRIVIALspace{U} }\ \subset\ \Re\left(\TRIVIALspace{U}\right)$ because if
		$a \ONESvect^T + \ONESvect b^T\ \ =\ \ \left[ a_i + b_j \right]_{i,j=1,\ldots,N}$ is real then there must be
		$\Im\left(a_1\right) = -\Im\left(b_1\right) = \ldots = -\Im\left(b_N\right) = \Im\left(a_2\right) = \ldots = \Im\left( a_N \right)\ =\ \phi$,
		that is\ \ \  $a\ \ =\ \ c\ +\ \Ii \phi \cdot \ONESvect$,\ \ \ $b\ \ =\ \ d\ -\  \Ii \phi \cdot \ONESvect$\ \ \ for $c,d,\phi$ real, and\ \
		$a \ONESvect^T + \ONESvect b^T \ =\ c \ONESvect^T + \ONESvect d^T \ =\ \Re\left( c \ONESvect^T + \ONESvect d^T \right)$.

		In a similar way one shows that  $\IMspaceOF{ \TRIVIALspace{U} }\ =\ \Ii \cdot \Im\left( \TRIVIALspace{U} \right)$:
		firstly $\Ii \cdot \Im\left( a \ONESvect^T + \ONESvect b^T \right) \ =\  (\Ii\cdot\Im(a)) \ONESvect^T + \ONESvect (\Ii \cdot \Im(b))^T$,
		secondly $a \ONESvect^T + \ONESvect b^T \in \IMspaceOF{ \TRIVIALspace{U} }$ is equal to some
		$(\Ii c) \ONESvect^T + \ONESvect (\Ii d)^T\ \ =\ \ \Ii \cdot \Im\left( (\Ii c) \ONESvect^T + \ONESvect (\Ii d)^T\right)$, 
		where $c,d$ are those real vectors for which (necessarily as above)\ \ \  $a\ \ =\ \ \phi \cdot \ONESvect \ +\ \Ii \cdot c$,\ \ \  
		$b\ \ =\ \ -\phi \cdot \ONESvect \ +\ \Ii \cdot d$\ \ \  for some real $\phi$.

		$\REspaceOF{ \TRIVIALspace{U} }$ and $\IMspaceOF{ \TRIVIALspace{U} }$ are real subspaces of $\TRIVIALspace{U}$ treated as a real space.
		Because $X \longrightarrow \Ii X$ is both a complex and real isomorphism from $\TRIVIALspace{U}$ onto $\TRIVIALspace{U}$, 
		$\IMspaceOF{ \TRIVIALspace{U} }\ =\ \Ii \cdot \REspaceOF{ \TRIVIALspace{U} }$.

		$\TRIVIALspace{U}\ =\ \REspaceOF{ \TRIVIALspace{U} } + \IMspaceOF{ \TRIVIALspace{U} }$, an algebraic sum (because
		$a \ONESvect^T + \ONESvect b^T  
					\ \ =\ \ 
			\left(\Re(a) \ONESvect^T + \ONESvect (\Re(b))^T \right)
				\ +\ 
			\left((\Ii \cdot \Im(a)) \ONESvect^T + \ONESvect (\Ii \cdot \Im(b))^T \right)$
		),
		and $\ZEROvect$ can be decomposed into the elements of $\REspaceOF{ \TRIVIALspace{U} },\ \IMspaceOF{ \TRIVIALspace{U} }$
		only in one way: $\ZEROvect\ =\ \ZEROvect + \ZEROvect$. 
		Hence $\TRIVIALspace{U}\ =\ \REspaceOF{ \TRIVIALspace{U} } \oplus \IMspaceOF{ \TRIVIALspace{U} }$.

	\item[d)]
		The first characterization in  (\ref{eq_V1_Im_characterization}) was copied from (\ref{eq_Vu_Im_definition}).

		The second characterization is a consequence of Lemma \ref{lem_properties_of_Cu_Du_Iu} {\bf e)}:\\
		$V \in \IMspaceOF{\Vu{U}{1}}
				\  \Longleftrightarrow\  
		\left(   V = -\CONJ{V}\ \mbox{and}\ \ V=\Iu{U}(V)   \right)
				 \ \stackrel{L.\ref{lem_properties_of_Cu_Du_Iu} \mathbf{e)}}{\Longleftrightarrow}\ \\
		\left(   \Cu{U}(V)\ \mbox{antiherm.  and}\ V=\Iu{U}(V)   \right)$.\\
		Obviously $\Cu{U}(V)$ can be replaced by $\Du{U}(V)$ on $\Vu{U}{1}$.

		$\IMspaceOF{\Vu{U}{1}}$, as the intersection of its first and second characterizing set in (\ref{eq_V1_Im_characterization}),
		is surely contained in the third one, also in the case when  $\Cu{U}$ is everywhere replaced with $\Du{U}$. 
		To justify the opposite inclusion we use Lemma \ref{lem_properties_of_Cu_Du_Iu} {\bf b)}:\\ 
		$\left(  V\ \mbox{imag. and}\ \Cu{U}(V) = -\hermTRANSPOSE{\Cu{U}(V)}   \right)
				\ \stackrel{L.\ref{lem_properties_of_Cu_Du_Iu} \mathbf{b)}}{\Longrightarrow}\ 
		\left(  V\ \mbox{imag. and}\ \Cu{U}(V) = -\Du{U}\left(\CONJ{V}\right)   \right)
				\ \Longrightarrow\
		\left(   V\ \mbox{imag. and}\ \Cu{U}(V) = \Du{U}(V)   \right)
				\ \Longrightarrow\ 
		V \in \IMspaceOF{\Vu{U}{1}}$\ \ , \\
		and similarly when $\Cu{U}$ is replaced with $\Du{U}$.

	\item[e)]
		We rewrite the last parametrization in  (\ref{eq_feasible_space_characterization}):
	\begin{eqnarray}
		\label{eq_parametrizing_Du_with_V1_Im}
		\FEASIBLEspace{U}
		&    =    & 
		\left\{
			\Ii R\ :\ \ R\ \mbox{real}\ \ \ \mbox{and}\ \ \ \Cu{U}(\Ii R)\ \mbox{antiherm.}
		\right\}
		\ \ \HADprod\ \ U                      \\
		\nonumber
		&   =   &
		\left\{
			V\ :\ \ V\ \mbox{imag.}\ \ \ \mbox{and}\ \ \ \Cu{U}(V)\ \mbox{antiherm.}
		\right\}
		\ \ \HADprod\ \ U                       \\
		\nonumber
		&    =    &
		\IMspaceOF{\Vu{U}{1}}
		\ \ \HADprod\ \ U                                  \ \ ,
	\end{eqnarray}
	where the last equality holds on account of the third characterization of $\IMspaceOF{\Vu{U}{1}}$ in 
	(\ref{eq_V1_Im_characterization}) at item {\bf d)} of the theorem.

	$\TANGENTspace{\ENPHASEDmatrices{U}}{U}$ is spanned by vectors (\ref{eq_enphasing_directions}), 
	that is it is equal to $\{\Ii (a \ONESvect^T + \ONESvect b^T) \HADprod U \ :\ \ a,b \in \REALcolumnVECTORS{N} \}$
	(we have proved this at the end of section \ref{sec_defect}). So:
	\begin{eqnarray}
		\label{eq_parametrizing_phasing_manifold_tangent_space_with_Tu}
		 \TANGENTspace{\ENPHASEDmatrices{U}}{U}
		&   =   &
		\left\{
			(\Ii a) \ONESvect^T\ +\ \ONESvect (\Ii b)^T \ :\ \ 
			a,b \in \REALcolumnVECTORS{N}
		\right\}
		\ \ \HADprod\ \  U                         \\
		\nonumber
		&   =   &
		\IMspaceOF{ \TRIVIALspace{U} }
		\ \ \HADprod\ \ U                           \ \ ,
	\end{eqnarray}
	because every element of $\IMspaceOF{ \TRIVIALspace{U} }$ can be written as $(\Ii a)\ONESvect^T + \ONESvect (\Ii b)^T$
	for some real column vectors $a,b$, as noticed in the third paragraph of the proof of item {\bf c)} after 'secondly'.

	The second parametrizations in  
			(\ref{eq_feasible_space_parametrization}), 
			(\ref{eq_enphasing_manifold_tangent_space_parametrization})
	are implied by equality\ \ $F \HADprod U \ \ =\ \ \Cu{U}(F) U$.  
		The conclusion (\ref{eq_undephased_defect_and_V1_dimension_equality}) is based on (\ref{eq_Vu_and_Vu_Re_and_Vu_Im_dimension}).

\end{description}
\PROOFend 

Can we easily extract  a basis of $\FEASIBLEspace{U}\ \ =\ \ \IMspaceOF{\Vu{U}{1}} \HADprod U$ from a known basis of $\Vu{U}{1}$ ?

Let $\left(   V_1,\ V_2,\ \ldots,\ V_{\undephasedDEFECT(U)} \right)$ be a basis of $\Vu{U}{1}$.
In general we can only say that 
		$\left( \Ii R_1,\ \Ii R_2,\ \ldots,\ \Ii R_{\undephasedDEFECT(U)} \right)$, 
a basis of the real space $\IMspaceOF{\Vu{U}{1}}$, can always be chosen from \\ 
		$\left(\ \  \   \Ii \cdot \Re\left( V_1 \right),\ \ \        \ldots,\ \ \            \Ii \cdot  \Re\left( V_{\undephasedDEFECT(U)} \right)         ,\ \ \ 
					\Ii \cdot \Im\left( V_1 \right),\ \ \        \ldots,\ \ \            \Ii \cdot \Im\left( V_{\undephasedDEFECT(U)} \right) \ \  \            \right)$,
an arrangement of vectors sitting, by Theorem \ref{theor_unitarity_and_eigenspaces_of_Iu} {\bf b)},   in $\IMspaceOF{\Vu{U}{1}}$.
It is because the following arrangements (think of them as matrices) are of equal complex rank $\undephasedDEFECT(U)$ since
they are obtained one from another by elementary operations and since 
	$V \in \Vu{U}{1}\ \Longrightarrow \CONJ{V} \in \Vu{U}{1}$ by Theorem \ref{theor_unitarity_and_eigenspaces_of_Iu} {\bf b)}:
\begin{eqnarray}
	\label{eq_V1_Im_basis_obtained_from_V1_basis}
	\lefteqn{
		\left(\ \ \       
			V_1,\ \ \
			\ldots,\ \ \  
			V_{\undephasedDEFECT(U)},\ \ \ 
			\CONJ{V}_1,\ \ \ 
			\ldots,\ \ \ 
			\CONJ{V}_{\undephasedDEFECT(U)}\ \ \      
		\right)
	}
	&   &   \\
	\nonumber
	&   \longrightarrow    &
	\left(\ \ \ 
		V_1 + \CONJ{V}_1,\ \ \ 
		\ldots,\ \ \ 
		V_{\undephasedDEFECT(U)}     +    \CONJ{V}_{\undephasedDEFECT(U)},\ \ \ 
		-2 \CONJ{V}_1,\ \ \ 
		\ldots,\ \ \ 
		-2 \CONJ{V}_{\undephasedDEFECT(U)}\ \ \ 
	\right)                                                            \\
	\nonumber
	&    \longrightarrow    &
	\left(\ \ \   
		V_1 + \CONJ{V}_1,\ \ \ 
		\ldots,\ \ \ 
		V_{\undephasedDEFECT(U)}     +    \CONJ{V}_{\undephasedDEFECT(U)},\ \ \ 
		V_1 - \CONJ{V}_1,\ \ \ 
		\ldots,\ \ \ 
		V_{\undephasedDEFECT(U)}     -    \CONJ{V}_{\undephasedDEFECT(U)}\ \ \ 
	\right)                      \\
	\nonumber
	&    \longrightarrow    &
	\left(\ \ \ 
		\Ii \cdot \Re\left( V_1 \right),\ \ \ 
		\ldots,\ \ \ 
		\Ii \cdot \Re\left( V_{\undephasedDEFECT(U)} \right),\ \ \ 
		\Ii \cdot \Im\left( V_1 \right),\ \ \ 
		\ldots,\ \ \ 
		\Ii \cdot \Im\left( V_{\undephasedDEFECT(U)} \right)\ \ \ 
	\right)                                                                                               \ \ ,
\end{eqnarray}
and the last arrangement contains also precisely $\undephasedDEFECT(U)$ vectors independent in the real sense because they are purely
imaginary.

\subsection{Equivalence}
	\label{subsec_equivalence}
	We return to the issue of equivalence of unitary matrices (see (\ref{eq_equivalence})) and describe 
how the operators $\Cu{U}$, $\Du{U}$ and the Berezin transform $\Iu{U}\ =\ \Cu{U}^{-1} \Du{U}$ change when
$U$ is transformed into a matrix equivalent to $U$. To this end we define  new objects:

\begin{definition}
	\label{def_pairs_action_stabilizer}
	\begin{description}
		\item[a)]
			$\enphasedPERMS$ is  the group of all $N \times N$ enphased permutation matrices, i.e. products of permutation
			and unitary diagonal matrices. $\enphasedPERMS \subset \UNITARY$.

		\item[b)]
			$\PERMS$ is the group of all $N \times N$ permutation matrices. $\PERMS \subset \enphasedPERMS$.

		\item[c)]
			$(S,T) X\ \ \stackrel{def}{=}\ \ S X T^{-1}\ \ =\ \ S X T^*$, where $(S,T) \in \PP$ and $X \in \COMPLEXmatricesNxN{N}$, 
			is an action of $\PP$ on $\COMPLEXmatricesNxN{N}$.\\
			(Consider the tensor product on $\columnVECTORS{N}$: $x \overline{\otimes} y \ \ \stackrel{def}{=}\ \ x y^T$.
			Let $A,B \in \COMPLEXmatricesNxN{N}$ and let $B$ be invertible.
			 Then $(Ax) \overline{\otimes} (B^{-T}y)\ \ =\ \ A(x \overline{\otimes} y) B^{-1}\ \ =\ \ (A,B)(x \overline{\otimes} y)$. This means that
			$X \longrightarrow (A,B)X$ is a tensor product of maps $x \rightarrow Ax$ and $y \rightarrow B^{-T}y$ with respect to the tensor product
			$\overline{\otimes}$ both on the side of arguments and results.)

		\item[d)]
			For group $\GROUP{E} \subset \PP$ let 
				\begin{itemize}
					\item
						$\STAB{\GROUP{E}}{X}\ \ \stackrel{def}{=}\ \ \left\{ (S,T) \in \GROUP{E}:\ (S,T)X = X \right\}$,
						the stabilizer of $X$ in group $\GROUP{E}$ 
						($\STAB{\GROUP{E}}{X}$ is a subgroup of $\GROUP{E}$),
					\item
						$\MAP{\GROUP{E}}{X}{Y}\ \ \stackrel{def}{=}\ \ \left\{ (S,T) \in \GROUP{E}:\ (S,T)X = Y \right\}$\\
						($\MAP{\GROUP{E}}{X}{Y}\ =\ (S,T) \STAB{\GROUP{E}}{X} \ =\ \STAB{\GROUP{E}}{Y} (S,T)$
						for any $(S,T) \in \MAP{\GROUP{E}}{X}{Y}$, if the last set is not empty).
					\end{itemize}

		\item[e)]
			For $(S,T) \in \PP$  let:
			\begin{itemize}
				\item
					$\left( \leftPERM{S}{T},\ \rightPERM{S}{T} \right)\ \in \PsPs$\ \  be the pair of 'skeleton' (underlying) permutation matrices
					such that\ \  $S\ =\ \leftDIAG{S}{T} \leftPERM{S}{T}$\ \  
					and\ \           $T\ =\ \rightDIAG{S}{T} \rightPERM{S}{T}$\ \ , where
				\item
					$\left( \leftDIAG{S}{T},\ \rightDIAG{S}{T} \right)$ is the pair of ('left') unitary diagonal matrices uniquely determined by $(S,T)$.
			\end{itemize}

		\item[f)]
			For any $A,B \in \COMPLEXmatricesNxN{N}$, where $B$ is invertible,  $\pairOPERATOR{A}{B}$ denotes the operator
			$ X \longrightarrow AXB^{-1}\ :\ \COMPLEXmatricesNxN{N} \longrightarrow \COMPLEXmatricesNxN{N}$.
			Note that $\pairOPERATOR{A}{B}$ is unitary wrt $\STinnerPRODUCT{}{}$ (defined in (\ref{eq_standard_inner_product})) if $A,\ B$ are unitary.
				
	\end{description}
\end{definition}

The next lemma gives some insight into Karabegov's results concerning his symmetry groups of a unitary $U$, to be presented later.
We will start from an observation. 
Let 
$\Delta \in \FEASIBLEspace{U}\ \ =\ \TANGENTspace{\UNITARY}{U}                  \ \ \cap\ \                   \TANGENTspace{\FIXEDmoduliMATRICES{U}}{U}$
(recall (\ref{eq_feasible_space})), 
that is $\Delta\ =\ \Ii R \HADprod U\ =\ EU$ for some real matrix $R$ and antihermitian $E$. 
$(S,T) \in \PP$ maps $U$ into an equivalent matrix $V\ =\ (S,T)U$.  
Because\  \  $(S,T)\UNITARY\ =\ \UNITARY$,\ \    $(S,T)\FIXEDmoduliMATRICES{U}\ =\ \FIXEDmoduliMATRICES{V}$\ \ 
and\ \ $X \longrightarrow (S,T)X$ is a real linear isomorphic map equal to its differential at $U$, then\ \ 
$(S,T)\TANGENTspace{\UNITARY}{U}\ =\ \TANGENTspace{\UNITARY}{V}$\ \   and \ \ 
$(S,T) \TANGENTspace{\FIXEDmoduliMATRICES{U}}{U}\ =\  \TANGENTspace{\FIXEDmoduliMATRICES{V}}{V}$.\ \ 
This means that 
			$(S,T) \FEASIBLEspace{U} \subset \FEASIBLEspace{V}
					\ \ =\ \ 
			\TANGENTspace{\UNITARY}{V}    \ \ \cap\ \        \TANGENTspace{\FIXEDmoduliMATRICES{V}}{V}$ 
(in fact $(S,T) \FEASIBLEspace{U} = \FEASIBLEspace{V}$ because $(S,T)^{-1} \FEASIBLEspace{V} \subset \FEASIBLEspace{U}$ on the same principle),
 in particular $(S,T)\Delta \in \FEASIBLEspace{V}$,
\ \ that is $(S,T)\Delta\ =\ \Ii Q \HADprod V\ =\ B V$ for some real matrix $Q$ and antihermitian $B$. 

We have used the notion of a map between manifolds (here $X \longrightarrow (S,T)X$) and its tangent map or differential ($X \longrightarrow (S,T)X$ again) between
their tangent spaces. However the fact that\ \  
	$(S,T)\Delta\ \ \in\ \  \TANGENTspace{\UNITARY}{V}    \ \ \cap\ \        \TANGENTspace{\FIXEDmoduliMATRICES{V}}{V}$\ \  
can be  checked directly using the description of both the tangent spaces (in section \ref{sec_defect})  and item {\bf a)}
of the lemma below.


\begin{lemma}
	\label{eq_ST_pair_acting_on_F_o_U_equal_CF_U}
	Let $U$ be a unitary matrix with no zero entries, and let $V\ \ =\ \ (S,T)U$ for some $(S,T) \in \PP$.
	\begin{description}
		\item[a)]
			Let $\Delta \in \COMPLEXmatricesNxN{N}$ be expressed as (which is always possible):
			\begin{equation}
				\label{eq_Delta_forms}
				\Delta\ \ =\ \ F \HADprod U      \ \ \stackrel{equivalently}{=}\ \   \Cu{U}(F) \cdot U                   \ \ .
			\end{equation}
			Then 
			\begin{equation}
				\label{eq_ST_pair_acting_on Delta}
				(S,T) \Delta   
						\ \ =\ \    
				\left( \left( \leftPERM{S}{T}, \rightPERM{S}{T} \right) F \right)     \HADprod     V
						\ \ =\ \ 
				\left( (S,S) \Cu{U}(F) \right) \cdot V
			\end{equation}
	\item[b)]
		For any $F \in \COMPLEXmatricesNxN{N}$ there holds:
		\begin{eqnarray}
			\label{eq_Cu_Cv_relation_by_ST_pair}
			\Cu{V}\left(  \left( \leftPERM{S}{T}, \rightPERM{S}{T} \right) F \right)
			&    =    &
			(S,S) \Cu{U}(F)                       \\
			\label{eq_Du_Dv_relation_by_ST_pair}
			\Du{V}\left(  \left( \leftPERM{S}{T}, \rightPERM{S}{T} \right) F \right)
			&    =    &
			(S,S) \Du{U}(F) 
		\end{eqnarray}

	\end{description}
\end{lemma}

\PROOFstart  
\begin{description}
	\item[a)]
		On the one hand $\Delta$ is transformed 
	(where $(S,T)\ =\ \left( \leftDIAG{S}{T}, \rightDIAG{S}{T} \right)\left( \leftPERM{S}{T}, \rightPERM{S}{T} \right)$\ ,
	and where we use the equality\ \  $\DIAGwith{a} X  \DIAGwith{b}\ \ =\ \ ab^T \HADprod X$\ \  
	for any $a,b \in \columnVECTORS{N}$ and $X \in \COMPLEXmatricesNxN{N}$,\ \ for $\DIAGwith{y}$ see Lemma \ref{lem_properties_of_Cu_Du_Iu} {\bf h)}):
	\begin{eqnarray}
		\label{eq_transforming_Delta}
		(S,T) \Delta
		&   =  & 
		S (F \HADprod U) T^*
					\ \ =\  \
		\leftDIAG{S}{T} \leftPERM{S}{T}  (F \HADprod U) \TRANSPOSE{\rightPERM{S}{T}} \CONJ{\rightDIAG{S}{T}}
					\ \ =\ \           \\
		\nonumber
		&  =  &
		\leftDIAG{S}{T} 
		\left( 
			\left(  \leftPERM{S}{T}       F \TRANSPOSE{\rightPERM{S}{T}} \right)
				\HADprod
			\left( \leftPERM{S}{T}    U    \TRANSPOSE{\rightPERM{S}{T}} \right)
		\right)
		\CONJ{\rightDIAG{S}{T}}                   \ \ =                 \\
		\nonumber
		&  =  &
		\left(  \leftPERM{S}{T}       F \TRANSPOSE{\rightPERM{S}{T}} \right)
			\HADprod
		\left(
			\leftDIAG{S}{T}  \leftPERM{S}{T}  U  \TRANSPOSE{\rightPERM{S}{T}}  \CONJ{\rightDIAG{S}{T}}
		\right)
					\ \ =\ \                                      \\
		\nonumber
		&  =  &
		\left(  \leftPERM{S}{T}       F \TRANSPOSE{\rightPERM{S}{T}} \right)            \HADprod       S U T^*
					\ \ =\ \ 
		\left(  \left( \leftPERM{S}{T}, \rightPERM{S}{T} \right)   F  \right)        \HADprod       V                \ \ .
	\end{eqnarray}
	On the other hand:
	\begin{equation}
		\label{eq_transforming_Delta_differently}
		(S,T) \Delta    \ \ =\ \    S\ \Cu{U}(F) \cdot U\ T^*    \ \ =\ \
		S\ \Cu{U}(F)\ S^*\ \ \cdot\ \ S U T^*     \ \ =\ \   \left( (S,S) \Cu{U}(F) \right)  \cdot V                \ \ .
	\end{equation}
	\item[b)]
		(\ref{eq_Cu_Cv_relation_by_ST_pair}) follows from the equality between the second and third expression in 
		(\ref{eq_ST_pair_acting_on Delta}). $F$ is arbitrary because $\Delta$ was. 
		To (\ref{eq_Du_Dv_relation_by_ST_pair}) we come using Lemma  \ref{lem_properties_of_Cu_Du_Iu}{\bf b)}
		and (\ref{eq_Cu_Cv_relation_by_ST_pair}):
		\begin{eqnarray}
			\label{eq_proving_Du_Dv_relation_by_ST_pair}
			(S,S)  \Du{U}(F)
					&  \stackrel{L.\ref{lem_properties_of_Cu_Du_Iu} \mathbf{b)}}{=}  &
			(S,S) \hermTRANSPOSE{ \Cu{U}\left( \CONJ{F} \right) }
					\  \ =\ \ 
			\hermTRANSPOSE{ (S,S) \Cu{U}\left( \CONJ{F} \right) }        
					\ \ \stackrel{(\ref{eq_Cu_Cv_relation_by_ST_pair})}{=}\ \                \\
			\nonumber
			&  =  &
			\hermTRANSPOSE{ \Cu{V}\left( \CONJ{  \left( \leftPERM{S}{T}, \rightPERM{S}{T} \right) F } \right) }
					\ \ \stackrel{L.\ref{lem_properties_of_Cu_Du_Iu} \mathbf{b)}}{=}\ \ 
			\Du{V}\left(  \left( \leftPERM{S}{T}, \rightPERM{S}{T} \right) F  \right)       \ \ .          
		\end{eqnarray}
		
\end{description}
\PROOFend 

Indeed, if $\Delta\ =\ F \HADprod U\ \ =\ \Cu{U}(F)U\ \ \in\ \ \FEASIBLEspace{U}$,  
that is  $F$ is imaginary and $\Cu{U}(F)$ is antihermitian, then 
	$\left( \leftPERM{S}{T}, \rightPERM{S}{T} \right)F$ is imaginary and $(S,S)\Cu{U}(F)\ =\ S \Cu{U}(F) S^*$ is antihermitian, 
so  $(S,T)\Delta$ given in (\ref{eq_ST_pair_acting_on Delta}) belongs to $\FEASIBLEspace{V}$.

Above (before the last lemma) we wrote that $(S,T) \FEASIBLEspace{U}\ \ =\ \ \FEASIBLEspace{V}$.
$X \longrightarrow (S,T)X$, denoted $\pairOPERATOR{S}{T}$, is thus an isomorphism between 
$\FEASIBLEspace{U}$ and $\FEASIBLEspace{V}$, where $U$ and $V$ are equivalent. In particular, 
for the corresponding dimensions (undephased defects), there holds $\undephasedDEFECT(U)\ \ =\ \ \undephasedDEFECT(V)$, 
if $U$ and $V$ are equivalent. This will also be formally stated in the corollary following the one which comes next.

Lemma \ref{eq_ST_pair_acting_on_F_o_U_equal_CF_U} has a very interesting consequence for the operator $\Iu{U}$.


\begin{corollary}
	\label{cor_Iu_Iv_operators_for_equivalent_U_V_are_similar}
	Let $U$ be a unitary matrix with no zero entries, and let $V\ \ =\ \ (S,T)U$ for some $(S,T) \in \PP$
	(Recall that $(S,T)\ =\ \left(\leftDIAG{S}{T}, \rightDIAG{S}{T} \right)  \left( \leftPERM{S}{T}, \rightPERM{S}{T} \right)$ . )

	\begin{description}
		\item[a)]
			Operators $\Iu{U}$ and $\Iu{V}$ are similar, and $\pairOPERATOR{\leftPERM{S}{T}}{\rightPERM{S}{T}}$ is the
			similarity operator:
			\begin{equation}
				\label{eq_Iu_Iv_operators_for_equivalent_U_V_are_similar}
				\Iu{V}
						\ \ =\ \ 
				\pairOPERATOR{\leftPERM{S}{T}}{\rightPERM{S}{T}}
				\Iu{U}
				\pairOPERATOR{\leftPERM{S}{T}}{\rightPERM{S}{T}}^{-1}
			\end{equation}

		\item[b)]
			Operators $\Ju{U}\ \ \stackrel{def}{=}\ \ \Cu{U} \Du{U}^{-1}$\ \ and $\Ju{V}$ are similar, and 
			$\pairOPERATOR{S}{S}$ is the similarity operator:
			\begin{equation}
				\label{eq_Ju_Jv_operators_for_equivalent_U_V_are_similar}
				\Ju{V}
						\ \ =\ \ 
				\pairOPERATOR{S}{S}     
				\Ju{U}
				\pairOPERATOR{S}{S}^{-1}
			\end{equation}
	\end{description}
\end{corollary}

\PROOFstart 
(\ref{eq_Cu_Cv_relation_by_ST_pair}) and (\ref{eq_Du_Dv_relation_by_ST_pair}) of 
Lemma \ref{eq_ST_pair_acting_on_F_o_U_equal_CF_U} {\bf b)}
can be rewritten as equalities between operators:
\begin{equation}
	\label{eq_Cu_Cv_and_Du_Dv_relation_by_ST_pair_operator_equalities}
	\Cu{V}      \pairOPERATOR{\leftPERM{S}{T}}{\rightPERM{S}{T}}
			\ \ =\ \ 
	\pairOPERATOR{S}{S}      \Cu{U}
					\ \ ,\ \ \ \ \ \ 
	\Du{V}      \pairOPERATOR{\leftPERM{S}{T}}{\rightPERM{S}{T}}
			\ \ =\ \ 
	\pairOPERATOR{S}{S}      \Du{U}                        \ \ .
\end{equation}
Hence 
\begin{equation}
	\label{eq_Cv_Dv_in_relation_to_Cu_Du_operator_equalities}
	\Cu{V}
			\ \ =\ \ 
	\pairOPERATOR{S}{S}      \Cu{U}     \pairOPERATOR{\leftPERM{S}{T}}{\rightPERM{S}{T}}^{-1}
					\ \ ,\ \ \ \ \ \ 
	\Du{V}
			\ \ =\ \ 
	\pairOPERATOR{S}{S}      \Du{U}     \pairOPERATOR{\leftPERM{S}{T}}{\rightPERM{S}{T}}^{-1}                    \ \ .
\end{equation}
where of course all the operators involved are invertible. The compositions of the above expressed $\Cu{V}$, $\Du{V}$ with 
their inverses:\ \ \ $\Cu{V}^{-1} \Du{V}$\ \ \  and\ \ \  $\Cu{V} \Du{V}^{-1}$\ \ \ lead to 
(\ref{eq_Iu_Iv_operators_for_equivalent_U_V_are_similar})  and (\ref{eq_Ju_Jv_operators_for_equivalent_U_V_are_similar}),
respectively.
\PROOFend   

Note that Corollary \ref{cor_Iu_Iv_operators_for_equivalent_U_V_are_similar} provides a necessary condition for
$U$ and $V$ being equivalent, namely that the  operators $\Iu{U}$ and $\Iu{V}$ (or $\Ju{U}$ and $\Ju{V}$)
must be similar in such a situation. In particular their characteristic polynomials must be the same.

Besides Corollary \ref{cor_Iu_Iv_operators_for_equivalent_U_V_are_similar} {\bf a)} justifies the earlier mentioned fact that the
undephased defect is invariant under equivalence operations on $U$, and other related facts:

\begin{corollary}
	\label{cor_DU_equal_DV_for_U_V_equivalent}
	Let $U$ and $V$ be as in Corollary \ref{cor_Iu_Iv_operators_for_equivalent_U_V_are_similar}, that is let them be equivalent.
	Let $(S,T) \in \MAP{\PP}{U}{V}$.

	Then $\undephasedDEFECT(U)\ \ =\ \ \undephasedDEFECT(V)$, i.e. the multiplicities of $1$ in the spectra of
	$\Iu{U}$ and $\Iu{V}$ are equal. Moreover, the associated eigenspaces are mapped one into the other:
	\begin{eqnarray}
		\label{eq_eigenspaces_for_1_for_U_and_V_are related_by similarity_operator}
		\Vu{V}{1}   
		&   =   &
		\pairOPERATOR{\leftPERM{S}{T}}{\rightPERM{S}{T}}\left(    \Vu{U}{1}    \right)      \ \ ,      \\
		\label{eq_Im_eigenspaces_for_1_for_U_and_V_are related_by similarity_operator}
		\IMspaceOF{\Vu{V}{1}} 
		&   =    &
		\pairOPERATOR{\leftPERM{S}{T}}{\rightPERM{S}{T}}\left(    \IMspaceOF{\Vu{U}{1}}    \right)     
	\end{eqnarray}
	(because the similarity operator $\pairOPERATOR{\leftPERM{S}{T}}{\rightPERM{S}{T}}$ maps real matrices into real ones),
	where the last equality corresponds to the transformation 
		$\FEASIBLEspace{V}\ \ =\ \ \pairOPERATOR{S}{T} \left( \FEASIBLEspace{U} \right)$
	by     Theorem \ref{theor_unitarity_and_eigenspaces_of_Iu} {\bf e)},
	        Lemma \ref{eq_ST_pair_acting_on_F_o_U_equal_CF_U} {\bf a)}
	and  (\ref{eq_Im_eigenspaces_for_1_for_U_and_V_are related_by similarity_operator}):
	\begin{eqnarray}
		\label{eq_transforming_Du_by_ST_pair}
		\lefteqn{\pairOPERATOR{S}{T} \left( \FEASIBLEspace{U} \right)  
				\ \ \stackrel{Th.\ref{theor_unitarity_and_eigenspaces_of_Iu} {\mathbf e)}}{=}\ \ }
		&  &  \\
		\nonumber
		&  &
		\pairOPERATOR{S}{T} \left( \IMspaceOF{\Vu{U}{1}} \HADprod U \right)
				\ \ \stackrel{Lem. \ref{eq_ST_pair_acting_on_F_o_U_equal_CF_U} {\mathbf a)}}{=}\ \ 
		\left( \pairOPERATOR{\leftPERM{S}{T}}{\rightPERM{S}{T}} \left(\IMspaceOF{\Vu{U}{1}} \right)    \right)     \HADprod     V
				\ \ \stackrel{(\ref{eq_Im_eigenspaces_for_1_for_U_and_V_are related_by similarity_operator})}{=}\ \ 
		\IMspaceOF{\Vu{V}{1}}   \HADprod    V   
				\ \ \stackrel{Th.\ref{theor_unitarity_and_eigenspaces_of_Iu} {\mathbf e)}}{=}\ \ 
		\FEASIBLEspace{V}             .
	\end{eqnarray}
\end{corollary}

Karabegov in his paper considers only the situation when $V\ =\ U$ and now will analyze this case in our own style.

First let us observe that for $(S,T) \in \STAB{\PP}{U}$ we have  that 
	$\pairOPERATOR{S}{T}\left( \FEASIBLEspace{U} \right)\ =\  \FEASIBLEspace{U}$, that is the stabilizer preserves
the feasible space, which is a consequence of (\ref{eq_transforming_Du_by_ST_pair}).
In this way the feasible space can be looked for  among the invariant real subspaces of the stabilizer,
treated as a collection of real isomorphisms, in which the methods of the representation theory could be applied.
We also want to say here that $\FEASIBLEspace{U}$ is to some degree determined by $\STAB{\PP}{U}$
if the latter is not trivial. Further,


\begin{corollary}
	\label{eq_Iu_Ju_commute_with_ops_related_to_ST_from_the_stabilizer_of_U}
	Let $U$ be a unitary matrix with no zero entries. For any $(S,T) \in \STAB{\PP}{U}$, that is $(S,T) \in \PP$ such that
	$(S,T)U\ =\ U$, there hold the following commutation relations:   
	\begin{description}
		\item[a)]
			\begin{equation}
				\label{eq_Iu_commutes_with_perm_op_related_to_ST_from_the_stabilizer_of_U}
				\Iu{U}  \pairOPERATOR{\leftPERM{S}{T}}{\rightPERM{S}{T}}
						\ \ =\ \ 
				\pairOPERATOR{\leftPERM{S}{T}}{\rightPERM{S}{T}}    \Iu{U}              \ \ ,
			\end{equation}
			where both the operators $\Iu{U}$ and $\pairOPERATOR{\leftPERM{S}{T}}{\rightPERM{S}{T}}$
			are unitary wrt $\UinnerPRODUCT{U}{}{}$ (defined in (\ref{eq_U_inner_product})).

		\item[b)]
			\begin{equation}
				\label{eq_Ju_commutes_with_SS_op_related_to_ST_from_the_stabilizer_of_U}
				\Ju{U}   \pairOPERATOR{S}{S}           \ \ =\ \      \pairOPERATOR{S}{S}   \Ju{U}           \ \ ,
			\end{equation}
			where both the operators $\Ju{U}$ and $\pairOPERATOR{S}{S}$ are unitary wrt $\STinnerPRODUCT{}{}$
			(where $\Ju{U}$ defined in Lemma \ref{cor_Iu_Iv_operators_for_equivalent_U_V_are_similar} {\bf b)},
			$\STinnerPRODUCT{}{}$ defined in (\ref{eq_standard_inner_product})).
	\end{description}
\end{corollary}

\PROOFstart 
Equality (\ref{eq_Iu_commutes_with_perm_op_related_to_ST_from_the_stabilizer_of_U}) at item {\bf a)}
is equivalent to 
	(\ref{eq_Iu_Iv_operators_for_equivalent_U_V_are_similar})
	of Corollary \ref{cor_Iu_Iv_operators_for_equivalent_U_V_are_similar} {\bf a)}
for $V\ \ =\ \ U\ \ =\ \ (S,T)U$.
In this special case also (\ref{eq_Ju_commutes_with_SS_op_related_to_ST_from_the_stabilizer_of_U})
at item {\bf b)} is equivalent to 
	(\ref{eq_Ju_Jv_operators_for_equivalent_U_V_are_similar}) 
	of Corollary \ref{cor_Iu_Iv_operators_for_equivalent_U_V_are_similar} {\bf b)}.

The unitarity of $\Iu{U}$ wrt to $\UinnerPRODUCT{U}{}{}$ is stated in Theorem \ref{theor_unitarity_and_eigenspaces_of_Iu} {\bf a)},
while the unitarity of $\pairOPERATOR{\leftPERM{S}{T}}{\rightPERM{S}{T}}$ has to be checked. To this end the unitarity of 
$\pairOPERATOR{S}{T}$ wrt $\STinnerPRODUCT{}{}$ and Lemma \ref{eq_ST_pair_acting_on_F_o_U_equal_CF_U} {\bf a)} will be used.
\begin{eqnarray}
	\label{eq_perm_op_related_to_ST_weighted_unitarity}
	\UinnerPRODUCT{U}{F}{G}
	&  =  &
	\STinnerPRODUCT{F \HADprod U}{G \HADprod U}
			\ \ =\ \ 
	\STinnerPRODUCT{(S,T)(F \HADprod U)}{(S,T)(G \HADprod U)}
			\ \ \stackrel{L.\ref{eq_ST_pair_acting_on_F_o_U_equal_CF_U} \mathbf{a)}}{=}                  \\
	\nonumber
	&   =   &
	\STinnerPRODUCT{    \left( \left(\leftPERM{S}{T},\rightPERM{S}{T}\right) F \right)   \HADprod  U    }
	                           {    \left( \left(\leftPERM{S}{T},\rightPERM{S}{T}\right) G \right)   \HADprod  U    }
			\ \ =\ \ 
	\UinnerPRODUCT{U}{   \left(\leftPERM{S}{T},\rightPERM{S}{T}\right) F   }
	                           {   \left(\leftPERM{S}{T},\rightPERM{S}{T}\right) G   }                 \ \ .
\end{eqnarray} 

The unitarity of $\Ju{U}$ wrt $\STinnerPRODUCT{}{}$ follows from Lemma \ref{lem_properties_of_Cu_Du_Iu} {\bf a)}:
\begin{equation}
	\label{eq_standard_unitarity_of_Ju}
	\STinnerPRODUCT{ \Cu{U} \Du{U}^{-1}(F) }{ \Cu{U} \Du{U}^{-1}(G) }
			\ \ \stackrel{L.\ref{lem_properties_of_Cu_Du_Iu} \mathbf{a)}}{=}\ \ 
	\UinnerPRODUCT{U}{  \Du{U}^{-1}(F) }{  \Du{U}^{-1}(G) }
			\ \ \stackrel{L.\ref{lem_properties_of_Cu_Du_Iu} \mathbf{a)}}{=}\ \ 
	\STinnerPRODUCT{F}{G}                 \ \ ,
\end{equation}
and the unitarity of $\pairOPERATOR{S}{S}$ is obvious.
\PROOFend  

Let us note that at item  {\bf a)} of Corollary \ref{eq_Iu_Ju_commute_with_ops_related_to_ST_from_the_stabilizer_of_U}
we deal with operators whose matrices are unitary in a fixed basis of $\COMPLEXmatricesNxN{N}$ orthonormal wrt 
$\UinnerPRODUCT{U}{}{}$. Because these matrices commute (as the operators commute) and they are normal, 
they are simultaneously unitarily diagonalizable (see \cite{Johnson}, Theorem 2.5.5). 
From this we conclude that $\Iu{U}$ and $\pairOPERATOR{\leftPERM{S}{T}}{\rightPERM{S}{T}}$ have a common eigenbasis
orthonormal wrt $\UinnerPRODUCT{U}{}{}$. It might even be possible to find a subset of $\STAB{\PP}{U}$ generating, 
through $(S,T) \longrightarrow \pairOPERATOR{\leftPERM{S}{T}}{\rightPERM{S}{T}}$, a set of operators which,
together with $\Iu{U}$, form a commuting family, in consequence a family with a common eigenbasis orthonormal wrt $\UinnerPRODUCT{U}{}{}$.
In this way an analysis of the stabilizer $\STAB{\PP}{U}$ can help in the spectral decomposition of $\Iu{U}$.
(At the end of this subsection we shall explain what other practical consequences in this respect 
	Corollary \ref{eq_Iu_Ju_commute_with_ops_related_to_ST_from_the_stabilizer_of_U}
has
	if a set of operators $\pairOPERATOR{\leftPERM{S}{T}}{\rightPERM{S}{T}}$ forming a finite group is considered as commuting with $\Iu{U}$.)

Similar remarks could be written  for the operators $\Ju{U}$ and $\pairOPERATOR{S}{S}$ at
	Corollary \ref{eq_Iu_Ju_commute_with_ops_related_to_ST_from_the_stabilizer_of_U}  {\bf b)}.
But why should we care about $\Ju{U}$ ?
Because the values of $\Du{U}$ on $\IMspaceOF{\Vu{U}{1}}$ 
	(which in turn is the domain of the  parametrization of $\FEASIBLEspace{U}$, the space of interest,
	provided by Theorem  \ref{theor_unitarity_and_eigenspaces_of_Iu} {\bf e)}) 
are  eigenvectors of $\Ju{U}\ =\ \Cu{U} \Du{U}^{-1}$  associated with the eigenvalue $1$:
\begin{equation}
	\label{eq_eigenvectors_of_Ju_for_eigenvalue_1}
	\Ii R \in \IMspaceOF{\Vu{U}{1}}
			\ \ \ \ \Longleftrightarrow\ \ \ \ 
	\Cu{U}(\Ii R)    \ \ =\ \      \Du{U}(\Ii R)
			\ \ \ \ \Longleftrightarrow\ \ \ \ 
	\Cu{U} \Du{U}^{-1} \left(   \Du{U}(\Ii R)  \right)       \ \ =\ \      \Du{U}(\Ii R)                   \ \ .
\end{equation}
Besides operators $\Iu{U}\ =\ \Cu{U}^{-1} \Du{U}$ and $\Du{U} \Cu{U}^{-1}$ have the same eigenvalues counting multiplicities
(see \cite{Johnson}, Theorem 1.3.20), and the eigenvalues of $\Ju{U}\ =\ \Cu{U} \Du{U}^{-1}\ =\ \left(\Du{U} \Cu{U}^{-1}\right)^{-1}$  
are the inverses of the eigenvalues of  $\Du{U} \Cu{U}^{-1}$.
However we have not analyzed $\Ju{U}$   in our efforts to characterize $\FEASIBLEspace{U}$ and $\undephasedDEFECT(U)$,
its dimension.

Karabegov in his article \cite{Karabegov} doesn't use the notion of a stabilizer. Instead he defines a symmetry group of a unitary matrix $U$.
Our formulation of this definition is different, but equivalent.


\begin{definition}
	\label{def_symmetry_group}
	$\GROUP{G}$ is a symmetry group (in the sense of Karabegov) of a unitary matrix $U$ if there exist monomorphisms
	$g \longrightarrow S(g),\ \ g \longrightarrow T(g):\ \ \GROUP{G} \rightarrow \enphasedPERMS$ satisfying
	\begin{equation}
		\label{eq_Smono_and_Tmono_equivalence}
		\forall g \in \GROUP{G}
				\ \ \ \ 
		S(g)\ \ =\ \ U  \cdot T(g)   \cdot  U^{-1}           \ \ .
	\end{equation}
\end{definition}

Fortunately we can identify a symmetry group of $U$ with a subgroup of the stabilizer $\STAB{\PP}{U}$:


\begin{lemma}
	\label{lem_identifying_symmetry_group_with_stabilizer_subgroup}
	\begin{description}
		\item[a)]
			If $\GROUP{G}$ is a symmetry group of $U$, for which the required monomorphism $g \rightarrow S(g)$, $g \rightarrow T(g)$
			were found, then $g \longrightarrow (S(g),T(g))$ is a monomorphism which maps $\GROUP{G}$ onto a subgroup
			$\GROUP{E}$ of the stabilizer $\STAB{\PP}{U}$.

		\item[b)]
			If $\GROUP{E}$ is a subgroup of the stabilizer $\STAB{\PP}{U}$, then it is a symmetry group of $U$ for which\ \ 
			$\GROUP{E} \ni (S,T) \longrightarrow S$,\ \  $\GROUP{E} \ni (S,T) \longrightarrow T$ are monomorphisms which 
			can be taken as the required $g \rightarrow S(g)$, $g \rightarrow T(g)$ in the definition of a symmetry group,
			respectively.

		\item[c)]
			$\GROUP{G}$, a group, is a symmetry group (in the sense of Karabegov) of $U$ if and only if $\GROUP{G}$ is isomorphic
			to some subgroup $\GROUP{E}$ of the stabilizer $\STAB{\PP}{U}$.
	\end{description}
\end{lemma}

\PROOFstart 
\begin{description}
	\item[a)]
		$g \longrightarrow (S(g),T(g))$ clearly is an injective homomorphism (between $\GROUP{G}$ and $\PP$), as built from monomorphisms.
		Its image, a subgroup of  $\PP$, is additionally contained in $\STAB{\PP}{U}$, also a subgroup of $\PP$. It is because
		from (\ref{eq_Smono_and_Tmono_equivalence}) we get that $U\ \ =\ \ S(g) U T(g)^{-1}\ \ =\ \ (S(g),T(g))U$,
		so $(S(g),T(g)) \in \STAB{\PP}{U}$, for any $g \in \GROUP{G}$. Hence this image is a subgroup of $\STAB{\PP}{U}$.  

	\item[b)]
		Let $\GROUP{E}$ be a subgroup of $\STAB{\PP}{U}$. Maps $(S,T) \longrightarrow S$ and $(S,T) \longrightarrow T$ 
		are obviously homomorphisms on $\GROUP{E}$. For any $(S,T) \in \GROUP{E} \subset \STAB{\PP}{U}$ there holds 
		$U\ =\ SUT^{-1}$, that is $S\ = UTU^{-1}$. 
		Thus condition (\ref{eq_Smono_and_Tmono_equivalence}) is satisfied for $(S,T) \longrightarrow S$ and $(S,T) \longrightarrow T$.
		If $(S_1,T_1) \neq (S_2,T_2)$ then $S_1 \neq S_2$ and $T_1 \neq T_2$ because $S$ uniquely determines $T$, and vice versa, for any 
		$(S,T) \in \STAB{\PP}{U}$. Hence the discussed maps are monomorphisms.

	\item[c)]
		In {\bf a)} we have shown that if $\GROUP{G}$ is a symmetry group of $U$ then it is isomorphic with its image in $\STAB{\PP}{U}$
		under\ \ $\GROUP{G}\ni g\ \ \longrightarrow\ \ (S(g), T(g))$.

		Now let $\GROUP{G}$ be a group, $\GROUP{E}$ be a subgroup of $\STAB{\PP}{U}$, and let $\sigma:\ \GROUP{G} \rightarrow \GROUP{E}$
		be an isomorphism between these groups. Then on account of what has been shown in the proof of {\bf b)} the compositions\ \ 
		$\left( (S,T) \rightarrow S \right) \sigma$\ \  and\ \   $\left( (S,T) \rightarrow T \right) \sigma$\ \ are monomorphisms from $\GROUP{G}$
		into $\enphasedPERMS$, which satisfy:
		\begin{equation}
			\label{eq_composed_monomorphisms_ST_equivalence}
			\forall g \in \GROUP{G}
						\ \ \ \ \ \ 
			\left( \left( (S,T) \rightarrow S \right) \sigma \right)(g)
					\ \ =\ \ 
			U         \cdot         \left( \left( (S,T) \rightarrow T \right) \sigma \right)(g)        \cdot         U^{-1}                 \ \ ,
		\end{equation}
		because 
			$\left(\ \    \left( \left( (S,T) \rightarrow S \right) \sigma \right)(g)\ \ ,\ \ 
						\left( \left( (S,T) \rightarrow T \right) \sigma \right)(g)\ \ \right)
			\ \ =\ \ 
			\sigma(g)\ \ \in \STAB{\PP}{U}$.
		We end with a conclusion that $\GROUP{G}$ is a symmetry group of $U$.
\end{description}
\PROOFend 

The results of 
	Lemma \ref{eq_ST_pair_acting_on_F_o_U_equal_CF_U} {\bf b)} 
for $U=V$ and 
	Corollary \ref{eq_Iu_Ju_commute_with_ops_related_to_ST_from_the_stabilizer_of_U}   {\bf a)} 
are formulated in Karabegov's work \cite{Karabegov} in terms of complex representations 
	(homomorphisms into the group of invertible operators on $\COMPLEXmatricesNxN{N}$)
of a symmetry group of $U$.
Having proved Lemma  \ref{lem_identifying_symmetry_group_with_stabilizer_subgroup} we instead speak about representations
of a subgroup $\GROUP{E}$ of the stabilizer $\STAB{\PP}{U}$. 
Before we provide a lemma  where these results are formulated in the Karabegov's style, 
let us first introduce the representations to be used there. 
The knowledge of their properties though is not needed  for  comprehension of this lemma.
The representations of interest are intoduced at items {\bf d)} and {\bf e)} of the lemma below.


\begin{lemma}
	\label{lem_representations_of_Stab_U}
	Let $U$ be a unitary matrix. Then:
	\begin{description}
		\item[a)]
			Map $(S,T) \longrightarrow (S,S)$ is a monomorphism from $\STAB{\PP}{U}$ into $\PP$

		\item[b)]
			Map $(S,T) \longrightarrow   \left( \leftPERM{S}{T}, \rightPERM{S}{T} \right)$ is a homomorphism from $\STAB{\PP}{U}$
			into $\PsPs$.

			If $U$ has no zero entries then its kernel is the set\ \ $\{ (zI,zI):\ \ABSOLUTEvalue{z}\ =\ 1 \}$.

		\item[c)]
			Map $(A,B) \longrightarrow \pairOPERATOR{A}{B}$ is a homomorphism from any of the groups listed below into
			the group of invertible linear operators on $\COMPLEXmatricesNxN{N}$.
			These homomorphisms have the respective kernels:
			\begin{description}
				\item[1.]
					The group of pairs of nonsingular complex matrices.\ \   $\ker= \{ (wI, wI):\ w\ \neq 0 \}$.

				\item[2.]
					$\UNITARY^2$ or\ \ $\PP$ or\ \ $\STAB{\PP}{U}$.\ \    $\ker\ \ =\ \ \{ (zI,zI):\ \ABSOLUTEvalue{z}\ =\ 1 \}$.

				\item[3.]
					$\PsPs$\ \ or\ \ $\STAB{\PsPs}{U}$.\ \ $\ker\ \ =\ \ \{(I,I)\}$\ \ (monomorphism).
			\end{description}

		\item[d)]
			Map $(S,T) \longrightarrow \pairOPERATOR{S}{S}$ is a homomorphism from $\STAB{\PP}{U}$ into
			the group of invertible linear operators on $\COMPLEXmatricesNxN{N}$. Its kernel is the set 
			$\{ (zI,zI):\ \ABSOLUTEvalue{z}\ =\ 1 \}$.

		\item[e)]
			Map $(S,T) \longrightarrow \pairOPERATOR{\leftPERM{S}{T}}{\rightPERM{S}{T}}$ is a homomorphism from $\STAB{\PP}{U}$
			into the group of invertible linear operators on $\COMPLEXmatricesNxN{N}$.

			If $U$ has no zero entries then its kernel is the set\ \ $\{ (zI,zI):\ \ABSOLUTEvalue{z}\ =\ 1 \}$.
	
	\end{description}
\end{lemma}

\PROOFstart 
\begin{description}
	\item[a)] 
		This is obviously a homomorphism. It is injective because $(S,T) \rightarrow S$ is injective 
	on $\STAB{\PP}{U}$ by Lemma \ref{lem_identifying_symmetry_group_with_stabilizer_subgroup} {\bf b)}.

	\item[b)]
		One can easily show that permutation patterns of $(S,T)$ and $(Q,R)$ get multiplied, when multiplying these pairs:
	$\left( \leftPERM{SQ}{TR}, \rightPERM{SQ}{TR} \right) \ \ =\ \ 
		\left( \leftPERM{S}{T}, \rightPERM{S}{T} \right)      \left( \leftPERM{Q}{R}, \rightPERM{Q}{R} \right)$,
	using the fact that any $(X,Y) \in \PP$ can be obtained from $\left( \leftPERM{X}{Y}, \rightPERM{X}{Y} \right)$
	by both left or right multiplying it by a pair of unitary diagonal matrices.

	Now about the kernel. If $(S,T)$ is in the set mentioned at {\bf b)} then it is mapped into $(I,I)$ as needed.
	Let us assume that $U$ has no zero entries. Let $(S,T) \in \STAB{\PP}{U}$ and let 
	$\left( \leftPERM{S}{T}, \rightPERM{S}{T} \right)\ \ =\ \ (I,I)$. 
	This means that $(S,T)\ =\ \left( \leftDIAG{S}{T}, \rightDIAG{S}{T} \right)$, a pair of unitary diagonal matrices. Then
	\begin{equation}
		\label{eq_DprimUDbiv_inv_equals_U}
		(S,T)U\ =\ \leftDIAG{S}{T} U \CONJ{\rightDIAG{S}{T}}\ \ \ =\ \ \ U
				\ \ \ \ \ \ \ \Longleftrightarrow\ \ \ \ \ \ \ 
		\left( a \CONJ{b}^T \right) \HADprod U \ \ \ =\ \ \ U
				\ \ \ \ \ \ \ \Longleftrightarrow\ \ \ \ \ \ \
		 a \CONJ{b}^T \ \ =\ \ \ONESvect \ONESvect^T       ,
	\end{equation}
	where $a,b \in \columnVECTORS{N}$ contain the diagonals of $\leftDIAG{S}{T}$ and $\rightDIAG{S}{T}$, respectively,
	and $\ONESvect$ is the all ones column vector. From the last equality we get \ \ 
		$a\ \ =\ \ z \ONESvect$,\ \ $b\ \ =\ \ z \ONESvect$,\ \ 
	where $\ABSOLUTEvalue{z}\ =\ 1$. This implies that $(S,T)\ \ =\ \ (zI,zI)$ with $z$ unimodular.

	\item[c)]
		Since $(A,B)(C,D)\ =\ (AC,BD)$ is maped into $\pairOPERATOR{AC}{BD}\ =\ \pairOPERATOR{A}{B}  \pairOPERATOR{C}{D}$,
	the considered map is a homomorphism. Its kernel on the group of pairs of invertible matrices is the set:
	\begin{equation}
		\label{eq_pair_to_pair_operator_map_kernel}
		K
				\ \ = \ \ 
		\left\{   (A,B)\ \mbox{inv.}:\ \pairOPERATOR{A}{B}\ =\ \IDENTITY{\COMPLEXmatricesNxN{N}}   \right\}
				\ \ =\ \ 
		\left\{   (A,B)\ \mbox{inv.}:\  \forall X \in \COMPLEXmatricesNxN{N}\ \ AXB^{-1} = X    \right\}      .
	\end{equation}

	If $(A,B) \in K$, then $AIB^{-1}\ =\ I$, so $A\ =\ B$.\ \ 
	If $(A,A) \in K$ then\ \
		$\forall X \in \COMPLEXmatricesNxN{N}\ \ AX = XA$.
	Let us take $X = \STbasis{i} \STbasis{i}^T$, having $1$ at the $i,i$th position and zeros elsewhere. Then:
	\begin{equation}
		\label{eq_AX_XA_equality_special_case}
		AX
				\ \ =\ \ 
		\COLUMNof{A}{i}    \STbasis{i}^T
				\ \ =\ \ 
		\STbasis{i}    \ROWof{A}{i}
				\ \ =\ \ 
		XA
							\ \ \ \ \ 
		\mbox{for}\ \ i\ =\ 1,2,\ldots,N                     \ \ ,
	\end{equation}
	where $\COLUMNof{A}{i}$ is the $i$th column and $\ROWof{A}{i}$ is the $i$th row of $A$. 
	From the above equalities we conclude that if $(A,A) \in K$ then $A$ is diagonal. In this situation also
	$A$ commutes with $\STbasis{1} \ONESvect^T$, having $1$'s in the first row and zeros elsewhere,
	so all diagonal elements of $A$ are equal. Thus we have shown that
	\begin{equation}
		\label{eq_kernel_at_case_1}
		K
				\ \ =\ \ 
		\left\{  (A,A):\ \forall X \in \COMPLEXmatricesNxN{N}\ \ AX = XA  \right\}
				\ \ =\ \ 
		\left\{   (wI,wI):\    w \neq 0    \right\}                         \ \ ,
	\end{equation}
	which is the kernel at item {\bf c)1.}. 
	(A comment: the set $\left\{ A \in \COMPLEXmatricesNxN{N} :\   \forall X \in \COMPLEXmatricesNxN{N}\ \ AX = XA  \right\}$
	is called the centre of algebra $\COMPLEXmatricesNxN{N}$ and because it is equal to $\{wI:\ w \in \COMPLEX\}$ 
	algebra $\COMPLEXmatricesNxN{N}$ is called central.)

	By restricting the domain of the considered homomorphism we restrict its kernel, so at item {\bf c)2.} the kernel is equal to
	\begin{equation}
		\label{eq_kernel_at_case_2}
		K\ \cap\ \UNITARY^2       \ \ =\ \        K\ \cap\ \PP      \ \ =\ \       K \cap \STAB{\PP}{U}     \ \ =\ \ 
		\{ (zI,zI):\ \ABSOLUTEvalue{z}\ =\ 1 \}                   \ \ ,
	\end{equation}
	and at item {\bf c)3.} the kernel is equal to:
	\begin{equation}
		\label{eq_kernel_at_case_3}
		K\ \cap\ \PsPs        \ \ =\ \        K\ \cap\ \STAB{\PsPs}{U}         \ \ =\ \                 \{(I,I)\}        \ \ .
	\end{equation}

	\item[d)]
		Map $(S,T) \longrightarrow \pairOPERATOR{S}{S}$ is a composition of
					the monomorphism $(S,T) \rightarrow (S,S)$ on $\STAB{\PP}{U}$ of item {\bf a)} 
		and 
					the homomorphism $(A,B) \rightarrow \pairOPERATOR{A}{B}$ on $\PP$ of item {\bf c)2.}.
		Therefore its kernel is the preimage of $\{ (zI,zI):\ \ABSOLUTEvalue{z}\ =\ 1 \}$ under this first monomorphism,
		which is equal to this set.

	\item[e)]
		Map $(S,T) \longrightarrow \pairOPERATOR{\leftPERM{S}{T}}{\rightPERM{S}{T}}$ is a composition of 
					the homomorphism $(S,T) \rightarrow  \left( \leftPERM{S}{T}, \rightPERM{S}{T} \right)$ on $\STAB{\PP}{U}$ of item {\bf b)}
		and
					the monomorphism $(A,B) \rightarrow \pairOPERATOR{A}{B}$ on $\PsPs$ of item {\bf c)3.}.
		Therefore its kernel is the kernel of this first homomorphism.

\end{description}
\PROOFend 

The additional assumption on $U$ at items {\bf b)} and {\bf e)} of the above lemma is important. Consider $U = I$ which has plenty of zeros as its entries.
In this case $\STAB{\PP}{I}\ \ =\ \ \left\{ (S,S):\ S \in \enphasedPERMS \right\}$, 
and the kernel of  $(S,T) \rightarrow  \left( \leftPERM{S}{T}, \rightPERM{S}{T} \right)$ on $\STAB{\PP}{I}$ is the set \\
$\left\{ (D,D):\ D\ \mbox{unitary diagonal} \right\}$, which is greater than $\{ (zI,zI):\ \ABSOLUTEvalue{z}\ =\ 1 \}$.
In general the kernel of  $(S,T) \rightarrow  \pairOPERATOR{\leftPERM{S}{T}}{\rightPERM{S}{T}}$ on $\STAB{\PP}{U}$
is equal to the kernel of  $(S,T) \rightarrow  \left( \leftPERM{S}{T}, \rightPERM{S}{T} \right)$ on $\STAB{\PP}{U}$ 
because $(A,B) \rightarrow \pairOPERATOR{A}{B}$ is monomorphic on $\PsPs$, the set of pairs of permutation matrices
(by Lemma \ref{lem_representations_of_Stab_U} {\bf c)}).

By restricting the homomorphisms on $\STAB{\PP}{U}$, considered in the above lemma, to a subgroup $\GROUP{E}$ of this stabilizer, 
we restrict their kernels by intersecting them with $\GROUP{E}$.

We are ready now to present the results of
Lemma \ref{eq_ST_pair_acting_on_F_o_U_equal_CF_U} {\bf b)} for $U=V$ and 
Corollary \ref{eq_Iu_Ju_commute_with_ops_related_to_ST_from_the_stabilizer_of_U}   {\bf a)}, 
plus some additional material on $\Ju{U}\ =\ \Cu{U} \Du{U}^{-1}$, in the Karabegov's way:


\begin{lemma}
	\label{eq_properties_of_representations_of_STAB_U}
	Let $U$ be  a unitary matrix with no zero entries and let $\GROUP{E}$ be a subgroup of $\STAB{\PP}{U}$. 
	\begin{description}
		\item[a)]
			Representation $(S,T) \longrightarrow \pairOPERATOR{S}{S}$ of $\GROUP{E}$ is unitary wrt $\STinnerPRODUCT{}{}$.
			Representation $(S,T) \rightarrow  \pairOPERATOR{\leftPERM{S}{T}}{\rightPERM{S}{T}}$ of $\GROUP{E}$
			is unitary wrt $\UinnerPRODUCT{U}{}{}$.

		\item[b)]
			Operators $\Cu{U}$ and $\Du{U}$ intertwine representations 
				$(S,T) \longrightarrow \pairOPERATOR{S}{S}$ and $(S,T) \rightarrow  \pairOPERATOR{\leftPERM{S}{T}}{\rightPERM{S}{T}}$
			of $\GROUP{E}$:
			\begin{equation}
				\label{eq_Cu_Du_intertwine_representations_of_STAB_U}
				\pairOPERATOR{S}{S}   
						\ \ =\ \ 
				\Cu{U}     \pairOPERATOR{\leftPERM{S}{T}}{\rightPERM{S}{T}}     \Cu{U}^{-1}
										\ \ \ \ ,\ \ \ \ \ \ \ \ 
				\pairOPERATOR{S}{S}   
						\ \ =\ \ 
				\Du{U}     \pairOPERATOR{\leftPERM{S}{T}}{\rightPERM{S}{T}}     \Du{U}^{-1}                   \ \ .
			\end{equation}

		\item[c)]
			Operator $\Iu{U}$ and $\Ju{U}$ commute with the respective representations of $\GROUP{E}$:
			\begin{equation}
				\label{eq_Iu_Ju_commute_with_representations_of_STAB_U}
				\Iu{U}  \pairOPERATOR{\leftPERM{S}{T}}{\rightPERM{S}{T}}
						\ \ =\ \ 
				\pairOPERATOR{\leftPERM{S}{T}}{\rightPERM{S}{T}}    \Iu{U}
									\ \ \ \ ,\ \ \ \ \ \ \ \ 
				\Ju{U}    \pairOPERATOR{S}{S}    
						\ \ =\ \ 
				\pairOPERATOR{S}{S}   \Ju{U}                    \ \ .
			\end{equation}

	\end{description} 
\end{lemma}

\PROOFstart 
\begin{description}
	\item[a)]
		The unitarity of $\pairOPERATOR{\leftPERM{S}{T}}{\rightPERM{S}{T}}$ and $\pairOPERATOR{S}{S}$
		has been shown in the proof Corollary \ref{eq_Iu_Ju_commute_with_ops_related_to_ST_from_the_stabilizer_of_U}
		{\bf a)} and {\bf b)}.

	\item[b)]
		This is a consequence of Lemma \ref{eq_ST_pair_acting_on_F_o_U_equal_CF_U} {\bf b)} for $V=U$.
		(See also the beginning of the proof of Corollary \ref{cor_Iu_Iv_operators_for_equivalent_U_V_are_similar}.)

	\item[c)]
		This follows from {\bf b)} or from 
		 Corollary \ref{eq_Iu_Ju_commute_with_ops_related_to_ST_from_the_stabilizer_of_U}.
\end{description}
\PROOFend 

Let us first note at this point that representation $(S,T)     \rightarrow     \pairOPERATOR{\leftPERM{S}{T}}{\rightPERM{S}{T}}$
is always finite (i.e. its image is finite), therefore we can use the theory of representations of finite groups to investigate its irreducible (minimal) 
invariant spaces. Secondly, the fact that $\Iu{U}$ commutes with all $\pairOPERATOR{\leftPERM{S}{T}}{\rightPERM{S}{T}}$'s
in this representation implies that every eigenspace $\Vu{U}{\lambda}$ of $\Iu{U}$ is invariant for this representation 
(the reader can easily check this), 
this generalizes (\ref{eq_eigenspaces_for_1_for_U_and_V_are related_by similarity_operator}) 
in Corollary \ref{cor_DU_equal_DV_for_U_V_equivalent} for $V\ =\ U$.
	(Because  $\pairOPERATOR{\leftPERM{S}{T}}{\rightPERM{S}{T}}$'s map real matrices into real ones,
	$\REspaceOF{\Vu{U}{\lambda}}$ and $\IMspaceOF{\Vu{U}{\lambda}}$ are invariant real subspaces for 
	$\pairOPERATOR{\leftPERM{S}{T}}{\rightPERM{S}{T}}$'s, treated as real isomorphisms, from the considered representation.
	This generalizes (\ref{eq_Im_eigenspaces_for_1_for_U_and_V_are related_by similarity_operator}) for $V = U$.)	
We can thus consider
restrictions $\RESTRICTEDto{\pairOPERATOR{\leftPERM{S}{T}}{\rightPERM{S}{T}}}{\Vu{U}{\lambda}}$ and further split
	$\Vu{U}{\lambda}\ =\    \left( \Vu{U}{\lambda} \right)'    \oplus     \left( \Vu{U}{\lambda} \right)''     \oplus     \ldots$
into irreducible invariant subspaces for the representation. 
    This is guaranteed by the Maschke's theorem in the representation theory.
Doing this we can make use of the operation of taking the orthogonal
complement wrt $\UinnerPRODUCT{U}{}{}$ within a subspace -- 
if 
	$\mathbbm{W} \subset \Vu{U}{\lambda}$ 
is invariant then
	$\UorthCOMPLEMENT{U}{\mathbbm{W}} \subset \Vu{U}{\lambda}$
is invariant too, because the representation (restricted to any invariant subspace) is unitary wrt $\UinnerPRODUCT{U}{}{}$ by 
	Lemma \ref{eq_properties_of_representations_of_STAB_U} {\bf a)}.
Of course not only $\Vu{U}{\lambda}$ but any invariant subspace of it can be split in this way as long as it is not irreducible.

Now, if in a  known (from some other source) decomposition of $\COMPLEXmatricesNxN{N}$ into irreducible invariant subspaces some 
of these subspaces generate 
(by restricting $\pairOPERATOR{\leftPERM{S}{T}}{\rightPERM{S}{T}}$'s to such a subspace) subrepresentations that are inequivalent
to all the other subrepresentations in the decomposition, then these selected subspaces are {\sl uniquely determined} and are a special type 
of the so called isotypic subspaces.  
This follows from the theorem on the uniqueness of the isotypical decomposition in the representation theory.
In such a situation a decomposition into irreducible invariant subspaces obtained from the eigenspaces of $\Iu{U}$:
	$\COMPLEXmatricesNxN{N}\ =\ 
			\left( \Vu{U}{\lambda} \right)'    \oplus     \left( \Vu{U}{\lambda} \right)''     \oplus     \ldots 
				\oplus
			\left( \Vu{U}{\mu} \right)'    \oplus     \left( \Vu{U}{\mu} \right)''     \oplus     \ldots $
must contain these isotypic irreducible components. 
In this way an analysis of the representation $(S,T)     \rightarrow     \pairOPERATOR{\leftPERM{S}{T}}{\rightPERM{S}{T}}$
(i.e. of its known, from some other source, irreducible subspaces)
may reveal subspaces over which $\Iu{U}$ is a scalar map and help in finding its eigenspaces !

Let an extremal example be the Fourier matrix $U = 1/\sqrt{N} \cdot F$ (see (\ref{eq_Fourier_Kron_prod}) in the next section), 
the stabilizer of which contains a subgroup which has a finite and abelian representation 
$(S,T)     \rightarrow     \pairOPERATOR{\leftPERM{S}{T}}{\rightPERM{S}{T}}$. As it is explained at the beginning of
section \ref{sec_fourier_defect} in proof {b)} of Theorem \ref{theor_eigenbasis_of_I_F}, 
this representation has $N^2$ independent isotypic (thus inequivalent) irreducible invariant $1$ dimensional subspaces in 
$\COMPLEXmatricesNxN{N}$, over which $\Iu{U}$ is a scalar map, hence the full eigenbasis of $\Iu{U}$ is given.

We should also remember that $\TRIVIALspace{U} \subset \Vu{U}{1}$ 
		(see Theorem \ref{theor_unitarity_and_eigenspaces_of_Iu} {\bf c)})
is always an invariant subspace for the discussed representation and so is 
	$\UorthCOMPLEMENT{U}{\TRIVIALspace{U}}$, the complement in   $\COMPLEXmatricesNxN{N}$. 
Therefore we can concentrate (in our search for other eigenspaces of $\Iu{U}$) 
on decomposing only $\UorthCOMPLEMENT{U}{\TRIVIALspace{U}}$  into irreducible invariant subspaces,
some of which, possibly, are simultaneously isotypic and thus $\Iu{U}$ is scalar on them.  This is analogous to the
presence of isotypic irreducible subspaces, over which $\Iu{U}$ is scalar, in the decomposition of $\COMPLEXmatricesNxN{N}$
in the paragraph preceding the one above.
	(Here is an additional explanation:
	$\Iu{U}$ is unitary wrt $\UinnerPRODUCT{U}{}{}$ and $\TRIVIALspace{U} \subset \Vu{U}{1}$ so both 
	$\TRIVIALspace{U}$ and $\UorthCOMPLEMENT{U}{\TRIVIALspace{U}}$ are invariant spaces of $\Iu{U}$. 
	$\RESTRICTEDto{\Iu{U}}{\UorthCOMPLEMENT{U}{\TRIVIALspace{U}}}$, unitary wrt $\UinnerPRODUCT{U}{}{}$,
	commutes with all  $\RESTRICTEDto{\pairOPERATOR{\leftPERM{S}{T}}{\rightPERM{S}{T}}}{ \UorthCOMPLEMENT{U}{\TRIVIALspace{U}} }$'s\ ,\ \ 
	so  $\UorthCOMPLEMENT{U}{\TRIVIALspace{U}}$ can be decomposed into eigensubspaces of 
	$\RESTRICTEDto{\Iu{U}}{\UorthCOMPLEMENT{U}{\TRIVIALspace{U}}}$ which are invariant irreducible for the representation
	restricted to $\UorthCOMPLEMENT{U}{\TRIVIALspace{U}}$. Just as it was in the case with $\COMPLEXmatricesNxN{N}$ 
	in the place of $\UorthCOMPLEMENT{U}{\TRIVIALspace{U}}$ above.)

This is what   allowes  Karabegov to fully decompose $\COMPLEXmatricesNxN{N}$ into the eigenspaces of $\Iu{U}$ in his Example 2
in \cite{Karabegov}. In the example a unitary $N \times N$ matrix\ \ $U\ =\ I_N + \frac{\theta -1}{N} J_N$ is considered, 
where $I_N$ is the identity matrix, $J_N$ is the all ones matrix, $\ABSOLUTEvalue{\theta}\ =\ 1$ but $\theta \neq \pm 1$,
and $N > 2$. Karabegov restricts his attention to $\GROUP{E}\ =\ \STAB{\PsPs}{U}$, a subgroup of $\STAB{\PP}{U}$,
being equal to $\left\{(P,P):\ \mbox{$P$ a permutation matrix}\right\}$ and to the representation of $\GROUP{E}$:
$(P,P)\ \longrightarrow\ \pairOPERATOR{P}{P}$.

\subsection{Kronecker products of unitary matrices}
	\label{subsec_Kronecker_products}

The rest of this section is devoted to the application of the Karabegov's theory to Kronecker products of unitary matrices.
It also supplements our paper \cite{KroneckerDefect} on the defect of such matrices.

In his papers \cite{Karabegov_old,Karabegov} Karabegov presented the spectral decomposition of $\Iu{F}$, the Berezin transform constructed 
for a Fourier matrix $F$ (see (\ref{eq_Fourier_Kron_prod}) in section \ref{sec_Fourier_matrices}). He discovered that the spectrum of
$\Iu{F}$ is formed by the entries of $F$. A consequence of this fact is that when we consider Fourier matrix $F \otimes G$ built from 
Fourier matrices $F$ and $G$, the eigenvalues of $\Iu{F \otimes G}$ are all the products $\lambda \cdot \mu$ where $\lambda$,$\mu$ are eigenvalues
of $\Iu{F}$,$\Iu{G}$, respectively. This clearly suggests a generalization, namely that $\Iu{U \otimes V}$ is a tensor product of 
operators $\Iu{U}$ and $\Iu{V}$, which is the subject of Lemma \ref{lem_Iu_for_Kron_prod_of_unitaries} below. The above property
concerning Fourier matrices also leads to the conditional multiplicativity of the undephased defect of a Fourier matrix $F$ with respect to Kronecker 
subproducts of $F$, see Lemma \ref{lem_multiplicativity_of_D_of_F} in section \ref{sec_fourier_defect}.


\begin{lemma}
	\label{lem_Iu_for_Kron_prod_of_unitaries}
	Let $U$ and $V$ be unitary matrices of size $N \times N$ and $M \times M$, respectively, 
	with nonzero entries. Let $F \in \COMPLEXmatricesNxN{N}$ and $G \in \COMPLEXmatricesNxN{M}$
	be arbitrary. Then (where $\otimes$ denotes the Kronecker product):
	\begin{description}
		\item[a)]  $\Cu{U \otimes V}( F \otimes G )      \ \ =\ \       \Cu{U}(F)  \otimes  \Cu{V}(G)$.
		\item[b)]  $\Du{U \otimes V}( F \otimes G )      \ \ =\ \       \Du{U}(F)  \otimes  \Du{V}(G)$.
		\item[c)]  $\Iu{U \otimes V}( F \otimes G )      \ \ =\ \         \Iu{U}(F)  \otimes   \Iu{V}(G)$.
	\end{description}
	Expressing it differently,\ \ \ $\Cu{U} \otimes \Cu{V} \ \ =\ \ \Cu{U \otimes V}$,\ \ \
								$\Du{U} \otimes \Du{V} \ \ =\ \ \Du{U \otimes V}$,\ \ \
								$\Iu{U} \otimes \Iu{V} \ \ =\ \ \Iu{U \otimes V}$,\ \ \
	where the tensor products of operators are taken with respect to the Kronecker products both on the side of arguments and results.

	The number of tensor (Kronecker)  factors in the above statements can be greater than 2 because the Kronecker product is
	associative.
\end{lemma}

\PROOFstart 
\begin{description}
	\item[a)]
		$ \Cu{U \otimes V}(F \otimes G)
				\ \  =\ \ 
		((F \otimes G) \HADprod (U \otimes V))  (U \otimes V)^*
				\ \  =\ \ 
		((U \HADprod F) \otimes (V \HADprod G))  \left(U^*  \otimes V^* \right)
				\ \ =\ \ 
		\left( (U \HADprod F)U^* \right)      \ \otimes\         \left( (V \HADprod G)V^* \right)
				\ \ =\ \ 
		\Cu{U}(F)     \ \otimes\      \Cu{V}(G)$.

	\item[b)]
		$\Du{U \otimes V}(F \otimes G)
				\ \ \stackrel{L.\ref{lem_properties_of_Cu_Du_Iu} \mathbf{b)}}{=}\ \ 
		\hermTRANSPOSE{ \Cu{U \otimes V}\left( \CONJ{F} \otimes \CONJ{G} \right)}
				\ \ \stackrel{\mathbf{a)}}{=}\ \ 
		\hermTRANSPOSE{ \Cu{U}\left( \CONJ{F} \right)    \otimes    \Cu{V}\left( \CONJ{G} \right) }
				\ \ =\ \     \\ 
		\hermTRANSPOSE{ \Cu{U}\left( \CONJ{F} \right) }        \otimes        \hermTRANSPOSE{ \Cu{V}\left( \CONJ{G} \right) }
				\ \ \stackrel{L.\ref{lem_properties_of_Cu_Du_Iu} \mathbf{b)}}{=}\ \ 
		\Du{U}(F)   \otimes   \Du{V}(G)$.

	\item[c)]
		In {\bf a)} and {\bf b)} we have proved that\ \ 
			$\Cu{U} \otimes \Cu{V} \ \ =\ \ \Cu{U \otimes V}$\ \    and\ \    $\Du{U} \otimes \Du{V} \ \ =\ \ \Du{U \otimes V}$.\ \ 
		Hence\ \ 
			$\Cu{U  \otimes V}^{-1}\ \ =\ \ \left(\Cu{U} \otimes \Cu{V}\right)^{-1}\ \ =\ \ \Cu{U}^{-1} \otimes \Cu{V}^{-1}$. \\
		So\ \ 
			$\Iu{U \otimes V}
					\ \ =\ \ 
			\Cu{U  \otimes V}^{-1}    \Du{U \otimes V}   
					\ \ =\ \   
			\left( \Cu{U}^{-1} \otimes \Cu{V}^{-1}  \right)       \left( \Du{U} \otimes \Du{V} \right)
					\ \ =\ \ 
			\left( \Cu{U}^{-1} \Du{U} \right)\    \otimes\      \left(\Cu{V}^{-1} \Du{V} \right)
					\ \ =\ \ 
			\Iu{U}\       \otimes\        \Iu{V}$.\ \ 
		We have used various identities concerning the tensor products of operators, but the reader unfamiliar with the multilinear algebra
		will surely accept the alternative reasoning, involving Kronecker products only:
		\begin{eqnarray}
			\label{eq_Iu_for_Kron_prod_of_unitaries_proof}
			&  &  \\
			\nonumber
			\Du{U}(F)  \otimes \Du{V}(G)       &    =    &     \Cu{U}\left(  \Iu{U}(F)  \right)      \otimes     \Cu{V}\left(  \Iu{V}(G)  \right)      \\
			\nonumber
										&       \Downarrow\ \ \mbox{using {\bf a)} and {\bf b)}}     &                           \\
			\nonumber
			\Du{U \otimes V}(F \otimes G)     &      =      &        \Cu{U \otimes V}\left(  \Iu{U}(F)  \otimes  \Iu{V}(G)  \right)                 \\
			\nonumber
			                                      &     \Downarrow                                                                &                            \\
			\nonumber
			\Iu{U \otimes V}(F \otimes G)
					\ \ =\ \ 
			\Cu{U \otimes V}^{-1}  \Du{U \otimes V}(F \otimes G)
								&   =    &
			\Iu{U}(F)  \otimes  \Iu{V}(G)
		\end{eqnarray}
		for any $F \in \COMPLEXmatricesNxN{N}$ and $G \in \COMPLEXmatricesNxN{M}$.
	
\end{description}
\PROOFend 

From the fact that 
	$\Iu{U \otimes V \otimes \ldots}(F \otimes G \otimes \ldots)\ \ =\ \   \Iu{U}(F) \otimes \Iu{V}(G) \otimes \ldots$
we conclude that if 
	$\left( \ithMATRIX{F}{i},\ \lambda_i \right)$, $\left( \ithMATRIX{G}{j},\ \mu_j \right)$, ...
are eigenpairs for
	$\Iu{U}$, $\Iu{V}$, ..., respectively,
then \\
	$\left( \ithMATRIX{F}{i} \otimes \ithMATRIX{G}{j} \otimes \ldots,\  \lambda_i \mu_j \ldots \right)$
are eigenpairs for 
	$\Iu{U \otimes V \otimes \ldots}$.
In this way an eigenbasis for $\Iu{U \otimes V \otimes \ldots}$ is generated by eigenbases for $\Iu{U}$, $\Iu{V}$, ...
(Kronecker product, as a tensor product, generates an independent set from sets of independent vectors).
In Theorem \ref{theor_unitarity_and_eigenspaces_of_Iu} {\bf a)} we said  that each such basis for $\Iu{X}$
can be chosen to be orthonormal wrt $\UinnerPRODUCT{X}{}{}$. Is then the Kronecker generated basis for
 $\Iu{U \otimes V \otimes \ldots}$ orthonormal wrt to $\UinnerPRODUCT{U \otimes V \otimes \ldots}{}{}$\ ?
It is, because:


\begin{lemma}
	\label{lem_weigthed_U_inner_product_induced}
	Let $U$ and $V$ be as in Lemma \ref{lem_Iu_for_Kron_prod_of_unitaries}. 	
	The inner product $\UinnerPRODUCT{U \otimes V}{}{}$ on 
		$\COMPLEXmatricesNxN{N} \otimes \COMPLEXmatricesNxN{M}\ \ =\ \ \COMPLEXmatricesNxN{MN}$
	is an induced inner product associated (and thus uniquely determined by) with the inner products
		$\UinnerPRODUCT{U}{}{}$ on $\COMPLEXmatricesNxN{N}$
	and
		$\UinnerPRODUCT{V}{}{}$ on $\COMPLEXmatricesNxN{M}$. 
	In other words, for any\ \ 
		$\ithMATRIX{F}{1}, \ithMATRIX{F}{2} \in \COMPLEXmatricesNxN{N}$\ \ 
	and\ \ 
		$\ithMATRIX{G}{1}, \ithMATRIX{G}{2} \in \COMPLEXmatricesNxN{M}$\ \ 
	there holds:
	\begin{equation}
		\label{eq_weigthed_U_inner_product_induced}
		\UinnerPRODUCT{U \otimes V}{\ithMATRIX{F}{1} \otimes \ithMATRIX{G}{1}}{\ithMATRIX{F}{2} \otimes \ithMATRIX{G}{2}}
				\ \ =\ \ 
		\UinnerPRODUCT{U}{\ithMATRIX{F}{1}}{\ithMATRIX{F}{2}}
			\cdot
		\UinnerPRODUCT{V}{\ithMATRIX{G}{1}}{\ithMATRIX{G}{2}}    \ \ .
	\end{equation}

	Because the Kronecker product is associative, the number of factors in the above statement can be greater than $2$.
\end{lemma}

\PROOFstart 
$\UinnerPRODUCT{U \otimes V}{\ithMATRIX{F}{1} \otimes \ithMATRIX{G}{1}}{\ithMATRIX{F}{2} \otimes \ithMATRIX{G}{2}}
		\ \ \stackrel{L.\ref{lem_properties_of_Cu_Du_Iu} \mathbf{a)}}{=}\ \ 
\STinnerPRODUCT{ \Cu{U \otimes V}\left( \ithMATRIX{F}{1} \otimes \ithMATRIX{G}{1} \right)}
					 { \Cu{U \otimes V}\left( \ithMATRIX{F}{2} \otimes \ithMATRIX{G}{2} \right)}
		\ \ \stackrel{L.\ref{lem_Iu_for_Kron_prod_of_unitaries} \mathbf{a)}}{=}\ \    \\
\STinnerPRODUCT{\Cu{U}\left( \ithMATRIX{F}{1} \right)  \otimes  \Cu{V}\left( \ithMATRIX{G}{1} \right)}
					 {\Cu{U}\left( \ithMATRIX{F}{2} \right)  \otimes  \Cu{V}\left( \ithMATRIX{G}{2} \right)}
		\ \ =\ \  \\
\STinnerPRODUCT{\Cu{U}\left( \ithMATRIX{F}{1} \right)}{\Cu{U}\left( \ithMATRIX{F}{2} \right)}
	\cdot
\STinnerPRODUCT{\Cu{V}\left( \ithMATRIX{G}{1} \right)}{\Cu{V}\left( \ithMATRIX{G}{2}\right)}
		\ \ \stackrel{L.\ref{lem_properties_of_Cu_Du_Iu} \mathbf{a)}}{=}\ \ 
\UinnerPRODUCT{U}{\ithMATRIX{F}{1}}{\ithMATRIX{F}{2}}
	\cdot
\UinnerPRODUCT{V}{\ithMATRIX{G}{1}}{\ithMATRIX{G}{2}}$. \\
The third equality is justified by the fact that the standard inner product on 
	$\COMPLEXmatricesNxN{N} \otimes \COMPLEXmatricesNxN{M}\ \ =\ \ \COMPLEXmatricesNxN{MN}$
is also induced, that is 
	$\STinnerPRODUCT{\ithMATRIX{F}{1} \otimes \ithMATRIX{G}{1}}{\ithMATRIX{F}{2} \otimes \ithMATRIX{G}{2}}
			\ \ =\ \ 
	\STinnerPRODUCT{\ithMATRIX{F}{1}}{\ithMATRIX{F}{2}}
		\cdot
	\STinnerPRODUCT{\ithMATRIX{G}{1}}{\ithMATRIX{G}{2}}$
for arbitrary arguments of the appropriate size.
\PROOFend 

Thus, if 
	$\left( \ithMATRIX{F}{i},\ i\ =\ 1,2,\ldots \right)$ is an eigenbasis for $\Iu{U}$ orthonormal wrt $\UinnerPRODUCT{U}{}{}$,\\
	$\left( \ithMATRIX{G}{j},\ j\ =\ 1,2,\ldots \right)$ is an eigenbasis for $\Iu{V}$ orthonormal wrt $\UinnerPRODUCT{V}{}{}$,
	...,\\
then
	$\UinnerPRODUCT{U \otimes V \otimes \ldots}{ \ithMATRIX{F}{i_1} \otimes \ithMATRIX{G}{j_1} \otimes \ldots }
														{ \ithMATRIX{F}{i_2} \otimes \ithMATRIX{G}{j_2} \otimes \ldots }
			\ \ =\ \ 
	\UinnerPRODUCT{U}{ \ithMATRIX{F}{i_1} }{  \ithMATRIX{F}{i_2} } 
		\cdot
	\UinnerPRODUCT{V}{ \ithMATRIX{G}{j_1} }{  \ithMATRIX{G}{j_2} }
		\cdot
	\ldots $ can take values only in the set $\{0,1\}$.

Besides, (always existing by Theorem \ref{theor_unitarity_and_eigenspaces_of_Iu} {\bf c)})  eigenpairs of $\Iu{U}$, $\Iu{V}$, ... :
	$\left( F,\ 1 \right)$, $\left( G,\ 1 \right)$, ..., respectively,
give rise to eigenvectors $F \otimes G \otimes \ldots$\ \  of\ \  $\Iu{U \otimes V \otimes \ldots}$\ \  associated with eigenvalue $1$, being
the product of the respective eigenvalues, $1$'s actually, associated with $F$, $G$, ... . Therefore we can write:


\begin{corollary}
	\label{cor_V1_spaces_Kron_multiplied_are_contained_in_V1_space}
	Let $U$ and $V$ be as in Lemma \ref{lem_Iu_for_Kron_prod_of_unitaries}. 
	Then  (the notation as in Theorem \ref{theor_unitarity_and_eigenspaces_of_Iu} {\bf b)}):
	\begin{description}
		\item[a)]
			$\Vu{U}{1} \otimes \Vu{V}{1}      \ \ \subset\ \        \Vu{U \otimes V}{1}$, \\
			where on the left we have the tensor product, associated with the Kronecker product, of spaces $\Vu{U}{1}$ and $\Vu{V}{1}$, 
			that is the space spanned by all the products $F \otimes G$ where $F \in \Vu{U}{1}$ and $G \in \Vu{V}{1}$.

		\item[b)]
			$\Ii\    \cdot\    \left(  \IMspaceOF{\Vu{U}{1}} \otimes  \IMspaceOF{\Vu{V}{1}}  \right)
					\ \ \subset\ \ 
			\IMspaceOF{\Vu{U \otimes V}{1}}$ \\
			for the corresponding real spaces.

		\item[c)]
			As a consequence, their dimensions satisfy\ \ 
				$\undephasedDEFECT(U) \cdot  \undephasedDEFECT(V)   \ \ \leq\ \     \undephasedDEFECT(U \otimes V)$.
	\end{description}

	The number of factors in the above statements can be greater than $2$.
	Then at {\bf b)} the multiplier $\Ii$ has to be replaced by $1$ when the number of factors is odd.
\end{corollary}
		
\PROOFstart 
In the proof we only remind the reader that if\ \  $\left( F_i \right)_{i=1..\undephasedDEFECT(U)}$,\ \ 
$\left( G_j \right)_{j=1..\undephasedDEFECT(V)}$\ \  are bases for $\Vu{U}{1}$,  $\Vu{V}{1}$, respectively,
then\ \ $\left( F_i \otimes G_j \right)_{(i,j) \in \{1..\undephasedDEFECT(U)\}\times\{1..\undephasedDEFECT(V)\}}$\ \ 
is a basis for $\Vu{U}{1} \otimes \Vu{V}{1}$. Therefore the dimensions of these spaces satisfy
	$\DIMC\left( \Vu{U}{1} \otimes \Vu{V}{1} \right) 
			\ \ =\ \ 
	\DIMC\left( \Vu{U}{1} \right)  \cdot  \DIMC\left( \Vu{V}{1} \right)
			\ \ =\ \ 
	\undephasedDEFECT(U) \cdot \undephasedDEFECT(V)$.
Analogously for the real spaces, the respective dimensions are the same by Theorem \ref{theor_unitarity_and_eigenspaces_of_Iu} {\bf e)}.
\PROOFend  

In our previous paper \cite{KroneckerDefect} we defined the undephased defect, $\undephasedDEFECT(U)$, 
just like in equation (\ref{eq_old_definition_of_undephased_defect}). Basing on (68),(69) from \cite{KroneckerDefect}
we recall this definition:
\begin{equation}
	\label{eq_old_definition_of_undephased_defect_repeated}
	\undephasedDEFECT(U)\ \ =\ \ \DIMR(\VuOLD{U})             \ \ ,
			\ \ \ \ \ \ \ \mbox{where}\ \ \ \ 
	\VuOLD{U}\ \ \stackrel{def}{=}\ \ \left\{ R\ \mbox{real}:\    (\Ii R \HADprod U )U^{*}\ \ \mbox{antihermitian} \right\}               \ \ .
\end{equation}
When we compare the third characterization of $\IMspaceOF{\Vu{U}{\lambda}}$ in
	(\ref{eq_Vu_Im_definition}) (Theorem \ref{theor_unitarity_and_eigenspaces_of_Iu} {\bf b)}),
the third characterization of $\IMspaceOF{\Vu{U}{1}}$ in 
	(\ref{eq_V1_Im_characterization}) (Theorem \ref{theor_unitarity_and_eigenspaces_of_Iu} {\bf d)})
and the above definition of  $\VuOLD{U}$ it will be clear for us that $\VuOLD{U}\ \ =\ \ \REspaceOF{\Vu{U}{1}}$
if the latter space can be defined, i.e. if $U$ has no zero entries. But at the time of writing \cite{KroneckerDefect}
we did not assume that $U$ had no zero entries. 
In section Introduction of \cite{KroneckerDefect}, in the paragraph preceding (1), we noted that in the case when $U$ had zero entries
for any $R \in \VuOLD{U}$ and the corresponding direction $\Ii R \HADprod U$ there existed infinitely many other real $R$'s 
producing this direction and thus belonging to $\VuOLD{U}$ ($R_{i,j}$ may vary if $U_{i,j}=0$). This means that if $U$ has zero entries, spaces\ \ \ 
$\FEASIBLEspace{U}
		\ \ \stackrel{(\ref{eq_feasible_space_characterization})}{=}\ \ 
\left\{
	\Ii R \HADprod U\ :\ \
	R\ \ \mbox{real}\ \ \ \ \ \ \mbox{and}\ \ \ \ \ \ \mbox{$(\Ii R \HADprod U) \hermTRANSPOSE{U}$\ antihermitian}
\right\}
		\ \ =\ \ 
\Ii \VuOLD{U}  \HADprod  U$\ \ \ 
and\ \ \ $\VuOLD{U}$\ \ \ are not isomorphic and $\DIMR(\VuOLD{U})$, the old $\undephasedDEFECT(U)$, is greater than
$\DIMR(\FEASIBLEspace{U})$, the new  $\undephasedDEFECT(U)$.  Otherwise the above spaces are isomorphic 
($X \longrightarrow \Ii X \HADprod U$ is a real isomorphism then) and both definitions of $\undephasedDEFECT(U)$ are equivalent.

Nevertheless it was trivially obtained in \cite{KroneckerDefect} that also in the case of the the old $\undephasedDEFECT(U)$ there holds 
	$\undephasedDEFECT(U \otimes V)\ \ \geq\ \ \undephasedDEFECT(U) \undephasedDEFECT(V)$ (Corollary 5.6, the number of factors can be greater than $2$)
for
	$\VuOLD{U} \otimes \VuOLD{V}\ \ \subset\ \ \VuOLD{U \otimes V}$ (follows from Lemma 5.5, where $\VuOLD{U} \otimes \VuOLD{V}$ is spanned by
Kronecker products $R \otimes S$ with $R \in \VuOLD{U}$ and $S \in \VuOLD{V}$), which resembles the facts in
Corollary \ref{cor_V1_spaces_Kron_multiplied_are_contained_in_V1_space}.

Besides, in \cite{KroneckerDefect} we found lower bounds for $\undephasedDEFECT(U)$ when $U$ is a Kronecker product of
unitary matrices (Theorem 5.4):
\begin{itemize}
	\item
		Let $U\ \ =\ \ \ithMATRIX{U}{1} \otimes \ldots \otimes \ithMATRIX{U}{r}$, where $\ithMATRIX{U}{k}$ is of size $N_k \times N_k$,
		$N_k \geq 2$, and there is only one $2$ in the sequence $\left( N_1,\ N_2,\ \ldots,\ N_r \right)$. Then
		\begin{equation}
			\label{eq_kron_prod_undephased_defect_bound_without_2s}
			\undephasedDEFECT(U)
					\ \ \geq\ \ 
			\left( 2N_1 - 1 \right)     \left( 2N_2 - 1 \right)           \cdot    \ldots    \cdot        \left( 2N_r - 1 \right)                \ \ .
		\end{equation}

	\item
		$U$ as above, but this time there are\ \  $x$\ \ $2$'s in the sequence   $\left( N_1,\ N_2,\ \ldots,\ N_r \right)$,
		where $x \geq 1$. Then
		\begin{equation}
			\label{eq_kron_prod_undephased_defect_bound_with_2s}
			\undephasedDEFECT(U)
					\ \ \geq\ \ 
			\left(   \prod_{N_k > 2}   \left( 2N_k - 1 \right)   \right)            \ \cdot\       2^{x-1} \left( 2^x  +   1 \right)            \ \ .
		\end{equation}
		This bound coincides with (\ref{eq_kron_prod_undephased_defect_bound_without_2s}) when $x=1$, otherwise it is greater 
		as explained in \cite{KroneckerDefect} in the paragraph containing (105),(106).
\end{itemize}
These bounds were obtained by tedious estimation of the rank of the matrix of the system 
((\ref{eq_iRU_Uhansposed_antihermitian_system_in_matrix_form}) or (\ref{eq_iRU_Uhansposed_antihermitian_system_in_entrywise_form})) 
defining space $\VuOLD{U}$. We remarked however (\cite{KroneckerDefect}, the paragraph containing (103),(104))
 that bound (\ref{eq_kron_prod_undephased_defect_bound_without_2s}) was a consequence of the multi factor version of
	$\undephasedDEFECT(U \otimes V)\ \ \geq\ \ \undephasedDEFECT(U) \undephasedDEFECT(V)$
and the fact that for each $k$\ \ \ $\undephasedDEFECT\left( \ithMATRIX{U}{k} \right)\ \geq\ 2N_k -1$
because  the space $\VuOLD{\ithMATRIX{U}{k}}$ contained the $2N_k -1$ dimensional space
	$\REspaceOF{\TRIVIALspace{\ithMATRIX{U}{k}}}  
			\ \ =\ \ 
	\left\{   a \ONESvect^T\ +\ \ONESvect b^T \ :\ \ a,b \in \REALcolumnVECTORS{N_k}     \right\}$
(see also Theorem \ref{theor_unitarity_and_eigenspaces_of_Iu} {\bf c)}).
Here we feel obliged to justify the bound (\ref{eq_kron_prod_undephased_defect_bound_with_2s}) in a comparably simple way.
(In the paragraph below, as in this one, we use the older definition (\ref{eq_old_definition_of_undephased_defect_repeated}) of the undephased defect,
and of course we do not assume that $U$ has no zero entries.
By comparing the characterization (16) of the dephased defect $\DEFECT(U)$ in our earlier paper \cite{Defect} with (20) therein we see
that $\DEFECT(U)\ =\   \DIMR(\VuOLD{U}) - (2N-1)\ =\    \undephasedDEFECT(U) - (2N-1)$, $N$ being  the size of $U$. In \cite{Defect}
we have shown, in Lemmas 3.1 and 3.2, the invariance of $\DEFECT(U)$ under equivalence operations, so 
$\undephasedDEFECT(U)$ also remains invariant under such operations.)

This is again based on inequality
	$\undephasedDEFECT(U \otimes V)\ \ \geq\ \ \undephasedDEFECT(U) \undephasedDEFECT(V)$,
also in the many factor case, and an observation that a $2 \times 2$ unitary matrix is equivalent to a real orthogonal matrix.
The last property of course extends to any Kronecker product of $2 \times 2$ unitary matrices:
if, using the notation from the previous subsection,\ \ $U_1 = \left(A_1, B_1\right) V_1$\ \  and\ \   $U_2 = \left(A_2, B_2\right) V_2$,\ \  where
$V_1, V_2$ are real orthogonal matrices and for any $k \in \{1,2\}$\ \ $\left(A_k,B_k\right)$ belongs to $\PP$ of the appropriate size,
then\ \  $U_1 \otimes U_2 \ =\ \left( A_1 \otimes A_2,  B_1 \otimes B_2 \right) \left( V_1 \otimes V_2 \right)$\ \  which means that $U_1 \otimes U_2$
is equivalent to the real orthogonal matrix  $V_1 \otimes V_2$.\ \ \
$U\ \ =\ \ \ithMATRIX{U}{1} \otimes \ldots \otimes \ithMATRIX{U}{r}$ can be permuted to become a product,
equivalent to $U$, with all $x$\ \ $2 \times 2$ factors on the left. These $x$\ \ $2 \times 2$ factors can be considered separately
as a single factor $U'$ equivalent to a real $2^x  \times 2^x$ orthogonal matrix $O'$, whose defect, 
	$\undephasedDEFECT\left( O' \right)\ =\ \undephasedDEFECT\left( U' \right)$, 
is not lower than 
	$2^x  \left( 2^x + 1 \right)/2$
which follows from our Lemma 3.7 (a) in \cite{Defect} 
(where in fact $\DEFECT(U)\ =\ \undephasedDEFECT(U) - (2N-1)$ is estimated from below by $(N-1)(N-2)/2$, $N$ being the size of $U$).
The remaining factors $\ithMATRIX{U}{k}$ in the permuted $U$, to the right of $U'$, have their sizes greater than $2$ and
we can only say that $\undephasedDEFECT\left( \ithMATRIX{U}{k} \right)\ \geq\ 2N_k - 1$ for each of them, as in the previous paragraph.
In this way
\begin{eqnarray}
	\label{eq_D_estimation_for_kron_prod_with_2s}
	\undephasedDEFECT(U)
	&    =    &
	\undephasedDEFECT(\mbox{permuted}\ U)
			\ \ =\ \ 
	\undephasedDEFECT\left(U' \otimes \bigotimes_{N_k > 2} \ithMATRIX{U}{k} \right)                                                  \\
	\nonumber
	&     \geq     &
	\undephasedDEFECT\left( U' \right)       \cdot       \prod_{N_k > 2} \undephasedDEFECT\left( \ithMATRIX{U}{k} \right)
			\ \ \geq\ \ 
	\frac{2^x \left( 2^x + 1 \right)}{2}     \cdot       \prod_{N_k > 2}   \left( 2N_k - 1 \right)           \ \ ,
\end{eqnarray}
which is  (\ref{eq_kron_prod_undephased_defect_bound_with_2s}).

In the above reasoning we used the bound on 
	$\DEFECT(U)\ \ =\ \ \undephasedDEFECT(U) - (2N-1)\ \ =\ \ \DIMR(\VuOLD{U}) - (2N-1)$, 
for a real orthogonal matrix $U$ of size $N$, which can be found in \cite{Defect}. 
Like it was the case with the bounds 
	(\ref{eq_kron_prod_undephased_defect_bound_without_2s}) and (\ref{eq_kron_prod_undephased_defect_bound_with_2s}),
this had been obtained by an estimation of the rank of the matrix of the system defining $\VuOLD{U}$. A closer inspection of this
matrix when an orthogonal $U$ has no zero entries leads to the conclusion that the bound is saturated then, that is\ \  
$\DEFECT(U)\ =\ (N-1)(N-2)/2$\ \ and\ \ $\undephasedDEFECT(U)\ =\ N(N+1)/2$.  This very result can be obtained on the ground
of Karabegov's theory, and it was observed by him in \cite{Karabegov_old} that:


\begin{lemma}
	\label{lem_Iu_selfadjoint_for_U_real_orthogonal}
	If $U$ is a real orthogonal $N \times N$ matrix with no zero entries, then 
		$\Iu{U}$ is  self adjoint wrt $\UinnerPRODUCT{U}{}{}$ 
	(equivalently its matrix is hermitian in a basis orthonormal wrt $\UinnerPRODUCT{U}{}{}$),
	so its spectrum consists of\ \ \ $\undephasedDEFECT(U)\ =\ N(N+1)/2$\ \ \ $1$'s\ \ \
	and\ \ \ $(N-1)N/2$\ \ \ $(-1)$'s.
\end{lemma}

\PROOFstart 
$\UinnerPRODUCT{U}{\Iu{U}(F)}{G}
		\ \ \stackrel{L.\ref{lem_properties_of_Cu_Du_Iu} \mathbf{a)}}{=}\ \ 
\STinnerPRODUCT{\Du{U}(F)}{\Cu{U}(G)}
		\ \ \stackrel{L.\ref{lem_properties_of_Cu_Du_Iu} \mathbf{b)}}{=}\ \ 
\STinnerPRODUCT{ \hermTRANSPOSE{ \Cu{U}\left( \CONJ{F} \right) } }
                           { \hermTRANSPOSE{ \Du{U}\left( \CONJ{G} \right) } }
		\ \ \stackrel{see\ below}{=}\ \ 
\STinnerPRODUCT{ \hermTRANSPOSE{ \CONJ{ \Cu{U}(F) } } }
                           { \hermTRANSPOSE{ \CONJ{ \Du{U}(G) } } }
		\ \ =\ \ 
\STinnerPRODUCT{ \TRANSPOSE{ \Cu{U}(F) } }{ \TRANSPOSE{ \Du{U}(G) } }
	\ \ =\ \ 
\STinnerPRODUCT{  \Cu{U}(F)  }{  \Du{U}(G) }
	\ \ \stackrel{L.\ref{lem_properties_of_Cu_Du_Iu} \mathbf{a)}}{=}\ \ 
\UinnerPRODUCT{U}{F}{\Iu{U}(G)}$ \\
for any $F$ and $G$. 
The third equality holds because for real $U$\ \   $\Cu{U}$ and $\Du{U}$ map real arguments into real results,
so, for example,\ \ $\Cu{U}\left(\CONJ{X + \Ii Y}\right)\ =\ \Cu{U}\left(X - \Ii Y\right) \ =\ \Cu{U}(X) - \Ii \Cu{U}(Y)\ =\ 
\CONJ{\Cu{U}(X) + \Ii \Cu{U}(Y)} \ =\ \CONJ{\Cu{U}(X + \Ii Y)}$,\ \ where $X,Y$ are real.

We have thus shown that $\Iu{U}$ is a self adjoint operator wrt $\UinnerPRODUCT{U}{}{}$.
Hence all its eigenvalues are real and of modulus $1$ because $\Iu{U}$ is unitary 
(Theorem \ref{theor_unitarity_and_eigenspaces_of_Iu} {\bf a)}). 
The proportion of $1$'s to $(-1)$'s in the spectrum can be established using
 the trace of $\Iu{U}$, which is $N$ by Lemma \ref{lem_properties_of_Cu_Du_Iu} {\bf g)}.
The sum of $N^2$ eigenvalues must be equal to $N$ which is possible only when
the multiplicity of $1$ is equal to $N(N+1)/2$. This is  the undephased defect 
$\undephasedDEFECT(U)$, by Theorem \ref{theor_unitarity_and_eigenspaces_of_Iu} {\bf e)},
in the case when $U$ is a real orthogonal matrix.
\PROOFend    

Let us repeat. If\ \  $U\  =\  \ithMATRIX{U}{1} \otimes \ldots \otimes \ithMATRIX{U}{x}$\ \  is a Kronecker product of 
unitary $2 \times 2$  matrices with no zero entries, then it is equivalent to a Kronecker product of  real 
orthogonal $2 \times 2$ matrices with no zero entries,\ \  
	$O\  =\  \ithMATRIX{O}{1} \otimes \ldots \otimes \ithMATRIX{O}{x}$.\ \ 
The defects of $U$ and $O$ are equal (Corollary \ref{cor_DU_equal_DV_for_U_V_equivalent}), 
$\undephasedDEFECT(U)\ =\ \undephasedDEFECT(O)\ =\ 2^x\left( 2^x + 1 \right)/2$,\ \ 
which is the value from Lemma \ref{lem_Iu_selfadjoint_for_U_real_orthogonal}.
This value can be used when estimating (Corollary \ref{cor_V1_spaces_Kron_multiplied_are_contained_in_V1_space} {\bf c)})
the undephased defect of a Kronecker product of unitary matrices with a subproduct consisting of $x$\ \ $2 \times 2$ matrices. 
It is however possible to obtain the spectral decomposition of $\Iu{U}$ for $U$ being such a subproduct without resorting to real orthogonal matrices, 
because the spectral decomposition of $\Iu{V}$, when $V$ is of size $2$, is easy.


\begin{lemma}
	\label{lem_2x2_unitary_Iu_spectral_decomposition}
	Let $U\ \ =\ \ \left[ \begin{array}{cc}  a  &  b  \\  c  &  d  \end{array} \right]$
	be a unitary matrix with no zero entries. Then $\Iu{U}$ (which acts on $\COMPLEXmatricesNxN{2}$)
	has  only two eigenspaces, orthogonal wrt $\UinnerPRODUCT{U}{}{}$, corresponding to eigenvalues
	$1$ and $-1$\ :
	\begin{eqnarray}
		\label{eq_2x2_unitary_1_eigenspace}
		\Vu{U}{1}
		&    =    &
		\SPANC\left(
			\Ii \left[ \begin{array}{cc}  1 & 1 \\ 0 & 0  \end{array} \right]     ,\  
			\Ii \left[ \begin{array}{cc}  0 & 0 \\ 1 & 1  \end{array} \right]     ,\ 
			\frac{1}{\sqrt{2} \ABSOLUTEvalue{a} \ABSOLUTEvalue{b}}  \Ii
				\left[
					\begin{array}{cc}
						\ABSOLUTEvalue{b}^2  &  -\ABSOLUTEvalue{a}^2  \\
						\ABSOLUTEvalue{a}^2  &  -\ABSOLUTEvalue{b}^2
					\end{array}
				\right]
		\right)
					 \ =
		\TRIVIALspace{U}                      ,                \\
		\label{eq_2x2_unitary_minus1_eigenspace}
		\Vu{U}{-1}
		&     =     &
		\SPANC\left(
			\frac{\ABSOLUTEvalue{a} \ABSOLUTEvalue{b}}{\sqrt{2}}  \Ii
				\left[
					\begin{array}{cc}
						\frac{1}{\ABSOLUTEvalue{a}^2}     &     -\frac{1}{\ABSOLUTEvalue{b}^2}      \\
						-\frac{1}{\ABSOLUTEvalue{b}^2}    &     \frac{1}{\ABSOLUTEvalue{a}^2}
					\end{array}
				\right]
		\right)                               \ \ ,
	\end{eqnarray}
	where the bases provided are orthonormal wrt $\UinnerPRODUCT{U}{}{}$  and $\TRIVIALspace{U}$ is the space
	necessarily contained in $\Vu{U}{1}$ (Theorem \ref{theor_unitarity_and_eigenspaces_of_Iu} {\bf c)}).

	The real spaces spanned by the vectors in  (\ref{eq_2x2_unitary_1_eigenspace}) and (\ref{eq_2x2_unitary_minus1_eigenspace})
	are $\IMspaceOF{\Vu{U}{1}} = \IMspaceOF{\TRIVIALspace{U}}$\ \ and\ \ $\IMspaceOF{\Vu{U}{-1}}$,\ \ respectively.
\end{lemma}

\PROOFstart 
Orthonormality of all the basis vectors listed in (\ref{eq_2x2_unitary_1_eigenspace})   and  (\ref{eq_2x2_unitary_minus1_eigenspace})
is trivial. Since
\begin{equation}
	\label{eq_3rd_basis_vector_in_Tu}
	\left[
		\begin{array}{cc}
			\ABSOLUTEvalue{b}^2  &  -\ABSOLUTEvalue{a}^2  \\
			\ABSOLUTEvalue{a}^2  &  -\ABSOLUTEvalue{b}^2
		\end{array}
	\right]
			\ \ =\ \ 
	\left[ 
		\begin{array}{cc}  
			\ABSOLUTEvalue{b}^2  -  \ABSOLUTEvalue{a}^2       &     \ABSOLUTEvalue{b}^2  -  \ABSOLUTEvalue{a}^2     \\
			                                     0                                           &                                         0
		\end{array}
	\right]
		\ +\ 
	\left[
		\begin{array}{cc}
			\ABSOLUTEvalue{a}^2     &    0      \\
			\ABSOLUTEvalue{a}^2     &    0      
		\end{array}
	\right]
		\ +\ 
	\left[
		\begin{array}{cc}
			0    &    -\ABSOLUTEvalue{b}^2        \\
			0    &    -\ABSOLUTEvalue{b}^2
		\end{array}
	\right]                              \ \ ,
\end{equation}
the third basis vector in (\ref{eq_2x2_unitary_1_eigenspace}) belongs to $\TRIVIALspace{U}$, like the first two.

Having the $3$ dimensional eigenspace $\Vu{U}{1}\ =\ \TRIVIALspace{U}$ we have almost no choice for the fourth eigenvector
-- it must be $\UinnerPRODUCT{U}{}{}$-orthogonal to $\Vu{U}{1}$  (according to Theorem \ref{theor_unitarity_and_eigenspaces_of_Iu} {\bf a)}),
as the corresponding eigenvalue $-1$, different from $1$, is fixed by the value $N=2$ of the
trace of $\Iu{U}$ (Lemma \ref{lem_properties_of_Cu_Du_Iu} {\bf g)}). 
    We will show directly that the eigenvalue is equal to $-1$ though. We start from
calculating $\Cu{U}$ at a vector proportional to the one in (\ref{eq_2x2_unitary_minus1_eigenspace})\ :
\begin{eqnarray}
	\nonumber
	\Cu{U}\left(
		\left[
			\begin{array}{cc}
				\frac{1}{\ABSOLUTEvalue{a}^2}     &     -\frac{1}{\ABSOLUTEvalue{b}^2}      \\
				-\frac{1}{\ABSOLUTEvalue{b}^2}    &     \frac{1}{\ABSOLUTEvalue{a}^2}
			\end{array}
		\right]
	\right)
	&      =      &         
	\left[
		\begin{array}{cc}
			\frac{a}{\ABSOLUTEvalue{a}^2}     &     -\frac{b}{\ABSOLUTEvalue{b}^2}      \\
			-\frac{c}{\ABSOLUTEvalue{b}^2}    &     \frac{d}{\ABSOLUTEvalue{a}^2}
		\end{array}
	\right]
	\left[
		\begin{array}{cc}
			\CONJ{a}         &       \CONJ{c}      \\
			\CONJ{b}         &       \CONJ{d}    
		\end{array}
	\right]
			\ \ =\ \ 
	\left[
		\begin{array}{cc}
			0        &       \left(  \frac{a \CONJ{c}}{\ABSOLUTEvalue{a}^2}     -      \frac{b \CONJ{d}}{\ABSOLUTEvalue{b}^2}    \right)       \\
			\left(  -\frac{\CONJ{a} c}{\ABSOLUTEvalue{b}^2}     +      \frac{\CONJ{b} d}{\ABSOLUTEvalue{a}^2}    \right)    &         0
		\end{array}
	\right]                          \\
	\label{eq_Cu_on_Vu_minus1_eigenvector}
	&      =       &
	\left[
		\begin{array}{cc}
			0      &     a \CONJ{c} \left(  \frac{1}{\ABSOLUTEvalue{a}^2}   +    \frac{1}{\ABSOLUTEvalue{b}^2}   \right)      \\
			-\CONJ{a} c \left(  \frac{1}{\ABSOLUTEvalue{a}^2}   +    \frac{1}{\ABSOLUTEvalue{b}^2}   \right)     &      0
		\end{array}
	\right]                                \ \ ,
\end{eqnarray}
to get an antihermitian matrix. The above argument of $\Cu{U}$, call it $G$, is real, so from 
Lemma \ref{lem_properties_of_Cu_Du_Iu} {\bf c)} and the fact that $\Cu{U}(G)$ is antihermitian it follows that
$\Iu{U}(G)\ =\ -G$, that is the eigenvalue corresponding to $G$ is equal to $-1$.

Finally, the real spaces spanned by the vectors in  (\ref{eq_2x2_unitary_1_eigenspace}) and (\ref{eq_2x2_unitary_minus1_eigenspace})
are obviously contained in $\IMspaceOF{\Vu{U}{1}} = \IMspaceOF{\TRIVIALspace{U}}$\ \ and\ \ $\IMspaceOF{\Vu{U}{-1}}$,\ \ respectively.
Their dimensions are  $3$ and $1$,  because the vectors in (\ref{eq_2x2_unitary_1_eigenspace}) are independent also in the real sense,
and these dimensions are equal to\ \ 
$\DIMR\left(  \IMspaceOF{\Vu{U}{1}} \right)\ =\ \DIMC\left( \Vu{U}{1} \right)\ =\ 3$\ \  and\ \  
$\DIMR\left(  \IMspaceOF{\Vu{U}{-1}} \right)\ =\ \DIMC\left( \Vu{U}{-1} \right)\ =\ 1$\ \   
(see Theorem \ref{theor_unitarity_and_eigenspaces_of_Iu} {\bf b)}),  respectively.
So the considered real spaces are equal to  $\IMspaceOF{\Vu{U}{1}}$ and $\IMspaceOF{\Vu{U}{-1}}$.
\PROOFend 

Getting back to Kronecker products of unitary matrices of size $2$, we derive the below corollary from:
\begin{itemize}
	\item 
	the commment  following immediately  the proof of Lemma \ref{lem_Iu_for_Kron_prod_of_unitaries}
	($(F,\lambda),(G,\mu)$, eigenpairs for $\Iu{U},\Iu{V}$, generate an eigenpair $(F \otimes G, \lambda \mu)$ for $\Iu{U \otimes V}$),

	\item
	the above Lemma \ref{lem_2x2_unitary_Iu_spectral_decomposition} (on eigenvectors of $\Iu{U}$ for $U$ of size $2$), 

	\item
	Lemma \ref{lem_weigthed_U_inner_product_induced} 
	(on the induced inner product  $\UinnerPRODUCT{U \otimes V}{F' \otimes G'}{F'' \otimes G''} \ =\ 
	\UinnerPRODUCT{U}{F'}{F''} \cdot \UinnerPRODUCT{V}{G'}{G''}$),

	\item  
	Corollary  \ref{cor_Iu_Iv_operators_for_equivalent_U_V_are_similar} {\bf a)}
	($\Iu{U}$ and $\Iu{V}$ are similar if $U$ and $V$ are equivalent),

	\item
	Lemma \ref{lem_Iu_selfadjoint_for_U_real_orthogonal} ($\Iu{U}$ is self adjoint for a real orthogonal $U$ and has $N(N+1)/2$\ \ $1$'s 
	in the spectrum).
\end{itemize}



\begin{corollary}
	\label{cor_spectral_decomposition_of_Iu_for_kron_prod_of_2x2_matrices}
	If $U$ of size $N =2^x$ is a Kronecker product of $x$  unitary $2 \times 2$ matrices with nonzero entries:
	\begin{equation}
		\label{eq_U_a_kron_prod_of_2x2_matrices}
		U\ \ =\ \ \ithMATRIX{U}{1} \otimes \ldots \otimes \ithMATRIX{U}{x}           \ \ ,
	\end{equation}
	then $\Iu{U}$ has only two eigenspaces $\Vu{U}{1}$ and $\Vu{U}{-1}$, where 
	$\COMPLEXmatricesNxN{N}\ =\ \Vu{U}{1}  \oplus\  \Vu{U}{-1}$, orthogonal wrt $\UinnerPRODUCT{U}{}{}$.
	An eigenbasis for $\Iu{U}$ orthonormal wrt $\UinnerPRODUCT{U}{}{}$ can be formed from Kronecker products
	\begin{equation}
		\label{eq_eigenvectors_for_U_a_kron_prod_of_2x2_matrices}
		\ithMATRIX{G}{1} \otimes  \ldots  \otimes \ithMATRIX{G}{x}                    \ \ ,
	\end{equation}
	where $\ithMATRIX{G}{k}$ is one of the four eigenvectors of $\Iu{\ithMATRIX{U}{k}}$ listed in 
	Lemma \ref{lem_2x2_unitary_Iu_spectral_decomposition}.  Denoting by $\lambda_{\ithMATRIX{G}{k}}$
	the corresponding eigenvalue, the multiplicity of $1$ in the spectrum of $\Iu{U}$ is the number
	of possible sequences $\left(  \lambda_{\ithMATRIX{G}{1}},\ \ldots,\ \lambda_{\ithMATRIX{G}{x}}   \right)$
	with an even number of $-1$'s. This multiplicity (i.e. $\undephasedDEFECT(U)$) must be equal to $2^x \left(2^x + 1 \right)/2$ because,
	as $U$ is equivalent to some real orthogonal matrix $O$, the operators $\Iu{U}$ and $\Iu{O}$ are similar.
\end{corollary}

It is clear now that it pays to analyze the whole of the spectrum  of $\Iu{U}$ when $U$ plays the role   
of a Kronecker factor of some larger unitary $V$ the defect of which needs to be calculated.

%
%

\section{Fourier matrices, symmetries and equivalence}
\label{sec_Fourier_matrices}

Let us introduce an alternative way of indexing rows and columns of an $N \times N$  complex matrix, namely one involving
the elements of $\ZN{N}\ =\ \{ \ZNclass{0},...,\ZNclass{N-1} \}\ =\ \INTEGER/ N \INTEGER$.  Let $i \ = \ZNclass{\REPR{i}} \in \ZN{N}$ designate
the ordinarily indexed $(\REPR{i}+1)$-th row (or column), where $\REPR{i} \in \{0,1,...,N-1\}$ represents
class $i \in \ZN{N}$. 

We shall adopt the following definition
of the $N \times N$ Fourier matrix $\FOURIER{N}$, which is the matrix of the discrete Fourier transform of an $N$ element 
complex vector:
\begin{equation}
	\label{eq_Fourier_matrix}
	\ELEMENTof{\FOURIER{N}}{i}{j} \ \ =\ \ \PHASE{\frac{2\pi}{N} \REPR{i} \REPR{j}},
\end{equation}
where $i,j \in \ZN{N}$ and $\REPR{i},\REPR{j}$ are arbitrary representatives of
$i,j$, respectively. $\ZN{n}$ in this context will be called the \DEFINED{indexing group of $\FOURIER{N}$} and denoted
$\Igroup{\FOURIER{N}}$.

To generalize the above concepts we consider the Kronecker product:
\begin{equation}
	\label{eq_Fourier_Kron_prod}
	F\ \ =\ \ \FOURIER{N_1} \otimes \FOURIER{N_2} \otimes \ldots \otimes \FOURIER{N_r}
\end{equation}
of size $N\ =\ N_1 N_2 \ldots N_r$. 
Such a matrix will be indexed alternatively with the elements of $\ZN{N_1} \times \ldots \times \ZN{N_r}$, called group indices. 
The new way of indexing, mandatory further, can be defined
in three equivalent ways presented in the list below. 
In the description we use $i,j \in \ZN{N_1} \times \ldots \times \ZN{N_r}$ and this group will further be called the
\DEFINED{indexing group of $F$} and denoted $\Igroup{F}$. 
$i,j$ thus decompose $i\ =\ \vectorINDEX{i_1,\ldots,i_r}$, $j\ =\ \vectorINDEX{j_1,\ldots,j_r}$, where of course $i_x,j_x \in \ZN{N_x}$.
Also $k,l \in \{1,\ldots,N\}$ are used in the role of ordinary indices corresponding to $i,j$.
\begin{itemize}
	\item
		Pair $i,j$ designates that $k,l$-th element of $F$ which is built as
		$\ELEMENTof{\FOURIER{N_1}}{i_1}{j_1} \cdot \ldots  \cdot \ELEMENTof{\FOURIER{N_r}}{i_r}{j_r}$
		in the Kronecker product (\ref{eq_Fourier_Kron_prod}) 
		(note that we index $\FOURIER{N_x}$ with the elements of its indexing group 
		$\Igroup{\FOURIER{N_x}}$ as described at the beginning of this section).

	\item
		$i,j$ correspond to those ordinary indices $k,l$ which are the positions of $i,j$, respectively, 
		in the indexing group $\Igroup{F}$ 
		ordered lexicographically: $\vectorINDEX{\ZNclass{0},\ldots,\ZNclass{0},\ZNclass{0}}$,
										$\vectorINDEX{\ZNclass{0},\ldots,\ZNclass{0},\ZNclass{1}}$,
										...,
										$\vectorINDEX{\ZNclass{0},\ldots,\ZNclass{0},\ZNclass{N_r - 1}}$,
										$\vectorINDEX{\ZNclass{0},\ldots,\ZNclass{1},\ZNclass{0}}$,
										...,
										$\vectorINDEX{\ZNclass{0},\ldots,\ZNclass{1},\ZNclass{N_r - 1}}$,
										...,
										$\vectorINDEX{\ZNclass{N_1 - 1},\ldots,\ZNclass{N_{r-1} - 1},\ZNclass{N_r - 1}}$. 

	\item
		the ordinary indices to which $i,j$ refer are given by the formulas:  
		\begin{eqnarray}
			\label{eq_ordinary_and_group_indices_relation}
			k  &  =  &
			\REPR{i_1} \left(N_2 \ldots N_r\right)\ +\ 
			\REPR{i_2} \left(N_3 \ldots N_r\right)\ +\ 
			\ldots\ +\ 
			\REPR{i_{r-1}} N_r\ +\ 
			\REPR{i_r}\ +\ 
			1,                         \\
			\nonumber
			l  &  =  &
			\REPR{j_1} \left(N_2 \ldots N_r\right)\ +\ 
			\REPR{j_2} \left(N_3 \ldots N_r\right)\ +\ 
			\ldots\ +\ 
			\REPR{j_{r-1}} N_r\ +\ 
			\REPR{j_r}\ +\ 
			1,         
		\end{eqnarray}
		where $\REPR{i_x}, \REPR{j_x} \in \{0,\ldots,N_x - 1\}$ are  representatives
		of classes $i_x, j_x \in \ZN{N_x}$, respectively.
\end{itemize}

Obviously, in the special case $F\ =\ \FOURIER{N}$ the indexing group is equal to $\ZN{N}$ and the above described way of indexing
with the use of group indices is just the way introduced at the beginning of this section. For example, for $F\ =\ \FOURIER{6}$ group
indices $\ZNclass{0}$, $\ZNclass{1}$, ..., $\ZNclass{5}$ correspond to ordinary indices $1$, $2$, ..., $6$, respectively.
In a more complicated example of 
$F\ =\ \FOURIER{2} \otimes \FOURIER{3}$ indexed by $\Igroup{F}\ =\ \ZN{2} \times \ZN{3}$, group indices 
$\vectorINDEX{\ZNclass{0},\ZNclass{0}}$, 
$\vectorINDEX{\ZNclass{0},\ZNclass{1}}$,
$\vectorINDEX{\ZNclass{0},\ZNclass{2}}$, 
$\vectorINDEX{\ZNclass{1},\ZNclass{0}}$,
$\vectorINDEX{\ZNclass{1},\ZNclass{1}}$, 
$\vectorINDEX{\ZNclass{1},\ZNclass{2}}$
should be interpreted as ordinary indices $1$, $2$, ..., $6$, respectively (note the lexicographic order of the group indices).

Let $F$ be the Kronecker product given in (\ref{eq_Fourier_Kron_prod}), of size $N$, indexed by the appropriate indexing group 
$\Igroup{F}$. Such matrices will further be called \DEFINED{Fourier matrices}. It is not difficult to show that for any
$i,j \in \Igroup{F}$:
\begin{equation}
	\label{eq_multiplying_rows_of_F}
	F_{i+j, :} \ \ =\ \ F_{i,:} \HADprod F_{j,:}     \ \ ,
\end{equation}
where on the right hand side we have the entrywise (Hadamard) product of the $i$-th and $j$-th row of $F$. Thus the map
$i \longrightarrow F_{i,:}$ is a monomorphism ($F$ has orthogonal hence distinct rows) from $\Igroup{F}$ into
the group of $N$ element rows with unimodular entries, with entrywise product. So the image of this monomorphism, 
the rows of $F$, forms a group isomorphic to $\Igroup{F}$, which  is called the character group of $\Igroup{F}$ in the representation theory. 
All that has just been said also applies to columns because $F$ is symmetric
($F_{x,y}\ =\ F_{y,x}$, both respect to ordinary and group indices which is a consequence of our definition of the Fourier matrices). 
Property (\ref{eq_multiplying_rows_of_F}) restricted to columns expresses the fact that each $j$th column (and symmetrically: each $i$th row), as
a function of the group row (column) index\ \  $i \longrightarrow F_{i,j}$\ \ ($j \longrightarrow F_{i,j}$),\ \  
is a homomorphism from $\Igroup{F}$ into the group of unimodular complex numbers
(in fact into the group of the roots of unity of order $\lcm\left(N_1,\ldots,N_r\right)$). 
As it is shown in the representation theory, columns of $F$ give us all such
homomorphism on $\Igroup{F}$. Since most of this section is based on that fact, we prove it below:


\begin{lemma}
	\label{lem_all_unimodular_morphisms_on_F_indexing_group}
	Let $F$ be a Fourier matrix of size $N$ indexed by group $\Igroup{F}$.\\
	If\ \  $\uniMORPHISM{\phi}:\ \left( \Igroup{F}, + \right) \longrightarrow \left( \COMPLEX \setminus \{0\}, \cdot \right)$ is a homomorphism, 
	i.e. $\uniMORPHISM{\phi}(i+j)\ =\ 
                 \uniMORPHISM{\phi}(i) \cdot \uniMORPHISM{\phi}(j)$, 
	then for some column index $j \in \Igroup{F}$ we have that\ \ 
	$\forall i \in \Igroup{F}\ \ \ \uniMORPHISM{\phi}(i)\ =\ F_{i,j}$.
\end{lemma}

\PROOFstart 
Let $x$ be an $N$ element column vector with entries (we use group indices related to $F$) $x_i \ = \uniMORPHISM{\phi}(i)$.
x can be expressed in a basis formed by the columns of $F$ (where $1/\sqrt{N} \cdot F$ is unitary): $x\ =\ F a$, so
$a\ =\ (1/N) \cdot F^* x\ =\ (1/N) \cdot \CONJ{F} x$ because $F\ =\ F^T$. 

Now, for $i \in \Igroup{F}$, using $\CONJ{F_{i,:}}\ =\ F_{-i,:}$ for rows, where the conjugation is the inverse operation in the
group of row of $F$ isomorphic to $\Igroup{F}$, let us calculate the square of the $i$-th component of $a$:
\begin{eqnarray}
	\label{eq_ai_squared}
	\lefteqn{a_i^2\ \ =}  &   &    \\
	\nonumber
	&  &
	\left( \frac{1}{N} \CONJ{F}_{i,:} x \right)^2 \ \ =\ \ 
	\left( \frac{1}{N} F_{-i,:} x \right)^2 \ \ = \ \
	\frac{1}{N^2} x^T \left(F_{-i,:} \right)^T  F_{-i,:}\,  x \ \ =    \\
	\nonumber
	&  &
	\frac{1}{N^2}  
	\sum_{k,l \in \Igroup{F}}  
		F_{-i,k} F_{-i,l}\,  \uniMORPHISM{\phi}(k) \uniMORPHISM{\phi}(l)
	\ \ =      \\
	\nonumber
	&  &
	\frac{1}{N^2}  
	\sum_{k,l \in \Igroup{F}} 
		F_{-i,k+l}\, \uniMORPHISM{\phi}(k+l)
	\ \ =\ \ 
	\frac{1}{N^2}  
	\sum_{k \in \Igroup{F}}
		\sum_{l \in \Igroup{F}}
			F_{-i,k+l}\, \uniMORPHISM{\phi}(k+l)
	\ \ =       \\
	\nonumber
	&  &
	\frac{1}{N^2}  
	\sum_{k \in \Igroup{F}}
		\sum_{m \in \Igroup{F}}
			F_{-i,m}\, \uniMORPHISM{\phi}(m)
	\ \ =\ \ 
	\frac{1}{N}  
	\sum_{k \in \Igroup{F}}
		\frac{1}{N} \CONJ{F}_{i,:} x
	\ \ =\ \ 
	\frac{1}{N}  
	\sum_{k \in \Igroup{F}}
		a_i
	\ \ =\ \ 
	a_i.
\end{eqnarray}
$a_i$ is thus equal either to $0$ or $1$. Let there be $n$ $1$'s in vector $a$. Then for $\GROUPzero{\Igroup{F}}$, the neutral element of 
$\Igroup{F}$,\ \ 
$\uniMORPHISM{\phi}\left(\GROUPzero{\Igroup{F}}\right)
	\ =\ 
x_{\GROUPzero{\Igroup{F}}}
	\ =\ 
F_{\GROUPzero{\Igroup{F}},:}\, a\ =\ n$\ \ 
 because $F_{\GROUPzero{\Igroup{F}},:}$ is an all $1$'s row, 
it is the neutral element of the group of rows of $F$. But $\uniMORPHISM{\phi}\left(\GROUPzero{\Igroup{F}}\right)\ =\ 1$,
 as $\uniMORPHISM{\phi}$ is a homomorphism,
so there is only one $1$ in $a$ which means that $x$ is a column of $F$.
\PROOFend 

We will occasionally use the group indices associated with $F$ to designate the elements of matrices other than $F$, but of the same size.
Special attention should be paid to the permutation matrices. 
Let $p:\{1 \ldots N\} \longrightarrow \{1 \ldots N\}$ be a permutation. By $\Pro{p}$ we understand an $N \times N$ permutation matrix
with $1$'s at the ordinarily indexed positions\ \ $k,p(k)$.\ \  Then, as one can easily check,\ \  $\Pro{p}^{-1}\ =\ \Pro{p^{-1}}$\ \ and\ \ 
$\Pro{p'} \Pro{p''} \ =\  \Pro{p'' p'}$. Further, let $c_F: \Igroup{F} \longrightarrow \{1 \ldots N\}$ be the bijection connecting 
(in the way described at the beginning of this section) the ordinary and 
group indices related to an $N \times N$ Fourier matrix $F$. Next, let for Fourier matrices $F,G$ of size $N$ and bijection\ \ 
$\varphi: \Igroup{F} \longrightarrow \Igroup{G}$\ \ $\Pro{\varphi}$ be defined as $\Pro{c_G \varphi c_F^{-1}}$. 
It has $1$'s at the positions\ \ $k, c_G \varphi c_F^{-1}(k)$,\ \  where $k \ =\ c_F(i)$ for some $i \in \Igroup{F}$, so there are $1$'s
at the positions\ \ $c_F(i), c_G(\varphi(i))$,\ \ or, speaking in terms of $\Igroup{F}$ for the rows and $\Igroup{G}$ for the columns, at the
positions\ \ $i, \varphi(i)$.\ \ Matrices $\Pro{\varphi}$ also satisfy the rules:
\begin{eqnarray}
	\label{eq_Pro_inversion_rule}
	\Pro{\varphi}^{-1} 
		\ \ =\ \ 
	\Pro{c_G \varphi c_F^{-1}}^{-1}
		\ \ =\ \
	\Pro{\left(c_G \varphi c_F^{-1}\right)^{-1}}
		\ \ =\ \
	\Pro{c_F \varphi^{-1} c_G^{-1}}
		\ \ =\ \ 
	\Pro{\varphi^{-1}}    \ \ ,
	   & &   \\
	\label{eq_Pro_multiplication_rule}
	\Pro{\varphi} \Pro{\psi} 
		\ \ =\ \ 
	\Pro{c_G \varphi c_F^{-1}}  \Pro{c_H \psi c_G^{-1}}
		\ \ =\ \ 
	\Pro{\left(c_H \psi c_G^{-1}\right)\left(c_G \varphi c_F^{-1}\right)}
		\ \ =\ \ 
	\Pro{c_H \psi \varphi c_F^{-1}}
		\ \ =\ \ 
	\Pro{\psi \varphi}     \ \ ,
	&  &  
\end{eqnarray}
where\ \ $\psi: \Igroup{G} \longrightarrow \Igroup{H}$\ \ is a bijection and $H$ is a Fourier matrix of the same size as $F,G$.

In the rest of this section we will use the objects introduced in Definition \ref{def_pairs_action_stabilizer}, 
namely $\PP$, the group of pairs of enphased permutation matrices of size $N$,
its subgroup $\PsPs$ consisting of all the pairs of permutation matrices,
and the action\ \ $(S,T)X\ =\ S X T^{-1}$\ \ of $\PP$ on the set of all complex matrices of size $N$.
We will be interested in a description of the stabilizers $\STAB{\PP}{F}$, $\STAB{\PsPs}{F}$ of a Fourier matrix $F$, 
as well as a description of $\MAP{\PP}{F}{G}$, $\MAP{\PsPs}{F}{G}$ for any Fourier matrices $F,G$ of the same size $N$.
In other words we ask which equivalence (permutation equivalence -- to be defined in a moment -- when $\PsPs$ is considered)
operations map $F$ back into $F$ and which, if any, map $F$ into $G$.
The equivalence has been defined in (\ref{eq_equivalence}) and here we apply it to rescaled unitary matrices $F,G$, which
are obviously equivalent if and only if $F\ =\ (S,T)G$ for some $(S,T) \in \PP$. Besides, we need:

\begin{definition}
	\label{def_isom_object}
	$\ISOMORPHISMS{F}{G}$ is the set  of all isomorphisms 
	$\ISOM{\chi}:\ \Igroup{F} \longrightarrow \Igroup{G}$ between the indexing groups of
	Fourier matrices $F$ and $G$ of the same size\\
	($\ISOMORPHISMS{F}{F}$ is a group. $\ISOMORPHISMS{F}{G}$ is nonempty if and only if $\Igroup{F}, \Igroup{G}$
	are isomorphic -- 
	which is not the case e.g. for $\Igroup{\FOURIER{4}} = \ZN{4}$, $\Igroup{\FOURIER{2} \otimes \FOURIER{2}} = \ZN{2} \times \ZN{2}$
	-- and then $\ISOMORPHISMS{F}{G} \ =\ \ISOM{\chi} \ISOMORPHISMS{F}{F} \ =\ \ISOMORPHISMS{G}{G} \ISOM{\chi}$
	for any $\ISOM{\chi} \in \ISOMORPHISMS{F}{G}$).
\end{definition}

The two lemmas and the corollary below address the issue of permutation equivalence of Fourier matrices. We regard matrices $X$ and $Y$ as 
\DEFINED{permutation equivalent} if $Y\ =\ (P', P'')X$ for some pair $(P',P'')$ of permutation matrices.
This is an equivalence relation in the usual sense.  
In the proofs below we use the properties (\ref{eq_Pro_inversion_rule}), (\ref{eq_Pro_multiplication_rule}).


\begin{lemma}
	\label{lem_permutation_equivalence_PFRequalsG_matrices_determined_by_isomorphisms}
	Let $F$, $G$ be Fourier matrices of the same size.
	Let $\chi',\ \chi'' :\ \Igroup{G} \longrightarrow \Igroup{F}$ be bijections.
	If for permutation matrices $\Pro{\chi'}$, $\Pro{\chi''}$ there holds:
	\begin{equation}
		\label{eq_perm_equivalence_with_some_matrices}
		\Pro{\chi'}  F  \Pro{\chi''}^T \ \ =\ \ G,
	\end{equation}
	then $\chi'$, $\chi''$ are isomorphisms between the indexing  groups $\Igroup{F}$ and $\Igroup{G}$.
\end{lemma}

\PROOFstart 
If (\ref{eq_perm_equivalence_with_some_matrices}) holds then the $j$-th column of $G$ ($j \in \Igroup{G}$) is the $\chi''(j)$-th column
of $F$ permuted with the use of $\Pro{\chi'}$:
\begin{equation}
	\label{eq_jth_column_of_G_as_a_permuted_column_of_F}
	G_{:,j}\ \ =\ \ 
	\ELEMENTof{ \Pro{\chi'} F }{:}{\chi''(j)} \ \ =\ \ 
	\Pro{\chi'}  F_{:,\chi''(j)}\ \ .
\end{equation}
Because $G_{i,j}$ is a homomorphism with respect to both $i$ and $j$ 
($i,j \in \Igroup{G}$, see also the comment and the lemma after (\ref{eq_multiplying_rows_of_F})), 
we have for the entrywise product of the $j_1$-th and $j_2$-th column of $G$ that:
\begin{equation}
	\label{eq_entrywise_product_of_columns_of_G}
	G_{:,j_1} \HADprod G_{:,j_2} \ \ =\ \ 
	G_{:,j_1 + j_2} \ \ =\ \ 
	\Pro{\chi'}  F_{:,\chi''(j_1 + j_2)}\ \ .
\end{equation}
The left hand side of the above is also equal to:
\begin{equation}
	\label{eq_entrywise_product_of_columns_of_G_otherwise}
	\left( \Pro{\chi'}  F_{:,\chi''(j_1)} \right)  \HADprod  \left( \Pro{\chi'}  F_{:,\chi''(j_2)}  \right)
	\ \ =\ \ 
	\Pro{\chi'}  \left( F_{:,\chi''(j_1)}   \HADprod  F_{:,\chi''(j_2)}   \right)
	\ \ =\ \ 
	\Pro{\chi'}  F_{:,\chi''(j_1) + \chi''(j_2)}
\end{equation}
as $F_{k,l}$ is also a homomorphism with respect to both $k$ and $l$   ($k,l \in \Igroup{F}$).
Since the left hand sides of 
(\ref{eq_entrywise_product_of_columns_of_G}) and 
(\ref{eq_entrywise_product_of_columns_of_G_otherwise}) 
are equal, we conclude that
\begin{equation}
	\label{eq_equality_of_columns_of_F_concluded}
	\Pro{\chi'}  F_{:,\chi''(j_1 + j_2)}\ \ =\ \   \Pro{\chi'}  F_{:,\chi''(j_1) + \chi''(j_2)}\ \ ,
\end{equation}
hence $\chi''(j_1 + j_2)\ =\ \chi''(j_1) + \chi''(j_2)$ because all the columns of $F$ are different. 
As $j_1$ and $j_2$ are arbitrary, we have shown that $\chi''$ is a homomorphism, apart from being a bijection.

Since, after transposition,  (\ref{eq_perm_equivalence_with_some_matrices})
converts to 
				$  \Pro{\chi''}  F  \Pro{\chi'}^T \ \ =\ \ G  $,
the above reasoning can be repeated for $\chi'$ to show that it also is an isomorphism.
\PROOFend  


\begin{lemma}
	\label{lem_P_determines_R_for_permutation_equivalence_PFRequalsG}
	Let $F$, $G$ be Fourier matrices indexed by isomorphic groups $\Igroup{F}$ and $\Igroup{G}$, and let
	$\ISOM{\chi'}:\ \Igroup{G} \longrightarrow \Igroup{F}$ be an  isomorphism between these groups.
	Then there exists the only bijection $\ISOM{\chi''}:\ \Igroup{G} \longrightarrow \Igroup{F}$,
	being an isomorphism, such that:
	\begin{equation}
		\label{eq_perm_equivalence_with_matrix_and_its_partner}
		\Pro{\chi'}  F  \Pro{\chi''}^T \ \ =\ \ G,
	\end{equation}
	which is equal to $\INVERSE{\hMAP{F}{G}{\ISOM{\chi'}}}$, 
	where $\hMAP{F}{G}{\ISOM{\chi'}}:\ \Igroup{F} \longrightarrow \Igroup{G}$ is the only bijection between
	these indexing groups such that:
	\begin{equation}
		\label{eq_h_defining property}
		F_{  \ISOM{\chi'}(i),\   j }\ \ =\ \ G_{i,\   \hMAP{F}{G}{\ISOM{\chi'}}(j) }
		\ \ \ \ 
		\mbox{for any}\ \ i \in \Igroup{G},\ \ j \in \Igroup{F}.
	\end{equation}
	By symmetry, $\ISOM{\chi'}\ =\ \INVERSE{ \hMAP{F}{G}{\ISOM{\chi''}}}$.
\end{lemma} 

\PROOFstart 
For $i \in \Igroup{G}$ the $i$-th row of $\Pro{\ISOM{\chi'}} F$ is the $\ISOM{\chi'}(i)$-th row of $F$:
\begin{equation}
	\label{eq_ith_row_of_PF}
	\ELEMENTof{ \Pro{\ISOM{\chi'}} F }{i}{:}
	\ \ =\ \ 
	F_{\ISOM{\chi'}(i),\ :}\ \ .
\end{equation}
Because $\ISOM{\chi'}:\ \Igroup{G} \longrightarrow \Igroup{F}$ is an isomorphism and $F_{k,l}$ is  a 
homomorphism with respect to both $k$ and $l$ ($k,l \in \Igroup{F}$), there holds:
\begin{eqnarray}
	\label{eq_i1_plus_i2_th_row_of_PF}
	\lefteqn{ \ELEMENTof{ \Pro{\ISOM{\chi'}} F }{i_1 + i_2}{:}\ \ =}
	&  &    \\
	\nonumber
	&  &
	F_{\ISOM{\chi'}(i_1 + i_2),\ :} \ \ =\ \ 
	F_{\ISOM{\chi'}(i_1) + \ISOM{\chi'}(i_2),\ :}\ \ =\ \ 
	F_{\ISOM{\chi'}(i_1),\ :}  \HADprod F_{\ISOM{\chi'}(i_2),\ :}\ \ =   \\
	\nonumber
	&  &
	\ELEMENTof{ \Pro{\ISOM{\chi'}} F }{i_1}{:}   \HADprod   \ELEMENTof{ \Pro{\ISOM{\chi'}} F }{i_2}{:}\ \ \ .
\end{eqnarray}
The equality between the leftmost and rightmost expressions means that every $j$th column of  $\Pro{\ISOM{\chi'}} F$  ($j \in \Igroup{F}$), 
treated as map $i \longrightarrow \ELEMENTof{ \Pro{\ISOM{\chi'}} F }{i}{j}$, is a homomorphism. 
From Lemma \ref{lem_all_unimodular_morphisms_on_F_indexing_group} we know that, because of the fact just mentioned,
all the columns of $\Pro{\ISOM{\chi'}} F$ must sit in $G$: the above function of $i$ is equal to $i \longrightarrow G_{i,l}$ for
some $l \in \Igroup{G}$ corresponding to $j$. 
Since both $\Pro{\ISOM{\chi'}} F$ and $G$ have all their columns different and an equal number of them, 
there is only one way to permute the columns of  $\Pro{\ISOM{\chi'}} F$ to convert it into $G$. That is to say, there exists
only one bijection $\hMAP{F}{G}{\ISOM{\chi'}}: \Igroup{F} \longrightarrow \Igroup{G}$ such that:
\begin{equation}
	\label{eq_jth_PF_column_as_some_G_column}
	\ELEMENTof{ \Pro{\ISOM{\chi'}} F }{:}{j} 
	\ \ =\ \ 
	G_{:, \hMAP{F}{G}{\ISOM{\chi'}}(j)} 
	\ \ \ \ 
	\mbox{for any}\ \   j \in \Igroup{F}\ \ ,
\end{equation}
or equivalently:
\begin{equation}
	\label{eq_PF_columns_as_G_columns}
	\Pro{\ISOM{\chi'}} F \ \ =\ \ G   \left( \Pro{\hMAP{F}{G}{\ISOM{\chi'}}}  \right)^T\ \ .
\end{equation}
Note that the above condition, in the form (\ref{eq_jth_PF_column_as_some_G_column}), 
expressed entrywise and combined with the general property (\ref{eq_ith_row_of_PF}) takes the
form (\ref{eq_h_defining property}) used in the lemma. 

Yet another shape of the condition will be obtained by trasforming equality (\ref{eq_PF_columns_as_G_columns}) :
\begin{equation}
	\label{eq_PF_columns_as_G_columns_transformed}
	\Pro{\ISOM{\chi'}}  F  \Pro{\hMAP{F}{G}{\ISOM{\chi'}}}
	\ \ =\ \ 
	\Pro{\ISOM{\chi'}}  F  \Pro{\INVERSE{\hMAP{F}{G}{\ISOM{\chi'}}}}^T
	\ \ =\ \ 
	G
\end{equation}
From Lemma \ref{lem_permutation_equivalence_PFRequalsG_matrices_determined_by_isomorphisms}
the bijection $\ISOM{\chi''}\ =\ \INVERSE{\hMAP{F}{G}{\ISOM{\chi'}}}$ determining the right 
transposed permutation matrix in the second product in (\ref{eq_PF_columns_as_G_columns_transformed})
is an isomorphism. It is unique because  $\hMAP{F}{G}{\ISOM{\chi'}}$ is unique.

Also, as\ \  $\Pro{\ISOM{\chi'}} F \Pro{\ISOM{\chi''}}^T\ =\ G$\ \  is equivalent, by transposition,
to\ \ $\Pro{\ISOM{\chi''}} F \Pro{\ISOM{\chi'}}^T\ =\ G$,\ \
$\ISOM{\chi'}$ is determined by $\ISOM{\chi''}$ precisely in the same way as 
$\ISOM{\chi''}$ is determined by $\ISOM{\chi'}$.
\PROOFend 
	
In our older article \cite{PermEqClasses} in Theorem 3.1 we described the permutation equivalence classes of Fourier matrices, 
called there the Kronecker products of Fourier matrices, without explicitly resorting to the underlying indexing groups.
The lemma below explains this classification.


\begin{corollary}
	\label{cor_Fouriers_perm_equivalent_iff_their_indexing_groups_isomorphic}
	Fourier matrices $F$ and $G$, indexed by the groups $\Igroup{F}$ and $\Igroup{G}$, are permutation equivalent, i.e.
	\begin{equation}
		\label{eq_perm_equivalence_of_Fouriers}
		\Pro{\chi'}  F  \Pro{\chi''}^T \ \ =\ \ G
	\ \ \ \ 
	\mbox{for some bijections}\ \ \chi', \chi'' : \Igroup{G} \longrightarrow \Igroup{F}
	\end{equation}
	if and only if groups $\Igroup{F}$ and $\Igroup{G}$ are isomorphic.
\end{corollary}

\PROOFstart 
By Lemma \ref{lem_permutation_equivalence_PFRequalsG_matrices_determined_by_isomorphisms},
(\ref{eq_perm_equivalence_of_Fouriers}) implies that $\chi':\ \Igroup{G} \longrightarrow \Igroup{F}$ is an isomorphism.
By Lemma \ref{lem_P_determines_R_for_permutation_equivalence_PFRequalsG}, 
if $\chi'$ is an isomorphism there exists  isomorphism $\chi'':\ \Igroup{G} \longrightarrow \Igroup{F}$ 
such that (\ref{eq_perm_equivalence_of_Fouriers}) holds.
\PROOFend  

At this point we remind the reader that any indexing group 
$\Igroup{F}\ =\ \ZN{N_1} \times \ZN{N_2} \times \ldots \times \ZN{N_r}$, as a finite abelian group, 
 is isomorphic to the unique direct product
\begin{equation}
	\label{eq_primary_ordered_direct_product_of_ZNs}
	\ZNprimePOWER{a}{1}{k}{1}{1} \times \ldots \times \ZNprimePOWER{a}{1}{k}{n_1}{1}
	\ \ \times \ \ \ldots\ \ \times\ \ 
	\ZNprimePOWER{a}{s}{k}{1}{s} \times \ldots \times \ZNprimePOWER{a}{s}{k}{n_s}{s}\ \ ,
\end{equation}
where $a_1 < \ldots < a_s$ are prime numbers and for any\ \  $t\ =\ 1 \ldots s$\ \  numbers
$k_1^{(t)} \leq \ldots \leq k_{n_t}^{(t)}$ are natural. 
The unique product (\ref{eq_primary_ordered_direct_product_of_ZNs}) corresponding to $\Igroup{F}$
is obtained from the original product $\Igroup{F}$ by splitting the orders of the factors into products of
relatively prime numbers, and finally by reordering the resulting direct product. 
For example we transform $\ZN{12} \times \ZN{6}$ into its isomorphic counterpart in this way:
$\ZN{12} \times \ZN{6}\ \ \rightarrow\ \ 
\ZN{4} \times \ZN{3} \times \ZN{6} \ \ \rightarrow\ \ 
\ZN{4} \times \ZN{3} \times \ZN{3} \times \ZN{2}\ \ \rightarrow\ \ 
\ZN{2} \times \ZN{4} \times \ZN{3} \times \ZN{3}$.

Two different direct products of the form (\ref{eq_primary_ordered_direct_product_of_ZNs}) are not isomorphic.
(The reason is that even if they both contain elements of a given order $x$, $x$ can be chosen so that there are different numbers
of these elements in both groups. See our paper \cite{PermEqClasses} where we count, though rather ineffectively, 
the number of rows of a given type, i.e. of a given order, in the group of rows of a Fourier matrix $F$ indexed by product 
(\ref{eq_primary_ordered_direct_product_of_ZNs}) being isomorphic to this group of rows. Nevertheless, we are talking here 
about a classical result in the group theory.)
For example $\ZN{8}$, $\ZN{2} \times \ZN{4}$ and $\ZN{2} \times \ZN{2} \times \ZN{2}$  are
not isomorphic, so the Fourier matrices indexed by these groups, $\FOURIER{8}$, $\FOURIER{2} \otimes \FOURIER{4}$
and $\FOURIER{2} \otimes \FOURIER{2} \otimes \FOURIER{2}$ are not permutation equivalent 
by Corollary \ref{cor_Fouriers_perm_equivalent_iff_their_indexing_groups_isomorphic}. 
Let another  example concern not isomorphic groups $\ZN{2} \times \ZN{6}\ \equiv\ \ZN{2} \times \ZN{2} \times \ZN{3}$ 
and $\ZN{12}\ \equiv\ \ZN{4} \times \ZN{3}$, leading  to
$\FOURIER{2} \otimes \FOURIER{6}$ and $\FOURIER{12}$ which are not permutation equivalent,
while $\FOURIER{2} \otimes \FOURIER{6}$ and $\FOURIER{2} \otimes \FOURIER{2} \otimes \FOURIER{3}$ are permutation equivalent,
and the same for the pair $\FOURIER{12},\ \FOURIER{4} \otimes \FOURIER{3}$. 
We have used $\equiv$ to denote an isomorphism between the groups.

The map $\ISOM{\chi} \longrightarrow \hMAP{F}{G}{\ISOM{\chi}}$ 
introduced in Lemma \ref{lem_P_determines_R_for_permutation_equivalence_PFRequalsG}   
 has a number of properties, used in calculations that come further.
These properties are collected in Lemma \ref{lem_h_FG_map_properties}. The corollary that follows deals with the special case when $F\ =\ G$.


\begin{lemma}
	\label{lem_h_FG_map_properties}
	Let the indexing groups $\Igroup{F}$, $\Igroup{G}$, $\Igroup{H}$ of Fourier matrices $F$, $G$, $H$ be isomorphic.
	Then:
	\begin{description}
		\item[(a)]
			$\hMAP{F}{G}{\ISOM{\chi}}:\ \Igroup{F} \longrightarrow \Igroup{G}$\ \  is an isomorphism, for any
			isomorphism\ \ $\ISOM{\chi}:\Igroup{G} \longrightarrow \Igroup{F}$.

		\item[(b)]
			If\ \ $\ISOM{\chi}:\Igroup{G} \longrightarrow \Igroup{F}$\ \ and\ \
			   $\ISOM{\eta}:\Igroup{H} \longrightarrow \Igroup{G}$\ \ are isomorphisms, then:
			\begin{equation}
				\label{eq_composing_hx_isomorphisms}
				\hMAP{F}{H}{\ISOM{\chi} \ISOM{\eta}}
				\ \ =\ \ 
				\hMAP{G}{H}{\ISOM{\eta}}
				\hMAP{F}{G}{\ISOM{\chi}}\ \ ,
			\end{equation}
			while for their inverses:
			\begin{equation}
				\label{eq_composing_inverses_of_hx_isomorphisms}
				\INVERSE{ \hMAP{F}{H}{\ISOM{\chi} \ISOM{\eta}} }
				\ \ =\ \ 
				\INVERSE{ \hMAP{F}{G}{\ISOM{\chi}} }
				\INVERSE{ \hMAP{G}{H}{\ISOM{\eta}} }\ \ .
			\end{equation}

		\item[(c)]
			If\ \ $\ISOM{\chi}:\Igroup{G} \longrightarrow \Igroup{F}$\ \ is an isomorphism, then:
			\begin{equation}
				\label{eq_h_map_of_inverted_isomorphism}
				\hMAP{G}{F}{\ISOM{\chi^{-1}}}
				\ \ =\ \ 
				\INVERSE{ \hMAP{F}{G}{\ISOM{\chi}} }\  \  .
			\end{equation}

		\item[(d)]
			Map (see Definition \ref{def_isom_object} for $\ISOMORPHISMS{\ldots}{\ldots}$)
			\begin{equation}
				\label{eq_x_hx_map}
				\ISOM{\chi} \longrightarrow  \hMAP{F}{G}{\ISOM{\chi}}\ :\ \ 
				\ISOMORPHISMS{G}{F}
				\longrightarrow
				\ISOMORPHISMS{F}{G}
			\end{equation}
			is a bijection, and its inverse is the map:
			\begin{equation}
				\label{eq_inverse_of_x_hx_map}
				\ISOM{\psi} \longrightarrow  \hMAP{G}{F}{\ISOM{\psi}}\ :\ \ 
				\ISOMORPHISMS{F}{G}
				\longrightarrow
				\ISOMORPHISMS{G}{F}\ \ .
			\end{equation}

		\item[(e)]
			Map
			\begin{equation}
				\label{eq_x_hx_isomorphism_inverted}
				\ISOM{\chi} \longrightarrow  \INVERSE{ \hMAP{F}{G}{\ISOM{\chi}} }\ :\ \ 
				\ISOMORPHISMS{G}{F}
				\longrightarrow
				\ISOMORPHISMS{G}{F}
			\end{equation}
			is an involution, i.e. it is the inverse of itself.
	\end{description}
\end{lemma}

\PROOFstart 
\begin{description}
	\item[(a)] 
		This is because $\INVERSE{ \hMAP{F}{G}{\ISOM{\chi}} }$ is an isomorphism by 
		Lemma \ref{lem_P_determines_R_for_permutation_equivalence_PFRequalsG}.

	\item[(b)]
		By Lemma \ref{lem_P_determines_R_for_permutation_equivalence_PFRequalsG} we have that\ \ \ 
		$\Pro{\chi} F \Pro{ \hMAP{F}{G}{\ISOM{\chi}} }\ =\ G$\ \ \   and\ \ \  
		$\Pro{\eta} G \Pro{ \hMAP{G}{H}{\ISOM{\eta}} }\ =\ H$,\ \ \ so
		\begin{eqnarray}
			\label{eq_from_F_through_G_to_H}
			\lefteqn{ H\ \ = }
			&  &    \\
			\nonumber
			&  &  
			\Pro{\ISOM{\eta}} \Pro{\ISOM{\chi}} 
			\cdot F \cdot 
			\Pro{ \hMAP{F}{G}{\ISOM{\chi}} } 
			\Pro{ \hMAP{G}{H}{\ISOM{\eta}} }
			\ \ =\ \ 
			\Pro{\ISOM{\chi} \ISOM{\eta}} 
			\cdot F \cdot
			\Pro{ \hMAP{G}{H}{\ISOM{\eta}}  \hMAP{F}{G}{\ISOM{\chi}} }
			\ \ =\ \ 
			\Pro{\ISOM{\chi} \ISOM{\eta}}
			\cdot F \cdot
			\Pro{ \INVERSE{  \hMAP{G}{H}{\ISOM{\eta}}  \hMAP{F}{G}{\ISOM{\chi}} } }^T\ \ ,
		\end{eqnarray}
		where $\ISOM{\chi} \ISOM{\eta}:\ \Igroup{H} \longrightarrow \Igroup{F}$\ \  and\ \  
		$\hMAP{G}{H}{\ISOM{\eta}}  \hMAP{F}{G}{\ISOM{\chi}}\ :\  \Igroup{F} \longrightarrow \Igroup{H}$
		are isomorphisms.

		By Lemma \ref{lem_P_determines_R_for_permutation_equivalence_PFRequalsG}
		the isomorphism $\INVERSE{  \hMAP{G}{H}{\ISOM{\eta}}  \hMAP{F}{G}{\ISOM{\chi}} }$
		designating the right transposed permutation matrix in (\ref{eq_from_F_through_G_to_H})
		is equal to the function $\INVERSE{  \hMAP{F}{H}{\ISOM{\chi} \ISOM{\eta}} }$. From this fact we get 
		(\ref{eq_composing_hx_isomorphisms}) and immediately (\ref{eq_composing_inverses_of_hx_isomorphisms}).

	\item[(c)]
		By Lemma \ref{lem_P_determines_R_for_permutation_equivalence_PFRequalsG}\ \ \ 
		$\Pro{\ISOM{\chi}} F \Pro{\hMAP{F}{G}{\chi}}\ =\ G$,\ \  which is equivalent to 
		\begin{equation}
			\label{eq_PFRequalsG_converted_to_PGRequalsF}
			\Pro{\ISOM{\chi}^{-1}} G \Pro{ \hMAP{F}{G}{\chi} }^T \ \ =\ \ F.
		\end{equation}
		Again, by Lemma \ref{lem_P_determines_R_for_permutation_equivalence_PFRequalsG}
		the isomorphism $\hMAP{F}{G}{\chi}$
		designating the right transposed permutation matrix must be equal to 
		$\INVERSE{ \hMAP{G}{F}{\ISOM{\chi}^{-1}} }$. This is equivalent to (\ref{eq_h_map_of_inverted_isomorphism}).

	\item[(d)]
		By Lemma \ref{lem_P_determines_R_for_permutation_equivalence_PFRequalsG}
		in the equality
		\begin{equation}
			\label{eq_PFRequalsG_with_right_isomorphism_considered}
			\Pro{\ISOM{\chi}} F \Pro{ \INVERSE{ \hMAP{F}{G}{\ISOM{\chi}} } }^T\ \ =\ \ G
		\end{equation}
		isomorphism $\ISOM{\chi'}\ =\ \ISOM{\chi} \in \ISOMORPHISMS{G}{F}$ is determined by 
		the isomorphism $\ISOM{\chi''}\ =\ \INVERSE{ \hMAP{F}{G}{\ISOM{\chi}} }$
		designating the right transposed permutation matrix:
		\begin{equation}
			\label{eq_x_determined_by_hx_isomorphism_inverted}
			\ISOM{\chi}\ \ =\ \ \INVERSE{   \hMAP{F}{G}{  \INVERSE{ \hMAP{F}{G}{\ISOM{\chi}} }  }   }\ \ ,
		\end{equation}
		so
		\begin{equation}
			\label{eq_x_inverted_determined_by_hx}
			\ISOM{\chi}^{-1}
			\ \ =\ \   
			\hMAP{F}{G}{  \INVERSE{ \hMAP{F}{G}{\ISOM{\chi}} }  }
			\ \ =\ \   
			\INVERSE{   \hMAP{G}{F}{   \hMAP{F}{G}{\ISOM{\chi}} }  }\ \ ,
		\end{equation}
		where property {\bf (c)} has been used. Inverting the above equality between the leftmost and rightmost 
		isomorphisms we get
		\begin{equation}
			\label{eq_x_equals_hGF_hFG_x}
			\ISOM{\chi}\ \ =\ \ \hMAP{G}{F}{   \hMAP{F}{G}{ \ISOM{\chi} }   }\ \ ,
		\end{equation}
		and precisely in the same way we show that for $\psi \in \ISOMORPHISMS{F}{G}$
		\begin{equation}
			\label{eq_y_equals_hFG_hGF_y}
			\ISOM{\psi}\ \ =\ \ \hMAP{F}{G}{   \hMAP{G}{F}{ \ISOM{\psi} }   }\ \ .
		\end{equation}
		We thus see that maps $\ISOM{\chi} \longrightarrow \hMAP{F}{G}{\ISOM{\chi}}$ and
		$\ISOM{\psi} \longrightarrow \hMAP{G}{F}{\ISOM{\psi}}$ are the inverses of each other, so 
		they are bijections.

	\item[(e)] 
		This statement is a consequence of (\ref{eq_x_determined_by_hx_isomorphism_inverted}),
		or can also be calculated as:
		\begin{equation}
			\label{eq_inverted_h_inverted_hx}
			\INVERSE{  \hMAP{F}{G}{  \INVERSE{  \hMAP{F}{G}{  \ISOM{\chi}  }  }  }  }
			\ \ =\ \ 
			\hMAP{G}{F}{  \hMAP{F}{G}{  \ISOM{\chi}  }  }
			\ \ =\ \ 
			\ISOM{\chi}\ \ ,
		\end{equation}
		where the first equality follows from {\bf (c)}, and the last one from {\bf (d)}.
\end{description}

\PROOFend  


\begin{corollary}
	\label{cor_h_FF_map_properties}
	Let $F$ be a Fourier matrix indexed by group $\Igroup{F}$. Then:
	\begin{description}
		\item[(a)]
			Map\ \  $\hMAP{F}{F}{\ISOM{\varrho}}:\ \Igroup{F} \longrightarrow \Igroup{F}$\ \  is an isomorphism,
			for any isomorphism\ \ $\ISOM{\varrho}:\ \Igroup{F} \longrightarrow \Igroup{F}$.

		\item[(b)]
			Map\ \  $\ISOM{\varrho} \longrightarrow \hMAP{F}{F}{\ISOM{\varrho}} :\ \ISOMORPHISMS{F}{F} \longrightarrow \ISOMORPHISMS{F}{F}$\ \  
			is an anti-isomorphism, hence
			\begin{equation}
				\label{eq_x_hFFx_map_is_antiisomorphism}
				\hMAP{F}{F}{\ISOM{\varrho} \ISOM{\sigma}}
				\ \ =\ \ 
				\hMAP{F}{F}{\ISOM{\sigma}}
				\hMAP{F}{F}{\ISOM{\varrho}}\  \ ,
			\end{equation}
			and it is an involution (the inverse of itself).

		\item[(c)]
			Map\ \  $\ISOM{\varrho} \longrightarrow \INVERSE{ \hMAP{F}{F}{\ISOM{\varrho}} }  
					   :\ \ISOMORPHISMS{F}{F} \longrightarrow \ISOMORPHISMS{F}{F}$ 
			is an isomorphism, hence
			\begin{equation}
				\label{eq_x_hFFx_isomorphism_inverted_map_is_isomorphism}
				\INVERSE{  \hMAP{F}{F}{\ISOM{\varrho} \ISOM{\sigma}}  }
				\ \ =\ \ 
				\INVERSE{  \hMAP{F}{F}{\ISOM{\varrho}}  }
				\INVERSE{  \hMAP{F}{F}{\ISOM{\sigma}}  }\  \ ,
			\end{equation}
			and it is an involution.

		\item[(d)]
			\begin{equation}
				\label{eq_x_hFFx_map_values_at_identity_and_x_inverted}
				\hMAP{F}{F}{\IDENTITY{\Igroup{F}}} \ \ =\ \  \IDENTITY{\Igroup{F}}
				\ \ \ \ ,\ \ 
				\hMAP{F}{F}{\ISOM{\varrho}^{-1}}  \ \ =\ \  \INVERSE{ \hMAP{F}{F}{\ISOM{\varrho}} }\ \ ,
			\end{equation}
			where $\IDENTITY{\Igroup{F}}$ is the identity on $\Igroup{F}$.
	\end{description}
\end{corollary}

\PROOFstart 
\begin{description}
	\item[(a)] 
		It follows from Lemma \ref{lem_h_FG_map_properties} {\bf (a)}.

	\item[(b)]
		Equality (\ref{eq_x_hFFx_map_is_antiisomorphism}) results from applying Lemma \ref{lem_h_FG_map_properties} {\bf (b)} 
		to\ \  $F$,\ \ $G=F$,\ \  $H = F$,\ \  	$\ISOM{\chi} = \ISOM{\varrho}$ and $\ISOM{\eta}  = \ISOM{\sigma}$. 
		Thus $\ISOM{\varrho} \longrightarrow \hMAP{F}{F}{\ISOM{\varrho}}$ is an anti-homomorphism 
		on $\ISOMORPHISMS{F}{F}$.

		Next, from Lemma \ref{lem_h_FG_map_properties} {\bf (d)} applied to matrices $F$ and $G = F$ we get that map 
		$\ISOM{\varrho} \longrightarrow \hMAP{F}{F}{\ISOM{\varrho}}$ on $\ISOMORPHISMS{F}{F}$ is a bijection and
		the very same map is its inverse. So the considered map is an involution and, because of being additionally an anti-homomorphism,
		it is an anti-isomorphism.

	\item[(c)]
		Inverting both sides of equality (\ref{eq_x_hFFx_map_is_antiisomorphism}) we get the equality in 
		(\ref{eq_x_hFFx_isomorphism_inverted_map_is_isomorphism}) which means that 
		map  $\ISOM{\varrho} \longrightarrow \INVERSE{ \hMAP{F}{F}{\ISOM{\varrho}} }$ is a homomorphism on
		$\ISOMORPHISMS{F}{F}$.

		From Lemma \ref{lem_h_FG_map_properties} {\bf (e)} applied to $F$ and $G=F$ it follows that 
		$\ISOM{\varrho} \longrightarrow \INVERSE{ \hMAP{F}{F}{\ISOM{\varrho}} }$ is an involution, that is, among other,
		it is a bijection. Because of being additionally a homomorphism it is an isomorphism.

	\item[(d)]
		It is an easy consequence of $\ISOM{\varrho} \longrightarrow \hMAP{F}{F}{\ISOM{\varrho}}$ being an 
		anti-homomorphism.
		
\end{description}

\PROOFend 

An immediate consequence of Lemmas 
\ref{lem_permutation_equivalence_PFRequalsG_matrices_determined_by_isomorphisms}
and 
\ref{lem_P_determines_R_for_permutation_equivalence_PFRequalsG}
is that $\MAP{\PsPs}{F}{G}$,
 and $\STAB{\PsPs}{F}\ =\ \MAP{\PsPs}{F}{F}$ as a special case 
(see Definition  \ref{def_pairs_action_stabilizer}),
can now be described in terms of permutation matrices designated by isomorphisms 
$\ISOM{\chi}:\ \Igroup{G} \longrightarrow \Igroup{F}$ and their images under 
$\ISOM{\chi} \longrightarrow \INVERSE{ \hMAP{F}{G}{\ISOM{\chi}} }$:


\begin{corollary}
	\label{cor_permutational_MAP_FG}
     \rule{0cm}{0cm}

	\begin{description}
		\item[(a)]
			\begin{eqnarray}
				\label{eq_permutational_STAB_F}
				\STAB{\PsPs}{F}
				&  =  &
				\left\{\ 
					\left(
						\Pro{\ISOM{\varrho}}
						\ ,\ 
						\Pro{  \INVERSE{ \hMAP{F}{F}{\ISOM{\varrho}} }  }
					\right)
					\ :\ \ 
					\ISOM{\varrho} \in \ISOMORPHISMS{F}{F}
					\ 
				\right\}     \\
				\nonumber
				&  =  &
				\left\{\ 
					\left(
						\Pro{  \INVERSE{ \hMAP{F}{F}{\ISOM{\varrho}} }  }
						\ ,\ 
						\Pro{\ISOM{\varrho}}
					\right)
					\ :\ \ 
					\ISOM{\varrho} \in \ISOMORPHISMS{F}{F}
					\ 
				\right\}\ \ .
			\end{eqnarray}

		\item[(b)]
			\begin{eqnarray}
				\label{eq_permutational_MAP_FG}
				\MAP{\PsPs}{F}{G}
				&  =  &
				\left\{\ 
					\left(
						\Pro{\ISOM{\chi}}
						\ ,\ 
						\Pro{  \INVERSE{ \hMAP{F}{G}{\ISOM{\chi}} }  }
					\right)
					\ :\ \ 
					\ISOM{\chi} \in \ISOMORPHISMS{G}{F}
					\ 
				\right\}     \\
				\nonumber
				&  =  &
				\left\{\ 
					\left(
						\Pro{  \INVERSE{ \hMAP{F}{G}{\ISOM{\chi}} }  }
						\ ,\ 
						\Pro{\ISOM{\chi}}
					\right)
					\ :\ \ 
					\ISOM{\chi} \in \ISOMORPHISMS{G}{F}
					\ 
				\right\}\ \ .
			\end{eqnarray}

	\end{description}
\end{corollary}

A little more can be said about the stabilizer of $F$ when this matrix is permutation equivalent to $\FOURIER{N}$, 
i.e. when $\Igroup{F}$ is isomorphic to $\Igroup{\FOURIER{N}} \ =\ \ZN{N}$
	(see Corollary \ref{cor_Fouriers_perm_equivalent_iff_their_indexing_groups_isomorphic}),
or, in other words,  when $\Igroup{F}$ is cyclic. In this situation
map $\ISOM{\varrho} \longrightarrow \hMAP{F}{F}{\ISOM{\varrho}}$ is just the identity on $\ISOMORPHISMS{F}{F}$.


\begin{lemma}
	\label{lem_permutational_stabilizer_of_F_when_I_F_is_cyclic}
	If the indexing group $\Igroup{F}$ of a Fourier matrix $F$ is cyclic then:
	\begin{equation}
		\label{eq_h_FF_is_identity_for_I_F_cyclic}
		\hMAP{F}{F}{\ISOM{\varrho}} \ \ =\ \ \ISOM{\varrho}
		\ \ \ \ \mbox{for any}\ \ 
		\ISOM{\varrho} \in \ISOMORPHISMS{F}{F}\ \ ,
	\end{equation}
	and 
	\begin{equation}
		\label{eq_permutational_stabilizer_of_F_when_I_F_is_cyclic}
		\STAB{\PsPs}{F} \ \ =\ \ 
		\left\{\ 
			\left(
				\Pro{\ISOM{\varrho}}\ ,\ \Pro{\ISOM{\varrho}}^T
			\right)
			\ :\ \ 
			\ISOM{\varrho} \in \ISOMORPHISMS{F}{F}\ 
		\right\}\ \ .
	\end{equation}
\end{lemma}

\PROOFstart  
Acording to Lemma \ref{lem_P_determines_R_for_permutation_equivalence_PFRequalsG}\ \ 
$\hMAP{F}{F}{\ISOM{\varrho}}:\ \Igroup{F} \longrightarrow \Igroup{F}$\ \  is the only bijection such that:
\begin{equation}
	\label{eq_h_FF_defining_property}
	\forall i,j \in \Igroup{F}\ \ \ \ \ \ F_{ \ISOM{\varrho}(i),\ j } \ \ =\ \ F_{ i,\ \hMAP{F}{F}{\ISOM{\varrho}}(j) }\ \ ,
\end{equation}
so we need to show that $\ISOM{\varrho}$ has the above property.

$\Igroup{F}$ is assumed to be cyclic, let $g$ denote one of its generators then. For some positive integers 
$\tilde{i}$, $\tilde{j}$ the elements of $\Igroup{F}$: $i$, $j$, $\ISOM{\varrho}(i)$ and $\ISOM{\varrho}(j)$ can be expressed
as multiplicities:\ \  $i\ =\ \tilde{i} g$,\ \  $j\ =\ \tilde{j} g$,\ \  $\ISOM{\varrho}(i)\ =\ \tilde{i} \ISOM{\varrho}(g)$\ \  and\ \ 
$\ISOM{\varrho}(j)\ =\ \tilde{j} \ISOM{\varrho}(g)$. 

Because $F_{k,l}$ is a homomorphism with respect to both $k$ and $l$  ($k,l  \in \Igroup{F}$), we check that (\ref{eq_h_FF_defining_property})
is true with $\hMAP{F}{F}{\ISOM{\varrho}}$ replaced with $\ISOM{\varrho}$\ :
\begin{eqnarray}
	\label{eq_hFF_defining_property_proof_with_hFFro_equal_ro}
	\lefteqn{  F_{\ISOM{\varrho}(i),\ j} \ \ =  }
	&  &   \\
	\nonumber
	&  &  
	F_{\tilde{i} \ISOM{\varrho}(g),\ \tilde{j} g} 
	\ \ =\ \ 
	\left(  F_{\tilde{i} \ISOM{\varrho}(g),\  g}  \right)^{\tilde{j}}
	\ \ =\ \ 
	\left(  F_{\ISOM{\varrho}(g),\  g}  \right)^{\tilde{i} \tilde{j}}
	\ \ =\ \ 
	\left(  F_{g,\  \ISOM{\varrho}(g)}  \right)^{\tilde{i} \tilde{j}}
	\ \ =\ \ 
	F_{\tilde{i} g,\ \tilde{j} \ISOM{\varrho}(g)}  
	\\
	\nonumber
	\lefteqn{ =\ \ F_{i,\ \ISOM{\varrho}(j)}\ \ , }
	&  &
\end{eqnarray}
where also the symmetry of $F$ has been used.

The last statement in the lemma is a consequence of Corollary \ref{cor_permutational_MAP_FG} {\bf (a)}.
\PROOFend 

It does not have to be so when $F$ is not permutation equivalent to $\FOURIER{N}$, as the below example presents.
Matrix $\FOURIER{2} \otimes \FOURIER{2}$ considered there is not permutation equivalent to $\FOURIER{4}$ because
groups $\ZN{2} \times \ZN{2}$ and $\ZN{4}$ are not isomorphic ($\ZN{2} \times \ZN{2}$ doesn't have a generator, i.e. an 
element of order $4$).

\begin{example}
	\label{ex_F2F2_hF_not_equal_to_identity}
	Let $F\ =\ \FOURIER{2} \otimes \FOURIER{2}$, then $\Igroup{F}\ =\ \ZN{2} \times \ZN{2}$. Let
	$\ISOM{\varrho}$ be an automorphism on $\Igroup{F}$ given by the table:

	\begin{tabular}{c||c|c|c|c|}
		$(i,j)$                            &     $(0,0)$    &    $(0,1)$    &    $(1,0)$     &    $(1,1)$      \\
		\hline 
		$\ISOM{\varrho}(i,j)$     &     $(0,0)$    &    $(1,1)$    &     $(1,0)$    &    $(0,1)$       
	\end{tabular} .

	Let $F'$ be defined, using the group indices from $\Igroup{F}$, by\ \   $F'_{i,j}\ =\ F_{\ISOM{\varrho}(i),\ j}$.
	Recall that the group indices $(0,0)$, $(0,1)$, $(1,0)$, $(1,1)$ correspond to the ordinary indices $1$, $2$, $3$, $4$.
	Below $F$, $F'$ and $F'$ in terms of the columns of $F$ are presented:
	\begin{equation}
		\label{eq_F_and_Fprim}
		F =
		\left[
			\begin{array}{rrrr}
				1    &   1   &   1   &   1     \\
				1    &  -1   &   1   &  -1     \\
				1    &   1   &  -1   &  -1     \\
				1    &  -1   &  -1   &   1
			\end{array}
		\right],
		\ \ \ \ 
		F' =
		\left[
			\begin{array}{rrrr}
				1    &   1   &   1   &   1     \\
				1    &  -1   &  -1   &  1     \\
				1    &   1   &  -1   &  -1     \\
				1    &  -1   &   1   &  -1
			\end{array}
		\right]
		\  =\ 
		\left[
			F_{:, (0,0)},\  
			F_{:, (0,1)},\ 
			F_{:, (1,1)},\ 
			F_{:, (1,0)}
		\right],
	\end{equation}
	so the values of $\hMAP{F}{F}{\ISOM{\varrho}}$
	(defined in Lemma \ref{lem_P_determines_R_for_permutation_equivalence_PFRequalsG})
	are:

	\begin{tabular}{c||c|c|c|c|}
		$(i,j)$                                                  &     $(0,0)$    &    $(0,1)$    &    $(1,0)$     &    $(1,1)$      \\
		\hline 
		$\hMAP{F}{F}{\ISOM{\varrho}}(i,j)$     &     $(0,0)$    &    $(0,1)$    &     $(1,1)$    &    $(1,0)$       
	\end{tabular} .

	In this case $\hMAP{F}{F}{\ISOM{\varrho}}\ \neq\ \ISOM{\varrho}$. The stabilizing pair, a member of 
	$\STAB{\PsPs}{F}$ according to Corollary \ref{cor_permutational_MAP_FG} {\bf (a)}, corresponding to 
	$\ISOM{\varrho}$ is the pair:
	\begin{equation}
		\label{eq_stabilizing_4x4_example_pair_for_ro}
		\left(\ 
			\Pro{  \ISOM{\varrho}  },\ 
			\Pro{  \INVERSE{ \hMAP{F}{F}{\ISOM{\varrho}} }  }\ 
		\right)
		\ \  \ =\ \  \
		\left(\ \
			\left[
				\begin{array}{cccc}
					1   &   0   &   0   &   0   \\
					0   &   0   &   0   &   1   \\
					0   &   0   &   1   &   0   \\
					0   &   1   &   0   &   0   
				\end{array}
			\right],
			\  \ 
			\left[
				\begin{array}{cccc}
					1   &   0   &   0   &   0   \\
					0   &   1   &   0   &   0   \\
					0   &   0   &   0   &   1   \\
					0   &   0   &   1   &   0   
				\end{array}
			\right]\ \ 
		\right)\ \  \ ,
	\end{equation}
	for which there holds\ \ 
	$\Pro{  \ISOM{\varrho}  }  \cdot F  \cdot \Pro{  \INVERSE{ \hMAP{F}{F}{\ISOM{\varrho}} }  }^T\ =\ F$.
\end{example}

Commenting the case of Fourier matrices associated with cyclic indexing groups let us note that applying
a permutational  stabilizing pair $\left( \Pro{\ISOM{\varrho}},  \Pro{\ISOM{\varrho}}^T \right)$ to $F$
($\ISOM{\varrho} \in \ISOMORPHISMS{F}{F}$):\ \ 
$\left( \Pro{\ISOM{\varrho}},  \Pro{\ISOM{\varrho}}^T \right) F 
		\ \ =\ \ 
	\Pro{\ISOM{\varrho}}  F  \Pro{\ISOM{\varrho}}$ is not in general a similarity transformation, because
$\Pro{\ISOM{\varrho}}$ is not usually equal to $\Pro{\ISOM{\varrho}}^{-1} \ =\ \Pro{\ISOM{\varrho}^{-1}}$.
For $\Pro{\ISOM{\varrho}} \ =\ \Pro{\ISOM{\varrho}^{-1}}$ to hold it is necessary and sufficient that $\ISOM{\varrho}\ =\ \ISOM{\varrho}^{-1}$
which in turn is equivalent to $\ISOM{\varrho}(g)\ =\ \ISOM{\varrho}^{-1}(g)$, where $g$ is any generator of $\Igroup{F}$,
because any $\ISOM{\chi} \in \ISOMORPHISMS{F}{F}$ is determined by its value on a generator.
$\ISOM{\varrho}(g)$ then is also a generator because 
	$\Igroup{F}\ =\ \ISOM{\varrho}\left( \Igroup{F} \right) \ =\ 
	\left\{ \ISOM{\varrho}(g), \ISOM{\varrho}(2g), \ldots , \ISOM{\varrho}\left(Ng = \GROUPzero{\Igroup{F}}\right) \right\} \ =\ 
	\left\{ \ISOM{\varrho}(g), 2 \ISOM{\varrho}(g), \ldots, N \ISOM{\varrho}(g) = \GROUPzero{\Igroup{F}} \right\}$, 
so $\ISOM{\varrho}(g) \ =\ \tilde{k} g$ for some $\tilde{k} \in \{ 1,\ldots,N-1\} \subset \INTEGER$ and 
$\gcd\left(\tilde{k},N\right)\ =\ 1$, where $N$ is the order of $\Igroup{F}$. 
Isomorphism $\ISOM{\chi}(g)$  satisfying $\ISOM{\chi}(g) \ =\ \tilde{n} g$ ($g$ fixed!) will further be denoted by
$\ISOM{\varrho}_{\tilde{n}}$ and then for any $\tilde{x} g \in \Igroup{F}$ there holds 
$\ISOM{\varrho}_{\tilde{n}}\left( \tilde{x} g \right) \ =\ \tilde{x} \ISOM{\varrho}_{\tilde{n}}\left(  g \right) \ =\ 
\tilde{x} \left( \tilde{n} g \right) \ =\ \left(\tilde{x} \tilde{n}\right) g$.

     Because $\ISOM{\varrho}_{\tilde{k}}\left( \ISOM{\varrho}_{\tilde{k}}(g) \right)  \ =\ \tilde{k}^2 g$
condition $\ISOM{\varrho}_{\tilde{k}}(g) \ =\ \ISOM{\varrho}_{\tilde{k}}^{-1}(g)$, 
equivalent to $\ISOM{\varrho}_{\tilde{k}}\left( \ISOM{\varrho}_{\tilde{k}}(g) \right) \ =\ g$, amounts to $\tilde{k}^2 g = g$, 
that is $\tilde{k}^2 \ =\ 1 \mod N$, where $\tilde{k}$ is relatively prime to $N$.

	Besides, the set of all generators of $\Igroup{F}$ forms an abelian group:\\ 
$\left( \left\{ \tilde{k} g:\ \tilde{k} \in \{1,\ldots,N-1\},\ \gcd\left( \tilde{k},N \right)=1 \right\},
	\ \ \tilde{k}g \bullet \tilde{l} g \ \stackrel{def}{=}\ \left( \tilde{k} \tilde{l} \right) g  \right)$,\ \ 
with $g$ playing the role of the neutral element. It is isomorphic to $\ISOMORPHISMS{F}{F}$ (thus abelian),
via $\tilde{k} g \longrightarrow \ISOM{\varrho}_{\tilde{k}}$, 
since 
$\ISOM{\varrho}_{\tilde{k}} \ISOM{\varrho}_{\tilde{l}} \ =\ 
 \ISOM{\varrho}_{\tilde{k}\tilde{l}}$\ \ 
as\ \  
$\ISOM{\varrho}_{\tilde{k}}  \ISOM{\varrho}_{\tilde{l}}(g)  \ =\ 
\left(\tilde{k}\tilde{l}\right) g$.

We end with two conclusions, both for a cyclic $\Igroup{F}$ of order $N$. The first one is that
$\ISOMORPHISMS{F}{F}\ = \ \left\{ \ISOM{\varrho}_{\tilde{k}}:\ \tilde{k} \in \{1,\ldots,N-1\},\ \gcd\left(\tilde{k},N\right) = 1 \right\}$,
a well known fact in the group theory. The second one is that 
$\left( \Pro{\ISOM{\varrho}_{\tilde{k}}},  \Pro{\ISOM{\varrho}_{\tilde{k}}}^T \right) F$ is a similarity transformation
if and only if $\tilde{k}^2 \ =\ 1 \mod N$.  This last condition is satisfied, for any $N$, by $\tilde{k}\ =\ 1, N-1$.
But for example for $N=9,10$ there are no more such $\tilde{k}$'s. On the other hand, for $N=12,24$ all $\tilde{k}$'s
from $\{1,\ldots,N-1\}$ relatively prime to $N$  satisfy the condition, that is every pair from $\STAB{\PsPs}{F}$  
produces a similarity transformation.

Here we end the digression on the permutational stabilizer of $F$ indexed by a cyclic  group. As for the general abelian case, 
we refer the reader to \cite{Rhea_automorphisms},\cite{Curran_automorphisms} for the description of $\ISOMORPHISMS{F}{F}$
(the group of automorphisms of $\Igroup{F}$) which in a sense
parametrizes $\STAB{\PsPs}{F}$ according to Corollary \ref{cor_permutational_MAP_FG} {\bf a)}.

We know a lot about pairs of permutation matrices, the elements of $\PsPs$, acting on Fourier matrices. To take the rest of the encompassing
group $\PP$ into consideration we must introduce the special matrices $\Wk{F}{k}$ and $\Zk{F}{k}$:


\begin{definition}
	\label{def_Wk_Zk_matrices}
	Let $F$ be a Fourier matrix of size $N$ indexed by group $\Igroup{F}$.
	\begin{description}
		\item[(a)]
			We define $\Wk{F}{k}$, where $k \in \Igroup{F}$, to be the $N \times N$ diagonal matrix 
			containing the $k$-th row (or column) of $F$ as its diagonal:
			\begin{equation}
				\label{eq_Wk_definition}
				\Wk{F}{k} \ \ =\ \ \diag\left(  F_{k,:} \right)
							 \ \ =\ \ \diag\left(  F_{:,k} \right)
				,\ \ \ \ \ \ \ \ \ k \in \Igroup{F}\ \ .
			\end{equation}
			
		\item[(b)]
			We define $\Zk{F}{k}$,   where $k \in \Igroup{F}$, to be the permutation matrix shifting every $j$-th
			entry of a vertical vector $v$ to the $(j-k)$-th position 
			when this vector is multiplied by $\Zk{F}{k}$:\ \  $v \longrightarrow \Zk{F}{k} \cdot v$, 
			\begin{equation}
				\label{eq_Zk_definition}
				\forall  i \in \Igroup{F}\ \ \ \ \ \ \ELEMENTof{\Zk{F}{k}}{i}{i+k}\ \ =\ \ 1\ \ .
		\end{equation}

	\end{description}
\end{definition}

In general $\Wk{F}{k}$ and $\Zk{F}{k}$ have a Kronecker product structure related to the
structure of $F$:


\begin{lemma}
	\label{lem_Wk_Zk_are _Kronecker_products}
	Let $F\ =\ \FOURIER{N_1} \otimes \ldots \otimes \FOURIER{N_r}$, indexed by 
	$\Igroup{F}\ =\ \ZN{N_1} \times \ldots \times \ZN{N_r}$. 
	Let $k\ =\ (k_1,\ldots,k_r) \in \Igroup{F}$. Then:
	\begin{description}
		\item[(a)]
			$\Wk{F}{k}\ \ =\ \ \Wk{\FOURIER{N_1}}{k_1} \otimes \ldots \otimes \Wk{\FOURIER{N_r}}{k_r}$\ \ ,

		\item[(b)]
			$\Zk{F}{k}\ \ =\ \ \Zk{\FOURIER{N_1}}{k_1} \otimes \ldots \otimes \Zk{\FOURIER{N_r}}{k_r}$\ \ ,

	\end{description}
	where, in accordance with Definition \ref{def_Wk_Zk_matrices}, 
	$\Wk{\FOURIER{N_x}}{k_x}$ is the $N_x \times N_x$ diagonal matrix containing the $k_x$-th row of $\FOURIER{N_x}$
	as its diagonal, and 
	$\Zk{\FOURIER{N_x}}{k_x}$ is the $N_x \times N_x$ circulant permutation matrix causing the cyclic shift of the coordinates 
	of a vertical vector $v$
	by $\REPR{k_x}$ positions upward in the multiplication $\Zk{\FOURIER{N_x}}{k_x} \cdot v$, where 
	$\REPR{k_x} \in \INTEGER$ represents class $k_x \in \ZN{N_x}$.
\end{lemma}

\PROOFstart   
	\begin{description}
		\item[(a)]
			It is a consequence of:
			\begin{equation}
				\label{eq_kth_row_of_F_as_Kronecker_product}
				F_{k,:} \ \ =\ \ 
				\ELEMENTof{\FOURIER{N_1}}{k_1}{:} \otimes \ldots \otimes \ELEMENTof{\FOURIER{N_r}}{k_r}{:}\ \ .
			\end{equation}

		\item[(b)]
			$\Zk{\FOURIER{N_x}}{k_x}$ has $1$'s at positions $i_x, (i_x + k_x)$ where $i_x \in \ZN{N_x}$. 
			
			The $i,j$-th entry of the Kronecker product at {\bf (b)}, where 
			$i\ =\ (i_1,\ldots,i_r),\ \ j\ =\ (j_1,\ldots,j_r)\ \ \in \Igroup{F}$, is equal to
			\begin{equation}
				\label{eq_ijth_element_of_Kronecker_product_of_Zks}
				\ELEMENTof{\Zk{\FOURIER{N_1}}{k_1}}{i_1}{j_1} 
				\cdot \ldots \cdot 
				\ELEMENTof{\Zk{\FOURIER{N_r}}{k_r}}{i_r}{j_r}
			\end{equation}
			and it is equal to $1$ if and only if all the factors are $1$'s, that is 
			$j_1\ =\ i_1\ +\ k_1$, ...,   $j_r\ =\ i_r\ +\ k_r$, or in a compact form $j\ =\ i + k$.
			But this just the definition of $\Zk{F}{k}$.

	\end{description}
\PROOFend  

Further, $\Wk{F}{k}$ and  $\Zk{F}{k}$ faithfully represent the indexing group $\Igroup{F}$, as it is stated in the next lemma. 
An easy proof  we ommit.


\begin{lemma}
	\label{lem_maps_k_Wk_and_k_Zk_are_monomorphisms}
	Let $F$ be an $N \times N$ Fourier matrix indexed by group $\Igroup{F}$. 

	Maps $k \longrightarrow \Wk{F}{k}$ and $k \longrightarrow \Zk{F}{k}$ 
	from $\Igroup{F}$ into the group of nonsingular complex $N \times N$ matrices are monomorphisms,
	so there holds:

	\begin{equation}
		\label{eq_Wk_homomorphism}
		\Wk{F}{k+l}\ \ =\ \ \Wk{F}{k} \cdot \Wk{F}{l}\ \ ,
		\ \ \ \ \ \ \ \ 
		\Wk{F}{-k}\ \ =\ \ \CONJ{\Wk{F}{k}}\ \ =\ \ \INVERSE{\Wk{F}{k}}\ \ ,
	\end{equation}

	\begin{equation}
		\label{eq_Zk_homomorphism}
		\Zk{F}{k+l}\ \ =\ \ \Zk{F}{k} \cdot \Zk{F}{l}\ \ ,
		\ \ \ \ \ \ \ \ 
		\Zk{F}{-k}\ \ =\ \ \TRANSPOSE{\Zk{F}{k}}\ \ =\ \ \INVERSE{\Zk{F}{k}}
	\end{equation}
	for any $k,l \in \Igroup{F}$.
\end{lemma}

Analyzing the stabilizer of a Fourier matrix $F$ we encounter a monomorphic image of 
group $\GsymmetryGROUP{F}$ defined in the lemma below (the proof of which is left as an exercise).
Group $\GsymmetryGROUP{F}$ and its image were considered by Alexander Karabegov
in \cite{Karabegov}, as a source of symmetries (stabilizing pairs) allowing a calculation of the defect of $F$ using
his characterization of this quantity \\
Recall in this context:
\begin{itemize}
	\item 
		equation (\ref{eq_undephased_defect_and_V1_dimension_equality}),
		i.e. 
		$\undephasedDEFECT(F) 
				\ \stackrel{(\ref{eq_undephased_defect_for_rescaled_unitary})}{=}\ 
		\undephasedDEFECT\left( \frac{1}{\sqrt{N}} F \right) 
				\ =\ 
		\DIMC\left( \Vu{\frac{1}{\sqrt{N}} F}{1} \right)$,
		the dimension of the $\lambda=1$ eigenspace of $\Iu{\frac{1}{\sqrt{N}} F}$,
	\item
	       Definition \ref{def_symmetry_group} of a symmetry group of  $\frac{1}{\sqrt{N}} F$,\ \  
		$\GsymmetryGROUP{F}$ belonging to this class as shall be seen later,
	\item
		which by Lemma \ref{lem_identifying_symmetry_group_with_stabilizer_subgroup} {\bf c)}
		is isomorphic to a subgroup of $\STAB{\PP}{\frac{1}{\sqrt{N}} F}$,
	\item
		which in turn is mapped by $(S,T) \longrightarrow \pairOPERATOR{\leftPERM{S}{T}}{\rightPERM{S}{T}}$
		of Lemma \ref{lem_representations_of_Stab_U} {\bf e)}
		into a subgroup of the group of invertible operators on $\COMPLEXmatricesNxN{N}$,
	\item 
		where this last subgroup commutes with $\Iu{\frac{1}{\sqrt{N}} F}$ by
		  Lemma \ref{eq_properties_of_representations_of_STAB_U} {\bf c)},
	\item
		hence, in the concrete case of $\GsymmetryGROUP{F}$ taken as a symmetry group, the $1$ dimensional invariant spaces of
		the last subgroup are $N^2$ independent eigenspaces of  $\Iu{\frac{1}{\sqrt{N}} F}$, which is explained in the comment
		after Lemma \ref{eq_properties_of_representations_of_STAB_U} and in the proof {\bf b)} of Theorem \ref{theor_eigenbasis_of_I_F};
		thus $\Vu{\frac{1}{\sqrt{N}} F}{1}$ is given.
\end{itemize}
We will return to this issue in the next section, in proof {\bf b)} of Theorem \ref{theor_eigenbasis_of_I_F}.


\begin{lemma}
	\label{lem_G_F_is_a_group}
	Let $F$ be a Fourier matrix, indexed by group $\Igroup{F}$.

	The set\ \  $\GsymmetryGROUP{F}\ =\ \Igroup{F} \times \Igroup{F} \times \{ z \in \COMPLEX:\ \ABSOLUTEvalue{z}=1\}$\ \ 
	 together with the product:
	\begin{equation}
		\label{eq_G_F_product}
		\left( s_1,\ t_1,\ \theta_1 \right) 
		\cdot
		\left( s_2,\ t_2,\ \theta_2 \right) 
		\  \  
		\stackrel{def}{=}
		\ \ 
		\left(  s_1 + s_2,\ \ t_1 + t_2,\ \   \theta_1 \cdot \theta_2 \cdot F_{s_2, t_1} \right)
	\end{equation}
	is a group.
\end{lemma}

It will be proved  in Lemma \ref{lem_G_F_is_a_symmetry_group_of_F} 
that $\GsymmetryGROUP{F}$ is indeed a symmetry group of the unitary matrix $1/\sqrt{N} \cdot F$ in the sense of Karabegov 
(Definition \ref{def_symmetry_group}) . But before, for that and other purposes,
we need this lemma of a computational nature:


\begin{lemma}
	\label{lem_computations_with_F_Wk_Zk}
	Let $F$ be a Fourier matrix indexed by group $\Igroup{F}$, to which all the indices used here belong. 

	\begin{description}
		\item[(a)]
			An entrywise multiplication of all the columns of $F$ by one of the columns amounts to permuting the columns.
			Analogously for the rows:
			\begin{eqnarray}
				\label{eq_Wk_F_equals_F_Zk_transposed}
				\Wk{F}{k} \cdot F    &   =   &   F \cdot \TRANSPOSE{\Zk{F}{k}} \  \ ,  \\
				\label{eq_F_Wk_equals_Zk_F}
				F \cdot \Wk{F}{k}    &   =   &   \Zk{F}{k} \cdot F\ \ . 
			\end{eqnarray}
			(A comment: Note that (\ref{eq_F_Wk_equals_Zk_F}) is nothing else but an eigenvalue-eigenvector decomposition of $\Zk{F}{k}$, 
			which is not a surprise because $\Zk{F}{k}$ is a Kronecker product of circulant matrices, 
			see Lemma \ref{lem_Wk_Zk_are _Kronecker_products} {\bf (b)},
			and the columns of $\FOURIER{N}$ form an eigenbasis for a circulant $N \times N$ matrix.
			Also note that   (\ref{eq_F_Wk_equals_Zk_F}) implies the equivalence of the faithful representations (monomorphisms) 
			of $\Igroup{F}$:\ $k \longrightarrow \Wk{F}{k}$ and $k \longrightarrow \Zk{F}{k}$, because $\Wk{F}{k}\ =\ F^{-1} \cdot \Zk{F}{k} \cdot F$.)

		\item[(b)]
			The commutation rule:
			\begin{equation}
				\label{eq_Zs_Wt_equals_Fst_Wt_Zs}
				\Zk{F}{s} \cdot \Wk{F}{t} \ \ =\ \ F_{s,t} \cdot \left( \Wk{F}{t} \cdot \Zk{F}{s} \right).
			\end{equation}                                                                                                              
			(A comment: Note that in an alternative form, 
			$\Zk{F}{s} \cdot \Wk{F}{t} \cdot \TRANSPOSE{\Zk{F}{s}} \ =\ F_{s,t} \cdot \Wk{F}{t}$, it states that appropriate 
			reshuffling of  the $t$-th row (column) of $F$
			amounts to multiplying it by a scalar.)    

		\item[(c)]
			An entrywise multiplication of all the rows of $F$ by one of the rows followed by an entrywise multiplication of all the columns of the result by
			a column of $F$ leads to a permuted $F$, multiplied by an entry of $F$\ :
			\begin{equation}
				\label{eq_Zs_F_Zt_transposed_equals_Fst_times_Wt_F_Ws}
				\Zk{F}{s} \cdot F \cdot \TRANSPOSE{\Zk{F}{t}}
				\ \ =\ \ 
				F_{s,t} \cdot \left( \Wk{F}{t} \cdot F \cdot \Wk{F}{s} \right)
			\end{equation}

	\end{description}
\end{lemma}

\PROOFstart   
\begin{description}
	\item[(a)]
			The columns of $F$ form a group under the entrywise product, $F_{:,k} \HADprod F_{:,l} \ =\  F_{:,\ k+l}$, where $k,l \in \Igroup{F}$.
			Consider the $j$-th column of (\ref{eq_Wk_F_equals_F_Zk_transposed}):
			\begin{equation}
				\label{eq_Wk_F_jth_column}
				\ELEMENTof{\Wk{F}{k} F}{:}{j} \ \ =\ \ 
				F_{:,k} \HADprod F_{:,j} \ \ =\ \ 
				F_{:,k+j} \ \ =\ \ 
				\ELEMENTof{ F \TRANSPOSE{\Zk{F}{k}} }{:}{j}\ \ ,
			\end{equation}
			since $\TRANSPOSE{\Zk{F}{k}}$ has $1$'s at the positions $(i+k),i$. Both sides of (\ref{eq_Wk_F_equals_F_Zk_transposed})
			are thus equal column by column. The neighbouring equality (\ref{eq_F_Wk_equals_Zk_F}) is the transposition 
			of (\ref{eq_Wk_F_equals_F_Zk_transposed}).

	\item[(b)]
		Consider expression $\Zk{F}{s} \cdot \Wk{F}{t} \cdot \TRANSPOSE{\Zk{F}{s}}$, in which the resulting diagonal matrix has
		column $\Zk{F}{s} \cdot F_{:,t}\ =\ \ELEMENTof{\Zk{F}{s} F}{:}{t}$ as its diagonal. 
		From (\ref{eq_F_Wk_equals_Zk_F}) at {\bf (a)} this column is equal to
		\begin{equation}
			\label{eq_F_Ws_tth_column}
			\ELEMENTof{ F \Wk{F}{s} }{:}{t}\ \ =\ \ 
			\ELEMENTof{ \Wk{F}{s} }{t}{t} \cdot F_{:,t} \ \ =\ \ 
			F_{s,t} \cdot F_{:,t}\ \ ,
		\end{equation}
		so
		\begin{equation}
			\label{eq_Zs_Wt_Zs_transposed_final_form}
			\Zk{F}{s} \cdot \Wk{F}{t} \cdot \TRANSPOSE{\Zk{F}{s}} \ \ =\ \ F_{s,t} \cdot \Wk{F}{t}\ \ .
		\end{equation}

	\item[(c)]
		Having thrown the matrices $\Wk{F}{s}$ and $\Wk{F}{t}$ to the left hand side we calculate using 
		{\bf (a)} and {\bf (b)}:
		\begin{eqnarray}
			\label{eq_WZFZW_product_transformation}
			\lefteqn{ \Wk{F}{-t}\ \Zk{F}{s} F \TRANSPOSE{\Zk{F}{t}} \Wk{F}{-s} \ \ \stackrel{(\ref{eq_F_Wk_equals_Zk_F})}{=}}
			&   &   \\
			\nonumber
			&  &
			\Wk{F}{-t}\ F \Wk{F}{s} \Zk{F}{-t}\ \Wk{F}{-s} 
					\ \ \stackrel{(\ref{eq_Zs_Wt_equals_Fst_Wt_Zs})}{=}\ \ 
			F_{-t,-s} \cdot \Wk{F}{-t}\ F \Wk{F}{s} \Wk{F}{-s}\ \Zk{F}{-t}
					 \ \ =\ \ 
			F_{-t,-s} \cdot \Wk{F}{-t}\ F \Zk{F}{-t} \ \ =    \\
			\nonumber
			&    &
			F_{-t,-s} \cdot \Wk{F}{-t}\ F \TRANSPOSE{\Zk{F}{t}}
					\ \ \stackrel{(\ref{eq_Wk_F_equals_F_Zk_transposed})}{=}\ \  
			F_{-t,-s} \cdot \Wk{F}{-t}\ \Wk{F}{t} F 
					\ \ =\ \ 
			F_{s,t} \cdot F\ \ ,
		\end{eqnarray}
		where $F_{-t,-s}\ =\ F_{-s,-t}\ =\ \CONJ{F_{s,-t}}\ =\ \CONJ{\CONJ{F_{s,t}}}\ =\ F_{s,t}$\ .
		The equality of the leftmost and rightmost expressions in (\ref{eq_WZFZW_product_transformation})
		is equivalent to (\ref{eq_Zs_F_Zt_transposed_equals_Fst_times_Wt_F_Ws}).

\end{description}
\PROOFend 

Returning to the main course:


\begin{lemma}
	\label{lem_G_F_is_a_symmetry_group_of_F}
	Let $F$ be a Fourier matrix of size $N$ indexed by group $\Igroup{F}$.

	Group $\GsymmetryGROUP{F}$ (introduced in Lemma \ref{lem_G_F_is_a_group}) is a symmetry group 
	of the unitary matrix  $1/\sqrt{N} \cdot F$ in the sense of Karabegov (Definition \ref{def_symmetry_group}). 
	In other words, there exist maps from $\GsymmetryGROUP{F}$ into the group of enphased permutation matrices 
	$\enphasedPERMS$:  
	\begin{eqnarray}
		\label{eq_S_F_map}
		(s,t, \theta)   &   \stackrel{\Smap{F}}{\longrightarrow}  &   \Sformula{F}{s}{t}{\theta}\ \ ,  \\
		\label{eq_T_F_map}
		(s,t, \theta)   &   \stackrel{\Tmap{F}}{\longrightarrow}   &   \Tformula{F}{s}{t}{\theta}\ \ ,
	\end{eqnarray}
	such that
	\begin{description}
		\item[(a)]
			$\Smap{F}$ and $\Tmap{F}$ are faithful representations (monomorphisms) of $\GsymmetryGROUP{F}$ in
			$\enphasedPERMS$.

		\item[(b)]
			Representations $\Smap{F}$ and $\Tmap{F}$ are equivalent, being intertwined by unitary $1/\sqrt{N} \cdot F$:
			\begin{equation}
				\label{eq_Smap_Tmap_intertwined}
				\Tmap{F}(s,t,\theta) \ \ =\ \ 
				\INVERSE{ \frac{1}{\sqrt{N}} F }
					\cdot 
				\Smap{F}(s,t,\theta) 
					\cdot 
				\left( \frac{1}{\sqrt{N}} F \right)     \ \ .
			\end{equation}
		
	\end{description}
\end{lemma}

\PROOFstart  
First we shall see that map $\Smap{F}$ is a homomorphism:
\begin{eqnarray}
	\label{eq_Smap_is_homomorphism}
	\lefteqn{\Smap{F}(s_1, t_1, \theta_1)\ \cdot\ \Smap{F}(s_2, t_2, \theta_2)\ \ =}
	&  &    \\
	\nonumber
	&  &
	\Sformula{F}{s_1}{t_1}{\theta_1}\ \cdot\ \Sformula{F}{s_2}{t_2}{\theta_2}
			\ \ \stackrel{(\ref{eq_Zs_Wt_equals_Fst_Wt_Zs})}{=}\ \  
	\theta_1 \theta_2 F_{t_1,s_2} \cdot \Wk{F}{s_1} \Wk{F}{s_2} \Zk{F}{t_1} \Zk{F}{t_2}
			\ \ \stackrel{L.\ref{lem_maps_k_Wk_and_k_Zk_are_monomorphisms}}{=}\ \ 
	\Sformula{F}{s_1 + s_2}{t_1 + t_2}{\theta_1 \theta_2 F_{s_2,t_1}} 
			\ \ =      \\
	\nonumber
	&    &
	\Smap{F}
		\left(
			s_1 + s_2,\ 
			t_1 + t_2,\ 
			\theta_1 \theta_2 F_{s_2,t_1}
		\right)
	\ \ =\ \ 
	\Smap{F}( (s_1, t_1, \theta_1) \cdot (s_2, t_2, \theta_2) )\ \ ,
\end{eqnarray}
where Lemma \ref{lem_maps_k_Wk_and_k_Zk_are_monomorphisms}, 
		 Lemma \ref{lem_computations_with_F_Wk_Zk} {\bf (b)} and
		the symmetry of $F$ have been used.

Homomorphism $\Smap{F}$, composed with the unitary similarity homomorphism 
\begin{equation}
	\label{eq_unitary_similarity_homomorphism_with_F}
	X\ \  \longrightarrow\ \  
	\INVERSE{\frac{1}{\sqrt{N}}  F}  
		\cdot 
	X 
		\cdot
	\left(  \frac{1}{\sqrt{N}}  F \right)
\end{equation}	
as in (\ref{eq_Smap_Tmap_intertwined}),
will be another homomorphism, $\Tmap{F}$\ :
\begin{eqnarray}
	\label{eq_converting_Smap_to_Tmap}
	\lefteqn{  
		\frac{1}{\sqrt{N}} F^{*} \cdot 
		\left(  \Sformula{F}{s}{t}{\theta} \right)  \cdot  
		\frac{1}{\sqrt{N}} F \ \ =
	}
	&    &    \\
	\nonumber
	&    &
	\frac{1}{N} \theta \cdot   F^* \Wk{F}{s} \Zk{F}{t} F 
			\ \ \stackrel{(\ref{eq_F_Wk_equals_Zk_F})}{=}\ \ 
	\frac{1}{N} \theta \cdot   F^* \Wk{F}{s}  F  \Wk{F}{t} 
			\ \ \stackrel{(\ref{eq_Wk_F_equals_F_Zk_transposed})}{=}\ \ 
	\frac{1}{N} \theta \cdot   F^* F  \TRANSPOSE{\Zk{F}{s}}  \Wk{F}{t}
			\ \ \stackrel{(\ref{eq_Zk_homomorphism}),(\ref{eq_Zs_Wt_equals_Fst_Wt_Zs})}{=}        \\
	\nonumber
	&     &
	\frac{1}{N} \theta \cdot 
	\left( N \cdot I_N \right) \cdot
	\left( F_{-s,t} \cdot \Wk{F}{t} \Zk{F}{-s} \right)
	\ \ =\ \ 
	\Tformula{F}{s}{t}{\theta}
	\ \ =\ \ 
	\ \ \Tmap{F}(s,t,\theta)\ \ ,
\end{eqnarray}
where Lemma \ref{lem_maps_k_Wk_and_k_Zk_are_monomorphisms}, 
		 Lemma \ref{lem_computations_with_F_Wk_Zk} {\bf (a,b)} and 
		the properties of $F$ have been used, $I_N$ is the identity matrix.
Thus we have explained item {\bf (b)} of the lemma.

Are $\Smap{F}$, $\Tmap{F}$ monomorphisms? 
	If $\Smap{F}(s_1, t_1, \theta_1) \ =\ \Smap{F}(s_2, t_2, \theta_2)$, 
that is $\Sformula{F}{s_1}{t_1}{\theta_1}\ =\ \Sformula{F}{s_2}{t_2}{\theta_2}$, 
then the underlying permutation matrices, $\Zk{F}{t_1}$ and $\Zk{F}{t_2}$, are equal, 
hence $t_1\ =\ t_2$ from Lemma \ref{lem_maps_k_Wk_and_k_Zk_are_monomorphisms}.
From this $\theta_1 \cdot \Wk{F}{s_1} \ =\ \theta_2 \cdot \Wk{F}{s_2}$ which is true only when
$s_1\ =\ s_2$, because $\theta_1, \theta_2 \neq 0$ and the rows of $F$ are independent.
Finally we conclude that also $\theta_1\ =\ \theta_2$.
	$\Smap{F}$ is thus a monomorphism and so is $\Tmap{F}$ as a composition of two monomorphisms,
$\Smap{F}$ and (\ref{eq_unitary_similarity_homomorphism_with_F}). The proof of item {\bf (a)}
has been completed.
\PROOFend   

A consequence of $\GsymmetryGROUP{F}$ being a symmetry group of $1/\sqrt{N} \cdot F$ is that 
a monomorphic image of $\GsymmetryGROUP{F}$ forms a part of the stabilizer of $F$ in $\PP$, 
which follows  from  Lemma \ref{lem_identifying_symmetry_group_with_stabilizer_subgroup} {\bf a)}.
This is also the subject of Corollary \ref{cor_ST_map_image_of_G_F_is_subgroup_of_PP_stabilizer_of_F} below.
 However, 
the reader could easily prove that the considered image of $\GsymmetryGROUP{F}$  is a subgroup of the
stabilizer using only formulas from Lemma \ref{lem_computations_with_F_Wk_Zk}, without considering 
$\GsymmetryGROUP{F}$ at all. A formal proof is ommited.


\begin{corollary}
	\label{cor_ST_map_image_of_G_F_is_subgroup_of_PP_stabilizer_of_F}
	Map
	\begin{equation}
		\label{eq_ST_map}
		(s,t,\theta)\ \ \longrightarrow\ \ 
		\left(\ \Smap{F}(s,t,\theta),\ \Tmap{F}(s,t,,\theta)\ \right) \ \ =\ \ 
		\left(\ \Sformula{F}{s}{t}{\theta},\ \Tformula{F}{s}{t}{\theta}\ \right)
	\end{equation}
	is a monomorphism between $\GsymmetryGROUP{F}$ and $\PP$.

	Its image is a subgroup of the stabilizer 
	$\STAB{\PP}{1/\sqrt{N} \cdot F}\ =\ \STAB{\PP}{F}$.
\end{corollary}

The subgroup of the stabilizer of a Fourier matrix $F$ in $\PP$ pointed at in 
Corollary \ref{cor_ST_map_image_of_G_F_is_subgroup_of_PP_stabilizer_of_F}
turns out to be the stabilizer of $F$ in group $\PzPz{F} \subset \PP$, 
where
\begin{equation}
	\label{eq_enphased_shifts_definition}
	\enphasedSHIFTS{F}\ \ =\ \ 
	\left\{ 
		D\cdot \Zk{F}{k}:\ 
		D\ \mbox{unitary diagonal},\ k \in \Igroup{F} 
	\right\},
\end{equation}
the group of all the enphased  matrices $\Zk{F}{k}$. The lemma below justifies this fact and the use of symbol
$\STAB{\PzPz{F}}{F}$ to denote the considered subgroup of $\STAB{\PP}{F}$.


\begin{lemma}
	\label{lem_G_F_image_under_STmap_is_PzPz_stabilizer_of_F}
	For a Fourier matrix $F$
	\begin{equation}
		\label{eq_G_F_image_under_STmap_is_PzPz_stabilizer_of_F}
		\STAB{\PzPz{F}}{F}\ \ =\ \ 
		\left\{\ 
			\left(\ \Smap{F}(s,t,\theta),\ \Tmap{F}(s,t,\theta)\  \right)\
			:\ \ 
			(s,t,\theta) \in \GsymmetryGROUP{F}\ 
		\right\}\ \ .
	\end{equation}
\end{lemma}

\PROOFstart 
Definitely the considered subgroup 
(the image of $\GsymmetryGROUP{F}$, on the right of (\ref{eq_G_F_image_under_STmap_is_PzPz_stabilizer_of_F})) 
is contained in $\STAB{\PzPz{F}}{F}$, 
because of the form of its elements and the last statement in Corollary \ref{cor_ST_map_image_of_G_F_is_subgroup_of_PP_stabilizer_of_F}.

On the other hand, if $\left( D' \Zk{F}{k},\   D'' \Zk{F}{l} \right)\ \in\  \STAB{\PzPz{F}}{F}$, where $D'$, $D''$ are  unitary
diagonal, then 
\begin{equation}
	\label{eq_F_shift_phased}
	F
		\ \ =\ \ 
	D' \Zk{F}{k}  \cdot F \cdot   \TRANSPOSE{\Zk{F}{l}}  \CONJ{D''}
		\ \ \stackrel{(\ref{eq_Zs_F_Zt_transposed_equals_Fst_times_Wt_F_Ws})}{=}\ \ 
	F_{k,l} \cdot  \left( D' \Wk{F}{l} \right)   F   \left( \Wk{F}{k} \CONJ{D''} \right)\ \ ,
\end{equation}
where Lemma \ref{lem_computations_with_F_Wk_Zk}  {\bf (c)} has been used.

$F$ has the $\GROUPzero{\Igroup{F}}$th row and column filled all with $1$'s, hence the same can be said about the right product in 
(\ref{eq_F_shift_phased}).  For this\ \    $D' \Wk{F}{l}$\  \  and\ \   $\Wk{F}{k} \CONJ{D''}$,\ \ which are unitary diagonal, must be scalar matrices,
that is\ \  
	$D' \Wk{F}{l}\ =\ \theta_1 \cdot I_N$
\ \ and\ \ 
	$\Wk{F}{k} \CONJ{D''}\ =\ \theta_2 \cdot I_N$, 
where
	$|\theta_1|\ =\ |\theta_2|\ =\ 1$.
Moreover, the scalars $\theta_1$, $\theta_2$ must satisfy $F_{k,l} \theta_1 \theta_2\ =\ 1$.
These requirements can be equivalently expressed as:
\begin{equation}
	\label{eq_Dprim_Dbis_necessary_forms}
	D'\ \ =\ \ \theta \cdot \Wk{F}{-l}\ \ ,
	\ \ \ \ \ \ \ \ 
	D''\ \ =\ \ \theta F_{k,l} \cdot \Wk{F}{k}
\end{equation}
for some unimodular $\theta$, so the considered stabilizing pair $\left( D' \Zk{F}{k},\ D'' \Zk{F}{l} \right)$ now reads:
\begin{equation}
	\label{eq_shift_phasing_stabilizing_pair_necessary_form}
	\left( D' \Zk{F}{k},\ D'' \Zk{F}{l} \right)
	\ \ =\ \ 
	\left( \theta \cdot \Wk{F}{-l} \Zk{F}{k},\  \theta F_{k,l} \cdot \Wk{F}{k} \Zk{F}{l} \right)
	\ \ =\ \ 
	\left(\ \Sformula{F}{(-l)}{k}{\theta},\ \Tformula{F}{(-l)}{k}{\theta}\ \right)\ \ ,
\end{equation}
which is $\left( \Smap{F}(-l,k,\theta),\ \Tmap{F}(-l,k,\theta) \right)$, an element of the examined 
subgroup of $\STAB{\PP}{F}$.
    Thus  we have shown that $\STAB{\PzPz{F}}{F}$  is contained in the considered subgroup.
\PROOFend  

The next two lemmas are mere calculations to be used in the proof of the third one, stating
that $\MAP{\PsPs}{F}{G} \STAB{\PzPz{F}}{F}\ =\ \STAB{\PzPz{G}}{G} \MAP{\PsPs}{F}{G}$ 
for any two permutation equivalent Fourier matrices $F$ and $G$
(recall Corollary \ref{cor_Fouriers_perm_equivalent_iff_their_indexing_groups_isomorphic} on the permutation equivalence
and Corollary \ref{cor_permutational_MAP_FG} {\bf (b)} on $\MAP{\PsPs}{F}{G}$).


\begin{lemma}
	\label{lem_permutation_similarity_of_Wk}
	Let $F$, $G$ be permutation equivalent Fourier matrices indexed by isomorphic groups 
	$\Igroup{F}$, $\Igroup{G}$. Let $\ISOM{\alpha}:\Igroup{G} \longrightarrow \Igroup{F}$
	be an isomorphism. Then for the permutation matrix $\Pro{\ISOM{\alpha}}$ there holds:
	\begin{equation}
		\label{eq_permutation_similarity_of_Wk}
		\Pro{\ISOM{\alpha}} \cdot \Wk{F}{s} \cdot \Pro{\ISOM{\alpha}}^T
			\ \ =\ \ 
		\Wk{G}{\hMAP{F}{G}{\ISOM{\alpha}}(s)}
			\ \ \ \ \ \ \mbox{for any}\ \ 
		s \in \Igroup{F}\ \ ,
	\end{equation}
	where $\hMAP{F}{G}{\ISOM{\alpha}}:\Igroup{F} \longrightarrow \Igroup{G}$ 
	is an isomorphism (by Lemma \ref{lem_h_FG_map_properties} {\bf (a)})
	defined in Lemma \ref{lem_P_determines_R_for_permutation_equivalence_PFRequalsG}.
\end{lemma}

\PROOFstart 
The diagonal of the left hand side of equality (\ref{eq_permutation_similarity_of_Wk}) is 
built from the consecutive entries of the column vector:
\begin{equation}
	\label{eq_diagonal_vector_of_PWPtransposed}
	\Pro{\ISOM{\alpha}} F_{:,s}
		\ \ =\ \ 
	\ELEMENTof{\Pro{\ISOM{\alpha}}  F}{:}{s}
		\ \ =\ \ 
	\ELEMENTof{G \Pro{\hMAP{F}{G}{\ISOM{\alpha}}}^T}{:}{s}
		\ \ =\ \ 
	G_{:,\ \hMAP{F}{G}{\ISOM{\alpha}}(s)}\ \ ,
\end{equation}
the $\hMAP{F}{G}{\ISOM{\alpha}}(s)$-th column of $G$, where we have used the equality
\begin{equation}
	\label{eq_PFRequalsG_equality_of_lemma2}
	\Pro{\ISOM{\alpha}}  \cdot F \cdot  \Pro{\INVERSE{\hMAP{F}{G}{\ISOM{\alpha}}}}^T
			\ =\ 
	\Pro{\ISOM{\alpha}}  \cdot F \cdot  \Pro{\hMAP{F}{G}{\ISOM{\alpha}}}
		\ \ \stackrel{this\ one}{=}\ \ 
	G
\end{equation}
guaranteed by Lemma \ref{lem_P_determines_R_for_permutation_equivalence_PFRequalsG}.
Therefore the resulting diagonal matrix is equal to $\Wk{G}{\hMAP{F}{G}{\ISOM{\alpha}}(s)}$ in
accordance with Definition \ref{def_Wk_Zk_matrices} {\bf (a)}.
\PROOFend  


\begin{lemma}
	\label{lem_permutation_similarity_of_Zt}
	Let $F$, $G$ be permutation equivalent Fourier matrices indexed by isomorphic 
	groups $\Igroup{F}$, $\Igroup{G}$. 
	Let $\ISOM{\alpha}, \ISOM{\beta} :\Igroup{G} \longrightarrow \Igroup{F}$ be isomorphisms.
	Then for the corresponding permutation matrices $\Pro{\ISOM{\alpha}}$, $\Pro{\ISOM{\beta}}$
	and given $t \in \Igroup{F}$ :
	\begin{eqnarray}
		\label{eq_permutation_equivalence_of_Zt}
		\Pro{\ISOM{\alpha}} \cdot \Zk{F}{t} \cdot \Pro{\ISOM{\beta}}^T
			\ \ =\ \ 
		\Zk{G}{u}
		&   \mbox{for some}   &
		u \in \Igroup{G}                        \\
		\nonumber
		&       \Updownarrow     &     \\
		\nonumber
		\alpha    &    =    &    \beta\ \ ,
	\end{eqnarray}
	and then $\Zk{F}{t}$ and  $\Zk{G}{\ISOM{\alpha}^{-1}(t)}$ are permutation similar:
	\begin{equation}
		\label{eq_permutation_similarity_of_Zt}
		\Pro{\ISOM{\alpha}} \cdot \Zk{F}{t} \cdot  \Pro{\ISOM{\alpha}}^T
			\ \ =\ \ 
		\Zk{G}{\ISOM{\alpha}^{-1}(t)}\ \ .
	\end{equation}
\end{lemma}

\PROOFstart  
	First we show that (\ref{eq_permutation_similarity_of_Zt}) holds, 
	using Definition \ref{def_Wk_Zk_matrices} {\bf (b)} of $\Zk{F}{t}$
	and property (\ref{eq_Pro_multiplication_rule}):
	\begin{eqnarray}
		\label{eq_perm_sim_of_Zt_transformation}
		\lefteqn{ \Pro{\ISOM{\alpha}} \cdot \Zk{F}{t} \cdot \Pro{\ISOM{\alpha}}^T \ \ =}
		&  &  \\
		\nonumber
		&  &
		\Pro{\ISOM{\alpha}} \cdot \Pro{(i \rightarrow i+t)} \cdot \Pro{\ISOM{\alpha}^{-1}}
			\ \ =\ \ 
		\Pro{\ISOM{\alpha}^{-1} (i \rightarrow i+t) \ISOM{\alpha}}
			\ \ =\ \ 
		\Pro{\left(j \rightarrow j + \ISOM{\alpha}^{-1}(t)\right)}
			\ \ =\ \ 
		\Zk{G}{\ISOM{\alpha}^{-1}(t)}\ \ ,
	\end{eqnarray}
	because for any isomorphism $\ISOM{\alpha}:\ \Igroup{G} \longrightarrow \Igroup{F}$
	\begin{equation}
		\label{eq_I_G_isom_I_F_property}
		\left(
			\ISOM{\alpha}^{-1} (i \rightarrow i+t) \ISOM{\alpha}
		\right)( j \in \Igroup{G} )
			\ \ =\ \ 
		\ISOM{\alpha}^{-1}\left(  \ISOM{\alpha}(j )+ t \right)
			\ \ =\ \ 
		j + \ISOM{\alpha}^{-1}(t)
			\ \ 
		\in \Igroup{G}\ \ .
	\end{equation}
In this way we have also proved implication $\Uparrow$ of the lemma.

On the other hand, let 
$\Pro{\ISOM{\alpha}} \cdot \Zk{F}{t} \cdot \Pro{\ISOM{\beta}}^T\  =\ \Zk{G}{u}$. 
Then using (\ref{eq_permutation_similarity_of_Zt})  and properties (\ref{eq_Pro_inversion_rule}), (\ref{eq_Pro_multiplication_rule}):
\begin{equation}
	\label{eq_PaZtPb_transformation}
	\left( 
		\Pro{\ISOM{\alpha}} \cdot \Zk{F}{t} \cdot \Pro{\ISOM{\alpha}}^T
	\right)
	\left(
		\Pro{\ISOM{\alpha}}
		\Pro{\ISOM{\beta}}^T
	\right)
		\ \ =\ \ 
	\Zk{G}{\ISOM{\alpha}^{-1}(t)}     \cdot     \Pro{\ISOM{\beta}^{-1} \ISOM{\alpha}}
		\ \ =\ \ 
	 \Zk{G}{\ISOM{\alpha}^{-1}(t)}    \cdot     \Pro{\ISOM{\eta}}
		\ \ =\ \ 
	\Zk{G}{u}\ \ ,
\end{equation}
from which  
\begin{equation}
	\label{eq_Peta_result}
	\Pro{\ISOM{\eta}}    \ \ =\ \     \Zk{G}{x}       \ \ ,
\end{equation}
where we have used Lemma \ref{lem_maps_k_Wk_and_k_Zk_are_monomorphisms}, and where
$\ISOM{\eta}\ =\ \ISOM{\beta}^{-1} \ISOM{\alpha} : \Igroup{G} \longrightarrow \Igroup{G}$ is an isomorphism, 
while $x\ =\ u - \ISOM{\alpha}^{-1} (t)\ \in \Igroup{G}$. 

The last equality (\ref{eq_Peta_result}) amounts 
to the equality between the permutations designating the involved permutation matrices:
$\eta\ =\ (j \rightarrow j + x)$.  
This causes that
$\GROUPzero{\Igroup{G}}
	\ =\ 
\ISOM{\eta}\left(\GROUPzero{\Igroup{G}}\right)
	\ =\ 
\GROUPzero{\Igroup{G}} + x$, so $x$ must be neutral in $\Igroup{G}$. 
For the remaining values of $\ISOM{\eta}$ we then
have $\ISOM{\eta}(j)\ =\ j+x\ =\ j$, thus $\ISOM{\eta}$ is just an identity on $\Igroup{G}$, 
that is $\ISOM{\alpha}\ =\ \ISOM{\beta}$. Implication $\Downarrow$ of the lemma has been explained.
\PROOFend  


\begin{lemma}
	\label{lem_MAP_STAB_equals_STAB_MAP}
	Let $F$ and $G$ be permutation equivalent Fourier matrices.

	Let $\ISOM{\psi'},\ \ISOM{\psi''}: \ \Igroup{G} \longrightarrow \Igroup{F}$ be isomorphisms such that
	$\ISOM{\psi''}\ =\ \INVERSE{\hMAP{F}{G}{\ISOM{\psi'}}}$, 
	$\ISOM{\psi'}\ =\ \INVERSE{\hMAP{F}{G}{\ISOM{\psi''}}}$   
	(map $\ISOM{\psi} \longrightarrow \hMAP{F}{G}{\ISOM{\psi}}$ is defined in 
	Lemma \ref{lem_P_determines_R_for_permutation_equivalence_PFRequalsG}, 
	see also Lemma \ref{lem_h_FG_map_properties} {\bf (a),(e)}).

	Let $s,t \in \Igroup{F}$ and let $\theta$ be a unimodular complex number.

	Consider the pairs:
	\begin{eqnarray}
		\nonumber   
		\left(     \Pro{\ISOM{\psi'}},\      \Pro{\ISOM{\psi''}}    \right) 
		&	\in\   &    \MAP{\PsPs}{F}{G}  \ \ ,  \\
		\nonumber     
		\left( \Smap{F}(s,t,\theta),\ \Tmap{F}(s,t,\theta) \right)
		\ \ =\ \ 
		\left( \Sformula{F}{s}{t}{\theta},\ \Tformula{F}{s}{t}{\theta} \right)
		&   \in   &   \STAB{\PzPz{F}}{F}   \ \ ,  \\
		\nonumber      
		\left( \Smap{G}(s',t',\theta),\ \Tmap{G}(s',t',\theta) \right)
		\ \ =\ \ 
		\left( \Sformula{G}{s'}{t'}{\theta},\ \Tformula{G}{s'}{t'}{\theta} \right)
		&   \in   &   \STAB{\PzPz{G}}{G}    \ \ ,
	\end{eqnarray}
	where the above inclusions follow from Corollary 
	\ref{cor_permutational_MAP_FG} {\bf (b)} and Lemma
	\ref{lem_G_F_image_under_STmap_is_PzPz_stabilizer_of_F}, and where 
	\begin{eqnarray}
		\label{eq_s_prim_definition}
		s'   &    =     &     \hMAP{F}{G}{\ISOM{\psi'}}(s)    \ \ =\ \    \INVERSE{\ISOM{\psi''}}(s)    \ \ ,   \\
		\label{eq_t_prim_definition}
		t'  &     =       &    \hMAP{F}{G}{\ISOM{\psi''}}(t)    \ \ =\ \    \INVERSE{\ISOM{\psi'}}(t)    \ \ .
	\end{eqnarray}
	Then the below element of $\MAP{\PP}{F}{G}$ can be written in two ways:
	\begin{equation}
		\label{eq_PR_ST_equals_SnewTnew_PR}
		\left( \Pro{\ISOM{\psi'}},\ \Pro{\ISOM{\psi''}} \right)
			\cdot
		\left( \Smap{F}(s,t,\theta),\ \Tmap{F}(s,t,\theta) \right)
				\ \ =\ \ 
		\left( \Smap{G}(s',t',\theta),\ \Tmap{G}(s',t',\theta) \right)
			\cdot
		\left( \Pro{\ISOM{\psi'}},\ \Pro{\ISOM{\psi''}} \right)
	\end{equation}
	and, as a consequence,
	\begin{equation}
		\label{eq_MAP_STAB_equals_STAB_MAP}
		\MAP{\PsPs}{F}{G}     \cdot     \STAB{\PzPz{F}}{F}
			\ \ =\ \ 
		\STAB{\PzPz{G}}{G}     \cdot     \MAP{\PsPs}{F}{G}
	\end{equation}
\end{lemma}

\PROOFstart  
To prove (\ref{eq_PR_ST_equals_SnewTnew_PR}) we calculate separately for the left and 
right component of the considered pair.
\begin{eqnarray}
	\label{eq_PR_ST_left_component}
	\lefteqn{   \Pro{\ISOM{\psi'}}    \cdot   \Smap{F}(s,t,\theta)    \ \ =}
	&    &    \\
	\nonumber
	&   &
	\Pro{\ISOM{\psi'}}    \cdot   \left( \Sformula{F}{s}{t}{\theta} \right)
			\ \ =\ \ 
	\theta
		\cdot
	\left(     \Pro{\ISOM{\psi'}}  \Wk{F}{s}  \Pro{\ISOM{\psi'}}^T     \right)
		\cdot
	\left(     \Pro{\ISOM{\psi'}}  \Zk{F}{t}  \Pro{\ISOM{\psi'}}^T     \right)
		\cdot
	\Pro{\ISOM{\psi'}}
			\ \ =\ \             \\
	\nonumber
	&    &
	\theta
		\cdot
	\Wk{G}{\hMAP{F}{G}{\ISOM{\psi'}}(s)}
		\cdot
	\Zk{G}{\INVERSE{\ISOM{\psi'}}(t)}
		\cdot
	\Pro{\ISOM{\psi'}} 
			\ \ =   
	\left(   \Sformula{G}{s'}{t'}{\theta}   \right)     \cdot     \Pro{\ISOM{\psi'}}
			\ \ =\ \ 
	\Smap{G}(s',t',\theta)   \cdot   \Pro{\ISOM{\psi'}}\ \ ,
\end{eqnarray}
where Lemmas \ref{lem_permutation_similarity_of_Wk}  and  \ref{lem_permutation_similarity_of_Zt}
have been used.  Now for the right component:
\begin{eqnarray}
	\label{eq_PR_ST_right_component}
	\lefteqn{   \Pro{\ISOM{\psi''}}    \cdot   \Tmap{F}(s,t,\theta)    \ \ =}
	&    &    \\
	\nonumber
	&   &
	\Pro{\ISOM{\psi''}}    \cdot   \left( \Tformula{F}{s}{t}{\theta} \right)
			\ \ =\ \ 
	\theta \CONJ{F_{s,t}}
		\cdot
	\left(     \Pro{\ISOM{\psi''}}  \Wk{F}{t}  \Pro{\ISOM{\psi''}}^T     \right)
		\cdot
	\left(     \Pro{\ISOM{\psi''}}  \Zk{F}{-s}  \Pro{\ISOM{\psi''}}^T     \right)
		\cdot
	\Pro{\ISOM{\psi''}}
			\ \ =\ \                          \\
	\nonumber
	&    &
	\theta \CONJ{F_{s,t}}
		\cdot
	\Wk{G}{\hMAP{F}{G}{\ISOM{\psi''}}(t)}
		\cdot
	\Zk{G}{\INVERSE{\ISOM{\psi''}}(-s)}
		\cdot
	\Pro{\ISOM{\psi''}}
			\ \ =
	\left(   \Tformula{G}{s'}{t'}{\theta}   \right)     \cdot     \Pro{\ISOM{\psi''}}
			\ \ =\ \ 
	\Tmap{G}(s',t',\theta)   \cdot   \Pro{\ISOM{\psi''}}\ \ ,
\end{eqnarray}
where we have taken the below transformation into account:
\begin{equation}
	\label{eq_G_st_prim_equals_F_st}
	G_{s',t'}    \ \ =\ \     G_{\INVERSE{\ISOM{\psi''}}(s),\ \hMAP{F}{G}{\ISOM{\psi''}}(t)}
				  \ \ =\ \     F_{\ISOM{\psi''}\left( \INVERSE{\ISOM{\psi''}}(s) \right),\ t}
				  \ \ =\ \     F_{s,t}
\end{equation}
which follows from the  property (\ref{eq_h_defining property})  defining $\hMAP{F}{G}{\ISOM{\psi''}}$, provided in 
Lemma \ref{lem_P_determines_R_for_permutation_equivalence_PFRequalsG}.
\PROOFend   

In the lemma below we show that for Fourier matrices $F$ and $G$ being permutation equivalent means the same as being equivalent.
In other words, and to remind the reader the notion of equivalence, this can be stated as follows:
\begin{equation}
	\label{eq_perm_equi_equivalent_to_equi_decribed}
	\exists (P,R) \in \PsPs\ \ \ (P,R)F\ =\ G
						\ \ \ \ \ \ \ \Longleftrightarrow\ \ \ \ \ \ \   
	\exists (S,T) \in \PP\ \ \ (S,T)F\ =\ G.
\end{equation}
This fact was also mentioned in section Conclusions of our older work \cite{PermEqClasses}.
What is more, we provide a factorization of $\MAP{\PP}{F}{G}$, 
the set of all pairs $(S,T) \in \PP$ transforming $F$ into $G$.


\begin{lemma}
	\label{lem_F_G_equivalent_iff_F_G_permutation_equivalent}
	Let $F$ and $G$ be Fourier matrices indexed by groups $\Igroup{F}$ and $\Igroup{G}$.

	\begin{description}
		\item[(a)]
			$F$ and $G$ are equivalent if and only if $F$ and $G$ are permutation equivalent,
			i.e. $\Igroup{F}$ and $\Igroup{G}$ are isomorphic (according to 
			Corollary  \ref{cor_Fouriers_perm_equivalent_iff_their_indexing_groups_isomorphic}).

		\item[(b)]
			If $F$ and $G$ are equivalent, then
			\begin{eqnarray}
				\label{eq_MAP_F_G_in_enph_perms_factorization}
				\MAP{\PP}{F}{G}     &     =      &      \MAP{\PsPs}{F}{G}       \cdot     \STAB{\PzPz{F}}{F}  \\ 
											&    =      &      \STAB{\PzPz{G}}{G}     \cdot     \MAP{\PsPs}{F}{G}\ \ ,
			\end{eqnarray}
			where the second equality follows from the final statement in Lemma \ref{lem_MAP_STAB_equals_STAB_MAP}.
	\end{description}
\end{lemma}

\PROOFstart  
If $F$ and $G$ are permutation equivalent, then obviously they are equivalent.

If $F$ and $G$ are equivalent, then any pair $(P'D',\ P''D'') \in \MAP{\PP}{F}{G}$, where $P',P''$ are permutation matrices
and $D',D''$ are unitary diagonal matrices, satisfies $G\ =\ P'D' \cdot F \cdot \CONJ{D''} \TRANSPOSE{P''}$. 
So $D'  F  \CONJ{D''}$ has a row and column filled all with $1$'s. 
  This means, that $D'  F$ has  a column filled all with $\theta_1$'s, where $|\theta_1|\ =\ 1$, which is only possible 
when $D'\ =\ \theta_1 \cdot \Wk{F}{-k}$  for some $k \in \Igroup{F}$ and applies to the $k$-th column then.
	This also means that $ F  \CONJ{D''}$ has a row filled all with $\theta_2$'s, where again $|\theta_2|\ =\ 1$,
which happens only when $\CONJ{D''}\ =\ \theta_2 \cdot \Wk{F}{-l}$ for some $l \in \Igroup{F}$ and applies to
the $l$-th row in this case. Then, by Lemma \ref{lem_computations_with_F_Wk_Zk} {\bf (c)} and due to the symmetry of $F$,
\begin{equation}
	\label{eq_DFD_is_shifted_F}
	D' F \CONJ{D''}     \ \ =\ \         \theta_1 \theta_2   \cdot    \Wk{F}{-k}  F  \Wk{F}{-l}
							\ \ =\ \         \theta_1 \theta_2 \CONJ{F_{-k,-l}}   \cdot   \Zk{F}{-l}  F  \TRANSPOSE{\Zk{F}{-k}}\ \ ,
\end{equation}
the only constant row and column of which are filled all with 
\begin{equation}
	\label{eq_DFD_constant_row_and_column_parameter}
	\theta_1 \theta_2 \CONJ{F_{-k,-l}}     \ \ =\ \      \theta_1 \theta_2 \CONJ{F_{k,l}}\ \ ,
	\ \ \ \ \ \ \ \ 
	\mbox{hence}\ \ \ \      \theta_1 \theta_2     \ \ =\ \      F_{k,l}\ \ 
\end{equation}
because, on the other hand, they are filled all with $1$'s as it has been said above.

Now since 
$P'D' \cdot F \cdot \CONJ{D''} \TRANSPOSE{P''}   
	\ \ =\ \     
P'  \Zk{F}{-l}  F  \TRANSPOSE{\Zk{F}{-k}}  \TRANSPOSE{P''}
	\ \ =\ \ 
G$,
the pair $\left( P' \Zk{F}{-l},\ P'' \Zk{F}{-k} \right)$ belongs to $\MAP{\PsPs}{F}{G}$, 
that is $F$ and $G$ are permutation  equivalent as claimed at item {\bf (a)} of the lemma.
Further, by Corollary \ref{cor_permutational_MAP_FG},
$\left( P' \Zk{F}{-l},\ P'' \Zk{F}{-k} \right)    \ \ =\ \     \left( \Pro{\ISOM{\psi'}},\ \Pro{\ISOM{\psi''}} \right)$
for some isomorphisms $\ISOM{\psi'}$,\ $\ISOM{\psi''}: \Igroup{G} \longrightarrow \Igroup{F}$ such that 
$\ISOM{\psi''}\ = \ \INVERSE{ \hMAP{F}{G}{\ISOM{\psi'}} }$  
(and consequently $\ISOM{\psi'}\ = \ \INVERSE{ \hMAP{F}{G}{\ISOM{\psi''}} }$).
Hence $( P',\ P'' )\ =\  \left( \Pro{\ISOM{\psi'}},\ \Pro{\ISOM{\psi''}} \right)    \cdot   \left( \Zk{F}{l},\ \Zk{F}{k} \right)$
and, from what has been established above about $D'$, $D''$,
\begin{eqnarray}
	\label{eq_factorization_of_PD_PD_pair_such_that_PD_PD_F_equals_G}
	\lefteqn{ \left(  P' D',\   P'' D'' \right)  \ \ = }
	&   &   \\
	\nonumber
	&    &   
	\left(    \Pro{\ISOM{\psi'}},\     \Pro{\ISOM{\psi''}}    \right)
	\left(    \Zk{F}{l},\ \Zk{F}{k}     \right)
	\left(    \theta_1   \cdot   \Wk{F}{-k},\         \CONJ{\theta_2}   \cdot   \Wk{F}{l}      \right)
		\ \ =        \\
	\nonumber
	&    &
	\left(    \Pro{\ISOM{\psi'}},\     \Pro{\ISOM{\psi''}}    \right)
	\left(    \theta_1  F_{l,-k}    \cdot    \Wk{F}{-k}  \Zk{F}{l},\      \CONJ{\theta_2}  F_{k,l}    \cdot    \Wk{F}{l}  \Zk{F}{k}    \right)
		\ \ =         \\
	\nonumber
	&    &
	\left(    \Pro{\ISOM{\psi'}},\     \Pro{\ISOM{\psi''}}    \right)
	\left(      \Sformula{F}{(-k)}{l}{\CONJ{\theta_2}},\      \Tformula{F}{(-k)}{l}{\CONJ{\theta_2}}      \right)
		\ \ =            \\
	\nonumber
	&       &
	\left(    \Pro{\ISOM{\psi'}},\     \Pro{\ISOM{\psi''}}    \right)
	\left(     \Smap{F}\left( -k,l, \CONJ{\theta_2} \right),\       \Tmap{F}\left( -k,l, \CONJ{\theta_2} \right)       \right)\ \ ,
\end{eqnarray}
where the second equality follows from Lemma \ref{lem_computations_with_F_Wk_Zk} {\bf (b)} and the third one from
(\ref{eq_DFD_constant_row_and_column_parameter}) and the symmetry of $F$.

We have thus shown that any $(P'D',\ P''D'') \in \MAP{\PP}{F}{G}$ belongs to $\MAP{\PsPs}{F}{G} \cdot \STAB{\PzPz{F}}{F}$. Obviously
$\MAP{\PsPs}{F}{G} \cdot \STAB{\PzPz{F}}{F}    \subset     \MAP{\PP}{F}{G}$. We conclude that item {\bf (b)}  of
the lemma is a true statement.
\PROOFend  

The factorizations of $\MAP{\PP}{F}{G}$ at item {\bf (b)} of the lemma just proved considered in the case of $F\ =\ G$
unveil the structure of the stabilizer of $F$ in $\PP$. In fact it is an internal semidirect product of its subgroups, the stabilizers
$\STAB{\PzPz{F}}{F}$ (the normal component) and $\STAB{\PsPs}{F}$. This fact is expressed at items
{\bf (a)}, {\bf (b)}, {\bf (d)} and {\bf (f)} of the Theorem \ref{theor_enph_perms_stabilizer_of_F_structure} below. 
	Before we go to a formal proof we will say shortly that 
{\bf (a)} follows from   Lemma  \ref{lem_F_G_equivalent_iff_F_G_permutation_equivalent} {\bf (b)}, while 
{\bf (d)} from {\bf (a)} and {\bf (c)} which in turn is a consequence of {\bf (a)}, {\bf (b)} and
(\ref{eq_PR_ST_equals_SnewTnew_PR})  in Lemma  \ref{lem_MAP_STAB_equals_STAB_MAP}.
{\bf (e)} is just a consequence of {\bf (a)}, {\bf (b)} and {\bf (d)}.


\begin{theorem}
	\label{theor_enph_perms_stabilizer_of_F_structure}
	Let $F$ be a Fourier matrix of size $N$ indexed by group $\Igroup{F}$.

	\begin{description}
		\item[(a)]
			\begin{eqnarray}
				\label{eq_STAB_F_as_product}
				\STAB{\PP}{F}
				&   =   &
				\STAB{\PsPs}{F}      \cdot     \STAB{\PzPz{F}}{F}        \\
				\nonumber
				&   =   &
				\STAB{\PzPz{F}}{F}         \cdot       \STAB{\PsPs}{F}     
			\end{eqnarray}

		\item[(b)]
			$\STAB{\PsPs}{F}\     \cap\      \STAB{\PzPz{F}}{F}     \ \ =\ \        \{ (I_N,\ I_N) \}$, 
			where $(I_N,\ I_N)$ is the pair of identity matrices, 
			the common neutral element of $\STAB{\PP}{F}$, $\STAB{\PsPs}{F}$ and $\STAB{\PzPz{F}}{F}$.

		\item[(c)]
			Any  $(S,T) \in \STAB{\PP}{F}$ can be decomposed uniquely in two ways:
			\begin{eqnarray}
				\label{eq_STAB_F_element_factorization}
				(S,T) 
				&    =     &
				\left(  \Pro{\ISOM{\varrho}'},\  \Pro{\ISOM{\varrho}''}  \right)
					\cdot
				\left(  \Sformula{F}{s}{t}{\theta},\   \Tformula{F}{s}{t}{\theta}   \right)          \\
				\nonumber
				&    =     &
				\left(  \Sformula{F}{s'}{t'}{\theta},\   \Tformula{F}{s'}{t'}{\theta}   \right)         
					\cdot
				\left(  \Pro{\ISOM{\varrho}'},\  \Pro{\ISOM{\varrho}''}  \right)
			\end{eqnarray}
			into factors belonging to $\STAB{\PsPs}{F}$ and $\STAB{\PzPz{F}}{F}$, in accordance with item {\bf (a)} above,
			where
			\begin{itemize}
				\item
					isomorphisms $\ISOM{\varrho}',\ISOM{\varrho}'':\ \Igroup{F} \longrightarrow \Igroup{F}$ satisfying
					$\ISOM{\varrho}''\ =\ \INVERSE{\hMAP{F}{F}{\ISOM{\varrho}'}}$\ \  and, equivalently,
					$\ISOM{\varrho}'\ =\ \INVERSE{\hMAP{F}{F}{\ISOM{\varrho}''}}$
					($\hMAP{F}{G}{\ISOM{\chi}}$ is defined in Lemma  \ref{lem_P_determines_R_for_permutation_equivalence_PFRequalsG},
					its properties are listed in Lemma \ref{lem_h_FG_map_properties} and Corollary \ref{cor_h_FF_map_properties}),

				\item
					$s,t \in \Igroup{F}$,\ \ and\ \ $\theta$ is a unimodular complex number,

				\item
					$s'\ =\ \INVERSE{\ISOM{\varrho}''}(s)\ =\ \hMAP{F}{F}{\ISOM{\varrho}'}(s)$\ \ \ \ and\ \ \ \ 
					$t'\ =\ \INVERSE{\ISOM{\varrho}'}(t)\ =\ \hMAP{F}{F}{\ISOM{\varrho}''}(t)$.

			\end{itemize}

		\item[(d)]
			$\STAB{\PzPz{F}}{F}$ is a normal subgroup of $\STAB{\PP}{F}$.

		\item[(e)]
			The map\ \ \ \  
			$\STAB{\PsPs}{F}    \ \ni\         \left(  \Pro{\ISOM{\varrho}'},\  \Pro{\ISOM{\varrho}''}  \right)
			\ \ \ \ \longrightarrow\ \ \ \  
			\left(  \Pro{\ISOM{\varrho}'},\  \Pro{\ISOM{\varrho}''}  \right)    \cdot    \STAB{\PzPz{F}}{F}\        
						\in\   \STAB{\PP}{F} /  \STAB{\PzPz{F}}{F}$
			\ \ \ \ 
			is an isomorphism.

		\item[(f)]
			$\STAB{\PsPs}{F}$ is a normal subgroup of  $\STAB{\PP}{F}$ only if $F = \FOURIER{2}$, that is
			$\Igroup{F}\ =\ \ZN{2}$, and then $\STAB{\PsPs}{F}\ =\ \{(I_2,\ I_2)\}$.

	\end{description}
\end{theorem}

\PROOFstart  
Let $\STAB{\PsPs}{F}$ be denoted by $\GROUP{H}$ and let $\STAB{\PzPz{F}}{F}$ be denoted by $\GROUP{G}$.

\begin{description}
	\item[(a)]
		Applying Lemma \ref{lem_F_G_equivalent_iff_F_G_permutation_equivalent}  {\bf (b)} to $F$ and $G\ =\ F$ we get
		precisely factorization (\ref{eq_STAB_F_as_product}) of $\STAB{\PP}{F}$, as $\MAP{\GROUP{E}}{F}{F}\ =\ \STAB{\GROUP{E}}{F}$.

	\item[(b)]
		Let $\GROUP{G} \ni \left( \Sformula{F}{s}{t}{\theta},\ \Tformula{F}{s}{t}{\theta} \right)$ also belong to  $\GROUP{H}$.
		It is then a pair of permutation matrices, so $\theta \cdot \Wk{F}{s}$ and $\theta \CONJ{F_{s,t}} \cdot \Wk{F}{t}$
		are  identity matrices. For that there must be $s\ =\ t\ =\ 0_{\Igroup{F}}$, the neutral element of $\Igroup{F}$, and 
		$\theta\ =\ 1$. ( It is because the $0_{\Igroup{F}}$-th  row of $F$, filled all with $1$'s,  is the only row of $F$ with constant
		phase.) But this means that the considered element of   $\GROUP{H} \cap \GROUP{G}$ is a pair of  identity matrices.

	\item[(c)]
		Factorizations (\ref{eq_STAB_F_element_factorization}) exist thanks to {\bf (a)}, the forms of the factors are forced by
		Corollary \ref{cor_permutational_MAP_FG} {\bf (a)} and Lemma \ref{lem_G_F_image_under_STmap_is_PzPz_stabilizer_of_F}.

		That the factorizations (\ref{eq_STAB_F_element_factorization}) are unique is a trivial consequence of
		$\GROUP{H} \cap \GROUP{G}\ =\ \left\{ \left( I_N, I_N \right) \right\}$ of item {\bf (b)}: if $hg\ =\ h'g'$ then 
		$\INVERSE{h'} h \ =\ g' g^{-1}\ =\ \left( I_N, I_N \right)$, so $g\ = g'$ and $h\ =\ h'$. Analogously for $gh\ =\ g'h'$.

		The equality between the two factorizations in (\ref{eq_STAB_F_element_factorization}) results from 
		equality (\ref{eq_PR_ST_equals_SnewTnew_PR}) in Lemma \ref{lem_MAP_STAB_equals_STAB_MAP} applied to
		Fourier matrices $F$ and $G\ =\ F$. 

	\item[(d)]
		To show that $\GROUP{G}$ is normal in $\GROUP{H} \cdot \GROUP{G}\ =\ \STAB{\PP}{F}$ we use 
		the following consequence of items {\bf (a)} and {\bf (c)}: 
		for any $h \in \GROUP{H}$ and $g \in \GROUP{G}$ there exists $g_h \in \GROUP{G}$ such that
		$hg\ =\ g_h  h$.  Then for any $h \in \GROUP{H}$ and $g \in \GROUP{G}$:
		\begin{itemize}
			\item
				$h g h^{-1} \ =\ g_h  h  h^{-1}\ \in \GROUP{G}$, from which $h \GROUP{G} h^{-1}   \subset   \GROUP{G}$,

			\item
				 $g\ =\ h  \left( h^{-1} g h \right)  h^{-1}\ =\ h  \left( g_{h^{-1}}  h^{-1} h \right)  h^{-1}\ =\ h g_{h^{-1}} h^{-1}$, 
				from which $\GROUP{G}  \subset   h \GROUP{G} h^{-1}$.

		\end{itemize}
		In this way we have shown that $x \longrightarrow h x h^{-1}$ maps $\GROUP{G}$ isomorphically onto itself, 
		$\GROUP{G}\ =\ h \GROUP{G} h^{-1}$.

		Now let us take an arbitrary $hg \in \GROUP{H} \GROUP{G}$. 
		Then $(hg) \GROUP{G} \INVERSE{hg}\ =\ h \left( g \GROUP{G} g^{-1} \right) h^{-1} \ =\ h \GROUP{G} h^{-1} \ =\ \GROUP{G}$,
		which means that $\GROUP{G}$ is normal in $\GROUP{H} \GROUP{G}$.

	\item[(e)]
		The considered map, 
		$\GROUP{H} \ni h\ \longrightarrow\ h\GROUP{G} \in \GROUP{H} \GROUP{G} / \GROUP{G}$ in our notation,
		is  a homomorphism: 
			$h_1 h_2 \GROUP{G} 
					\ =\ 
			h_1 \left( h_2 \GROUP{G} \right) \GROUP{G} 
					\ =\ 
			h_1 \left( \GROUP{G} h_2  \right) \GROUP{G}
					\ =\ 
			\left( h_1 \GROUP{G} \right) \left( h_2 \GROUP{G} \right)$,
		as $\GROUP{G}$ is normal. 
		It is a monomorphism because of statement {\bf (b)} of the lemma:
		if $h_1 \GROUP{G} \ =\ h_2 \GROUP{G}$ for $h_1,h_2 \in \GROUP{H}$ then $h_2^{-1} h_1 \in \GROUP{G}$, 
		so $h_2^{-1} h_1 \in \GROUP{H} \cap \GROUP{G}$, that is $h_2^{-1} h_1\ =\ (I_N, I_N)$, hence $h_1\ =\ h_2$.
		The map is an epimorphism because for any $hg \in \GROUP{H} \GROUP{G}$ the coset $(hg) \GROUP{G}$ is
		equal to $h(g \GROUP{G})\ =\ h \GROUP{G}$, the value of the considered map at $h \in \GROUP{H}$.

	\item[(f)]
		Assume that  $\GROUP{H}$ is normal in $\GROUP{H} \GROUP{G}$. Then $g \GROUP{H} \ =\ \GROUP{H} g$
		for any $g \in \GROUP{G}$. This means that for any $g \in \GROUP{G}$ and $h \in \GROUP{H}$ there exists $h_g \in \GROUP{H}$
		such that $gh\ =\ h_g  g$.  

		On the other hand, as we have said in the proof of statement {\bf (d)}, $h g \ =\ g_h h$ (due to {\bf (a)} and {\bf (c)}). 
		On account of the remark in the paragraph above,
		$g_h h \ =\ h_{g_h}  g_h$. Together we have $hg \ =\  h_{g_h}  g_h$, so $g\ =\ g_h$ because of the uniqueness of factorizations 
		(item {\bf (c)}) , hence $hg\  =\ g_h h\ =\ gh$, that is simply $hg\ =\ gh$. Here $h$ and $g$ are arbitrary.

		The fact that $g\ =\ g_h$ in equality $hg\ =\ g_h h$ for a fixed $h$ and any $g$, expressed in terms of symbols used in expression
		(\ref{eq_STAB_F_element_factorization}), translates into the following: 
			for $\ISOM{\varrho}'$, $\ISOM{\varrho}''$ fixed
			the longer bracket in the upper product in (\ref{eq_STAB_F_element_factorization})
					is equal
			to the longer bracket in the lower product, for any $s$, $t$, $\theta$.
		Thus for any $s,t$ the underlying permutation  matrices are equal: $\Zk{F}{t}\ =\ \Zk{F}{t'}$ and $\Zk{F}{-s}\ =\ \Zk{F}{-s'}$.
		This holds, by Lemma \ref{lem_maps_k_Wk_and_k_Zk_are_monomorphisms}, 
		only when $t\ =\ t'$ and $s\ =\ s'$ for any $s,t$, which implies that isomorphisms
		$\ISOM{\varrho}'$, $\ISOM{\varrho}''$ are identities on $\Igroup{F}$. 
		Or that $\left( \Pro{\ISOM{\varrho}'},\ \Pro{\ISOM{\varrho}''} \right)\ =\ (I_N, I_N)$

		The element $h$ fixed in the above reasoning, 
		 representing $\left( \Pro{\ISOM{\varrho}'},\ \Pro{\ISOM{\varrho}''} \right)$, has been chosen arbitrarily. 
		We conclude that every $h \in \GROUP{H}$, a group assumed to be normal in $\GROUP{H} \GROUP{G}$, 
		is equal to $(I_N, I_N)$: $\STAB{\PsPs}{F}\ =\ \GROUP{H}\ =\ \{ (I_N, I_N)\}$.

		This is only possible when $F\ =\ \FOURIER{2}$, that is $\Igroup{F}\ =\ \ZN{2}$, because only then the group of 
		automorphisms on $\Igroup{F}$, $\ISOMORPHISMS{F}{F}$, to which $\ISOM{\varrho}'$, $\ISOM{\varrho}''$ parametrizing
		the elements of $\GROUP{H}$ belong, contains only the identity on $\Igroup{F}$. And indeed such a trivial 
		$\GROUP{H}\ =\ \{ (I_N, I_N)\}$ is a normal subgroup of  $\GROUP{H} \GROUP{G}$.

		If $\Igroup{F}\ =\ \ZN{N_1} \times \ldots \times \ZN{N_r}$ has at least one $k$-th factor of order $N_k > 2$ (we do not 
		consider $\ZN{1}$ corresponding to $\FOURIER{1}\ =\ [1]$ here), then an automorphism on $\Igroup{F}$ different from the identity
		can be defined in this way:
		\begin{equation}
			\label{eq_isom_not_identity_for_some_Zk_of_order_greater_than_2}
			\ISOM{\varphi}\left(   x_1,\ x_2,\ \ldots,\ x_r   \right)
				\ \ =\ \ 
			\left(   x_1,\ \ldots,\ x_{k-1},\ -x_k,\ x_{k+1},\ \ldots,\ x_r   \right)\ \ .
		\end{equation}
		If $\Igroup{F}\ =\ \ZN{2} \times \ldots \times \ZN{2}$ has more than one factor, then an automorphism with this property
		can be taken to be:
		\begin{equation}
			\label{eq_isom_not_identity_for_more_than_one_Z2_only}
			\ISOM{\varphi}\left(   x_1,\ x_2,\ \ldots,\ x_r   \right)
				\ \ =\ \ 
			\left(  x_r,\ x_1,\ x_2,\ \ldots,\ x_{r-1} \right)\ \ .
		\end{equation}

\end{description}
\PROOFend   

Conversely, knowing the stabilizer of a Fourier matrix $F$ (as we define it), and knowing 
that $G$ is permutation equivalent (hence equivalent) to	$F$: 
					$\MAP{\PsPs}{F}{G}\ \neq \emptyset$, 
and that $H$ is only equivalent to $F$:
					$\MAP{\PP}{F}{H}\ \neq \emptyset$,
one can construct these nonempty sets as:
\begin{eqnarray}
	\label{eq_PERMS_MAP_FG_conctructed_from_PERMS_STAB_F}
	\MAP{\PsPs}{F}{G}   &   =   &   \left( P', P'' \right)    \cdot    \STAB{\PsPs}{F}\ \ ,     \\
	\label{eq_ENPHASED_PERMS_MAP_FH_constructed_from_ENPHASED_PERMS_STAB_F}
	\MAP{\PP}{F}{H}   &   =   &   \left( S, T \right)    \cdot    \STAB{\PP}{F}\ \ ,
\end{eqnarray}
where $\left( P', P'' \right)$ is an arbitrary element of $\MAP{\PsPs}{F}{G}$ and $(S,T)$ is an element of $\MAP{\PP}{F}{H}$.
These are elementary consequences of our dealing here with the action ($(X,Y)A\ \stackrel{def}{=} XAY^{-1}$) of 
$\PP$ or $\PsPs$ on the set of square matrices (where all matrices are of a fixed size $N$). For completeness let us add to the above
constructions a number of other ones:
\begin{eqnarray}
	\label{eq_PERMS_MAP_FG_conctructed_from_PERMS_STAB_G}
	\MAP{\PsPs}{F}{G}   &   =   &           \STAB{\PsPs}{G}     \cdot       \left( P', P'' \right)\ \ ,     \\
	\label{eq_ENPHASED_PERMS_MAP_FH_constructed_from_ENPHASED_PERMS_STAB_H}
	\MAP{\PP}{F}{H}   &   =   &           \STAB{\PP}{H}      \cdot       \left( S, T \right)\ \ ,                          \\
	\label{eq_PERMS_STAB_G_conctructed_from_PERMS_STAB_F}
	\STAB{\PsPs}{G}   &    =     &      \left(  P',P''  \right)      \cdot    \STAB{\PsPs}{F}      \cdot      \left(  P',P''  \right)^{-1} \ \ ,    \\
	\label{eq_ENPHASED_PERMS_STAB_H_conctructed_from_ENPHASED_PERMS_STAB_F}
	\STAB{\PP}{H}   &    =     &      \left(  S,T  \right)      \cdot    \STAB{\PP}{F}      \cdot      \left(  S,T  \right)^{-1} \ \ ,
\end{eqnarray}
where $ \left( P', P'' \right)$ and $(S,T)$ as above.

If $G$, for example, is a Fourier matrix more generally defined  
as the one corresponding to an ordered finite abelian group $\GROUP{G}\ =\ \{g_1,g_2,\ldots,g_N \}$
according to the formula:
\begin{equation}
	\label{eq_general_Fourier_matrix_based_on_characters}
	G_{k,l}\  =\ G_{g_k,g_l}\ \    \stackrel{def}{=}\ \          \tau_{g_k}\left(  g_l  \right)\ \ ,
\end{equation}
where
\begin{itemize}
	\item
		 $k,l   \in    \{1,\ldots,N\}$ ordinary indices,\ \ $g_k,g_l \in \GROUP{G}$ group indices, 
	
	\item
		$g \longrightarrow \tau_g$ is an arbitrary isomorphism between $\GROUP{G}$ and its character group
		(group of homomorphisms from $\GROUP{G}$ into $\COMPLEX \setminus \{0\}$ with pointwise multiplication),
\end{itemize}
then $G$ is permutation equivalent (hence equivalent) to some Fourier matrix $F$ of ours, for which its indexing group $\Igroup{F}$
is isomorphic to $\GROUP{G}$. 
	($\GROUP{G}$, as a finite abelian group, is isomorphic to some $\ZN{N_1} \times \ldots \times \ZN{N_r}$.
	To justify that $G$ and $\FOURIER{N_1} \otimes \ldots \otimes \FOURIER{N_r}$ are permutation equivalent then
	we argue in a way similar to that in the proof of Lemma \ref{lem_P_determines_R_for_permutation_equivalence_PFRequalsG}.)
For such $G$ we can use constructions 
	(\ref{eq_PERMS_MAP_FG_conctructed_from_PERMS_STAB_F}) - 
	(\ref{eq_ENPHASED_PERMS_STAB_H_conctructed_from_ENPHASED_PERMS_STAB_F})
(with $H$ replaced by $G$)
to 'translate' the stabilizers from $F$ to $G$ and to build sets $\MAP{...}{F}{G}$ from the stabilizers.

Now let $F$ and $G$ be again Fourier matrices as we define them. 
Let $F$ and $G$ be equivalent 
(which is the same as being permutation equivalent by Lemma \ref{lem_F_G_equivalent_iff_F_G_permutation_equivalent} {\bf (a)}).
Then $\MAP{\PsPs}{F}{G}$ contained in $\MAP{\PP}{F}{G}$ is nonempty. 
Using (\ref{eq_ENPHASED_PERMS_MAP_FH_constructed_from_ENPHASED_PERMS_STAB_F}) with $H$ replaced by $G$ and 
$(S,T)\ =\ \left( P',P'' \right) \in \MAP{\PsPs}{F}{G}$ we get 
(where the second equality holds on account of Theorem \ref{theor_enph_perms_stabilizer_of_F_structure} {\bf (a)},
 while the third due to (\ref{eq_PERMS_MAP_FG_conctructed_from_PERMS_STAB_F}) ) :
\begin{eqnarray}
	\label{eq_ENPHASED_PERMS_MAP_FG_constructed_from_PERMS_MAP_FG_and_ENPHASED_SHIFTS_STAB_F_}
	\lefteqn{ \MAP{\PP}{F}{G}\ \ = }    &   &    \\
	\nonumber
	&   &
	\left( P', P'' \right)    \cdot     \STAB{\PP}{F}
		\ \ =\ \ 
	\left( P', P'' \right)    \cdot     \STAB{\PsPs}{F}    \cdot    \STAB{\PzPz{F}}{F}
		\ \ =\ \           \\
	\nonumber
	&    &
	\MAP{\PsPs}{F}{G}     \cdot      \STAB{\PzPz{F}}{F}\ \ ,
\end{eqnarray}
which is also stated in Lemma \ref{lem_F_G_equivalent_iff_F_G_permutation_equivalent} {\bf (b)}, the first
equality in (\ref{eq_MAP_F_G_in_enph_perms_factorization}). However, we cannot assume that
equivalent $F$ and $G$ are simultaneously permutation equivalent. This is only guaranteed by 
Lemma \ref{lem_F_G_equivalent_iff_F_G_permutation_equivalent} {\bf (a)}. Having proved that and having aquired,
possibly independently, knowledge about the stabilizer of $F$, we are allowed to carry out the calculation
(\ref{eq_ENPHASED_PERMS_MAP_FG_constructed_from_PERMS_MAP_FG_and_ENPHASED_SHIFTS_STAB_F_}).

It is reasonable to suppose that the 'translated' stabilizers 
(\ref{eq_PERMS_STAB_G_conctructed_from_PERMS_STAB_F}),
(\ref{eq_ENPHASED_PERMS_STAB_H_conctructed_from_ENPHASED_PERMS_STAB_F})
are also internal semidirect products just like their right hand side counterparts, the stabilizers of $F$. The more that it is so when
$G$ and $H$ are Fourier matrices as we define them 
(see Theorem \ref{theor_enph_perms_stabilizer_of_F_structure} and the preceding paragraph). 
In fact it is true in general thanks to the below simple fact, a consequence of\ \ $g \longrightarrow a^{-1} g a$\ \ being
an automorphism of some group $\GROUP{A}$:


\begin{lemma}
	\label{lem_translating_semidirect_product}
	Let group $\GROUP{B}$, contained in a larger group $\GROUP{A}$, be an internal semidirect product
	of its subgroups $\GROUP{H},\GROUP{G} \subset \GROUP{B}$. 
	That is to say, 
				$\GROUP{B}\ =\ \GROUP{H} \GROUP{G}$,
				$\GROUP{H} \cap \GROUP{G}$ contains only the neutral element,
				$\GROUP{H}$ is not normal in $\GROUP{B}$, while $\GROUP{G}$ is normal in $\GROUP{B}$.

	Then, for any $a \in \GROUP{A}$,\ \ \  $a^{-1} \GROUP{B} a$ is an internal semidirect product of its subgroups
	$a^{-1} \GROUP{H} a$ and $a^{-1} \GROUP{G} a$, where $a^{-1} \GROUP{H} a$ is not normal in $a^{-1} \GROUP{B} a$, 
	while $a^{-1} \GROUP{G} a$ is normal in $a^{-1} \GROUP{B} a$.
\end{lemma}

%
%

\section{The defect of a Fourier matrix}
	\label{sec_fourier_defect}

Our aim in this section is to characterize and calculate the undephased defect  (\ref{eq_undepased_defect}) 
of a unitary Fourier matrix of size $N$,\ \ $1/\sqrt{N} \cdot F$,\ \ 
where $F$ is given in (\ref{eq_Fourier_Kron_prod}). The notation adopted in the previous sections will be used. Our calculations will be based on 
the characterization of the undephased defect of $U$, $\undephasedDEFECT(U)$, as the multiplicity of $1$ in the spectrum of the Berezin transform
$\Iu{U}$ introduced by Karabegov in \cite{Karabegov}, \cite{Karabegov_old} (see Theorem \ref{theor_unitarity_and_eigenspaces_of_Iu} {\bf e)}).
The spectrum and the eigenspaces of $\Iu{1/\sqrt{N} \cdot F}$ are described in the theorem below. We will prove it in three ways:
\begin{description}
	\item[a)]   We check directly that the provided vectors are indeed eigenvectors of $\Iu{1/\sqrt{N} \cdot F}$.

	\item[b)]   We show that the eigenspaces of $\Iu{1/\sqrt{N} \cdot F}$ are determined by the fact that 
				$\Iu{1/\sqrt{N} \cdot F}$ commutes with all the elements
				of the image of the representation $(S,T) \rightarrow \pairOPERATOR{\leftPERM{S}{T}}{\rightPERM{S}{T}}$ of some
				subgroup of the stabilizer $\STAB{\PP}{1/\sqrt{N} \cdot F}\ =\ \STAB{\PP}{F}$, 
				according to Lemma \ref{eq_properties_of_representations_of_STAB_U} {\bf c)}.  This is a realization of the ideas
				contained in the comments following Lemma \ref{eq_properties_of_representations_of_STAB_U}.

	\item[c)]   We make use of the fact that $\Iu{1/\sqrt{N} \cdot F}$ and the mentioned above elements of the image of the representation
				form a commuting family, whose common eigenbasis is determined uniquely (up to scaling factors) by this family deprived
				of $\Iu{1/\sqrt{N} \cdot F}$.  Here we follow the procedure outlined just after 
				Corollary \ref{eq_Iu_Ju_commute_with_ops_related_to_ST_from_the_stabilizer_of_U}.

\end{description}


\begin{theorem}
	\label{theor_eigenbasis_of_I_F}
	Let $U$ be the unitary $N \times N$  Fourier matrix indexed by group 
	$\Igroup{F}\ \ =\ \ \ZN{N_1} \times \ldots \times \ZN{N_r}$\ :
	\begin{equation}
		\label{eq_U_being_Fourier_matrix}
		U               \ \ =\ \              \frac{1}{\sqrt{N}}   \cdot    F             
					\ \ =\ \                \frac{1}{\sqrt{N}}  \cdot    \FOURIER{N_1} \otimes \ldots \otimes \FOURIER{N_r}      \ \ .
	\end{equation}
	Then the $N^2$ pairs  (where $F_{:,k}$ is the $k$th column of $F$, $k \in \Igroup{F}$):
	\begin{equation}
		\label{eq_eigenvector_eigenvalue_pair}
		\left( F_{:,s}  F_{:,t}^T\ \ ,\ \ F_{s,t} \right) \ \ \ \ \ \ \ \  \mbox{where}\ \ s,t \in \Igroup{F}\ \ ,
	\end{equation}
	are the (eigevector, eigenvalue) pairs of $\Iu{U}$, producing a complete eigenbasis for $\Iu{U}$, orthogonal wrt 
	$\STinnerPRODUCT{}{}$ and $\UinnerPRODUCT{U}{}{}$\ .
\end{theorem}

\PROOFstart 
\begin{description}
\item[a)]
	In this proof we just check that $F_{:,s}  F_{:,t}^T$ is an eigenvector of $\Iu{U}$, using various lemmas and definitions indicated over the
$=$ signs:
\begin{eqnarray}
	\label{eq_I_F_acting_on_eigenvector}
	\lefteqn{\Iu{U}\left(  F_{:,s}  F_{:,t}^T  \right)     \ \ = }     &    &           \\
	\nonumber
	&        &
	\Cu{U}^{-1} \Du{U}\left(  F_{:,s}  F_{:,t}^T  \right)
			\ \ \stackrel{L.\ref{lem_properties_of_Cu_Du_Iu}\mathbf{h)}}{=}\ \ 
	\Cu{U}^{-1}\left(   \frac{1}{\sqrt{N}}F   \DIAGwith{F_{:,t}}  \hermTRANSPOSE{\frac{1}{\sqrt{N}}F}   \DIAGwith{F_{:,s}}   \right)
			\ \ \stackrel{D.\ref{def_Wk_Zk_matrices}\mathbf{a)}}{=}\ \ 
	                                         \\
	\nonumber
	&    &
	\Cu{U}^{-1}\left(   \frac{1}{\sqrt{N}}F   \Wk{F}{t}   \hermTRANSPOSE{\frac{1}{\sqrt{N}}F}    \Wk{F}{s}    \right)
			\ \ \stackrel{L.\ref{lem_computations_with_F_Wk_Zk}\mathbf{a)}}{=}\ \     
	\Cu{U}^{-1}\left(  \Zk{F}{t} \cdot \frac{1}{\sqrt{N}}F  \hermTRANSPOSE{\frac{1}{\sqrt{N}}F}   \Wk{F}{s}   \right)
			\ \ =\ \                                                                    \\
	\nonumber
	&   &
	\Cu{U}^{-1}\left(  \Zk{F}{t}    \Wk{F}{s}   \right)
			\ \ \stackrel{L.\ref{lem_computations_with_F_Wk_Zk}\mathbf{b)}}{=}\ \ 
	\Cu{U}^{-1}\left(   F_{t,s} \cdot   \Wk{F}{s}  \Zk{F}{t}   \right)
			\ \ \stackrel{F=F^T}{=}\ \                                                           \\
	\nonumber
	&      &
	F_{s,t} \cdot \Cu{U}^{-1}\left(     \Wk{F}{s}  \Zk{F}{t}  \cdot  \frac{1}{\sqrt{N}}F  \hermTRANSPOSE{\frac{1}{\sqrt{N}}F}  \right)
			\ \ \stackrel{L.\ref{lem_computations_with_F_Wk_Zk}\mathbf{a)}}{=}\ \ 
	F_{s,t} \cdot \Cu{U}^{-1}\left(     \Wk{F}{s} \cdot  \frac{1}{\sqrt{N}}F   \Wk{F}{t}    \hermTRANSPOSE{\frac{1}{\sqrt{N}}F}  \right)
			                 \\
	\nonumber
	&       &
			\ \ \stackrel{D.\ref{def_Wk_Zk_matrices}\mathbf{a)}}{=}\ \  
	F_{s,t} \cdot \Cu{U}^{-1}\left(   \DIAGwith{F_{:,s}}  \cdot  \frac{1}{\sqrt{N}}F   \DIAGwith{F_{:,t}}  \hermTRANSPOSE{\frac{1}{\sqrt{N}}F}   \right)
			\ \ \stackrel{L.\ref{lem_properties_of_Cu_Du_Iu}\mathbf{h)}}{=}\ \ 
	F_{s,t} \cdot \left( F_{:,s}  F_{:,t}^T  \right)              \ \ ,
\end{eqnarray}
and indeed $F_{s,t}$ is the corresponding eigenvalue.  Moreover, all these eigenvectors of $\Iu{U}$ are pairwise orthogonal wrt $\STinnerPRODUCT{}{}$ and 
$\UinnerPRODUCT{U}{}{}$ because
$	N \UinnerPRODUCT{U}{ F_{:,s}  F_{:,t}^T}{ F_{:,k}  F_{:,l}^T}
			\ =\ 
	\STinnerPRODUCT{ F_{:,s}  F_{:,t}^T}{ F_{:,k}  F_{:,l}^T}
			\ =\ 
	\trace \hermTRANSPOSE{ F_{:,k}  F_{:,l}^T  } \left( F_{:,s}  F_{:,t}^T \right)
			\ =\ 
	\trace     \CONJ{F}_{:.l}     F^{*}_{:,k}     F_{:,s}  F^T_{:,t}
			\ =\ 
	\trace  F^{*}_{:,k}     F_{:,s}     F^T_{:,t}   \CONJ{F}_{:.l}
			\ =\ 
	\left( F^{*}_{:,k}     F_{:,s} \right)    \left(  F^{*}_{:,l}     F_{:,t}   \right)$, 
which is nonzero only when $k=s$ and $l=t$.


\item[b)]
	This proof is based on the representation theory of finite groups. 
The stabilizer of $\frac{1}{\sqrt{N}} F$ in $\PP$:
$\STAB{\PP}{\frac{1}{\sqrt{N}} F}\ \ =\ \ \STAB{\PP}{F}$ is described in Theorem \ref{theor_enph_perms_stabilizer_of_F_structure}, 
but its normal factor  $\STAB{\PzPz{F}}{F}$ is already presented in 
	Lemma \ref{lem_G_F_image_under_STmap_is_PzPz_stabilizer_of_F}
and
	Corollary \ref{cor_ST_map_image_of_G_F_is_subgroup_of_PP_stabilizer_of_F}.

The representation $(S,T) \rightarrow \pairOPERATOR{\leftPERM{S}{T}}{\rightPERM{S}{T}}$, 
introduced in Lemma \ref{lem_representations_of_Stab_U} {\bf e)} and
used in Lemma  \ref{eq_properties_of_representations_of_STAB_U} {\bf c)}, 
on the elements of  $\STAB{\PzPz{F}}{F}$ is expressed by the formula:
\begin{equation}
	\label{eq_perm_repr_of_STAB_Z_F_stabilizer}
	\left(\Smap{F}(-s,t,\theta),\ \Tmap{F}(-s,t,\theta) \right)
		 \ =\  
	\left(  
		\theta \cdot    \Wk{F}{-s}  \Zk{F}{t}     \ ,\ \            
		\theta  \CONJ{F}_{-s,t}  \cdot      \Wk{F}{t}  \Zk{F}{s}   
	\right)
					\longrightarrow
	\pairOPERATOR{  \Zk{F}{t}  }{  \Zk{F}{s}  }      \ \ ,
\end{equation}
so its image on $\STAB{\PzPz{F}}{F}$ is finite.   The image can also be regarded as the image of the faithful (monomorphic)
representation\ \ $(t,s)  \longrightarrow  \pairOPERATOR{  \Zk{F}{t}  }{  \Zk{F}{s}  }$ of the finite group $\Igroup{F} \times \Igroup{F}$.
The image is of course abelian, so by the Schur's lemma in the theory its irreducible invariant spaces are $1$ dimensional. 
If $\SPACE{W} \subset \COMPLEXmatricesNxN{N}$ were invariant irreducible and of dimension greater than $1$, the fact
that $\RESTRICTEDto{\pairOPERATOR{  \Zk{F}{x}  }{  \Zk{F}{y}  }}{\SPACE{W}}$ commutes with all the  
$\RESTRICTEDto{\pairOPERATOR{  \Zk{F}{t}  }{  \Zk{F}{s}  }}{\SPACE{W}}$'s
would imply (Schur) that 
$\RESTRICTEDto{\pairOPERATOR{  \Zk{F}{x}  }{  \Zk{F}{y}  }}{\SPACE{W}}
		\ =\ 
\alpha_{x,y} \RESTRICTEDto{\pairOPERATOR{ I_N }{ I_N }}{\SPACE{W}}$ (a scalar map, $I_N$ is the identity matrix), for all $(x,y) \in \Igroup{F} \times \Igroup{F}$,
so every subspace of $\SPACE{W}$ would then be invariant  which contradicts our assumption on the irreducibility of $\SPACE{W}$.

Let us consider the invariant irreducible spaces $\SPANC\left(F_{:,k}  F_{:,l}^T\right)$ and check that the respective subrepresentations
$(t,s)  \longrightarrow    \RESTRICTEDto{ \pairOPERATOR{  \Zk{F}{t}  }{  \Zk{F}{s}  } }{   \SPANC\left(F_{:,k}  F_{:,l}^T\right) }$
are inequivalent:
\begin{eqnarray}
	\label{eq_ZtZs_operator_acting_on_FkFl_vector}
	\pairOPERATOR{  \Zk{F}{t}  }{  \Zk{F}{s}  }\left(F_{:,k}  F_{:,l}^T\right)
	&     =     &
	\Zk{F}{t}  F_{:,k}  F_{:,l}^T  \TRANSPOSE{\Zk{F}{s}}
			\ \ =\ \ 
	\Zk{F}{t}  F_{:,k}  \TRANSPOSE{\Zk{F}{s}  F_{:,l}}             \\
	\nonumber
	&      \stackrel{L.\ref{lem_computations_with_F_Wk_Zk}\mathbf{a)}}{=}      &
	\left(  F_{t,k} \cdot  F_{:,k}  \right)    \TRANSPOSE{  F_{s,l}  \cdot   F_{:,l} }  
			\ \ =\ \ 
	\left( F_{t,k} F_{s,l} \right)     \cdot    \left(F_{:,k}  F_{:,l}^T\right)                   \ \ ,
\end{eqnarray}
where scalar  $F_{t,k} F_{s,l}$ is the $1 \times 1$ matrix of 
$\RESTRICTEDto{ \pairOPERATOR{  \Zk{F}{t}  }{  \Zk{F}{s}  } }{   \SPANC\left(F_{:,k}  F_{:,l}^T\right) }$. 
Further $(t,s) \longrightarrow F_{t,k} F_{s,l}$ is the matrix representation corresponding to  
$(t,s)  \longrightarrow    \RESTRICTEDto{ \pairOPERATOR{  \Zk{F}{t}  }{  \Zk{F}{s}  } }{   \SPANC\left(F_{:,k}  F_{:,l}^T\right) }$.
This matrix representation is nothing else but a homomorphism  from $\Igroup{F} \times \Igroup{F}$ into 
$\COMPLEX \setminus \{0\}$, or the $(k,l)$th column of $F \otimes F$. Since all these columns are different, 
all subrepresentations 
$(t,s)  \longrightarrow    \RESTRICTEDto{ \pairOPERATOR{  \Zk{F}{t}  }{  \Zk{F}{s}  } }{   \SPANC\left(F_{:,k}  F_{:,l}^T\right) }$
are inequivalent.
Any other irreducible  representation of $\Igroup{F} \times \Igroup{F}$ on a $1$ dimensional (by the Schur's lemma) space $\SPACE{W}$ 
is equivalent to one of those 
$(t,s)  \longrightarrow    \RESTRICTEDto{ \pairOPERATOR{  \Zk{F}{t}  }{  \Zk{F}{s}  } }{   \SPANC\left(F_{:,k}  F_{:,l}^T\right) }$
because, by Lemma \ref{lem_all_unimodular_morphisms_on_F_indexing_group} applied to $\Igroup{F}  \times \Igroup{F}$ indexing
$F \otimes F$, the matrix representation of $\Igroup{F} \times \Igroup{F}$ corresponding to this other representation 
forms a column of $F \otimes F$.

But irreducible invariant subspaces of the image of  $(t,s)  \longrightarrow  \pairOPERATOR{  \Zk{F}{t}  }{  \Zk{F}{s}  }$
can be found in a different way. 
According to  Lemma  \ref{eq_properties_of_representations_of_STAB_U} {\bf c)} $\Iu{U}$ commutes with all the 
$\pairOPERATOR{  \Zk{F}{t}  }{  \Zk{F}{s}  }$'s, so the eigenspaces of $\Iu{U}$ are invariant for the considered representation.
These eigenspaces can further be split (by the Maschke's theorem in the theory) 
into $1$ dimensional irreducible invariant eigensubspaces $\SPACE{W}_i$ 
(see also the comments following Lemma \ref{eq_properties_of_representations_of_STAB_U}) and, as it has been explained
in the last paragraph above, every 
$(t,s)  \longrightarrow    \RESTRICTEDto{ \pairOPERATOR{  \Zk{F}{t}  }{  \Zk{F}{s}  } }{   \SPACE{W}_i }$
is equivalent to some
$(t,s)  \longrightarrow    \RESTRICTEDto{ \pairOPERATOR{  \Zk{F}{t}  }{  \Zk{F}{s}  } }{   \SPANC\left(F_{:,k}  F_{:,l}^T\right) }$.

The theory says (the theorem on the uniqueness of the isotypical decomposition) 
that the direct sums of irreducible invariant subspaces associated with equivalent subrepresentations,
so called isotypic subspaces, are uniquely determined only by the considered representation.
In our case it is  $(t,s)  \longrightarrow  \pairOPERATOR{  \Zk{F}{t}  }{  \Zk{F}{s}  }$ and the above property means that
all the $N^2$ eigensubspaces $\SPACE{W}_i$ are just all the $\SPANC\left(F_{:,k}  F_{:,l}^T\right)$'s. Thus we have found
a full eigenbasis for $\Iu{U}$. Other facts in the theorem have been justified in proof {\bf a)}.

\item[c)]
	The family of operators\ \ 
		$\left\{ \pairOPERATOR{  \Zk{F}{t}  }{  \Zk{F}{s}  }:\ s,t \in \Igroup{F} \right\}
			\ \cup\ 
		\left\{ \Iu{U} \right\}$\ \ 
is a commuting family, because $\left\{ \pairOPERATOR{  \Zk{F}{t}  }{  \Zk{F}{s}  }:\ s,t \in \Igroup{F} \right\}$ is a commuting family
(a consequence of  Lemma \ref{lem_maps_k_Wk_and_k_Zk_are_monomorphisms} and $\Igroup{F}$ being abelian)
and all its members commute with $\Iu{U}$ by  
	Corollary \ref{eq_Iu_Ju_commute_with_ops_related_to_ST_from_the_stabilizer_of_U} {\bf a)}.
Here $\pairOPERATOR{  \Zk{F}{t}  }{  \Zk{F}{s}  }$ is the $\pairOPERATOR{\leftPERM{S}{T}}{\rightPERM{S}{T}}$ operator 
corresponding to \\
	$\left(  
		\theta \cdot    \Wk{F}{-s}  \Zk{F}{t}     \ ,\ \            
		\theta  \CONJ{F}_{-s,t}  \cdot      \Wk{F}{t}  \Zk{F}{s}   
	\right)
			\ \in\ 
	 \STAB{\PzPz{F}}{F}   \ \subset\ \STAB{\PP}{F}\ =\ \STAB{\PP}{\frac{1}{\sqrt{N}}F}$ 
(see the first paragraph of proof {\bf b)}).

The family\ \ 
		$\left\{ \pairOPERATOR{  \Zk{F}{t}  }{  \Zk{F}{s}  }:\ s,t \in \Igroup{F} \right\}
			\ \cup\ 
		\left\{ \Iu{U} \right\}$\ \ 
is a family of operators unitary wrt $\UinnerPRODUCT{U}{}{}$ 
		(see Corollary \ref{eq_Iu_Ju_commute_with_ops_related_to_ST_from_the_stabilizer_of_U} {\bf a)})
and wrt $\STinnerPRODUCT{}{}$ because for $U=\frac{1}{\sqrt{N}}F$ these inner products are proportional. 
Therefore, as a family of operators normal wrt $\STinnerPRODUCT{}{}$, this family has a common eigenbasis orthonormal wrt
$\STinnerPRODUCT{}{}$ (see for example \cite{Johnson}, Theorem 2.5.5).

We know that $\left\{ F_{:,k} F_{:,l}^T :\ k,l \in \Igroup{F} \right\}$ is a common eigenbasis, orthogonal wrt $\STinnerPRODUCT{}{}$,
for  the family $\left\{ \pairOPERATOR{  \Zk{F}{t}  }{  \Zk{F}{s}  }:\ s,t \in \Igroup{F} \right\}$. The orthogonality has been checked at the end
of proof {\bf a)}, while the eigenvector property in (\ref{eq_ZtZs_operator_acting_on_FkFl_vector}). We will prove that every eigenvector
common for this smaller family is proportional to some $F_{:,k} F_{:,l}^T$.

Consider the orthogonal basis  $\left\{ F_{:,k} F_{:,l}^T :\ k,l \in \Igroup{F} \right\}$ of $\COMPLEXmatricesNxN{N}$. Let $\SUPP(V)$, the support of 
$V \in \COMPLEXmatricesNxN{N}$, be defined as
\begin{equation}
	\label{eq_support_wrt_FkFl_basis}
	\SUPP(V)       
			\ \ \stackrel{def}{=}\ \ 
	\left\{ 
			(x,y)\  \in\  \Igroup{F} \times \Igroup{F}             
				\ :\ \            
			\mbox{the coordinate of $V$ at $F_{:,x} F_{:,y}^T$ is nonzero}
	\right\}                \ \ .
\end{equation}
Let $V$ be an eigenvector of all $\pairOPERATOR{  \Zk{F}{t}  }{  \Zk{F}{s}  }$'s, then it is a common eigenvector of all 
$\pairOPERATOR{  \Zk{F}{t}  }{  \Zk{F}{0_{\Igroup{F}}  }}$'s.  Let us fix $t \in \Igroup{F}$. 
According to (\ref{eq_ZtZs_operator_acting_on_FkFl_vector})  $V$ is associated with the  eigenvalue $F_{t,k_t}$ for some $k_t$,
because $F_{0_{\Igroup{F}},l} \ =\ 1$ for any $l$.  Because the eigenspace of $\pairOPERATOR{  \Zk{F}{t}  }{  \Zk{F}{0_{\Igroup{F}}  }}$
associated with this eigenvalue is spanned 
	-- see (\ref{eq_ZtZs_operator_acting_on_FkFl_vector}) -- 
by those $ F_{:,k'} F_{:,l'}^T$ for which $F_{t,k'}\ =\ F_{t,k_t}$, there holds 
	$\SUPP(V) \subset \left\{k':\ F_{t,k'}\ =\ F_{t,k_t}\right\} \times \Igroup{F}$. 
This inclusion can be written for any $t$, so 
	$\SUPP(V) \subset  \bigcap_{t \in \Igroup{F}}  \left(  \left\{k':\ F_{t,k'}\ =\ F_{t,k_t}\right\} \times \Igroup{F}  \right)$, 
where the intersection on the right is not empty because we assume $V \neq \ZEROvect$. So the intersection contains
the set $\{k\} \times \Igroup{F}$ for some $k$. This means that 
	$\forall_{t \in \Igroup{F}}\ F_{t,k} \ =\ F_{t,k_t}$,
that is the $t$-th element of the $k$-th column of $F$ is the eigenvalue of $\pairOPERATOR{  \Zk{F}{t}  }{  \Zk{F}{0_{\Igroup{F}}  }}$
corresponding to $V$. Now, if the above intersection is greater than $\{k\} \times \Igroup{F}$, it contains $\{m\} \times \Igroup{F}$
where $m\ \neq\ k$. Then again column $F_{:,m}$ contains the identical sequence of the eigenvalues $F_{t,k_t}$ of 
$\pairOPERATOR{  \Zk{F}{t}  }{  \Zk{F}{0_{\Igroup{F}}  }}$'s associated with $V$. But it is impossible since columns of $F$ are different
and we wouldn't get the same eigenvalues for every $t$. Hence the above intersection is confined to $\{k\} \times \Igroup{F}$.

In a similar way and using operators $\pairOPERATOR{  \Zk{F}{0_{\Igroup{F}}  }  }{  \Zk{F}{s}}$ we can show that 
	$\SUPP(V) \subset \Igroup{F} \times \{l\}$ for some $l$.
As $V\ \neq\ \ZEROvect$ and, as a result of the above considerations, $\SUPP(V) \subset \{(k,l)\}$, there must be 
$V \ =\ \alpha    \cdot    F_{:,k} F_{:,l}^T$\ \  for some $k,l \in \Igroup{F}$ and $\alpha \neq 0$. 
So the only possible common eigenbasis both for 
	$\left\{ \pairOPERATOR{  \Zk{F}{t}  }{  \Zk{F}{s}  }:\ s,t \in \Igroup{F} \right\}$
and\ \  
	$\left\{ \pairOPERATOR{  \Zk{F}{t}  }{  \Zk{F}{s}  }:\ s,t \in \Igroup{F} \right\}   \  \cup\     \Iu{U}$\ \ 
is that consisting of vectors proportional to $F_{:,k} F_{:,l}^T$'s.
Other facts in the theorem have been justified in proof {\bf a)}.
\end{description}
\PROOFend 

On account of the above Theorem \ref{theor_eigenbasis_of_I_F} as well as Theorem \ref{theor_unitarity_and_eigenspaces_of_Iu} {\bf e)}
a calculation of $\undephasedDEFECT(1/\sqrt{N} \cdot F)$ amounts to counting the number of $1$'s in the spectrum of 
$\Iu{1/\sqrt{N} \cdot F}$ formed by all the entries of $F$. Before we perform this calculation in full generality, 
we shall solve a simpler taks of counting $1$'s in a \DEFINED{primary} Kronecker product of the form:
\begin{equation}
	\label{eq_primary_Fourier_Kronecker_product}
	F\ \ =\ \ 
	\FOURIER{a^{k_1}}      \otimes      \FOURIER{a^{k_2}}    \otimes    \ldots    \otimes   \FOURIER{a^{k_s}}        \ \ ,         
\end{equation}
indexed by group $\Igroup{F}\ =\ \ZN{a^{k_1}} \times \ZN{a^{k_2}} \times \ldots \times \ZN{a^{k_s}}$, where
\begin{itemize}
	\item    $a$\ \ is a prime number,
	\item    $s$\ \ is a number of factors,\ \ $s \geq 1$
	\item    $k_x$'s\ \ are nonnegative integer exponents about which we will most often assume that one of them is nonzero; 
			the case $k_1 =\ldots = k_s = 0$ will be treated separately
			as $\undephasedDEFECT(F)\ =\ \undephasedDEFECT([1])\ =\ 1$ then,
	\item     the size of a primary $F$ shall again be denoted by $N$,\ \  $N \ =\ a^{k_1} \cdot \ldots \cdot a^{k_s}$,
	\item  \begin{equation}
				\label{eq_k_max}
				k_{max} \ \ =\ \ \max_{x \in \{1,\ldots,s\}}  k_x
			\end{equation}
			is the associated maximal exponent,
	\item  \begin{equation}
				\label{eq_number_of_kx_below_value}
				\kxBELOW{l}  \ \ \stackrel{def}{=}\ \ 
				\mbox{the number of $x$ in $\left\{1,  \ldots,  s \right\}$ for which $k_x < l$}    \ \ ,
			\end{equation}
			where $l$ is an integer number ; note that quantities $\kxBELOW{1}$, $\kxBELOW{2}$, ... 
			determine how many $0$'s, $1$'s , ... are among $k_x$'s.
			The quantities $\kxBELOW{l}$ originate from the quantities $\tilde{n}_l \ =\ \kxBELOW{l} + 1$, though defined 
			in a different way, from \cite{PermEqClasses} (Definition 2.16), where they were called 'introduction indices'.
			They fully determine a primary Kronecker product of Fourier matrices up to the order of factors, which has also been stated in 
			Theorem 2.20 in \cite{PermEqClasses}, similarly as it was the case with $\kxBELOW{l}$'s determining $k_x$'s.
\end{itemize}

Each row of $\FOURIER{a^{k_x}}$ contains all the roots of unity of order $a^k$ in equal quantities, where $k$ depends on the row taken and 
belongs to $\{0,1,\ldots,k_x\}$, and all the values of $k$ occur among the rows.  
  As a consequence, a similar property applies to $F$ with the only difference that $k \in \{0,1,\ldots,k_{max}\}$. 
It was written in the introductory part of section  \ref{sec_Fourier_matrices} 
	(see (\ref{eq_multiplying_rows_of_F}) and the surrounding paragraph)
that the group of rows of $F$ is isomorphic to $\Igroup{F}$:  $(i+j) \longrightarrow F_{i+j,:}\ =\ F_{i,:} \HADprod F_{j,:}$. 
Thus $i$ has order $a^m$ in group $\Igroup{F}$  if and only if $F_{i,:}$ has order $a^m$ in the group of rows, which in turn
is equivalent to $F_{i,:}$ being composed of all the types of roots of unity of order $a^m$. Note that the values $a^m$ are the only ones
allowed as possible orders of $i \in \Igroup{F}$ because these orders must divide the number of elements in $\Igroup{F}$, 
$N\ =\ a^{k_1 + \ldots + k_s}$, hence $m \in \{0,1,\ldots,k_{max}\}$. Knowing the number of $i$'s of order $a^m$ in $\Igroup{F}$
we know the number of rows of the corresponding type in $F$, and each such row contains $N/(a^m)$\ \  $1$'s, as well as the same number
of identical roots of unity of order $a^m$ of any other type. These data are sufficient to count the total number of $1$'s in $F$ and get 
the undephased defect of $1/\sqrt{N} \cdot F$.

We start from a lemma which will allow us to provide alternative  formulas for the defect $\undephasedDEFECT(1/\sqrt{N} \cdot F)$.


\begin{lemma}
	\label{lem_sum_of_min_of_kx_and_m}
	Let $(k_1,\ k_2,\ \ldots,\ k_s)$ be a nonempty ($s \geq 1$) sequence of nonnegative integers $k_x$,
	such as the one used to define $F$ in (\ref{eq_primary_Fourier_Kronecker_product}),
	and let $m \geq 1$ be an integer. 
		Then, for $k_{max}$ defined in  (\ref{eq_k_max}) and $\kxBELOW{l}$ defined in (\ref{eq_number_of_kx_below_value}),
	there holds:
	\begin{description}
		\item[a)]
			\begin{equation}
				\label{eq_sum_of_min_of_kx_and_m}
				\sum_{x=1}^s  \min\left(k_x, m \right)
						\ \ =\ \ 
				s \cdot m\ -\ \sum_{l=1}^m  \kxBELOW{l}           \ \ ,
			\end{equation}

		\item[b)]
			\begin{equation}
				\label{eq_sum_of_all_kx}
				\sum_{x=1}^s  k_x
						\ \ =\ \ 
				s \cdot k_{max} \ -\ \sum_{l=1}^{k_{max}}  \kxBELOW{l}         \ \ ,
			\end{equation}
	\end{description}
	where we do not asssume that $k_{max} > 0$.
\end{lemma}

\PROOFstart  
\begin{description}
\item[a)]
	The calculation below is based on the assumption that all the sums exist.
We leave  it to the comment that follows what to do if at some step it is not the case. 
Conditions $k_x < (m-l)$, where $l=0,\ldots,m-1$, should be read as $\{x:\ k_x < (m-l)\}$.
\begin{eqnarray}
	\label{eq_calculating_sum_of_min_of_kx_and_m}
	\lefteqn{  \sum_{x=1}^s  m  \ -\  \sum_{x=1}^s  \min\left(k_x, m\right)
						\ \ =\ \ 
				\sum_{x=1}^s  \left( m    \ -\     \min\left(k_x, m\right) \right)     }
	&   &   \\
	\nonumber
	&  =  &
	\sum_{k_x < m} \left( m          \ -\            k_x \right)
			\ \ =\ \ 
	\sum_{k_x < m}  1    \ +\     \sum_{k_x < m} \left( (m-1)  \ -\  k_x  \right)     \\
	\nonumber
	&  =  &
	\sum_{k_x < m}  \!\!\! 1    \ +\     \sum_{k_x < (m-1)} \!\!\!\!\!\!\! \left( (m-1)  \ -\  k_x  \right)
			\ \ =\ \ 
	\sum_{k_x < m}  \!\!\! 1    \ +\     \sum_{k_x < (m-1)} \!\!\!\!\!\!\! 1     \ +\  \sum_{k_x < (m-1)} \!\!\!\!\!\!\! \left( (m-2)  \ -\  k_x  \right)       \\
	\nonumber
	&   =   &
	\sum_{k_x < m}  \!\!\!  1    \ + \!\!\!\!     \sum_{k_x < (m-1)}  \!\!\!\!\!\!\! 1     \ +\!\!    \sum_{k_x < (m-2)}  \!\!\!\!\!\!\! \left( (m-2)   -   k_x  \right)
			\ =\  
	\sum_{k_x < m}  \!\!\!  1    \ +\!\!\!\!      \sum_{k_x < (m-1)} \!\!\!\!\!\!\! 1     \ +\!\!     \sum_{k_x < (m-2)}\!\!\!\!\!\!\!  1     \ +\!\!    
	 \sum_{k_x < (m-2)} \!\!\!\!\!\!\!  \left( (m-3)    -    k_x  \right)     \\
	\nonumber
	&    =    &    \ldots     \\
	\nonumber
	&   =   &
	\sum_{k_x < m}  1    \ +\  \ldots   \ +\     \sum_{k_x < 2}  1     \ +\       \sum_{k_x < 1}  \left( 1  \ -\  k_x  \right)
			\ \ =\ \ 
	\sum_{k_x < m}  1    \ +\  \ldots   \ +\     \sum_{k_x < 2}  1     \ +\       \sum_{k_x < 1}  1    \\
	\nonumber
	&    =    &
	\kxBELOW{m}  \ +\   \ldots    \ +\   \kxBELOW{2}   \ +\  \kxBELOW{1}                    \ \ .
\end{eqnarray}
This proves (\ref{eq_sum_of_min_of_kx_and_m}) provided there are no problems with the existence of the sums in the above expressions.
If at some moment  the extracted sum ($l \geq 0$):
\begin{equation}
	\label{eq_extracted_sum}
	\sum_{k_x < (m-l)}   \!\!\!\!\!\!\!  \left( (m-l)   -   k_x  \right)
\end{equation}
doesn't  exist because $0\ =\ \kxBELOW{m-l} \geq \ldots \geq \kxBELOW{1} \geq 0$   then we can immediately replace it by 
$\kxBELOW{m-l}\ +\ \kxBELOW{m-l-1}\ +\ \ldots\ +\ \kxBELOW{1}$.

\item[b)]
	This we obtain by replacing $m$ with $k_{max}$ in (\ref{eq_sum_of_min_of_kx_and_m}).

\end{description}

\PROOFend 

Next, as promised, we count the number of elements of a given order in the indexing group $\Igroup{F}$.
We did it also in Lemma 2.19 in \cite{PermEqClasses} in terms of the 'introduction indices' $\tilde{n}_l \ =\ \kxBELOW{l} + 1$.
But we performed these calculations in a highly inefficient way, using Lemma 2.15 in \cite{PermEqClasses} on how
the quantities $\ROWStype{a}{m}$  (introduced in Lemma \ref{lem_number_of_rows_of_type} below) change when a new factor
is added to an existing primary Kronecker product of Fourier matrices.


\begin{lemma}
	\label{lem_number_of_rows_of_type}
	The number of elements of order $a^m$ in the indexing group 
		$\Igroup{F}\ =\ \ZN{a^{k_1}} \times  \ldots \times \ZN{a^{k_s}}$, 
	i.e. the number $\ROWStype{a}{m}$ of rows in\ \ 
		$F\  =\ \FOURIER{a^{k_1}}     \otimes    \ldots    \otimes   \FOURIER{a^{k_s}}$\ \
	composed of the $a^m$th  roots of unity $\PHASE{\frac{2\pi}{a^m} n}$, reads:
	\begin{itemize} 
		\item 
			If $m=0$\ :
			\begin{equation}
				\label{eq_number_of_type_1_rows}
				\ROWStype{a}{0}
						\ \ \ =\ \ \
				a^{\sum_{x=1}^{s} \min\left( k_x, 0 \right)}
						\ \ \ =\ \ \
				a^{s \cdot 0\ -\ \sum_{l=1}^{0} \kxBELOW{l}}
						\ \ \ =\ \ \
				1   \ \ ,
			\end{equation}
		\item
			if $m \in \{ 1, \ldots, k_{max} \}$\ :
			\begin{eqnarray}
				\label{eq_number_of_type_a_to_m_rows}
				\ROWStype{a}{m}
						&    =    &
				a^{\sum_{x=1}^{s} \min\left( k_x, m \right)}
					\ \ -\ \ 
				a^{\sum_{x=1}^{s} \min\left( k_x, m-1 \right)}                     \\
				\nonumber
						&     =     &
				a^{s \cdot m\ -\ \sum_{l=1}^{m} \kxBELOW{l}}
						\ \ -\ \ 
				a^{s \cdot (m-1)\ -\ \sum_{l=1}^{(m-1)} \kxBELOW{l}}       \ \ ,
			\end{eqnarray}
	\end{itemize}
	where we assume that $k_{max} > 0$
	(if $k_{max}=0$ then there is only $\ROWStype{a}{0} = 1$ element in $\Igroup{F} = \ZN{1} \times \ldots \times \ZN{1}$
		(row in $\FOURIER{1} \otimes \ldots \otimes \FOURIER{1}\ =\ [1]$), of order $a^0 = 1$).

	In other words:
	\begin{eqnarray}
		\label{number_of_type_certain_rows}
		\ROWStype{a}{0}    &    =    &     1            \\
		\nonumber   &   &    \\
		\nonumber
		\ROWStype{a}{1}    &     =    &     a^{\sum_{x=1}^{s} \min\left( k_x, 1 \right)}   \ \ -\ \  1           \\
		\nonumber
								&    =    &     a^{s \cdot 1\ -\  \kxBELOW{1}}   \ \ -\ \ 1               \\
		\nonumber   &   &   \\
		\nonumber
		\ROWStype{a}{2}    &     =    &      a^{\sum_{x=1}^{s} \min\left( k_x, 2 \right)}	      \ \ -\ \     a^{\sum_{x=1}^{s} \min\left( k_x, 1 \right)}    \\
		\nonumber
								&     =    &      a^{s \cdot 2\ -\ \left(\kxBELOW{1} + \kxBELOW{2}\right)}    \ \ -\ \     a^{s \cdot 1\ -\ \kxBELOW{1}}       \\
		\nonumber
		&      \vdots      &       \\
		\nonumber
		\ROWStype{a}{k_{max}}   &    =    &      a^{\sum_{x=1}^{s} \min\left( k_x, k_{max} \right)}     
															\ \ -\ \     
														a^{\sum_{x=1}^{s} \min\left( k_x, k_{max}-1 \right)}         \\
		\nonumber
		&    =    &
														a^{s \cdot k_{max}\ -\ \sum_{l=1}^{k_{max}} \kxBELOW{l}}
															\ \ -\ \ 
														a^{s \cdot (k_{max}-1)\ -\ \sum_{l=1}^{(k_{max}-1)} \kxBELOW{l}}
	\end{eqnarray}	

	Thus all the possible orders have been exhausted and all the elements in $\Igroup{F}$ (rows in $F$) counted.
\end{lemma}

\PROOFstart 
	Let $\ORDERinGROUP{\GROUP{G}}{g}$ denote the order of $g$ in group $\GROUP{G}$ with the neutral element $\GROUPzero{\GROUP{G}}$. 
Recall that  $\ORDERinGROUP{\GROUP{G}}{g} \ \ =\ \ \min_{g^k=\GROUPzero{\GROUP{G}},\ k\in\{1,2,\ldots\}} k$, 
though in the considered case of the abelian group $\Igroup{F}$  the additive notation is used: $kg$ instead of $g^k$.

Let us count the number of elements in $\Igroup{F}$ with order not exceeding $a^m$, that is belonging to 
$\{1,a,a^2,\ldots,a^m\}$, where $m \leq k_{max}$, because\\ 
	$\ORDERinGROUP{\Igroup{F}}{\vectorINDEX{i_1,\ldots,i_s}}
		\  \ =\ \ 
	\lcm\left(  \ORDERinGROUP{\ZN{a^{k_1}}}{i_1},\ \ldots,\  \ORDERinGROUP{\ZN{a^{k_s}}}{i_s}   \right)$ 
\ \ and\ \ 
	$\ORDERinGROUP{\ZN{a^{k_x}}}{i_x}  \in \left\{ 1,a,\ldots,a^{k_x} \right\}$.

If $m=0$ this number is equal to $1$ and it corresponds to the neutral element 
$\GROUPzero{\Igroup{F}}\ =\ \vectorINDEX{\GROUPzero{\ZN{a^{k_1}}},\ldots,\GROUPzero{\ZN{a^{k_s}}}}$.

If $m>0$  we count the number of those $\vectorINDEX{i_1,\ldots,i_s}\in \Igroup{F}$ which satisfy
$a^m \vectorINDEX{i_1,\ldots,i_s}\ =\ \GROUPzero{\Igroup{F}}$, 
equivalently $\forall x\ \ a^m   i_x \  =\ \GROUPzero{\ZN{a^{k_x}}}$. To construct an element of $\Igroup{F}$ with this property,
we choose $i_x \in \ZN{a^{k_x}}$\ :
\begin{itemize}
	\item
		If $m \geq k_x$ then for any $i_x \in \ZN{a^{k_x}}\ \ a^m i_x\ =\ a^{m-k_x} a^{k_x} i_x \ =\ \GROUPzero{\ZN{a^{k_x}}}$,
		so there are $a^{k_x}$ possible choices.

	\item
		If $m < k_x$ then $a^m i_x \ =\ \GROUPzero{\ZN{a^{k_x}}}$ is equivalent to \\
		$i_x \in \left\{ \ZNclass{0},\ \ZNclass{a^{k_x-m}},\ \ZNclass{2 \cdot a^{k_x-m}},\ \ldots,\ \ZNclass{\left(a^m-1\right) \cdot a^{k_x-m}} \right\}
		\subset \ZN{a^{k_x}}$,
		because a representative of $i_x$ must be divided by $a^{k_x - m}$,
		so there are $a^m$ possible choices
		(here $\ZNclass{b}$ denotes $b + a^{k_x} \INTEGER \ \ \in\ \ \ZN{a^{k_x}} = \INTEGER/a^{k_x} \INTEGER$).
\end{itemize}

In effect, we can choose\ \  $\vectorINDEX{i_1,\ldots,i_s}$\ \  in\ \  
	$a^{\min\left( k_1, m \right)} \cdot  a^{\min\left( k_2, m \right)}  \cdot \ldots \cdot a^{\min\left( k_s, m \right)}$\ \ 
ways so that its order does not exceed $a^m$. Note that the formula $a^{\sum_{x=1}^{s} \min\left( k_x,\ m \right)}$ works
also for $m=0$ and when some or all $k_x$'s are equal to zero.

Now the number of elements of $\Igroup{F}$ with a fixed order $a^m$ is equal to the number of elements with order $\leq a^m$
minus the number of element with order $\leq a^{m-1}$, which gives us the formulas in the lemma. Their alternative versions result from 
Lemma \ref{lem_sum_of_min_of_kx_and_m}.
\PROOFend 

Everything has been prepared to calculate the defect of a primary Kronecker product of Fourier matrices
$F$ in (\ref{eq_primary_Fourier_Kronecker_product}), $\undephasedDEFECT(F)$, 
which is defined analogously to $\undephasedDEFECT(1/\sqrt{N} \cdot F)$ for the unitary $1/\sqrt{N} \cdot F$, and both these values are equal
(see section \ref{sec_defect} starting from the paragraph containing (\ref{eq_feasible_space_for_rescaled_unitary})).
The effect of our calculation is in the theorem below.
		
We offer the reader a number of formulas so that he can choose the one that suits him best, either for further research or for
numerical implementation, but we have not analyzed their numerical efficiency. Preparing item {\bf d)} we cared about the 
possibility to implement the partial formulas in a single way for the two cases considered there.


\begin{theorem}
	\label{theor_defect_of_primary_Kron_prod_of_Fouriers}
	The undephased defect (defined in (\ref{eq_undephased_defect_for_rescaled_unitary})) $\undephasedDEFECT(F)$ of a primary Kronecker product
	of Fourier matrices \ \ 
		$F\ \ =\ \ \FOURIER{a^{k_1}}     \otimes    \ldots    \otimes   \FOURIER{a^{k_s}}$\ \ 
	(where $a$ is prime, $s \geq 1$, $k_x \geq 0$ integer) of size $N\ =\ a^{k_1} \cdot \ldots \cdot a^{k_s}$ can be expressed as:
	\begin{description}
		\item[a)]  
			\begin{equation}
				\label{eq_F_primary_defect_formula_A_basic}
				\frac{N}{1} \cdot \ROWStype{a}{0} 
					\ +\ 
				\frac{N}{a} \cdot \ROWStype{a}{1}
					\ +\ 
				\frac{N}{a^2} \cdot \ROWStype{a}{2}
					\ +\ 
				\ldots
					\ +\ 
				\frac{N}{a^{k_{max}}} \cdot \ROWStype{a}{k_{max}}                \ \ ,
			\end{equation}
			if $k_{max} > 0$, otherwise it is $1 = \frac{N}{1} \ROWStype{a}{0}$,
			where $\ROWStype{a}{m}$ is the number of rows in $F$ containing the $a^m$th roots of unity, introduced in
			Lemma \ref{lem_number_of_rows_of_type}, and $k_{max}$ is defined in (\ref{eq_k_max}).

		\item[b)]  
			\begin{eqnarray}
				\label{eq_F_primary_defect_formula_B_factored_with_min_kx_m}
				\lefteqn{
					\frac{N}{a}
					\left(    \rule{0cm}{0.8cm}
						(a-1)
						\left(   \rule{0cm}{0.5cm}
							1
								\ +\ 
							a^{-1\ +\  \sum_{x=1}^{s} \min\left( k_x,1 \right)}
								\ +\ 
							a^{-2\  +\ \sum_{x=1}^{s} \min\left( k_x,2 \right)}
							\ +\ 
							\ldots       
						\right.
					\right.
				}
				&   &   \\
				\nonumber
				&   &   
				\left.   \rule{0cm}{0.8cm}
					\left.    \rule{0cm}{0.5cm}
							\ +\ 
						a^{-k_{max}\ +\ \sum_{x=1}^{s} \min\left( k_x,k_{max} \right)}
					\right)
						\ \ +\ \ 
					a^{-k_{max}\ +\ \sum_{x=1}^{s} \min\left( k_x,k_{max} \right)}
				\right)        
			\end{eqnarray}
			if $k_{max} > 0$, otherwise it is\ \ $1\ = \ \frac{N}{a}\left( (a-1)(1)\ +\ 1 \right)$. 
			Note that\\
				$-k_{max}\ +\ \sum_{x=1}^{s} \min\left( k_x,k_{max} \right)
					\ \ =\ \ 
				-k_{max}\  +\ \sum_{x=1}^{s} k_x$.

		\item[c)] 
			\begin{eqnarray}
				\label{eq_F_primary_defect_formula_C_factored_with_sk_minus_sum_kxbelows}
				\lefteqn{
					\frac{N}{a}
					\left(    \rule{0cm}{0.8cm}
						(a-1)
						\left(   \rule{0cm}{0.5cm}
							1
								\ +\ 
							a^{s \cdot 1\ -\  \left( \kxBELOW{1} + 1 \right)}
								\ +\ 
							a^{s \cdot 2\  -\  \left( \kxBELOW{1} + 1 \right)  \ -\  \left( \kxBELOW{2} + 1 \right)}
							\ +\ 
							\ldots       
						\right.
					\right.
				}
				&   &   \\
				\nonumber
				&   &   
				\left.   \rule{0cm}{0.8cm}
					\left.    \rule{0cm}{0.5cm}
							\ +\ 
						a^{s \cdot k_{max}\ -\ \sum_{l=1}^{k_{max}} \left( \kxBELOW{l} + 1 \right)}
					\right)
						\ \ +\ \ 
					a^{s \cdot k_{max}\ -\ \sum_{l=1}^{k_{max}} \left( \kxBELOW{l} + 1 \right)}
				\right)        
			\end{eqnarray}
			if $k_{max} > 0$, otherwise it is\ \ $1\ = \ \frac{N}{a}\left( (a-1)(1)\ +\ 1 \right)$,
			and where $\kxBELOW{l}$ defined in (\ref{eq_number_of_kx_below_value}).
			Also\ \ 
				$s \cdot k_{max}\ -\ \sum_{l=1}^{k_{max}} \left( \kxBELOW{l} + 1 \right)
					\ \ =\ \ 
				-k_{max}\  +\ \sum_{x=1}^{s} k_x$.

		\item[d)]   
			If $k_{max} > 0$, let $b_1, b_2, \ldots, b_r \in \{1,\ldots, s\}$  be distinct numbers,
			designating $0 < k_{b_1} < \ldots < k_{b_r}\ =\ k_{max}$, 
			such that\ \  $\bigcup_{x \in \{1,\ldots,s\},\ k_x >0} \left\{ k_x \right\}  \  =\  \left\{ k_{b_1}, k_{b_2}, \ldots, k_{b_r} \right\}$,\ \ 
			that is each nonzero $k_x$ is represented in sequence $\left( k_{b_1}, \ldots, k_{b_r} \right)$ only once.
			$\undephasedDEFECT(F)$ then reads:
			\begin{equation}
				\label{eq_F_primary_defect_formula_D_factored_with_poly_geometric_series}
				\frac{N}{a}
				\left(  \rule{0cm}{0.8cm}
					(a-1)
					\left(  \rule{0cm}{0.5cm}
						1     \ +\ 
						\defectPOLY{0}{k_{b_1}}   \ +\   \defectPOLY{k_{b_1}}{k_{b_2}} \ +\ \ldots\  \ +\  
						\defectPOLY{k_{b_{r-1}}}{k_{b_r}}
					\right)
						\ \ +\ \ 
					a^{-k_{max}\  +\ \sum_{x=1}^{s} k_x}
				\right)
			\end{equation}
			where the internal sum is equal to\ \ $1\ +\ \defectPOLY{0}{k_{b_1}}$\ \  if $\left(k_1,k_2,\ldots,k_s\right)$ contains only $0$'s and $k_{b_1}$'s,
			and where the component polynomials $\defectPOLY{k_{\alpha}}{k_{\beta}}$  are defined below
			(where $\left( k_{\alpha}, k_{\beta} \right) \ = \ \left( k_{b_y} , k_{b_{y+1}} \right)$ or 
			           $\left( k_{\alpha}, k_{\beta} \right) \ = \ \left( 0 , k_{b_{1}} \right)$):
			\begin{itemize}
				\item
					In the situations when: 
					\begin{enumerate}
						\item  $\left( k_{\alpha}, k_{\beta} \right)\ =\ \left( k_{b_{r-1}} , k_{b_r} \right)$ when $r \geq 2$ ($s \geq 2$ then) and
								$k_{b_r} = k_{max}$ occurs only once in $\left(k_1,\ldots,k_s\right)$, or

						\item  $\left( k_{\alpha}, k_{\beta} \right)\ =\ \left( 0 , k_{b_r} \right)$ when $r=1$ and
								$k_{b_1} = k_{max}$ occurs only once in $\left(k_1,\ldots,k_s\right)$,
					\end{enumerate}
					-- then $\kxBELOW{k_{\beta}}\ =\ s-1$  (see (\ref{eq_number_of_kx_below_value})) -- 
					the internal sum of (\ref{eq_F_primary_defect_formula_D_factored_with_poly_geometric_series})
					has the last component equal to:
					\begin{eqnarray}
						\label{eq_Poly_ka_kb_when_kmax_only_once_among_kx}
						\lefteqn{\defectPOLY{k_{\alpha}}{k_{\beta}}    \ \ =}
						&       &                            \\
						\nonumber
						&   &
						a^{\sum_{k_x \leq k_{\alpha}} k_x}      \cdot      \left( k_{\beta} \ -\ k_{\alpha} \right)
								\ \ =                   \\
						\nonumber
						&    &
						a^{\sum_{k_x \leq k_{\alpha}} k_x  \ +\   ( s - \kxBELOW{k_{\beta}} - 1 )    \left( k_{\alpha} + 1 \right)}      
							\cdot      \left( k_{\beta} \ -\ k_{\alpha} \right)     
								\ \ =\ \      \\
						\nonumber
						&     &
						a^{s(k_{\alpha} + 1)\ -\ \sum_{l=1}^{k_{\alpha} + 1} \left( \kxBELOW{l} + 1 \right)}      \cdot      \left( k_{\beta} \ -\ k_{\alpha} \right)
								\ \ =            \\
						\nonumber
						&    &
						a^{s k_{\alpha}  \ -\ \sum_{l=1}^{k_{\alpha}} \left( \kxBELOW{l} + 1 \right)}     \cdot      \left( k_{\beta} \ -\ k_{\alpha} \right)  
								\ \ =              \\
						\nonumber
						&   &
						a^{s k_{\alpha}  \ -\ \sum_{l=1}^{k_{\alpha}} \left( \kxBELOW{l} + 1 \right)     \ +\   (s - \kxBELOW{k_{\beta}} - 1 )}     
						\cdot      \left( k_{\beta} \ -\ k_{\alpha} \right)      \ \ .
					\end{eqnarray}
					In the special case (item 2. above) when $\left( k_1, k_2, \ldots, k_s \right)$ contains only one nonzero entry $k_{max}\ =\ k_{\beta}$
					(then $\kxBELOW{k_{\beta}} = \kxBELOW{1} = s-1$)
					-- also when $s=1$  --
					the internal sum 
					of (\ref{eq_F_primary_defect_formula_D_factored_with_poly_geometric_series})
					contains only
					\begin{equation}
						\label{eq_Poly_ka_kb_when_kmax_only_once_among_kx_which_are_all_zeros_except_kmax}
						1\ +\ \defectPOLY{0}{k_{\beta}}  \ \ =\ \ 1\ +\ k_{\beta}      \ \ ,
					\end{equation}
					where $\defectPOLY{0}{k_{\beta}}$ can be derived from any of the formulas in (\ref{eq_Poly_ka_kb_when_kmax_only_once_among_kx})
					by replacing $k_{\alpha}$ with $0$
					(not existing sums are treated as zeros).

				\item
					In other situations (requiring that $s \geq 2$): 
					\begin{enumerate}
						\item  $\left( k_{\alpha}, k_{\beta} \right)\ =\ \left(k_{b_{r-y-1}}, k_{b_{r-y}}\right)$ for $1 \leq y \leq r-2$ when
								$r \geq 3$ ($s \geq 3$ then), or

						\item $\left( k_{\alpha}, k_{\beta} \right)\ =\ \left(0, k_{b_1}\right)$ when $r \geq 2$, then $k_{b_1} < k_{b_r}=k_{max}$, or

						\item $\left( k_{\alpha}, k_{\beta} \right)\ =\ \left(k_{b_{r-1}}, k_{b_r}\right)$ and $k_{b_r} = k_{max}$ occurs in 
								$\left(k_1,\ldots,k_s\right)$ more than once, or

						\item $\left( k_{\alpha}, k_{\beta} \right)\ =\ \left(0, k_{b_r}\right)$ when $r=1$ and $k_{b_1} = k_{max}$ occurs in
								$\left(k_1,\ldots,k_s\right)$ more than once,
					\end{enumerate}
					for which $\kxBELOW{k_{\beta}} < s-1$:
					\begin{eqnarray}
						\label{eq_kb_is_not_the_only_one_the_greatest_among_kx}
						\lefteqn{\defectPOLY{k_{\alpha}}{k_{\beta}}    \ \ =}
						&       &                            \\
						\nonumber
						&     &
						a^{\sum_{k_x \leq k_{\alpha}} k_x  \ +\   ( s - \kxBELOW{k_{\beta}} - 1 )    \left( k_{\alpha} + 1 \right)}      
							\cdot 
						\frac{1   \ -\   \left( a^{ s - \kxBELOW{k_{\beta}} - 1} \right)^{k_{\beta} - k_{\alpha}}}
							  {1   \ -\   \left( a^{ s - \kxBELOW{k_{\beta}} - 1} \right)}
								\ \ =
						\\
						\nonumber
						&   &  
						a^{s \left( k_{\alpha} + 1 \right)  \ -\ \sum_{l=1}^{k_{\alpha} + 1} \left( \kxBELOW{l} + 1 \right)}
							\cdot 
						\frac{1   \ -\   \left( a^{ s - \kxBELOW{k_{\beta}} - 1} \right)^{k_{\beta} - k_{\alpha}}}
							  {1   \ -\   \left( a^{ s - \kxBELOW{k_{\beta}} - 1} \right)}
								\ \ =
						\\
						\nonumber
						&    &
						a^{s k_{\alpha}     \ +\   ( s - \kxBELOW{k_{\beta}} - 1 )    \ -\   \sum_{l=1}^{k_{\alpha}} \left( \kxBELOW{l} + 1 \right)}
							\cdot 
						\frac{1   \ -\   \left( a^{ s - \kxBELOW{k_{\beta}} - 1} \right)^{k_{\beta} - k_{\alpha}}}
							  {1   \ -\   \left( a^{ s - \kxBELOW{k_{\beta}} - 1} \right)}
					\end{eqnarray}
					In the special case (items 2. and 4. above) when $k_{\alpha}\ =\ 0$ and $k_{\beta}\ =\ k_{b_1}\ =\ \min_{k_x > 0} k_x$  
					the internal sum of (\ref{eq_F_primary_defect_formula_D_factored_with_poly_geometric_series})
					begins with
					\begin{equation}
						\label{eq_kb_is_not_the_only_one_the_greatest_among_kx_and_it_is_the_smallest_nonzero_kx}
						1 \ +\ \defectPOLY{0}{k_{\beta}} 
								\ \ \ =\ \ \ 
						1 \ \ +\ \ 
						a^{ s - \kxBELOW{k_{\beta}} - 1 }  
							\cdot
						\frac{1   \ -\   \left( a^{ s - \kxBELOW{k_{\beta}} - 1} \right)^{k_{\beta}}}
							  {1   \ -\   \left( a^{ s - \kxBELOW{k_{\beta}} - 1} \right)}                        \ \ ,
					\end{equation}
					and there are no other polynomials in the internal sum if the only nonzero entries (more than one of these) 
					in $\left( k_1, k_2, \ldots, k_s \right)$ are those equal to $k_{b_1} = k_{max}$ (item 4. above). 
					$\defectPOLY{0}{k_{\beta}}$ of (\ref{eq_kb_is_not_the_only_one_the_greatest_among_kx_and_it_is_the_smallest_nonzero_kx})
					can  be obtained from any of the formulas in (\ref{eq_kb_is_not_the_only_one_the_greatest_among_kx})
					by replacing $k_{\alpha}$ with $0$.
			\end{itemize}

			The orders of the consecutive polynomials $1$, $\defectPOLY{0}{k_{b_1}}$, $\defectPOLY{k_{b_1}}{k_{b_2}}$, ..., $\defectPOLY{k_{b_{r-1}}}{k_{b_r}}$
			($1$, $\defectPOLY{0}{k_{max}}$ if $r=1$) form an increasing sequence except that the last two entries in this sequence of orders are equal if and only if
			$k_{max}$ occurs once in $\left(k_1,\ldots,k_s\right)$.
		
		\item[e)] 
			\begin{eqnarray}
				\label{eq_F_primary_defect_formula_D_sum_product_version_of_C}
				\lefteqn{
					\frac{N}{a}
					\left(    \rule{0cm}{0.8cm}
						(a-1)
						\left(   \rule{0cm}{0.5cm}
							1
								\ +\ 
							a^{s \ -\  \left( \kxBELOW{1} + 1 \right)}
							\left(   \rule{0cm}{0.5cm}
								1
									\ +\ 
								a^{s \ -\  \left( \kxBELOW{2} + 1 \right)}
								\left(   \rule{0cm}{0.5cm}
									1
										\ +\ 
									\ldots
								\right.
							\right.
						\right.
					\right.
				}
				&   &   \\
				\nonumber
				&  &
				\left.    \rule{0cm}{0.8cm}
					\left.    \rule{0cm}{0.5cm}
						\left.     \rule{0cm}{0.5cm}
							\left.     \rule{0cm}{0.5cm}
								\ldots
									\ +\ 
								a^{s \ -\  \left( \kxBELOW{(k_{max}-1)} + 1 \right)}
								\left(  \rule{0cm}{0.5cm}
									1
										\ +\ 
									a^{s \ -\  \left( \kxBELOW{k_{max}} + 1 \right)}
								\right)
								\ldots
							\right)
						\right)
					\right)
						\ \ +\ \ 
					a^{s \cdot k_{max}\ -\ \sum_{l=1}^{k_{max}} \left( \kxBELOW{l} + 1 \right)}   
				\right)  
			\end{eqnarray}
			if $k_{max} > 0$.
			Using the quantities $k_{b_1}$, ..., $k_{b_r}$ introduced at item {\bf d)}:\ \ $s - \left(\kxBELOW{m} + 1\right)$ is 
			nonnegative, and it is equal to zero if and only if
			\begin{itemize}
				\item
					$k_{b_{r-1}} + 1\ \leq m\ \leq\ k_{b_r}$\ \ for\ \ $r \geq 2$,

				\item
					$1\ \leq\ m\ \leq\ k_{b_1}=k_{max}$, i.e. all $m$, if\ \ $r=1$,
			\end{itemize}
			and $k_{max}$ occurs in $\left(k_1,\ldots,k_s\right)$ only once.
				
	\end{description}
\end{theorem}

\PROOFstart 
\begin{description}
	\item[a)]  
		$\undephasedDEFECT(F)\ =\ \undephasedDEFECT(1/\sqrt{N} \cdot F)$ is the number of $1$'s in matrix $F$
(see the paragraph containing (\ref{eq_primary_Fourier_Kronecker_product})). 
A row with $a^m$th roots of unity , of length $N$, contains 
$N/a^m$ $1$'s as well as $N/a^m$\ \ $\PHASE{\frac{2\pi}{a^m}n}$'s for any $n=1,\ldots,(a^m-1)$. 
The number of such rows (which are of order $a^m$ in the group of rows of $F$, isomorphic with $\Igroup{F}$), 
	$\ROWStype{a}{m}$,
is given in Lemma \ref{lem_number_of_rows_of_type}, where $m \in \left\{0,1,\ldots,k_{max}\right\}$ by this lemma. 
Thus all the rows of $F$ are exhausted.

	\item[b)]  
		This formula is obtained from (\ref{eq_F_primary_defect_formula_A_basic}) at {\bf a)} by substituting 
$a^{\sum_{x=1}^{s} \min\left( k_x, m \right)}
					\ \ -\ \ 
a^{\sum_{x=1}^{s} \min\left( k_x, m-1 \right)}$
for $\ROWStype{a}{m}$ in accordance with Lemma \ref{lem_number_of_rows_of_type}. 
	Going in the opposite direction, in expression (\ref{eq_F_primary_defect_formula_B_factored_with_min_kx_m}), of the form
$(N/a)((a-1)(1\ +\ \ldots)\ +\ \ldots)$, product $(N/a) a \cdot 1$ produces $(N/1) \ROWStype{a}{0}\ =\ N$, the number of
$1$'s in the $\GROUPzero{\Igroup{F}}$th row of $F$. Further
\begin{eqnarray}
	\label{eq_A_into_B_transformation_leading_to_summand_for_exponent_1}
	\lefteqn{ \frac{N}{a}\left(a \cdot a^{-1\ +\  \sum_{x=1}^{s} \min\left( k_x,1 \right)}    \ -\ 1 \cdot 1 \right)   \ \ =}
	&  &  \\
	\nonumber
	&   &
	\frac{N}{a}\left(   a^{\sum_{x=1}^{s} \min\left( k_x,1 \right)}  \ -\ 1 \right)  
			\ \ =\ \ 
	\frac{N}{a} \cdot \ROWStype{a}{1}    \ \ , 
\end{eqnarray}

\begin{eqnarray}
	\label{eq_A_into_B_transformation_leading_to_summand_for_exponent_2}
	\lefteqn{ \frac{N}{a}\left(a \cdot a^{-2\ +\  \sum_{x=1}^{s} \min\left( k_x,2 \right)}    
					\ -\ 
				1 \cdot a^{-1\ +\  \sum_{x=1}^{s} \min\left( k_x,1 \right)} \right)   \ \ =}
	&  &  \\
	\nonumber
	&  &
	\frac{N}{a^2}\left(   a^{\sum_{x=1}^{s} \min\left( k_x,2 \right)}    \ -\     a^{\sum_{x=1}^{s} \min\left( k_x,1 \right)}   \right)
			\ \ =\ \ 
	\frac{N}{a^2} \cdot \ROWStype{a}{2}         \ \ ,
\end{eqnarray}
and so on until
\begin{eqnarray}
	\label{eq_A_into_B_transformation_leading_to_summand_for_exponent_kmax}
	\lefteqn{ \frac{N}{a}\left(a \cdot a^{-k_{max}\ +\  \sum_{x=1}^{s} \min\left( k_x, k_{max} \right)}    
					\ -\ 
				1 \cdot a^{-\left( k_{max} - 1 \right) \ +\  \sum_{x=1}^{s} \min\left( k_x, k_{max}-1 \right)} \right)   \ \ =}
	&  &  \\
	\nonumber
	&  &
	\frac{N}{a^{k_{max}}}\left(   a^{\sum_{x=1}^{s} \min\left( k_x, k_{max} \right)}    \ -\     a^{\sum_{x=1}^{s} \min\left( k_x, k_{max}-1 \right)}   \right)
			\ \ =\ \ 
	\frac{N}{a^{k_{max}}} \cdot \ROWStype{a}{k_{max}}        \ \ .
\end{eqnarray}
The remaining product, $-1 \cdot  a^{-k_{max}\ +\  \sum_{x=1}^{s} \min\left( k_x, k_{max} \right)}$, is compensated for by the 
the rightmost power of $a$ in the outer bracket in  (\ref{eq_F_primary_defect_formula_B_factored_with_min_kx_m}).
We easily check that this formula works  when $k_{max}=1$:
\begin{eqnarray}
	\label{eq_B_formula_for_kmax_equal_1}
	\undephasedDEFECT(F)   
			&   =   &
	\frac{N}{1} \ROWStype{a}{0}  \ +\ \frac{N}{a}  \ROWStype{a}{1}
			\ \ =\ \
	\frac{N}{1} \cdot 1  \ +\     \frac{N}{a} \cdot \left(   a^{\sum_{x=1}^{s} \min\left( k_x,1 \right)}  \ -\ 1 \right)    \\
	\nonumber
	&    =    &
	\frac{N}{a}
	\left( 
		(a-1) 
		\left( 1 \ +\   a^{-1  \ +\ \sum_{x=1}^{s} \min\left( k_x,1 \right)}  \right)
			\ +\ 
		a^{-1  \ +\ \sum_{x=1}^{s} \min\left( k_x,1 \right)}
	\right)
\end{eqnarray}

	\item[c)]  
		This is the same formula as the one at {\bf b)}, the only difference is that every\ \ $\sum_{x=1}^{s} \min\left( k_x,m \right)$\ \ 
has been replaced by\ \  $s \cdot m\ -\ \sum_{l=1}^m  \kxBELOW{l}$,\ \  which is allowed thanks to Lemma \ref{lem_sum_of_min_of_kx_and_m}.

	\item[d)]  
		Here we further develop formula (\ref{eq_F_primary_defect_formula_B_factored_with_min_kx_m}) at {\bf b)} using the fact that
there can be groups of equal numbers among the exponents $k_1$, ..., $k_s$. Let $k_{\alpha}\ <\ k_{\beta}$ represent two such neighbouring groups
so that there is no $k_x$ between them. We also have to consider a special case when $k_{\alpha} = 0$ and $k_{\beta} = \min_{k_x > 0} k_x$,
though we do not assume that there are zeros among $k_1$, ..., $k_s$.  For each such pair $k_{\alpha},\ k_{\beta}$ we add up the 
neighbouring components in the internal sum of (\ref{eq_F_primary_defect_formula_B_factored_with_min_kx_m}), related to this pair, to form
polynomial $\defectPOLY{k_{\alpha}}{k_{\beta}}$. All those polynomials in a sum with $1$ will form the internal sum of 
(\ref{eq_F_primary_defect_formula_B_factored_with_min_kx_m}). $\defectPOLY{k_{\alpha}}{k_{\beta}}$ thus reads (where the sum can be a one 
summand sum if $k_{\alpha} + 1 \ =\ k_{\beta}$):
\begin{eqnarray}
	\label{eq_ka_kb_sum_of_neighbouring_monomials}
	\lefteqn{ \defectPOLY{k_{\alpha}}{k_{\beta}} \ \ =}
	&    &     \\
	\nonumber
	&   &
	a^{-(k_{\alpha} + 1)  \ +\  \sum_{x=1}^{s} \min\left( k_x,  k_{\alpha} + 1 \right)}
		\ +\ 
	a^{-(k_{\alpha} + 2)  \ +\  \sum_{x=1}^{s} \min\left( k_x,  k_{\alpha} + 2 \right)}
		\ +
	\ldots
		+\ 
	a^{-k_{\beta} \ +\  \sum_{x=1}^{s} \min\left( k_x,  k_{\beta} \right)}  \ \ ,
\end{eqnarray}
which is equal to:
\begin{eqnarray}
	\label{eq_defect_poly_transformation}
	\sum_{d=1}^{k_{\beta} - k_{\alpha}}
		a^{-(k_{\alpha} + d)  \ +\  \sum_{x=1}^{s} \min\left( k_x,  k_{\alpha} + d \right)}
			&  =  &
	\sum_{d=1}^{k_{\beta} - k_{\alpha}}
		a^{-(k_{\alpha} + d) \ +\  \sum_{k_x \leq k_{\alpha}} k_x    \ +\  \sum_{k_x \geq k_{\beta}} (k_{\alpha} + d)}   \ \ =     \\
	\nonumber
	\sum_{d=1}^{k_{\beta} - k_{\alpha}}
		a^{(s - \kxBELOW{k_{\beta}} - 1)(k_{\alpha} + d)  \ +\ \sum_{k_x \leq k_{\alpha}} k_x}
			& = &
	a^{(s - \kxBELOW{k_{\beta}} - 1)(k_{\alpha} + 1)   \ +\    \sum_{k_x \leq k_{\alpha}} k_x}
		\cdot
	\sum_{d=0}^{k_{\beta} - k_{\alpha} - 1}
		\left( a^{(s - \kxBELOW{k_{\beta}} - 1)}  \right)^d       \ \ ,
\end{eqnarray}
where on the right we have a sum of a geometric series. Note that $\sum_{k_x \leq k_{\alpha}} k_x$ may not exist 
(when $k_{\alpha}=0$ is absent from $\left(k_1,\ldots,k_s\right)$), it is $0$ then, 
and that $\sum_{k_x \geq k_{\beta}} (k_{\alpha} + d)$ exists even if $k_{\beta}\ =\ k_{max}$.\\ 
$s - \kxBELOW{k_{\beta}}$\ \  is the number of $k_x$'s greater or equal to $k_{\beta}$, according to the definition of $\kxBELOW{k_{\beta}}$
in (\ref{eq_number_of_kx_below_value}).

If\ \ $s - \kxBELOW{k_{\beta}} - 1\ =\ 0$\ \ the above expression reduces to
\begin{equation}
	\label{eq_defect_poly_when_kb_maximal_occurs_once}
	a^{\sum_{k_x \leq k_{\alpha}} k_x}     \cdot     \left(   k_{\beta}  \ -\     k_{\alpha}  \right)
			\ \ =\ \ 
	a^{\sum_{k_x \leq k_{\alpha}} k_x    \ +\     (s - \kxBELOW{k_{\beta}} - 1)(k_{\alpha} + 1)}     \cdot     \left(   k_{\beta}  \ -\     k_{\alpha}  \right)
	\ \ .
\end{equation}
To an existing sum $\sum_{k_x \leq k_{\alpha}} k_x$  (when $k_{\alpha}$ occurs in $\left(k_1,\ldots,k_s\right)$) one can apply 
Lemma \ref{lem_sum_of_min_of_kx_and_m}  {\bf b)} with $\left(k_1,\ldots,k_s\right)$ replaced with its subsequence composed of all $k_x$'s not greater 
than $k_{\alpha}$, $s$ replaced with $\kxBELOW{k_{\beta}}$ (this many of them there are), $k_{max}$ replaced with $k_{\alpha}$, so\ \ 
	$\sum_{k_x \leq k_{\alpha}} k_x  \ \ =\ \ \kxBELOW{k_{\beta}} \cdot k_{\alpha} \ -\ \sum_{l=1}^{k_{\alpha}} \kxBELOW{l}$,\ \ 
where the $\kxBELOW{l}$'s used have the same value  wrt to both $\left(k_1,\ldots,k_s\right)$ and the subsequence.  This works also for $k_{\alpha}\ =\ 0$,
irrespective of whether the sum on the left exists, that is whether $\left(k_1,\ldots,k_s\right)$ contains $0$'s  (nonexistent sums are treated as $0$'s).
Thus the exponent in (\ref{eq_defect_poly_when_kb_maximal_occurs_once}) can be further written as 
(where in any case $\kxBELOW{k_{\alpha}+1} \ =\ \kxBELOW{k_{\beta}}$):
\begin{eqnarray}
	\label{eq_defect_poly_multiplier_exponent_in_any_case}
	\lefteqn{
		\kxBELOW{k_{\beta}} \cdot k_{\alpha} \ -\ \sum_{l=1}^{k_{\alpha}} \kxBELOW{l}
			\ +\ 
		\left( s - \kxBELOW{k_{\beta}} - 1 \right) \left( k_{\alpha} + 1 \right)     \ \ =
	}
	&     &     \\
	\nonumber
	&    &
	\kxBELOW{k_{\beta}} \left( k_{\alpha} + 1 \right)    \ -\  \sum_{l=1}^{k_{\alpha} + 1} \kxBELOW{l}
		\ +\ 
	\left( s - \kxBELOW{k_{\beta}} - 1 \right) \left( k_{\alpha} + 1 \right)       \ \ =                   \\
	\nonumber
	&     &
	(s-1) \left( k_{\alpha} + 1 \right)              \ -\         \sum_{l=1}^{k_{\alpha} + 1} \kxBELOW{l}
			\  \ =\ \ 
	s \left( k_{\alpha} + 1 \right)        \ -\       \sum_{l=1}^{k_{\alpha} + 1} \left( \kxBELOW{l}  +  1\right)      \ \ .
\end{eqnarray}
In the considered case when\ \  $s - \kxBELOW{k_{\beta}} - 1 \ =\ 0$,\ \  i.e.\ \   $\kxBELOW{k_{\beta}}\ =\ \kxBELOW{k_{\alpha}+1} \ =\ s-1$,\ \ 
this final expression is equal to $s k_{\alpha} \ -\ \sum_{l=1}^{k_{\alpha}} \left( \kxBELOW{l}  +  1\right)$, hence 
(\ref{eq_defect_poly_when_kb_maximal_occurs_once}) has the shape:
\begin{eqnarray}
	\label{eq_defect_poly_when_kb_maximal_occurs_once_remaining_formulas}
	a^{s ( k_{\alpha} + 1 )        \ -\       \sum_{l=1}^{k_{\alpha} + 1} ( \kxBELOW{l}  +  1 )}         \left(   k_{\beta}  \ -\     k_{\alpha}  \right)
			&    =    &
	a^{s k_{\alpha} \ -\ \sum_{l=1}^{k_{\alpha}} ( \kxBELOW{l}  +  1 )}       \left(   k_{\beta}  \ -\     k_{\alpha}  \right)                   \\
	\nonumber
	&    =     &
	a^{s k_{\alpha} \ -\ \sum_{l=1}^{k_{\alpha}} ( \kxBELOW{l}  +  1 )   \ +\   ( s - \kxBELOW{k_{\beta}} - 1 )}    
	(   k_{\beta}  \ -\     k_{\alpha}  )       \ \ .
\end{eqnarray}
In this way all the formulas in (\ref{eq_Poly_ka_kb_when_kmax_only_once_among_kx}) have been explained, 
as well as in (\ref{eq_Poly_ka_kb_when_kmax_only_once_among_kx_which_are_all_zeros_except_kmax})
because our considerations include the case when\ \  $k_{\alpha}\ =\ 0$\ \  and\ \  $k_{\beta}\ =\ \min_{k_x > 0} k_x$,\ \ that is when
(\ref{eq_ka_kb_sum_of_neighbouring_monomials}) is the leftmost part (neglecting $1$) of the internal sum of
(\ref{eq_F_primary_defect_formula_B_factored_with_min_kx_m}).
\bigskip

If\ \  $s - \kxBELOW{k_{\beta}} - 1 \ >\ 0$,\ \  in the last expression in (\ref{eq_defect_poly_transformation})  we have a sum of a geometric   
series, which gives us:
\begin{equation}
	\label{eq_defect_poly_weighted_geometric_sum}
	a^{(s - \kxBELOW{k_{\beta}} - 1)(k_{\alpha} + 1)   \ +\    \sum_{k_x \leq k_{\alpha}} k_x}      
		\cdot 
	\frac{1   \ -\   \left( a^{ s - \kxBELOW{k_{\beta}} - 1} \right)^{k_{\beta} - k_{\alpha}}}
	{1   \ -\   \left( a^{ s - \kxBELOW{k_{\beta}} - 1} \right)}      \ \ .
\end{equation}
The exponent on the left can be transformed as in (\ref{eq_defect_poly_multiplier_exponent_in_any_case}) 
(see the paragraph preceding (\ref{eq_defect_poly_multiplier_exponent_in_any_case})) and the final expression there is equal to
\ \ $s k_{\alpha} \ -\ \sum_{l=1}^{k_{\alpha}} \left( \kxBELOW{l}  +  1 \right)   \ +\   \left( s - \kxBELOW{k_{\beta}} - 1 \right)$
because (in any case) $\kxBELOW{k_{\beta}} \ =\ \kxBELOW{k_{\alpha} + 1}$.  Hence the above (\ref{eq_defect_poly_weighted_geometric_sum})
also reads:
\begin{eqnarray}
	\label{eq_defect_poly_weighted_geometric_sum_remaining_formulas}
	a^{  s ( k_{\alpha} + 1 )        \ -\       \sum_{l=1}^{k_{\alpha} + 1} ( \kxBELOW{l}  +  1 )  }
		\cdot 
	\frac{1   \ -\   \left( a^{ s - \kxBELOW{k_{\beta}} - 1} \right)^{k_{\beta} - k_{\alpha}}}
	{1   \ -\   \left( a^{ s - \kxBELOW{k_{\beta}} - 1} \right)}  \ \ =
	&   &     \\
	\nonumber
	a^{  s k_{\alpha} \ -\ \sum_{l=1}^{k_{\alpha}} ( \kxBELOW{l}  +  1  )   \ +\    ( s - \kxBELOW{k_{\beta}} - 1 )  }
		\cdot 
	\frac{1   \ -\   \left( a^{ s - \kxBELOW{k_{\beta}} - 1} \right)^{k_{\beta} - k_{\alpha}}}
	{1   \ -\   \left( a^{ s - \kxBELOW{k_{\beta}} - 1} \right)}    \ \ .
	&    &   
\end{eqnarray}
We have justified the formulas sitting in (\ref{eq_kb_is_not_the_only_one_the_greatest_among_kx}) and 
(\ref{eq_kb_is_not_the_only_one_the_greatest_among_kx_and_it_is_the_smallest_nonzero_kx}).
\bigskip

The polynomials $\defectPOLY{0}{k_{b_1}}$, $\defectPOLY{k_{b_1}}{k_{b_2}}$, ..., $\defectPOLY{k_{b_{r-1}}}{k_{b_r}}$
are of nondecreasing orders because they group neighbouring summands in the internal sum of 
	(\ref{eq_F_primary_defect_formula_B_factored_with_min_kx_m})
	(or  \ref{eq_F_primary_defect_formula_C_factored_with_sk_minus_sum_kxbelows}), 
and these summands have nondecreasing exponents not smaller than $0$\ :
\begin{itemize}
	\item   Looking at (\ref{eq_F_primary_defect_formula_C_factored_with_sk_minus_sum_kxbelows}),
			the first exponent $s - \kxBELOW{1} - 1  \geq 0$ because $s \geq 1$ and $\kxBELOW{1} \leq \kxBELOW{k_{max}} \leq (s-1)$,
			where $k_{max} > 0$ as assumed at item {\bf d)}.

	\item  Further, consider the difference between the neighbouring exponents  
			of the $m$th and the $(m-1)$th summand (where $1 \leq m \leq k_{max}$), 
			including the difference between the first exponent and the $0$th exponent (equal to $0$)
			associated with the leading $1$:
			\begin{equation}
				\label{eq_neighbouring_exponents_difference}
				\left( sm - \sum_{l=1}^{m} \left( \kxBELOW{l} + 1 \right) \right)
					\ -\ 
				\left( s(m-1) - \sum_{l=1}^{m-1} \left( \kxBELOW{l} + 1 \right) \right)
						\ \ =\ \ 
				s - \kxBELOW{m} - 1
						\ \ \geq\ \ 
				0\ \ ,
			\end{equation}
			because $\kxBELOW{m} \leq \kxBELOW{k_{max}} \leq (s-1)$, where $k_{max} > 0$.
\end{itemize}

The difference (\ref{eq_neighbouring_exponents_difference}) is zero only when:
\begin{enumerate}
	\item
		$m \in \left\{ k_{b_{r-1}} +1,\  k_{b_{r-1}} +2,\ \ldots,\ k_{b_r} = k_{max} \right\}$, 
		in the case when $r \geq 2$ (at least two different nonzeros among $k_x$'s),

	\item
		$m \in \left\{ 1,\ 2,\ \ldots,\ k_{b_1}=k_{max} \right\}$ when $r=1$,
\end{enumerate}
and, simultaneously, there is a single $k_{max}$ in $\left( k_1,\ldots,k_s \right)$.

Let us consider the polynomials $\defectPOLY{k_{\alpha}}{k_{\beta}}$ formed by adding up the internal summands of
	(\ref{eq_F_primary_defect_formula_C_factored_with_sk_minus_sum_kxbelows})
corresponding to these $m$ listed above in the case 1.:\ \ $\defectPOLY{k_{b_{r-1}}}{k_{b_r}}$\ \ 
											and the case 2.:\ \ $\defectPOLY{0}{k_{b_1}}$,\ \ 
and let us also assume that there is only one $k_{max}$ in $\left( k_1,\ldots,k_s \right)$.
The $r$-th polynomial $\defectPOLY{k_{b_{r-1}}}{k_{b_r}}$
	\ \ (for $r=1$ we consider $\defectPOLY{0}{k_{b_1}}$)\ \   
is built from the $\left( k_{b_{r-1}} + 1\right)$-th, ..., $\left(k_{b_r} \right)$-th summands of the internal sum
of (\ref{eq_F_primary_defect_formula_C_factored_with_sk_minus_sum_kxbelows})
	\ \ (summands $1$-th, ..., $\left( k_{b_1}\right)$-th when $r=1$), where the leading $1$ in this internal sum is treated as the
$0$-th summand.  
Their orders are all equal to the order of
the  $\left(k_{b_{r-1}}\right)$-th summand sitting in $\defectPOLY{k_{b_{r-2}}\mbox{\ or\ $0$}}{k_{b_{r-1}}}$
	\ \ (the leading $1$ if $r=1$)\ , because, as it has been said above, the differences between the orders of these particular 
summands (corresponding to the values of $m$ in case 1. and 2. for the value of $r$) are all zeros, including the difference
between the orders of the  $\left(k_{b_{r-1}} + 1\right)$-th and  $\left(k_{b_{r-1}}\right)$-th summands 
			\ \ (between the orders of the $1$-st summand and the leading $1$, equal to $s - \kxBELOW{1} - 1$, for $r=1$).
We now see that  the orders of $\defectPOLY{k_{b_{r-1}}}{k_{b_r}}$ and $\defectPOLY{k_{b_{r-2}}\mbox{\ or\ $0$}}{k_{b_{r-1}}}$
			\ \ ($\defectPOLY{0}{k_{b_1}}$ and the leading $1$, when $r=1$)
are equal.
However, if $r \geq 2$ the order of $\defectPOLY{k_{b_{p-1}}}{k_{b_p}}$ 
is greater than  the order of  $\defectPOLY{k_{b_{p-2}}}{k_{b_{p-1}}}$ for $1 \leq p < r$. (where 
for  $p=2$ we have to consider the pair $\defectPOLY{k_{b_1}}{k_{b_2}}$, $\defectPOLY{0}{k_{b_1}}$,
for $p=1$ we take the pair $\defectPOLY{0}{k_{b_1}}$, the leading $1$).
  It is because the orders of $a^{sm - \sum_{l=1}^{m} \left( \kxBELOW{l} + 1 \right)}$ for $m=1,2,\ldots,k_{b_{r-1}}$ form
an increasing sequence, greater than $0$ as $s-\kxBELOW{1}-1 > 0$  for $r \geq 2$. 
Simply, the differences (\ref{eq_neighbouring_exponents_difference}) are greater than $0$ for $m \leq k_{b_{r-1}}$
since $\kxBELOW{m} \leq s-2$ then as there are at least two $k_x$'s not smaller than $m$: $k_{b_{r-1}}$ and $k_{b_r}$.

If we assume that there is more than one $k_{max}$ in $\left( k_1,\ldots,k_s \right)$ then again the differences 
(\ref{eq_neighbouring_exponents_difference}) are all greater than $0$, this time for all $m$ in $\left\{1,\ldots,k_{max}\right\}$,
since $\kxBELOW{m} \leq s-2$ as there are at least two $k_x$'s not smaller than $m$: these very $k_{max}$'s. 
That is why the orders of summands in the internal sum of 
	(\ref{eq_F_primary_defect_formula_C_factored_with_sk_minus_sum_kxbelows})
form an increasing sequence greater than $0$ (ommiting the leading $1$), and so do the orders of the polynomials
$\defectPOLY{0}{k_{b_1}}$, ..., $\defectPOLY{k_{b_{r-1}}}{k_{b_r}}$.

	\item[e)]  
This is another form of the formula 
	(\ref{eq_F_primary_defect_formula_C_factored_with_sk_minus_sum_kxbelows})
at item {\bf c)}. The statement under it is obvious and we have met it in the the proof of item {\bf d)}
-- see the inequality (\ref{eq_neighbouring_exponents_difference}) and the paragraph after it.
\end{description}

\PROOFend 

Analyzing the formulas in Theorem \ref{theor_defect_of_primary_Kron_prod_of_Fouriers}, we can also say something about the
dephased defect (definition (\ref{eq_dephased_defect_for_rescaled_unitary}))
$\DEFECT(F)\ =\ \undephasedDEFECT(F) - (2N-1)$ of an $N \times N$ primary Kronecker product $F$  !


\begin{corollary}
	\label{cor_D_and_d_factors_for_primary_kron_prod_F}
	For a primary Kronecker product of Fourier matrices $F\ \ =\ \ \FOURIER{a^{k_1}}     \otimes    \ldots    \otimes   \FOURIER{a^{k_s}}$,
	where $k_{max} > 0$\ :
	\begin{description}
		\item[a)]
			$\undephasedDEFECT(F)$ is divided by 
			$N/a
					\ \ =\ \ 
			a^{-1\ +\ \sum_{x=1}^s k_x} 
					\ \ =\ \  
			a^{-1\ +\ s \cdot k_{max} \ -\ \sum_{l=1}^{k_{max}}  \kxBELOW{l}}$.

		\item[b)]
			$\DEFECT(F)$ is divided by $(a-1)^2$.
	\end{description}
\end{corollary}

\PROOFstart  
\begin{description}

\item[a)]
	$\undephasedDEFECT(F)$ has factor $N/a$ in formula 
		(\ref{eq_F_primary_defect_formula_C_factored_with_sk_minus_sum_kxbelows}),
while the expression in the outer bracket in this formula is also an integer number. It is because every exponent in the inner bracket
is nonnegative, as $s - \kxBELOW{l} - 1\ \geq\ 0$ for $l \in \left\{ 1,\ldots,k_{max} \right\}$, where $k_{max} > 0$.

\item[b)]
	Let us recall: let $p(a)\ =\ (a-1)w(a)$ be a polynomial divided by $(a-1)$. Now if $p'(1)\ =\ w(1)\ \ \stackrel{assume!}{=} 0$,
for the derivative of the corresponding polynomial function, then $w$ is also divided by $(a-1)$ and we have that
$p(a)\ =\ (a-1)^2 v(a)$.

	Let us return to formula (\ref{eq_F_primary_defect_formula_B_factored_with_min_kx_m}) at 
		Theorem \ref{theor_defect_of_primary_Kron_prod_of_Fouriers} {\bf b)} 
for $\undephasedDEFECT(F)$ and write shortly the dephased defect $\DEFECT(F)$ as:
\begin{eqnarray}
	\nonumber	
	\DEFECT(F)    &    =    & 
	\undephasedDEFECT(F)\ -\ (2N - 1) 
			\ \ =\ \ 
	\frac{N}{a}\left(\ \   (a-1) \left(\ 1\ +\ \ldots\ \right)  \ \ +\ \ a^{\ldots}   \ \ \right)   \ \ \ -\ \  \ (2N-1)        \\
	\label{eq_dephased_defect_of_primary_Fourier_kron_prod}
	&   =   &
	\frac{1}{a}\left(\ \ \ N \left(\ \ (a-1)\left(\ 1\ +\ \ldots\ \right)   \ \ +\ \ a^{\ldots}  \ \ -\ \ 2a \ \ \right)   \ \ \ +\ \ \ a \ \ \ \right)
	\\
	\nonumber
	&   =    &       \frac{1}{a}  p(a) \ \ .
\end{eqnarray}
Our $p(a)$ is the polynomial in the outer bracket in (\ref{eq_dephased_defect_of_primary_Fourier_kron_prod}) 
(where $N = a^{\sum_{x=1}^s k_x}$) and it is a polynomial because every exponent there is nonnegative,
as explained in the proof of item {\bf a)}. 
We check that $p(1)\ =\ 0$ and next, using the fact that $\VALUEat{\left(a^k\right)'}{1}$, the value of the
derivative of  $a^k$ at $1$,  equals $k$ for $k = 0,1,2,\ldots$, we calculate:
\begin{eqnarray}
	\label{eq_p_poly_derivative_at_1}
	\lefteqn{ p'(1) \ \ = }    &     &     \\
	\nonumber
	&   =   &
	\left(\sum_{x=1}^s k_x\right)  \left(\ \ 0\left(\ 1\ +\ \ldots\ \right)  \ \ +\ \ 1\ \ -\ \ 2\ \ \right)  \\
	\nonumber
	&    &
	+\ \ \
	1\left(\ \ 1\left(\ 1\ +\ k_{max}\ \right)  \ \ +\ \ 0\VALUEat{\left(\ 1\ +\ \ldots\ \right)'}{1} \ \ + \ \ \left(-k_{max} + \sum_{x=1}^s k_x\right)   \ \ -\ \ 2\ \ \right)
		\ \ \ +\ \ \
	1                            \\
	\nonumber
	&    =   &
	\left(\sum_{x=1}^s k_x \right)\left(\ \ -1\ \ \right) \ \ \ +\ \ \ \left(\ \ \left(\sum_{x=1}^s k_x\right)\ \ -\ \ 1\ \ \right)   \ \ \ +\ \ \ 1
		\ \ \ =\ \ \ 
	0                         \ \ .
\end{eqnarray}
We thus obtain that $\DEFECT(F)\ \ =\ \ (1/a) p(a)\ \ =\ \ (1/a) (a-1)^2 v(a)$, 
where polynomial $v(a)$ has integer coefficients because $p(a)$ has integer coefficients. 
Let $a$ be the only prime divisor of the size of $F$. 
Because $a \cdot \DEFECT(F)\ \ =\ \ (a-1)^2 v(a)$, where $\DEFECT(F)$ is an integer number, 
$a$ must divide the integer number $v(a)$ as $a$ doesn't divide $(a-1)^2$.
That is $v(a)/a$ is an integer number resulting from dividing  $\DEFECT(F)$ by $(a-1)^2$.
\end{description}

\PROOFend  

In a table that follows we have collected the values of $\undephasedDEFECT(F)$ deprived of factor $N/a$ in the second column, and the values
of $\DEFECT(F)$ deprived of factor $(a-1)^2$ in the third column, for various primary Kronecker products
	$F\ \ =\ \ \FOURIER{a^{k_1}}     \otimes    \ldots    \otimes   \FOURIER{a^{k_s}}$.
$F$ in such a case is represented by partition $\left[ k_1,\ k_2,\ \ldots,\ k_s \right]$, a non-increasing  sequence of nonzero integers adding up to
$\log_a N$. Every primary product can be reshuffled to have non-increasing exponents running from left to right, and  a change in the order of Kronecker 
factors of $F$ amounts to permuting $F$ which changes neither $\undephasedDEFECT(F)$ nor $\DEFECT(F)$  
(Corollary \ref{cor_DU_equal_DV_for_U_V_equivalent} and the fact that  $\undephasedDEFECT(F)\ =\ \undephasedDEFECT(1/\sqrt{N} \cdot F)$,
see (\ref{eq_undephased_defect_for_rescaled_unitary})). 
Thus, for example, the row of the below table containing partition $[2\ 2\ 1]$ features
$(a/N)\undephasedDEFECT(F)$ and $\DEFECT(F)/\left((a-1)^2\right)$ for $F\ =\ \FOURIER{a^2} \otimes \FOURIER{a^2} \otimes \FOURIER{a}$, but also
for $F\ =\ \FOURIER{a^2} \otimes \FOURIER{a} \otimes \FOURIER{a^2}$ and for $F\ =\ \FOURIER{a} \otimes \FOURIER{a^2} \otimes \FOURIER{a^2}$. 
The partitions associated with a given value of $\log_a N$ occur in neighbouring  rows of the table and they are ordered lexicographically, while $\log_a N$
increases down the table. The polynomials in the second and third columns are factored over the ring of integers, so in some cases we can see additional
integer factors of  $\undephasedDEFECT(F)$ or $\DEFECT(F)$.

{\center Table 1. Factored undephased ($\undephasedDEFECT$) and dephased ($\DEFECT$) defects of inequivalent primary Kronecker products of Fourier matrices.}
\begin{longtable}{|l|l|l}    
\ensuremath{\mbox{partition}}   &    \ensuremath{\frac{a}{N} \undephasedDEFECT(F)}       &     \ensuremath{\frac{1}{(a-1)^2} \DEFECT(F)} \\
------------------  & --------------------------- & ------------------------------------------------------------
\endfirsthead
\ensuremath{\mbox{partition}}   &    \ensuremath{\frac{a}{N} \undephasedDEFECT(F)}       &     \ensuremath{\frac{1}{(a-1)^2} \DEFECT(F)}  \\
------------------  & --------------------------- & ------------------------------------------------------------
\endhead
		\ensuremath{[1]}     &     \ensuremath{2a-1}      &      \ensuremath{0}                                                \\
		\ensuremath{[2]}     &     \ensuremath{3a-2}      &       \ensuremath{1}                                               \\
		\ensuremath{[1\ \ 1]}    &   \ensuremath{a^2 + a - 1}        &          \ensuremath{a+1}                        \\
		\ensuremath{[3]}      &     \ensuremath{4a - 3}    &       \ensuremath{2a + 1}                                       \\
		\ensuremath{[2\ \ 1]}   &   \ensuremath{2a^2 -  1}         &      \ensuremath{2a^2  +   2a   +   1}            \\
		\ensuremath{[1\ \ 1\ \ 1]}
		 &   \ensuremath{a^3  + a  - 1}       &      \ensuremath{(a+1) \left( a^2 + a + 1 \right)}       \\
		\ensuremath{[4]}   &     \ensuremath{5a - 4}       &     \ensuremath{3 a^2  +  2 a   +   1}                 \\
		\ensuremath{[3\ \ 1]}     &    \ensuremath{3 a^2   - a  - 1}      &      \ensuremath{ 3 a^3  + 3 a^2  + 2 a + 1}                 \\
		\ensuremath{[2\ \ 2]}     &    \ensuremath{a^3 + a^2  - 1}     &     \ensuremath{(a+1)  \left( a^3  +  2 a^2  + a + 1 \right)}             \\
		\ensuremath{[2\ \ 1\ \ 1]}     &    \ensuremath{2 a^3 - a^2 + a - 1}       &     \ensuremath{2 a^4   +  3 a^3  +  3 a^2  + 2 a  + 1}     \\
		\ensuremath{[1\ \ 1\ \ 1\ \ 1]}
		&      \ensuremath{a^4 +  a  -  1}     
		&    \ensuremath{(a+1)\left(a^2 + 1\right)\left(a^2 + a + 1 \right)}             \\
		\ensuremath{[5]}   &     \ensuremath{6 a - 5}      &      \ensuremath{4 a^3  +  3 a^2  + 2 a  + 1}                      \\
		\ensuremath{[4\ \ 1]}     &    \ensuremath{4 a^2  -  2 a  -  1}    &     \ensuremath{4 a^4  +  4 a^3  +  3 a^2  +  2 a  +  1}        \\
		\ensuremath{[3\ \ 2]}      &    \ensuremath{2 a^3  - 1}     &      \ensuremath{(a+1)  \left(  2 a^4  +  2 a^3  +  2 a^2  +  a  +  1  \right)}       \\
		\ensuremath{[3\ \ 1\ \ 1]}    &     \ensuremath{3 a^3  -  2 a^2  +  a  -  1}       &     \ensuremath{3 a^5  +  4 a^4  + 4 a^3  +  3 a^2  + 2 a  + 1}       \\
		\ensuremath{[2\ \ 2\ \ 1]}
		&      \ensuremath{a^4  +  a^3  -  a^2  +  a  -  1}  
		&  \ensuremath{(a+1)\left(a^2 + a + 1\right)\left(a^3  +  a^2  +  1\right)}    \\
		 \ensuremath{[2\  \ 1\ \ 1\ \ 1]}     &     \ensuremath{2 a^4  -  a^3  +  a  -  1}       &        \ensuremath{2a^6 + 3a^5 + 4a^4 + 4a^3 + 3a^2 + 2a + 1}   \\
		\ensuremath{[1\ \ 1\ \ 1\ \ 1\ \ 1]}
		&    \ensuremath{\left(a^2 \!-\! a\! +\! 1\right) \left(a^3\! +\! a^2\! -\! 1\right)}\!\!\!     
		&  \ensuremath{(a+1)\left(a^2 + 1\right)\left(a^4 + a^3 + a^2 + a + 1\right)}                                                   \\
		\ensuremath{[6]}   &    \ensuremath{7a - 6}     &     \ensuremath{5a^4 + 4a^3 + 3a^2 + 2a + 1}                             \\
		\ensuremath{[5\ \ 1]}    &    \ensuremath{5a^2 - 3a - 1}    &     \ensuremath{5a^5 + 5a^4 + 4a^3 + 3a^2 + 2a + 1}    \\
		\ensuremath{[4\ \ 2]}   &   \ensuremath{3a^3  -  a^2  - 1}      &     \ensuremath{3a^6 +  5a^5 + 5a^4 + 4a^3 + 3a^2 + 2a  + 1}      \\
		\ensuremath{[4\ \ 1\ \ 1]}   &   \ensuremath{4a^3 - 3a^2 + a - 1}    &    \ensuremath{4a^6 + 5a^5 + 5a^4 + 4a^3 + 3a^2 + 2a + 1}    \\
		\ensuremath{[3\ \ 3]}
		&    \ensuremath{a^4 + a^3 - 1}       &    \ensuremath{(a+1)\left(a^2+a+1\right)\left(a^4+a^3+a^2+1\right)}   \\
		\ensuremath{[3\ \ 2\ \ 1]}   &   \ensuremath{2a^4 - a^2 + a - 1}    &    \ensuremath{(a+1)\left(2a^6 + 2a^5 + 3a^4  + 2a^3 + 2a^2 + a + 1\right)}   \\
		\ensuremath{[3\ \ 1\ \ 1\ \ 1]}   &   \ensuremath{3a^4  -  2a^3  +  a  -  1}    &   \ensuremath{3a^7 + 4a^6 + 5a^5 + 5a^4 + 4a^3 + 3a^2 + 2a + 1} \\
		\ensuremath{[2\ \ 2\ \ 2]}
		&   \ensuremath{a^5 + a^3 - a^2 + a - 1}    
		&  \ensuremath{\left(a^2+a+1\right)\left(a^6+a^5+2a^4+2a^3+a^2+a+1\right)} \\
		\ensuremath{[2\ \ 2\ \ 1\ \ 1]}  &  \ensuremath{a^5 + a^4 - a^3 + a - 1}  
		&  \ensuremath{(a+1)\left(a^2 + 1\right)\left(a^5 + 2a^4 + a^3 + a^2 + a + 1\right)}  \\
		\ensuremath{[2\ \ 1\ \ 1\ \ 1\ \ 1]}  &  \ensuremath{2a^5 - a^4 + a - 1}  &  \ensuremath{2a^8 + 3a^7 + 4a^6 + 5a^5 + 5a^4 + 4a^3 + 3a^2 + 2a + 1} \\
		\ensuremath{[1\ \ 1\ \ 1\ \ 1\ \ 1\ \ 1]}
		&   \ensuremath{a^6 + a - 1}  
		&  \ensuremath{(a+1)\left(a^2+a+1\right)\left(a^2-a+1\right)\left(a^4+a^3+a^2+a+1\right)}    \\

\end{longtable}

Having analyzed  primary Kronecker products of Fourier matrices we move on  
to the most general case, the Fourier matrices being arbitrary products  $\FOURIER{N_1} \otimes \ldots \otimes \FOURIER{N_r}$ of
smaller Fourier matrices associated with the cyclic groups $\ZN{N_1}$, ..., $\ZN{N_r}$. Recall that the indexing group $\Igroup{F}$ of $F$ is
$\ZN{N_1} \times \ldots \times \ZN{N_r}$ then (see the beginning of section \ref{sec_Fourier_matrices}). $F$ and $G$ are (permutation) equivalent
if and only if their indexing groups are isomorphic 
	(see Corollary \ref{cor_Fouriers_perm_equivalent_iff_their_indexing_groups_isomorphic} and
	 Lemma \ref {lem_F_G_equivalent_iff_F_G_permutation_equivalent} {\bf (a)}), 
so in such a case not only\ \  $\undephasedDEFECT(F)\ =\ \undephasedDEFECT(G)$\ \  and\ \ $\DEFECT(F)\ =\ \DEFECT(G)$
	(by Corollary \ref{cor_DU_equal_DV_for_U_V_equivalent} taken together with (\ref{eq_undephased_defect_for_rescaled_unitary})), 
but also the spectra of the accompanying Berezin transforms $\Iu{1/\sqrt{N} \cdot F}$ and $\Iu{1/\sqrt{N} \cdot  G}$ are identical as composed of the entries of
$F$ and $G$, respectively (Theorem \ref{theor_eigenbasis_of_I_F}).

In the paragraph following the proof of Corollary \ref{cor_Fouriers_perm_equivalent_iff_their_indexing_groups_isomorphic}
we outlined the procedure of transforming an indexing group $\Igroup{F}$ into an isomorphic counterpart which, as a direct product of groups,
has factors of the orders being powers of primes only. Analogous algorithm converts the underlying $F$ into its permutation equivalent counterpart
	$G\ =\ G_1 \otimes \ldots \otimes G_p$,
where every $G_i$ is a primary Kronecker product (i.e. of the type discussed earlier in this section) of size $a_i^{m_i}$. Assuming that 
$G_i$'s are of maximal length, that is that $a_i$'s are distinct primes, there holds 
	$\undephasedDEFECT(F)\ =\ \undephasedDEFECT(G)\ =\ \undephasedDEFECT\left(G_1\right) \cdot \ldots \cdot \undephasedDEFECT\left(G_p\right)$,
which reflects the conditional multiplicativity of $\undephasedDEFECT$ with respect to Kronecker subproducts (factors) holding for Fourier matrices:


\begin{lemma}
	\label{lem_multiplicativity_of_D_of_F}
	Let $F$ of size $N$ and $G$ of size $M$ be Fourier matrices (like in (\ref{eq_Fourier_Kron_prod})).
	\begin{description}
		\item[a)]
			If $\gcd(M,N)\ =\ 1$ then the number of $1$'s in $F \otimes G$ (in the spectrum of $\Iu{\frac{1}{\sqrt{MN}} F \otimes G}$) is equal to the 
		product of the number of $1$'s in $F$ (in the spectrum of $\Iu{\frac{1}{\sqrt{N}} F}$) and the number of $1$'s in $G$ 
		(in the spectrum of $\Iu{\frac{1}{\sqrt{M}} G}$).
		That is to say,\ \ \ \ $\undephasedDEFECT(F \otimes G) \ \ =\ \ \undephasedDEFECT(F) \cdot \undephasedDEFECT(G)$.

		\item[b)]
			If $\gcd(M,N)\ >\ 1$ then\ \ \ \ $\undephasedDEFECT(F \otimes G) \ \ >\ \ \undephasedDEFECT(F) \cdot \undephasedDEFECT(G)$.
	\end{description}
\end{lemma}

\PROOFstart 
\begin{description}
	\item[a)]
	The sizes $N_i$ of the Kronecker factors  $F_{N_i}$ of $F$ are divisors of $N$, the size of $F$, so the entries of $F$ belong to the group
of the $N$th roots of unity. Similarly, the entries of $G$ belong to the group of the $M$-th roots of unity. Let us take an entry from $F$, 
$\PHASE{2\pi \frac{k}{N}}$, and an entry from $G$, $\PHASE{2\pi \frac{l}{M}}$. If the corresponding (their product) entry of $F \otimes G$
is equal to $1$, $k/N\ +\ l/M\ \ =\ \ (kM + lN)/(MN)$ is an integer number, so $N$ divides $k$ and $M$ divides $l$ because $M,N$ are relatively prime.
Thus the considered entries of $F$ and $G$ are $1$'s, so all $1$'s in $F \otimes G$ appear as products of $1$'s from $F$ and $G$. 
The bracket version and the last statement are consequences of   Theorem \ref{theor_eigenbasis_of_I_F}, equality 
(\ref{eq_undephased_defect_and_V1_dimension_equality}) in Theorem \ref{theor_unitarity_and_eigenspaces_of_Iu} {\bf e)}
and finally (\ref{eq_undephased_defect_for_rescaled_unitary}), because for example:
\begin{equation}
	\label{eq_series_of_equal_defects}
	\undephasedDEFECT(F)
			\ \ \stackrel{(\ref{eq_undephased_defect_for_rescaled_unitary})}{=}\ \ 
	\undephasedDEFECT\left( \frac{1}{\sqrt{N}} F \right)
			\ \ \stackrel{(\ref{eq_undephased_defect_and_V1_dimension_equality})}{=}\ \ 
	\DIMC\left( \Vu{ \frac{1}{\sqrt{N}} F }{1} \right)
			 \ \ \stackrel{Th. \ref{theor_eigenbasis_of_I_F}}{=}\ \ 
	\mbox{the number of $1$'s in $F$}     \ \ .
\end{equation}

	\item[b)]
	We already know from Corollary \ref{cor_V1_spaces_Kron_multiplied_are_contained_in_V1_space} {\bf c)} that 
(using also (\ref{eq_undephased_defect_for_rescaled_unitary}))\ \ 
		$\undephasedDEFECT(F \otimes G) \ \geq\ \undephasedDEFECT(F) \cdot \undephasedDEFECT(G)$.

A row of $F$, as a function of the  group column  index, is a homomorphism on the indexing group $\Igroup{F}$  
(from (\ref{eq_multiplying_rows_of_F}) and $F\ =\ F^T$). Therefore its image is a group, the  set of all the roots of $1$ of the order
associated with this row. Each of the roots occurs the same number of times, for column indices belonging to an appropriate coset
of the kernel of the considered homomorphism.

Let 
	$F\ \ =\ \ \FOURIER{N_1}     \otimes    \ldots    \otimes   \FOURIER{N_r}$, it is indexed by $\ZN{N_1} \times \ldots \times \ZN{N_r}$.
The group of rows of $F$, under entrywise product $\HADprod$, is isomorphic to $\Igroup{F}$
	(explained in the paragraph containing (\ref{eq_multiplying_rows_of_F})).
Let us take $\left(i_1,i_2,\ldots,i_r\right) \in \Igroup{F}$ such that each $i_k$ is a generator in $\ZN{N_k}$, so its order is equal to $N_k$.
Then the order of $\left(i_1,i_2,\ldots,i_r\right)$ in $\Igroup{F}$ is equal to $x\ =\ \lcm\left(N_1, N_2, \ldots, N_r\right)$. 
Thus also $F$ contains a row of this order, i.e. a row containing all the $x$-th roots of unity, and each root occurs the same number of times as it has been said above.
Analogously,
	$G\ \ =\ \ \FOURIER{M_1}     \otimes    \ldots    \otimes   \FOURIER{M_s}$
has a row with all the $\lcm\left(M_1, M_2, \ldots, M_s\right)$th roots of $1$.

Now, if $\gcd\left( N\ =\ N_1 \cdot \ldots \cdot N_r,\ M\ =\ M_1 \cdot \ldots \cdot M_s\right)\ >\ 1$ then there exists a prime number $a>1$ 
that divides both $\lcm\left(N_1, N_2, \ldots, N_r\right)$ and $\lcm\left(M_1, M_2, \ldots, M_s\right)$, so the chosen above rows of $F$ and $G$
contain (in equal quantities), among other values, all the $a$th roots of unity (as subgroups of the groups 'contained' in these rows). Thus, in a Kronecker product of
these rows, being a row in $F \otimes G$, we always get occasions to produce $1$'s out of numbers not equal to $1$. 
That is why the number of $1$'s in $F \otimes G$ exceeds $\undephasedDEFECT(F) \cdot \undephasedDEFECT(G)$.
\end{description}
\PROOFend    

The example in the paragraph following the proof of Corollary \ref{cor_Fouriers_perm_equivalent_iff_their_indexing_groups_isomorphic}
features the (isomorphic) indexing groups of (permutation equivalent) Fourier matrices 
	$F_{12} \otimes F_6$,  $F_4 \otimes F_3 \otimes F_6$,  $F_4 \otimes F_3 \otimes F_3 \otimes F_2$, $F_2 \otimes F_4 \otimes F_3 \otimes F_3$.
All these matrices have equal defects (by Corollary \ref{cor_DU_equal_DV_for_U_V_equivalent} and (\ref{eq_undephased_defect_for_rescaled_unitary})) 
and according to Lemma \ref{lem_multiplicativity_of_D_of_F} their defect is equal to
	$\undephasedDEFECT\left( F_2 \otimes F_4 \otimes F_3 \otimes F_3\right) \ \ =\ \ 
		\undephasedDEFECT\left( F_2 \otimes F_4\right) \cdot \undephasedDEFECT\left(  F_3 \otimes F_3\right)$, that is it can be calculated as a product
of the defects of the maximal primary Kronecker subproducts extracted in a  matrix transformation process if needed.

Lemma \ref{lem_multiplicativity_of_D_of_F} {\bf a)}, together with Theorem \ref{theor_defect_of_primary_Kron_prod_of_Fouriers},
 allows to explain in an easy way the formula from Theorem 5.3 {\bf (b)} in \cite{Defect},
expressing the dephased defect of $\FOURIER{N}$, which was proved there with the use of number theoretic tools (see also Corollary 5.7 in \cite{Banica_def_gen_four_matr}).
The corresponding undephased defect of $\FOURIER{N}$, where $N\ =\ a_1^{n_1} \cdot \ldots \cdot a_r^{n_r}$ for distinct primes $a_1,\ldots,a_r$ and where
$\FOURIER{N}$ is then (permutation) equivalent to $\FOURIER{a_1^{n_1}} \otimes \ldots \otimes \FOURIER{a_r^{n_r}}$, reads:
\begin{eqnarray}
	\label{eq_undephased_defect_for_F_N}
	\undephasedDEFECT\left(   \FOURIER{a_1^{n_1}} \otimes \ldots \otimes \FOURIER{a_r^{n_r}}  \right)
	& \stackrel{Lem.\ref{lem_multiplicativity_of_D_of_F} \mathbf{a)}}{=}   &
	\undephasedDEFECT\left( \FOURIER{a_1^{n_1}} \right) \cdot \ldots \cdot \undephasedDEFECT\left( \FOURIER{a_r^{n_r}} \right)     \\
	\nonumber
	&		\stackrel{Th.\ref{theor_defect_of_primary_Kron_prod_of_Fouriers}}{=}      &
	\prod_{k=1}^{r} \left( \left(n_k + 1\right) a_k^{n_k}   \ -\       n_k   a_k^{n_k -1} \right)        \\
	\nonumber
	&    =      &
	\prod_{k=1}^{r}   \left(  a_k^{n_k} \left(  n_k  +  1  -  \frac{n_k}{a_k} \right) \right)      
			\ \ =\ \ 
	N \prod_{k=1}^{r}    \left(  n_k  +  1  -  \frac{n_k}{a_k} \right)\ \ .     
\end{eqnarray}

Another proof of the above conditional multiplicativity of $\undephasedDEFECT(F)$ can be found in Banica's work \cite{Banica_def_gen_four_matr} where
the author derives, in the proof of Theorem 4.2, the formula ($N$ is the size of $F$):
\begin{equation}
	\label{eq_undephased_defect_indexing_group_formula}
	\undephasedDEFECT(F)\ \ =\ \  \sum_{i \in \Igroup{F}} \frac{N}{\ORDERinGROUP{\Igroup{F}}{i}}\ \ ,\ \ \ \ 
	\mbox{where}\ \ \ORDERinGROUP{\Igroup{F}}{i}\ \ \mbox{is the order of $i$ in $\Igroup{F}$}\ \ ,
\end{equation}
generalizing, to arbitrary Fourier matrices, the concept of the $M \times M$ PCM matrices from \cite{Defect} (Definition 5.1) parametrizing the feasible space
for  $F_M$:\ \ \\
$\FEASIBLEspace{\FOURIER{M}}\ =\ \left\{  \Ii R \HADprod \FOURIER{M}\ :\ \ R = P \FOURIER{M}\ \ \mbox{and}\ \ P\ \mbox{is PCM} \right\}$
	(this is (\ref{eq_feasible_space_for_rescaled_unitary}) for $V = \FOURIER{M}$ and it has been justified in the proof of Theorem 5.2 in \cite{Defect},
	with the use of an explicit form of the system defining $\FEASIBLEspace{V}$ like that in (\ref{eq_iRU_Uhansposed_antihermitian_system_in_entrywise_form}),
	or, in a more general setting of arbitrary Fourier matrices, in the proof of Proposition 4.1 in \cite{Banica_def_gen_four_matr}).
Similar results obtained independently by us were presented in \cite{Mattriad2011}. Namely, that for any Fourier matrix $F$  indexed by $\Igroup{F}$,
\begin{equation}
	\label{eq_parametrization_of_feasible_space_with_PCM_matrices}
	\FEASIBLEspace{F}
			\ \ =\ \ 
	\left\{
		\Ii R \HADprod F\ :\ \ 
		R\ =\ PF\ \ \ \ \mbox{and}\ \ \ \ \left(\forall_{i,j \in \Igroup{F}}\ \ P_{i+j,j}\ =\ P_{i,j}\ \ \mbox{and}\ \ P_{i,-j}\ =\ \CONJ{P}_{i,j}\right)
	\right\}\ \ ,
\end{equation}
which is a bit differently formulated in Proposition 4.1 in \cite{Banica_def_gen_four_matr}. The proof is easy and it is based on property 
(\ref{eq_multiplying_rows_of_F}). Calculating the number of real parameters in matrix $P$ leads to (\ref{eq_undephased_defect_indexing_group_formula}),
as in Theorem 4.2 in \cite{Banica_def_gen_four_matr}. Note that $R\ =\ PF$ is real as it is required by (\ref{eq_feasible_space_for_rescaled_unitary}) !
Also note that (\ref{eq_undephased_defect_indexing_group_formula}) is a generalization of (\ref{eq_F_primary_defect_formula_A_basic}).

\begin{example}
	\label{ex_PCM_feasible_space_parametrization_for_F2xF4}
	Consider:
	\begin{equation}
		\label{eq_F2_and_F4_matrices}
		\FOURIER{2}
				\ \ =\ \ 
		\left[ \begin{array}{cc}  1  &  1  \\   1  &  -1  \end{array} \right]
							\ \ \ \ \ \ \mbox{and}\ \ \ \ \ \ 
		\FOURIER{4}
				\ \ =\ \ 
		\left[
			\begin{array}{cccc}
				1 & 1 & 1 & 1 \\
				1 & \Ii & -1 & -\Ii \\
				1 & -1 & 1 & -1 \\
				1 & -\Ii & -1 & \Ii 
			\end{array}
		\right]                       \ \ .
	\end{equation}
	$\FOURIER{2} \otimes \FOURIER{4}$ is indexed by $\Igroup{\FOURIER{2} \otimes \FOURIER{4}}\ =\ \ZN{2} \times \ZN{4}$. Further,
	\begin{equation}
		\label{eq_PCM_feasible_space_parametrization_for_F2xF4}
		\FEASIBLEspace{\FOURIER{2} \otimes \FOURIER{4}}
				\ \ =\ \ 
		\left\{ 
			\Ii \left(P \left(\FOURIER{2} \otimes \FOURIER{4}\right)\right) \HADprod \left(\FOURIER{2} \otimes \FOURIER{4}\right)\ :\ \ 
			P \in \SPACE{P}_{\FOURIER{2} \otimes \FOURIER{4}}
		\right\}                        \ \ ,
	\end{equation}
	where
	\begin{equation}
		\label{eq_PCM_space_for_F2xF4}
		\SPACE{P}_{\FOURIER{2} \otimes \FOURIER{4}}
				\ \ =\ \ 
		\left\{\ \ \ 
			\begin{array}{c}
				\begin{array}{cccccccc}
					{\scriptscriptstyle 00}\      &
					{\scriptscriptstyle 01}\      &
					{\scriptscriptstyle 02} \      &
					{\scriptscriptstyle 03} \      &
					{\scriptscriptstyle 10} \     &
					{\scriptscriptstyle 11} \      &
					{\scriptscriptstyle 12}  \     &
					{\scriptscriptstyle 13}    
				\end{array}                                \\
				\left[
					\begin{array}{cccccccc}  
						\alpha_0  \!\!\!\!\!\!\!\!\!\!\!\!\!\!\!\!{\scriptscriptstyle 00}\ \ \ \ \ \ \ \
						u_0     &     \beta_0     &    \CONJ{u}_0    &    \gamma_0    &    w_0    &    \delta_0    &   \CONJ{w}_0       \\
						\alpha_1  \!\!\!\!\!\!\!\!\!\!\!\!\!\!\!\!{\scriptscriptstyle 01}\ \ \ \ \ \ \ \
						u_0     &     \beta_1     &    \CONJ{u}_0    &    \gamma_1    &   w_1    &    \delta_1    &    \CONJ{w}_1       \\
						\alpha_2  \!\!\!\!\!\!\!\!\!\!\!\!\!\!\!\!{\scriptscriptstyle 02}\ \ \ \ \ \ \ \
						u_0     &     \beta_0     &    \CONJ{u}_0    &    \gamma_2    &   w_0    &    \delta_2    &    \CONJ{w}_0       \\
						\alpha_3  \!\!\!\!\!\!\!\!\!\!\!\!\!\!\!\!{\scriptscriptstyle 03}\ \ \ \ \ \ \ \
						u_0     &     \beta_1     &    \CONJ{u}_0    &    \gamma_3    &   w_1    &    \delta_3    &    \CONJ{w}_1       \\
						\alpha_4  \!\!\!\!\!\!\!\!\!\!\!\!\!\!\!\!{\scriptscriptstyle 10}\ \ \ \ \ \ \ \
						u_1     &     \beta_2     &    \CONJ{u}_1    &    \gamma_0    &    w_1    &    \delta_2    &   \CONJ{w}_1       \\
						\alpha_5  \!\!\!\!\!\!\!\!\!\!\!\!\!\!\!\!{\scriptscriptstyle 11}\ \ \ \ \ \ \ \
						u_1     &     \beta_3     &    \CONJ{u}_1    &    \gamma_1    &   w_0    &    \delta_3    &    \CONJ{w}_0       \\
						\alpha_6  \!\!\!\!\!\!\!\!\!\!\!\!\!\!\!\!{\scriptscriptstyle 12}\ \ \ \ \ \ \ \
						u_1     &     \beta_2     &    \CONJ{u}_1    &    \gamma_2    &   w_1    &    \delta_0    &    \CONJ{w}_1       \\
						\alpha_7  \!\!\!\!\!\!\!\!\!\!\!\!\!\!\!\!{\scriptscriptstyle 13}\ \ \ \ \ \ \ \
						u_1     &     \beta_3     &    \CONJ{u}_1    &    \gamma_3    &   w_0    &    \delta_1    &    \CONJ{w}_0       \\
					\end{array}
				\right]
			\end{array}
			\ \ :\ \ 
			\begin{array}{c}
				\alpha_{\ldots}, \beta_{\ldots}, \gamma_{\ldots}, \delta_{\ldots} \ \in\ \REALS    \\
				u_{\ldots}, w_{\ldots} \ \in\ \COMPLEX
			\end{array}\ \ \ 
		\right\}                          \ \ .
	\end{equation}
	$\SPACE{P}_{\FOURIER{2} \otimes \FOURIER{4}}$ is parametrized by $28$ real parameters, so 
	$\undephasedDEFECT\left( 	\FOURIER{2} \otimes \FOURIER{4} \right)\ =\ 28$.
\end{example}

Any single real parameter sitting in $P_{i,j}$ occurs in the $j$-th column at all the row positions belonging to $i + \langle j \rangle$, a coset of
the subgroup of $\Igroup{F}$ generated by $j$, and there are $|\Igroup{F}|/|	\langle j \rangle | \ =\ N/|\langle j \rangle |$ such cosets in $\Igroup{F}$,
where $|H|$ denotes the number of elements in $H$.
Also $\langle j \rangle \ =\ \langle -j \rangle$, hence when $j \neq -j$ and thus the corresponding columns contain complex parameters
we calculate the number of real parameters as if each column of $P$ contained independent real parameters satisfying $P_{i+j,j}\ =\ P_{i,j}$.

In another closely related work by Banica \cite{Banica_first_order_deform}
the author describes a simple parametrization of the space of those real $R$'s which generate, through $R \longrightarrow \Ii R \HADprod \FOURIER{M}$,
the feasible space $\FEASIBLEspace{\FOURIER{M}}$ (see Theorem 4.10 in \cite{Banica_first_order_deform}). 
In both the cited papers by Banica the reader can find a different formula expressing the undephased
defect of a primary Kronecker product of Fourier matrices than those of ours in Theorem \ref{theor_defect_of_primary_Kron_prod_of_Fouriers}.
See for example Proposition 5.5 in \cite{Banica_def_gen_four_matr} which, for an $N \times N$ primary Kronecker product $F$ indexed by $\Igroup{F}$,
features a formula for the quantity $\delta\left( \Igroup{F} \right)$ leading to 
$\undephasedDEFECT(F) \ =\ \left| \Igroup{F} \right| \cdot \delta\left( \Igroup{F} \right) \ =\ N \cdot \delta\left( \Igroup{F} \right)$.

In the next table we have collected inequivalent Fourier matrices with primary Kronecker subproducts, representing (by permutation equivalence) 
all Fourier matrices of the size up to $100$.
They are accompanied by their undephased ($\undephasedDEFECT$) and dephased ($\DEFECT$) defects. 
To represent the matrices we again use partitions as in the previous table,
so for example $F_4 \otimes F_2 \otimes F_3 \otimes F_3$ is represented by $2^{[2\ 1]} 3^{[1\ 1]}$. 
The accompanying values of $\undephasedDEFECT$ and $\DEFECT$ for this particular matrix
concern also  the permutation equivalent $F_{12} \otimes F_6$, $F_4 \otimes F_3 \otimes F_6$,  $F_4 \otimes F_3 \otimes F_3 \otimes F_2$ from our example 
in the paragraph just after the proof of Lemma \ref{lem_multiplicativity_of_D_of_F}. 
In the table $N$ denotes the size of matrices in the rows that follow.

{\center Table 2. Undephased ($\undephasedDEFECT$) and dephased ($\DEFECT$) defects of  inequivalent Fourier matrices of the size up to 100.
}\footnote{The extended table of defects for $N$ ranging between $1$ and $10000$ is available at \cite{Extended_defect_table}.}

\begin{longtable}{|l|l|l|}
\ensuremath{$F$}   &    \ensuremath{\undephasedDEFECT(F)}       &     \ensuremath{\DEFECT(F)} \\
---------------------------------------------  & ------------------ & ------------------
\endfirsthead
\ensuremath{$F$}   &    \ensuremath{\undephasedDEFECT(F)}       &     \ensuremath{\DEFECT(F)}  \\
---------------------------------------------  & ------------------ & ------------------
\endhead

....................$N = 1....................$ & ....................  &  ....................    \\
$1^{[1]}$   &   $1$   &   $0$   \\
....................$N = 2....................$ & ....................  &  ....................    \\
$2^{[1]}$   &   $3$   &   $0$   \\
....................$N = 3....................$ & ....................  &  ....................    \\
$3^{[1]}$   &   $5$   &   $0$   \\
....................$N = 4....................$ & ....................  &  ....................    \\
$2^{[2]}$   &   $8$   &   $1$   \\
$2^{[1\ 1]}$   &   $10$   &   $3$   \\
....................$N = 5....................$ & ....................  &  ....................    \\
$5^{[1]}$   &   $9$   &   $0$   \\
....................$N = 6....................$ & ....................  &  ....................    \\
$2^{[1]}$$3^{[1]}$   &   $15$   &   $4$   \\
....................$N = 7....................$ & ....................  &  ....................    \\
$7^{[1]}$   &   $13$   &   $0$   \\
....................$N = 8....................$ & ....................  &  ....................    \\
$2^{[3]}$   &   $20$   &   $5$   \\
$2^{[2\ 1]}$   &   $28$   &   $13$   \\
$2^{[1\ 1\ 1]}$   &   $36$   &   $21$   \\
....................$N = 9....................$ & ....................  &  ....................    \\
$3^{[2]}$   &   $21$   &   $4$   \\
$3^{[1\ 1]}$   &   $33$   &   $16$   \\
....................$N = 10....................$ & ....................  &  ....................    \\
$2^{[1]}$$5^{[1]}$   &   $27$   &   $8$   \\
....................$N = 11....................$ & ....................  &  ....................    \\
$11^{[1]}$   &   $21$   &   $0$   \\
....................$N = 12....................$ & ....................  &  ....................    \\
$2^{[2]}$$3^{[1]}$   &   $40$   &   $17$   \\
$2^{[1\ 1]}$$3^{[1]}$   &   $50$   &   $27$   \\
....................$N = 13....................$ & ....................  &  ....................    \\
$13^{[1]}$   &   $25$   &   $0$   \\
....................$N = 14....................$ & ....................  &  ....................    \\
$2^{[1]}$$7^{[1]}$   &   $39$   &   $12$   \\
....................$N = 15....................$ & ....................  &  ....................    \\
$3^{[1]}$$5^{[1]}$   &   $45$   &   $16$   \\
....................$N = 16....................$ & ....................  &  ....................    \\
$2^{[4]}$   &   $48$   &   $17$   \\
$2^{[3\ 1]}$   &   $72$   &   $41$   \\
$2^{[2\ 2]}$   &   $88$   &   $57$   \\
$2^{[2\ 1\ 1]}$   &   $104$   &   $73$   \\
$2^{[1\ 1\ 1\ 1]}$   &   $136$   &   $105$   \\
....................$N = 17....................$ & ....................  &  ....................    \\
$17^{[1]}$   &   $33$   &   $0$   \\
....................$N = 18....................$ & ....................  &  ....................    \\
$2^{[1]}$$3^{[2]}$   &   $63$   &   $28$   \\
$2^{[1]}$$3^{[1\ 1]}$   &   $99$   &   $64$   \\
....................$N = 19....................$ & ....................  &  ....................    \\
$19^{[1]}$   &   $37$   &   $0$   \\
....................$N = 20....................$ & ....................  &  ....................    \\
$2^{[2]}$$5^{[1]}$   &   $72$   &   $33$   \\
$2^{[1\ 1]}$$5^{[1]}$   &   $90$   &   $51$   \\
....................$N = 21....................$ & ....................  &  ....................    \\
$3^{[1]}$$7^{[1]}$   &   $65$   &   $24$   \\
....................$N = 22....................$ & ....................  &  ....................    \\
$2^{[1]}$$11^{[1]}$   &   $63$   &   $20$   \\
....................$N = 23....................$ & ....................  &  ....................    \\
$23^{[1]}$   &   $45$   &   $0$   \\
....................$N = 24....................$ & ....................  &  ....................    \\
$2^{[3]}$$3^{[1]}$   &   $100$   &   $53$   \\
$2^{[2\ 1]}$$3^{[1]}$   &   $140$   &   $93$   \\
$2^{[1\ 1\ 1]}$$3^{[1]}$   &   $180$   &   $133$   \\
....................$N = 25....................$ & ....................  &  ....................    \\
$5^{[2]}$   &   $65$   &   $16$   \\
$5^{[1\ 1]}$   &   $145$   &   $96$   \\
....................$N = 26....................$ & ....................  &  ....................    \\
$2^{[1]}$$13^{[1]}$   &   $75$   &   $24$   \\
....................$N = 27....................$ & ....................  &  ....................    \\
$3^{[3]}$   &   $81$   &   $28$   \\
$3^{[2\ 1]}$   &   $153$   &   $100$   \\
$3^{[1\ 1\ 1]}$   &   $261$   &   $208$   \\
....................$N = 28....................$ & ....................  &  ....................    \\
$2^{[2]}$$7^{[1]}$   &   $104$   &   $49$   \\
$2^{[1\ 1]}$$7^{[1]}$   &   $130$   &   $75$   \\
....................$N = 29....................$ & ....................  &  ....................    \\
$29^{[1]}$   &   $57$   &   $0$   \\
....................$N = 30....................$ & ....................  &  ....................    \\
$2^{[1]}$$3^{[1]}$$5^{[1]}$   &   $135$   &   $76$   \\
....................$N = 31....................$ & ....................  &  ....................    \\
$31^{[1]}$   &   $61$   &   $0$   \\
....................$N = 32....................$ & ....................  &  ....................    \\
$2^{[5]}$   &   $112$   &   $49$   \\
$2^{[4\ 1]}$   &   $176$   &   $113$   \\
$2^{[3\ 2]}$   &   $240$   &   $177$   \\
$2^{[3\ 1\ 1]}$   &   $272$   &   $209$   \\
$2^{[2\ 2\ 1]}$   &   $336$   &   $273$   \\
$2^{[2\ 1\ 1\ 1]}$   &   $400$   &   $337$   \\
$2^{[1\ 1\ 1\ 1\ 1]}$   &   $528$   &   $465$   \\
....................$N = 33....................$ & ....................  &  ....................    \\
$3^{[1]}$$11^{[1]}$   &   $105$   &   $40$   \\
....................$N = 34....................$ & ....................  &  ....................    \\
$2^{[1]}$$17^{[1]}$   &   $99$   &   $32$   \\
....................$N = 35....................$ & ....................  &  ....................    \\
$5^{[1]}$$7^{[1]}$   &   $117$   &   $48$   \\
....................$N = 36....................$ & ....................  &  ....................    \\
$2^{[2]}$$3^{[2]}$   &   $168$   &   $97$   \\
$2^{[2]}$$3^{[1\ 1]}$   &   $264$   &   $193$   \\
$2^{[1\ 1]}$$3^{[2]}$   &   $210$   &   $139$   \\
$2^{[1\ 1]}$$3^{[1\ 1]}$   &   $330$   &   $259$   \\
....................$N = 37....................$ & ....................  &  ....................    \\
$37^{[1]}$   &   $73$   &   $0$   \\
....................$N = 38....................$ & ....................  &  ....................    \\
$2^{[1]}$$19^{[1]}$   &   $111$   &   $36$   \\
....................$N = 39....................$ & ....................  &  ....................    \\
$3^{[1]}$$13^{[1]}$   &   $125$   &   $48$   \\
....................$N = 40....................$ & ....................  &  ....................    \\
$2^{[3]}$$5^{[1]}$   &   $180$   &   $101$   \\
$2^{[2\ 1]}$$5^{[1]}$   &   $252$   &   $173$   \\
$2^{[1\ 1\ 1]}$$5^{[1]}$   &   $324$   &   $245$   \\
....................$N = 41....................$ & ....................  &  ....................    \\
$41^{[1]}$   &   $81$   &   $0$   \\
....................$N = 42....................$ & ....................  &  ....................    \\
$2^{[1]}$$3^{[1]}$$7^{[1]}$   &   $195$   &   $112$   \\
....................$N = 43....................$ & ....................  &  ....................    \\
$43^{[1]}$   &   $85$   &   $0$   \\
....................$N = 44....................$ & ....................  &  ....................    \\
$2^{[2]}$$11^{[1]}$   &   $168$   &   $81$   \\
$2^{[1\ 1]}$$11^{[1]}$   &   $210$   &   $123$   \\
....................$N = 45....................$ & ....................  &  ....................    \\
$3^{[2]}$$5^{[1]}$   &   $189$   &   $100$   \\
$3^{[1\ 1]}$$5^{[1]}$   &   $297$   &   $208$   \\
....................$N = 46....................$ & ....................  &  ....................    \\
$2^{[1]}$$23^{[1]}$   &   $135$   &   $44$   \\
....................$N = 47....................$ & ....................  &  ....................    \\
$47^{[1]}$   &   $93$   &   $0$   \\
....................$N = 48....................$ & ....................  &  ....................    \\
$2^{[4]}$$3^{[1]}$   &   $240$   &   $145$   \\
$2^{[3\ 1]}$$3^{[1]}$   &   $360$   &   $265$   \\
$2^{[2\ 2]}$$3^{[1]}$   &   $440$   &   $345$   \\
$2^{[2\ 1\ 1]}$$3^{[1]}$   &   $520$   &   $425$   \\
$2^{[1\ 1\ 1\ 1]}$$3^{[1]}$   &   $680$   &   $585$   \\
....................$N = 49....................$ & ....................  &  ....................    \\
$7^{[2]}$   &   $133$   &   $36$   \\
$7^{[1\ 1]}$   &   $385$   &   $288$   \\
....................$N = 50....................$ & ....................  &  ....................    \\
$2^{[1]}$$5^{[2]}$   &   $195$   &   $96$   \\
$2^{[1]}$$5^{[1\ 1]}$   &   $435$   &   $336$   \\
....................$N = 51....................$ & ....................  &  ....................    \\
$3^{[1]}$$17^{[1]}$   &   $165$   &   $64$   \\
....................$N = 52....................$ & ....................  &  ....................    \\
$2^{[2]}$$13^{[1]}$   &   $200$   &   $97$   \\
$2^{[1\ 1]}$$13^{[1]}$   &   $250$   &   $147$   \\
....................$N = 53....................$ & ....................  &  ....................    \\
$53^{[1]}$   &   $105$   &   $0$   \\
....................$N = 54....................$ & ....................  &  ....................    \\
$2^{[1]}$$3^{[3]}$   &   $243$   &   $136$   \\
$2^{[1]}$$3^{[2\ 1]}$   &   $459$   &   $352$   \\
$2^{[1]}$$3^{[1\ 1\ 1]}$   &   $783$   &   $676$   \\
....................$N = 55....................$ & ....................  &  ....................    \\
$5^{[1]}$$11^{[1]}$   &   $189$   &   $80$   \\
....................$N = 56....................$ & ....................  &  ....................    \\
$2^{[3]}$$7^{[1]}$   &   $260$   &   $149$   \\
$2^{[2\ 1]}$$7^{[1]}$   &   $364$   &   $253$   \\
$2^{[1\ 1\ 1]}$$7^{[1]}$   &   $468$   &   $357$   \\
....................$N = 57....................$ & ....................  &  ....................    \\
$3^{[1]}$$19^{[1]}$   &   $185$   &   $72$   \\
....................$N = 58....................$ & ....................  &  ....................    \\
$2^{[1]}$$29^{[1]}$   &   $171$   &   $56$   \\
....................$N = 59....................$ & ....................  &  ....................    \\
$59^{[1]}$   &   $117$   &   $0$   \\
....................$N = 60....................$ & ....................  &  ....................    \\
$2^{[2]}$$3^{[1]}$$5^{[1]}$   &   $360$   &   $241$   \\
$2^{[1\ 1]}$$3^{[1]}$$5^{[1]}$   &   $450$   &   $331$   \\
....................$N = 61....................$ & ....................  &  ....................    \\
$61^{[1]}$   &   $121$   &   $0$   \\
....................$N = 62....................$ & ....................  &  ....................    \\
$2^{[1]}$$31^{[1]}$   &   $183$   &   $60$   \\
....................$N = 63....................$ & ....................  &  ....................    \\
$3^{[2]}$$7^{[1]}$   &   $273$   &   $148$   \\
$3^{[1\ 1]}$$7^{[1]}$   &   $429$   &   $304$   \\
....................$N = 64....................$ & ....................  &  ....................    \\
$2^{[6]}$   &   $256$   &   $129$   \\
$2^{[5\ 1]}$   &   $416$   &   $289$   \\
$2^{[4\ 2]}$   &   $608$   &   $481$   \\
$2^{[4\ 1\ 1]}$   &   $672$   &   $545$   \\
$2^{[3\ 3]}$   &   $736$   &   $609$   \\
$2^{[3\ 2\ 1]}$   &   $928$   &   $801$   \\
$2^{[3\ 1\ 1\ 1]}$   &   $1056$   &   $929$   \\
$2^{[2\ 2\ 2]}$   &   $1184$   &   $1057$   \\
$2^{[2\ 2\ 1\ 1]}$   &   $1312$   &   $1185$   \\
$2^{[2\ 1\ 1\ 1\ 1]}$   &   $1568$   &   $1441$   \\
$2^{[1\ 1\ 1\ 1\ 1\ 1]}$   &   $2080$   &   $1953$   \\
....................$N = 65....................$ & ....................  &  ....................    \\
$5^{[1]}$$13^{[1]}$   &   $225$   &   $96$   \\
....................$N = 66....................$ & ....................  &  ....................    \\
$2^{[1]}$$3^{[1]}$$11^{[1]}$   &   $315$   &   $184$   \\
....................$N = 67....................$ & ....................  &  ....................    \\
$67^{[1]}$   &   $133$   &   $0$   \\
....................$N = 68....................$ & ....................  &  ....................    \\
$2^{[2]}$$17^{[1]}$   &   $264$   &   $129$   \\
$2^{[1\ 1]}$$17^{[1]}$   &   $330$   &   $195$   \\
....................$N = 69....................$ & ....................  &  ....................    \\
$3^{[1]}$$23^{[1]}$   &   $225$   &   $88$   \\
....................$N = 70....................$ & ....................  &  ....................    \\
$2^{[1]}$$5^{[1]}$$7^{[1]}$   &   $351$   &   $212$   \\
....................$N = 71....................$ & ....................  &  ....................    \\
$71^{[1]}$   &   $141$   &   $0$   \\
....................$N = 72....................$ & ....................  &  ....................    \\
$2^{[3]}$$3^{[2]}$   &   $420$   &   $277$   \\
$2^{[3]}$$3^{[1\ 1]}$   &   $660$   &   $517$   \\
$2^{[2\ 1]}$$3^{[2]}$   &   $588$   &   $445$   \\
$2^{[2\ 1]}$$3^{[1\ 1]}$   &   $924$   &   $781$   \\
$2^{[1\ 1\ 1]}$$3^{[2]}$   &   $756$   &   $613$   \\
$2^{[1\ 1\ 1]}$$3^{[1\ 1]}$   &   $1188$   &   $1045$   \\
....................$N = 73....................$ & ....................  &  ....................    \\
$73^{[1]}$   &   $145$   &   $0$   \\
....................$N = 74....................$ & ....................  &  ....................    \\
$2^{[1]}$$37^{[1]}$   &   $219$   &   $72$   \\
....................$N = 75....................$ & ....................  &  ....................    \\
$3^{[1]}$$5^{[2]}$   &   $325$   &   $176$   \\
$3^{[1]}$$5^{[1\ 1]}$   &   $725$   &   $576$   \\
....................$N = 76....................$ & ....................  &  ....................    \\
$2^{[2]}$$19^{[1]}$   &   $296$   &   $145$   \\
$2^{[1\ 1]}$$19^{[1]}$   &   $370$   &   $219$   \\
....................$N = 77....................$ & ....................  &  ....................    \\
$7^{[1]}$$11^{[1]}$   &   $273$   &   $120$   \\
....................$N = 78....................$ & ....................  &  ....................    \\
$2^{[1]}$$3^{[1]}$$13^{[1]}$   &   $375$   &   $220$   \\
....................$N = 79....................$ & ....................  &  ....................    \\
$79^{[1]}$   &   $157$   &   $0$   \\
....................$N = 80....................$ & ....................  &  ....................    \\
$2^{[4]}$$5^{[1]}$   &   $432$   &   $273$   \\
$2^{[3\ 1]}$$5^{[1]}$   &   $648$   &   $489$   \\
$2^{[2\ 2]}$$5^{[1]}$   &   $792$   &   $633$   \\
$2^{[2\ 1\ 1]}$$5^{[1]}$   &   $936$   &   $777$   \\
$2^{[1\ 1\ 1\ 1]}$$5^{[1]}$   &   $1224$   &   $1065$   \\
....................$N = 81....................$ & ....................  &  ....................    \\
$3^{[4]}$   &   $297$   &   $136$   \\
$3^{[3\ 1]}$   &   $621$   &   $460$   \\
$3^{[2\ 2]}$   &   $945$   &   $784$   \\
$3^{[2\ 1\ 1]}$   &   $1269$   &   $1108$   \\
$3^{[1\ 1\ 1\ 1]}$   &   $2241$   &   $2080$   \\
....................$N = 82....................$ & ....................  &  ....................    \\
$2^{[1]}$$41^{[1]}$   &   $243$   &   $80$   \\
....................$N = 83....................$ & ....................  &  ....................    \\
$83^{[1]}$   &   $165$   &   $0$   \\
....................$N = 84....................$ & ....................  &  ....................    \\
$2^{[2]}$$3^{[1]}$$7^{[1]}$   &   $520$   &   $353$   \\
$2^{[1\ 1]}$$3^{[1]}$$7^{[1]}$   &   $650$   &   $483$   \\
....................$N = 85....................$ & ....................  &  ....................    \\
$5^{[1]}$$17^{[1]}$   &   $297$   &   $128$   \\
....................$N = 86....................$ & ....................  &  ....................    \\
$2^{[1]}$$43^{[1]}$   &   $255$   &   $84$   \\
....................$N = 87....................$ & ....................  &  ....................    \\
$3^{[1]}$$29^{[1]}$   &   $285$   &   $112$   \\
....................$N = 88....................$ & ....................  &  ....................    \\
$2^{[3]}$$11^{[1]}$   &   $420$   &   $245$   \\
$2^{[2\ 1]}$$11^{[1]}$   &   $588$   &   $413$   \\
$2^{[1\ 1\ 1]}$$11^{[1]}$   &   $756$   &   $581$   \\
....................$N = 89....................$ & ....................  &  ....................    \\
$89^{[1]}$   &   $177$   &   $0$   \\
....................$N = 90....................$ & ....................  &  ....................    \\
$2^{[1]}$$3^{[2]}$$5^{[1]}$   &   $567$   &   $388$   \\
$2^{[1]}$$3^{[1\ 1]}$$5^{[1]}$   &   $891$   &   $712$   \\
....................$N = 91....................$ & ....................  &  ....................    \\
$7^{[1]}$$13^{[1]}$   &   $325$   &   $144$   \\
....................$N = 92....................$ & ....................  &  ....................    \\
$2^{[2]}$$23^{[1]}$   &   $360$   &   $177$   \\
$2^{[1\ 1]}$$23^{[1]}$   &   $450$   &   $267$   \\
....................$N = 93....................$ & ....................  &  ....................    \\
$3^{[1]}$$31^{[1]}$   &   $305$   &   $120$   \\
....................$N = 94....................$ & ....................  &  ....................    \\
$2^{[1]}$$47^{[1]}$   &   $279$   &   $92$   \\
....................$N = 95....................$ & ....................  &  ....................    \\
$5^{[1]}$$19^{[1]}$   &   $333$   &   $144$   \\
....................$N = 96....................$ & ....................  &  ....................    \\
$2^{[5]}$$3^{[1]}$   &   $560$   &   $369$   \\
$2^{[4\ 1]}$$3^{[1]}$   &   $880$   &   $689$   \\
$2^{[3\ 2]}$$3^{[1]}$   &   $1200$   &   $1009$   \\
$2^{[3\ 1\ 1]}$$3^{[1]}$   &   $1360$   &   $1169$   \\
$2^{[2\ 2\ 1]}$$3^{[1]}$   &   $1680$   &   $1489$   \\
$2^{[2\ 1\ 1\ 1]}$$3^{[1]}$   &   $2000$   &   $1809$   \\
$2^{[1\ 1\ 1\ 1\ 1]}$$3^{[1]}$   &   $2640$   &   $2449$   \\
....................$N = 97....................$ & ....................  &  ....................    \\
$97^{[1]}$   &   $193$   &   $0$   \\
....................$N = 98....................$ & ....................  &  ....................    \\
$2^{[1]}$$7^{[2]}$   &   $399$   &   $204$   \\
$2^{[1]}$$7^{[1\ 1]}$   &   $1155$   &   $960$   \\
....................$N = 99....................$ & ....................  &  ....................    \\
$3^{[2]}$$11^{[1]}$   &   $441$   &   $244$   \\
$3^{[1\ 1]}$$11^{[1]}$   &   $693$   &   $496$   \\
....................$N = 100....................$ & ....................  &  ....................    \\
$2^{[2]}$$5^{[2]}$   &   $520$   &   $321$   \\
$2^{[2]}$$5^{[1\ 1]}$   &   $1160$   &   $961$   \\
$2^{[1\ 1]}$$5^{[2]}$   &   $650$   &   $451$   \\
$2^{[1\ 1]}$$5^{[1\ 1]}$   &   $1450$   &   $1251$   \\

\end{longtable}

%
%

\section{Conclusions}
	\label{sec_conclusions}
	To complete this long work, let us recap the main facts established in it.
For the consistency and existence of all the characterizations of defect in this article, as well as for the correctness of 
the bounds on the dimensions of manifolds involving the defect,  we assume throughout that the considered
unitary matrices have no zero entries.

The undephased defect $\undephasedDEFECT(U)$ of an $N \times N$ unitary matrix $U$ quantifies by how much the dimension
of the space $\TANGENTspace{\UNITARY}{U}$ (tangent to the $N^2$ dimensional unitary manifold $\UNITARY$ at $U$) is reduced 
when it is mapped by the differential $\DIFFERENTIAL{f}{U}$ into $\TANGENTspace{\BISTOCHSPACE}{f(U)}\ =\ \BISTOCHSPACE$
(the $(N-1)^2$ doubly stochastic hyperplane), where $f$ squares the moduli of the entries of its argument: 
	$f\left( \left[ V_{i,j} \right]_{i,j=1\ldots N} \right) \ =\ \left[ \ABSOLUTEvalue{V_{i,j}}^2 \right]_{i,j=1\ldots N}$  
(see (\ref{eq_dim_ker_image_sum_equality_for_differential_of_f})). The dephased defect $\DEFECT(U)\ =\ \undephasedDEFECT(U) - (2N-1)$ 
measures this reduction  with respect to the level of 
$(N-1)^2\ =\ \DIMR\left(\BISTOCHSPACE\right)$ instead of $N^2\ =\ \DIMR\left( \UNITARY \right)$
(see (\ref{eq_defect_original_definition})).

$\undephasedDEFECT(U)$ also is the dimension of the intersection $\FEASIBLEspace{U}$ of the tangent space $\TANGENTspace{\UNITARY}{U}$
and $\TANGENTspace{\FIXEDmoduliMATRICES{U}}{U}$  (see (\ref{eq_undepased_defect})), 
where $\FIXEDmoduliMATRICES{U}$ is the $N^2$ dimensional manifold
containing $U$ and composed of all the matrices with the moduli of their entries identical with those in $U$. 
The defect $\undephasedDEFECT(U)$ thus serves as an upper bound (see (\ref{eq_undephased_manifold_dim_estimation}))
on the dimensions of manifolds stemming from $U$ and contained in
$\UNITARY \cap \FIXEDmoduliMATRICES{U}$ (e.g. complex Hadamard matrices), while $\DEFECT(U)$ plays a similar role
(see (\ref{eq_dephased_manifold_tangent_space_dim_estimation})) when the dephased 
(e.g. with the $1$st row and column of their members fixed) manifolds are considered.
    $\undephasedDEFECT(U)$ and $\DEFECT(U)$ are well suited for rescaled unitary matrices such as Fourier matrices
(see (\ref{eq_undephased_defect_for_rescaled_unitary}),(\ref{eq_dephased_defect_for_rescaled_unitary})).

$\undephasedDEFECT(U)$ is the dimension of the complex eigenspace $\Vu{U}{1}$, associated with $\lambda = 1$,
of the $\UinnerPRODUCT{U}{}{}$-unitary operator $\Iu{U}$  on $\COMPLEXmatricesNxN{N}$
(see Theorem \ref{theor_unitarity_and_eigenspaces_of_Iu} {\bf a),e)})  , 
where
$\Iu{U}(X)\ =\ U^{\ominus 1}  	\HADprod     \left(      U      \left(    \CONJ{X} \HADprod U    \right)^{*}    U  \right)$
and $\HADprod$ is the entrywise product while $U^{\ominus 1}$ is the entrywise inverse of $U$, 
 with the trace equal to $N$ (see Lemma \ref{lem_properties_of_Cu_Du_Iu} {\bf g)}).
The above intersection $\FEASIBLEspace{U}$ of the tangent spaces is parametrized by the $\undephasedDEFECT(U)$ dimensional
real subspace $\IMspaceOF{\Vu{U}{1}}$ of $\Vu{U}{1}$: $\FEASIBLEspace{U}\ =\ \IMspaceOF{\Vu{U}{1}}  \HADprod  U$
(see Theorem \ref{theor_unitarity_and_eigenspaces_of_Iu} {\bf e)}).

Operator $\Iu{U}$ satisfies 
$\Iu{\ithMATRIX{U}{1} \otimes \ldots \otimes \ithMATRIX{U}{r}}\left( \ithMATRIX{X}{1} \otimes \ldots \otimes \ithMATRIX{X}{r} \right)
\ =\ \Iu{\ithMATRIX{U}{1}}\left( \ithMATRIX{X}{1} \right) \otimes \ldots \otimes \Iu{\ithMATRIX{U}{r}}\left( \ithMATRIX{X}{r} \right)$,
where $\ithMATRIX{U}{k}, \ithMATRIX{X}{k}, \Iu{\ithMATRIX{U}{k}}\left( \ithMATRIX{X}{k} \right)\ \in\ \COMPLEXmatricesNxN{N_k}$,
that is $\Iu{\ithMATRIX{U}{1} \otimes \ldots \otimes \ithMATRIX{U}{r}}$ is the tensor product of $\Iu{\ithMATRIX{U}{k}}$'s
(see Lemma \ref{lem_Iu_for_Kron_prod_of_unitaries} {\bf c)}), hence its eigenvalues are products of the eigenvalues of $\Iu{\ithMATRIX{U}{k}}$'s
and the corresponding eigenvectors (matrices) are Kronecker products of the eigenvectors of $\Iu{\ithMATRIX{U}{k}}$'s
corresponding to these eigenvalues of $\Iu{\ithMATRIX{U}{k}}$'s (see the comment following the proof of Lemma \ref{lem_Iu_for_Kron_prod_of_unitaries}).
Besides, from $\UinnerPRODUCT{\ithMATRIX{U}{k}}{}{}$-orthonormal eigenbases  for $\Iu{\ithMATRIX{U}{k}}$'s we can construct
an $\UinnerPRODUCT{\ithMATRIX{U}{1} \otimes \ldots \otimes \ithMATRIX{U}{r}}{}{}$-orthonormal eigenbasis for 
$\Iu{\ithMATRIX{U}{1} \otimes \ldots \otimes \ithMATRIX{U}{r}}$ (see Lemma \ref{lem_weigthed_U_inner_product_induced}) using
the Kronecker product $\otimes$. 
    $\undephasedDEFECT\left( \ithMATRIX{U}{1} \otimes \ldots \otimes \ithMATRIX{U}{r} \right) 
			\geq
	\undephasedDEFECT\left( \ithMATRIX{U}{1} \right) \cdot \ldots \cdot \undephasedDEFECT\left( \ithMATRIX{U}{r} \right)$
because 
	$\Vu{ \ithMATRIX{U}{1} \otimes \ldots \otimes \ithMATRIX{U}{r} }{1} 
			\supset 
	\Vu{\ithMATRIX{U}{1}}{1} \otimes \ldots \otimes  \Vu{\ithMATRIX{U}{r}}{1}$,
where on the right we have the space spanned by all the  Kronecker products of the elements of $\Vu{\ithMATRIX{U}{k}}{1}$'s
of the form $\bigotimes_{k=1\ldots r} X_k$, where $X_k \in \Vu{\ithMATRIX{U}{k}}{1}$
(see Corollary \ref{cor_V1_spaces_Kron_multiplied_are_contained_in_V1_space} {\bf a),c)}).

If $U$ is an $N \times N$ real orthogonal matrix then $\Iu{U}$ is $\UinnerPRODUCT{U}{}{}$-self adjoint and 
$\undephasedDEFECT(U)\ =\ N(N+1)/2$ (see Lemma \ref{lem_Iu_selfadjoint_for_U_real_orthogonal}), and
so is the value of undephased defect for any $U$ equivalent to a real orthogonal matrix, in particular for $U$ being a Kronecker product
of $2 \times 2$ unitary matrices when $N$ is even (see Corollary \ref{cor_spectral_decomposition_of_Iu_for_kron_prod_of_2x2_matrices}).

For an $N \times N$ unitary $U$ the eigenspace $\Vu{U}{1}$ contains at least the $(2N-1)$ dimensional complex space
$\TRIVIALspace{U}\ 
	=\ 
\left\{a \ONESvect^T + \ONESvect b^T\ :\ \ a,b \in \columnVECTORS{N},\ \   \ONESvect\ =\ [1\ 1\ \ldots\ 1]^T \ \right\}$, 
	and its real subspace $\IMspaceOF{\Vu{U}{1}}$ contains the $(2N-1)$ dimensional real space 
$\IMspaceOF{\TRIVIALspace{U}}
	\ =\ 
\{\Ii (a \ONESvect^T + \ONESvect b^T)\ :\ \ a,b \in \REALcolumnVECTORS{N} \}$ (see Theorem \ref{theor_unitarity_and_eigenspaces_of_Iu} {\bf c)}).
Therefore $\FEASIBLEspace{U}\ =\ \IMspaceOF{\Vu{U}{1}} \HADprod U$ contains the real $(2N-1)$ dimensional space
$\TANGENTspace{\ENPHASEDmatrices{U}}{U}
		\ =\   
\{\Ii (a \ONESvect^T + \ONESvect b^T) \HADprod U \ :\ \ a,b \in \REALcolumnVECTORS{N} \}$ (see Theorem \ref{theor_unitarity_and_eigenspaces_of_Iu} {\bf e)}),
the space tangent at $U$ to the enphasing manifold
$\ENPHASEDmatrices{U}\ =\ \left\{  D_r U D_c\ :\ \ D_r,\ D_c\ \ \ \mbox{unitary diagonal} \right\}$ 
(see the paragraph containing (\ref{eq_enphasing_manifold})).
For this and other reasons
$\undephasedDEFECT\left( \ithMATRIX{U}{1} \otimes \ldots \otimes \ithMATRIX{U}{r} \right)
		\geq
\left(   \prod_{N_k > 2}   \left( 2N_k - 1 \right)   \right)            \ \cdot\       2^{x-1} \left( 2^x  +   1 \right)$,
where $N_k$ is the size of $\ithMATRIX{U}{k}$,  $r\in\{1,2,\ldots\}$ and  $x\in\{0,1,\ldots\}$ is the number of the $2 \times 2$ Kronecker factors 
(see (\ref{eq_kron_prod_undephased_defect_bound_without_2s}),(\ref{eq_kron_prod_undephased_defect_bound_with_2s})).

Next we consider the action $(S,T)X\ =\ SXT^{-1}$ of the group $\PP$ of pairs of enphased permutation matrices on the
$N \times N$ complex matrices (see Definition \ref{def_pairs_action_stabilizer}). 
Unitaries $U$ and $V$ satisfying $V\ =\ (S,T)U$ for some $(S,T) \in \PP$ are called equivalent
(see (\ref{eq_equivalence})), or permutation equivalent if for some such $(S,T)$ it is equal to $\left( \leftPERM{S}{T},\ \rightPERM{S}{T} \right)$, the 
underlying pair of permutation matrices (see the paragraph before Lemma \ref{lem_permutation_equivalence_PFRequalsG_matrices_determined_by_isomorphisms}). 
Operators $\Iu{U}$ and $\Iu{V=(S,T)U}$ are similar 
(see Corollary \ref{cor_Iu_Iv_operators_for_equivalent_U_V_are_similar} {\bf a)}), 
besides $\pairOPERATOR{S}{T}\left( \FEASIBLEspace{U} \right)\ =\ \FEASIBLEspace{V}$,\ \ 
$\pairOPERATOR{\leftPERM{S}{T}}{\rightPERM{S}{T}}\left( \Vu{U}{1} \right)\ =\ \Vu{V}{1}$,\ \ 
$\pairOPERATOR{\leftPERM{S}{T}}{\rightPERM{S}{T}}\left( \IMspaceOF{\Vu{U}{1}} \right)\ =\ \IMspaceOF{\Vu{V}{1}}$,
and thus $\undephasedDEFECT(U)\ =\ \undephasedDEFECT(V)$ and $\DEFECT(U)\ =\ \DEFECT(V)$ 
(see  Corollary \ref{cor_DU_equal_DV_for_U_V_equivalent}), where 
$\pairOPERATOR{S}{T}:\ X \longrightarrow (S,T)X$ is an isomorphism on $\COMPLEXmatricesNxN{N}$ for $(S,T)\in\PP$.
$\Iu{U}$ and $\pairOPERATOR{\leftPERM{S}{T}}{\rightPERM{S}{T}}$ commute for any $(S,T) \in \STAB{\PP}{U}$,
the stabilizer of $U$ in $\PP$ under the considered action 
(see Corollary \ref{eq_Iu_Ju_commute_with_ops_related_to_ST_from_the_stabilizer_of_U} {\bf a)}
  and Lemma \ref{eq_properties_of_representations_of_STAB_U} {\bf c)} ),
so the eigenspaces of $\Iu{U}$ can be searched for among the irreducible invariant subspaces of the complex
representation $\STAB{\PP}{U} \supset E \ni (S,T) \longrightarrow \pairOPERATOR{S}{T}$, whose image for any subgroup $E$ of $\STAB{\PP}{U}$
is finite (see the comments following Lemma \ref{eq_properties_of_representations_of_STAB_U}).

Further we concentrate on the Fourier matrix $F\ =\ \FOURIER{N_1} \otimes \ldots \otimes \FOURIER{N_r}$, 
being the character table of  the finite abelian group $\ZN{N_1} \times \ldots \times \ZN{N_r}$, called here the indexing group of $F$ and denoted $\Igroup{F}$
(see the paragraphs containing (\ref{eq_Fourier_matrix}),(\ref{eq_Fourier_Kron_prod})). 
We describe $\STAB{\PsPs}{F}$, the stabilizer of $F$
in $\PsPs \subset \PP$ consisting of all the pairs of permutation matrices  (see Corollary \ref{cor_permutational_MAP_FG} {\bf a)}), 
$\STAB{\PzPz{F}}{F}$, the stabilizer of $F$ in $\PzPz{F} \subset \PP$ consisting of all the enphased shift matrices 
(see Lemma \ref{lem_G_F_image_under_STmap_is_PzPz_stabilizer_of_F}).
Then $\STAB{\PP}{F}\ =\STAB{\PsPs}{F}   \cdot   \STAB{\PzPz{F}}{F}$, it is an internal semidirect product of the factors on the right,
where the second factor is normal (see Theorem \ref{theor_enph_perms_stabilizer_of_F_structure}). As far as the equivalence of 
Fourier matrices is concerned, we have established  that Fourier matrices are equivalent if and only if they are permutation equivalent
(see Lemma \ref{lem_F_G_equivalent_iff_F_G_permutation_equivalent} {\bf a)}), but the latter property holds if and only if the underlying
indexing groups $\Igroup{F}$ and $\Igroup{G}$ are isomorphic (see Corollary \ref{cor_Fouriers_perm_equivalent_iff_their_indexing_groups_isomorphic}).
We have also presented a factorization of $\MAP{\PP}{F}{G}$, composed of all $(S,T) \in \PP$ such that $(S,T)F\ =\ G$. Namely
$\MAP{\PP}{F}{G}
		\ =\ 
\MAP{\PsPs}{F}{G} \cdot \STAB{\PzPz{F}}{F}
		\ =\ 
\STAB{\PzPz{G}}{G}     \cdot     \MAP{\PsPs}{F}{G}$
(see Lemma \ref{lem_F_G_equivalent_iff_F_G_permutation_equivalent} {\bf b)}).

Getting back to the defect, pairs $\left( F_{:,s} F_{:,t}^T\ \mbox{i.e. column times column transposed},\ \ F_{s,t} \right)$,
for $s,t \in \Igroup{F}$ (the indexing group of $F$) are (eigenvector, eigenvalue) pairs for operator $\Iu{\frac{1}{\sqrt{N}} F}$, 
where $N$ is the size of $F$, i.e. the order of $\Igroup{F}$ (see Theorem \ref{theor_eigenbasis_of_I_F}). 
Hence the quantity $\undephasedDEFECT(F)\ =\ \undephasedDEFECT\left( \frac{1}{\sqrt{N}} F \right)$ 
(see (\ref{eq_undephased_defect_for_rescaled_unitary})) is equal to the number of entries of $F$ equal to $1$.
We have used this characterization to derived formulas expressing $\undephasedDEFECT(F)$ for
$F\ =\ \FOURIER{a^{k_1}} \otimes \ldots \otimes \FOURIER{a^{k_s}}$, where $a$ is a prime number 
(see Theorem \ref{theor_defect_of_primary_Kron_prod_of_Fouriers}). $F$ is called primary in such a case and then
$\DEFECT(F)$ is divided by $(a-1)^2$  -- see Corollary \ref{cor_D_and_d_factors_for_primary_kron_prod_F} {\bf b)}.

Every Fourier matrix $F$ (as we define it) can be transformed, by a permutation equivalence operation, into a Kronecker product
$F'\ =\ \ithMATRIX{G}{1} \otimes \ldots \otimes \ithMATRIX{G}{r}$ with $\ithMATRIX{G}{k}$'s being primary Fourier matrices associated with distinct primes
(see the paragraph preceding Lemma \ref{lem_multiplicativity_of_D_of_F}), and then 
$\undephasedDEFECT(F)\ \ =\ \undephasedDEFECT(F')\ =\ 
\undephasedDEFECT\left( \ithMATRIX{G}{1} \right) \cdot \ldots \cdot \undephasedDEFECT\left( \ithMATRIX{G}{r} \right)$.
This last equality holds not only when  $\ithMATRIX{G}{k}$'s are primary associated with distinct primes, but always when they are
Kronecker factors of $F'$  of the sizes which are pairwise relatively prime -- see Lemma \ref{lem_multiplicativity_of_D_of_F} {\bf a)}.

We have included a table of $\undephasedDEFECT(F)$'s and $\DEFECT(F)$'s, as polynomials of the variable $a$, for the primary Fourier matrices $F$ associated with
a prime number $a$ and of the size $N$ (a power of $a$) satisfying  $a \leq N \leq a^6$. The table -- see Table 1. -- is located after Corollary \ref{cor_D_and_d_factors_for_primary_kron_prod_F}.
Another table is the one featuring $\undephasedDEFECT(F)$'s and $\DEFECT(F)$'s
of the representatives $F$ of all the (permutation) equivalence classes among the Fourier matrices of the size $N$ between $1$ and $100$.
This table -- see Table 2. -- ends section \ref{sec_fourier_defect}. The extended table, for $N$ ranging between $1$ and $10000$, is
available at \cite{Extended_defect_table}.

For the comments on the connections between the contents of this article with other people's work see section \ref{sec_introduction}
and other relevant sections. We hope that we have created a sound basis for any further research in the area of defect of unitary matrices,
especially complex Hadamard matrices, and the smooth families of these whose members have the moduli of their entries fixed.

%
%

\section{Acknowledgments}
The author would like to thank Karol \.Zyczkowski  from
Center for Theoretical Physics (Polish Academy of Sciences) in Warsaw,
for numerous discussions on the subject and reading this article.

The author is also greateful to Alexander Karabegov for the comments on his papers.

Financial support by the grant number N N202 090239 of
Polish Ministry of Science and Higher Education is gratefully acknowledged.

%
%


\begin{thebibliography}{99}

\bibitem{Defect}
      W. Tadej, K. \.Zyczkowski,
      Defect of a unitary matrix,
      {\sl Linear Algebra and its Applications} {\bf vol. 429 issues 2-3},
      pp. 447-481  (2008).\\
      Preprint {\sl http://www.arxiv.org/abs/math.RA/0702510}

\bibitem{Catalogue_printed}
      W. Tadej and K. \.Zyczkowski,
      A concise guide to complex Hadamard matrices,
      {\sl Open Systems \& Infor. Dynamics}  {\bf 13}, 
      pp. 133-177 (2006).

\bibitem{Catalogue}
	The online catalogue of complex Hadamard matrices.\\
      {\sl http://chaos.if.uj.edu.pl/${\sim}$karol/hadamard/}

\bibitem{KroneckerDefect}
	W. Tadej, 
      Defect of a Kronecker product of unitary matrices,
      {\sl Linear Algebra and its Applications} {\bf vol. 436 issue 7},
      pp. 1924-1959  (2012).\\
      Preprint {\sl http://arxiv.org/abs/1009.4037}

\bibitem{Karabegov}
      A. Karabegov,
      A mapping from the unitary to doubly stochastic
      matrices and symbols on a finite set,
      {\sl AIP Conf. Proc.} {\bf vol. 1079}, 
      pp. 39-50 (2008).\\
      Preprint {\sl  arXiv:0806.2357v1 [math.QA]}

\bibitem{Karabegov_old}
	A. Karabegov,
	The reconstruction of a unitary matrix from the moduli of its elements and symbols on a finite phase space,
	Preprint  YERPHI-1194(71)-89, 14 pp.

\bibitem{Johnson}
	Roger A. Horn, Charles R. Johnson,
	Matrix analysis,
	Cambridge University Press 1999.

\bibitem{PermEqClasses}
	W. Tadej,
      Permutation equivalence classes of Kronecker products
      of unitary Fourier matrices, 
      {\sl Linear Algebra and its Applications} {\bf vol. 418  issues 2-3}, pp. 719-736  (2006).\\
      Preprint {\sl http://www.arxiv.org/abs/math.RA/0501233}

\bibitem{Rhea_automorphisms}
	Christopher J. Hillar and Darren L. Rhea,
	Automorphisms of Finite Abelian Groups,
	{\sl The American Mathematical Monthly} {\bf vol. 114 no. 10},
	pp. 917-923 (2007).

\bibitem{Curran_automorphisms}
	J.N.S. Bidwell and M.J. Curran,
	Automorphisms of Finite Abelian Groups,
	{\sl Mathematical Proceedings of the Royal Irish Academy} {vol. 110 issue 1},
	pp. 57-71  (2011).

\bibitem{Mattriad2011}
	Wojciech Tadej,
	Defect of a Kronecker product of Fourier matrices,
	presentation at MATTRIAD 2011, 12-16 July, Tomar, Portugal,
	{\sl http://www.mattriad2011.ipt.pt/download/MT11$\_$Program(04Jul).pdf}


\bibitem{Banica_def_gen_four_matr}
	Teodor Banica,
	The defect of generalized Fourier matrices,
	{\sl Linear Algebra and its Applications}  {\bf vol. 438 issue 9},
	pp. 3667-3688 (2013).\\
	Preprint {\sl arXiv:1210.2556 [math.CO]}



\bibitem{Banica_first_order_deform}
	Teodor Banica,
	First order deformations of the Fourier matrix,
	{\sl Journal of Mathematical Physics}  {\bf vol. 55 issue 1 id. 012201}   (2014).\\
	Preprint {\sl arxiv: 1302.4153 [math.CO]}

\bibitem{Extended_defect_table}
	The Extended table of defects of inequivalent Fourier matrices of the size $N$ between $1$ and $10000$.\\
	{\sl http://chaos.if.uj.edu.pl/$\sim$karol/hadamard/chm\underline{\ }definitions/defect\underline{\ }10k.pdf}
\end{thebibliography}
\end{document}